\documentclass[a4paper,openany,12pt,oneside]{report}
\usepackage{psfrag}
\usepackage{epsfig}
\usepackage{graphicx}
\usepackage[all]{xy}
\xyoption{poly}
\usepackage{amsmath}
\usepackage{amssymb}
\usepackage{latexsym}
\usepackage{fancyhdr}
\usepackage{pb-diagram}
\usepackage{amsfonts}
\usepackage{eucal}
\usepackage{textcomp}
\usepackage[latin1]{inputenc}
\usepackage{graphicx}
\newtheorem{cor}{Corollary}[chapter]
\newtheorem{defin}{Definition}[chapter]
\newtheorem{teor}{Theorem}[chapter]
\newtheorem{prop}{Proposition}[chapter]
\newtheorem{lem}{Lemma}[chapter]

\newcommand{\uno}{\boldsymbol{u}}
\newcommand{\uu}{(u_{2\alpha-1} , u_{2\alpha} )}
\newcommand{\uua}{(u_{2\alpha-1} , u_{2\alpha} )}
\newcommand{\un}{u_{2\alpha-1}}
\newcommand{\um}{ u_{2\alpha} }
\newcommand{\barun}{{\overline{u}}_{2\alpha-1}}
\newcommand{\barum}{{\overline{u}}_{2\alpha}}

\newcommand{\sq}{{}^{1\pm 2}\square_{1\pm 4}^{1\pm 3}}
\newcommand{\sqp}{{}^{1+ 2}\square_{1+ 4}^{1+ 3}}

\newcommand{\sqmpp}{{}^{1- 2}\square_{1+ 4}^{1+ 3}}
\newcommand{\sqpm}{{}^{1\pm 2}\square_{1\pm 4}^{1\pm 3}}

\newcommand{\cost}{cos{\, \theta}}
\newcommand{\costq}{cos^2{\, \theta}}
\newcommand{\costd}{cos{\, 2\theta}}
\newcommand{\sit}{sin{\, \theta}}
\newcommand{\sitd}{sin{\, 2\theta}}
\newcommand{\sitq}{sin^2{\, \theta}}

\newcommand{\bfl}{\begin{flushleft}}
\newcommand{\efl}{\end{flushleft}}
\newcommand{\ka}{\mathrm{K\ddot{\mathrm{a}}hler}}

\newcommand{\grm}{Gr_4 (\mathbb{R}^8)}
\newcommand{\grms}{Gr_4 (\mathbb{R}^7)}

\newcommand{\mal}{M_{\alpha} }
\newcommand{\zw}{(\underline{z}, \underline{w})}
\newcommand{\za}{z_{\alpha}}
\newcommand{\bza}{\overline{z}_{\alpha}}
\newcommand{\wa}{w_{\alpha}}
\newcommand{\bwa}{\overline{w}_{\alpha}}
\newcommand{\zap}{z_{2\alpha}}
\newcommand{\zad}{z_{2\alpha-1}}
\newcommand{\zed}{z_{2\alpha+1}}
\newcommand{\bzap}{\overline{z}_{2\alpha}}
\newcommand{\bzad}{\overline{z}_{2\alpha-1}}
\newcommand{\bzed}{\overline{z}_{2\alpha+1}}
\newcommand{\wap}{w_{2\alpha}}
\newcommand{\wad}{w_{2\alpha-1}}
\newcommand{\bwap}{\overline{w}_{2\alpha}}
\newcommand{\bwad}{\overline{w}_{2\alpha-1}}

\begin{document}
\thispagestyle{empty}

{}

\vspace{3 cm}

\center{ Dottorato di Ricerca in Matematica}

{}

 \vspace{2 cm}

 \section*{Positive\,\, Self\,\, Dual\,\, Einstein\,\ Orbifolds\,\,
 with\,\, One-Dimensional\,\, Isometry\,\, Group}
\pagenumbering*
\vspace{2 cm}

\center{ $\mathbf{ Luca\,\, Bisconti}$ }


{}

\vspace{10 cm}

{}

\center{Universit$\grave{\mathrm{a}}$ degli Studi di Roma}

\center{"Tor Vergata"  }

\center{Dipartimento di  Matematica}

\newpage
\thispagestyle{empty}

\center{$\mathbf{XVIII\,\, Ciclo\,\, del\,\, Dottorato\,\, di\,\,
 Ricerca\,\, in\,\, Matematica}$ }

\center{$\mathbf{Universit\grave{\mathrm{a}}\,\, degli\,\, Studi\,\,
di \,\, Roma\,\, "Tor\,\,Vergata"}$  }


\vspace{4 cm}

\[
\left. \begin{array}{lcccl}
\mathrm{ Relatore\,\, di\,\, Tesi:\,\,    }  &\hspace{0,5 cm} & &\hspace{0,5 cm} &      \mathrm{Coordinatore\,\, del\,\, Dottorato:} \\
\mathrm{Prof.\,\, Paolo\,\, Piccinni}  &\hspace{0,5 cm} & &\hspace{0,5 cm} &  \mathrm{ Prof.\,\, Mauro\,\, Nacinovich} \\
\mathrm{(Dipartimento\,\, di\,\, Matematica,}  &\hspace{0,5 cm} & &\hspace{0,5 cm} & \mathrm{(Dipartimento\,\, di\,\, Matematica,}  \\
\mathrm{Universit\grave{a}\,\, di\,\, Roma\,\, "La\,\, Sapienza")}   &\hspace{0,5 cm} & &\hspace{0,5 cm} &   \mathrm{Universit\grave{a}\,\, di\,\, Roma\,\, "Tor\,\, Vergata")}  \\
\end{array}\right.
\]

\vspace{8 cm}

 \[
 \left. \begin{array}{c}
\mathrm{Tesi\,\, di \,\, Dottorato\,\, preparata\,\, per \,\, conseguire}  \\
\mathrm{  il\,\, titolo\,\, di\,\, Dottore\,\, di\,\, Ricerca \,\, in\,\, Matematica} \\
\end{array} \right.
\]

\pagenumbering{Roman}
\setcounter{page}{0}
 \tableofcontents

\chapter*{Introduction}  \markboth{Introduction}{}
\addcontentsline{toc}{chapter}{Introduction}

\bfl Oriented Riemannian geometry in dimension 4 has some special features
 coming from the action of the Hodge operator. The consequent decomposition
 of the Weyl conformal curvature tensor gives rise to the notion of self-dual
 Riemannian 4-manifolds, and this 
allows to establish the Penrose twistor
construction with its remarkable link with 3-dimensional complex geometry \cite{ah}.
The twistor correspondence has been extended to quaternion and
  quaternion $\ka$ manifolds $\mathcal{M}^{4n}$, $n\geq 2$ \cite{sa}, \cite{le},
and several
 motivations lead today to recognize that self-dual Einstein metrics
(SDE) are for $n=1$ the most appropriate choice for the
 quaternion $\ka$ condition.
One of these motivations comes from the Galicki - Lawson quotient
 construction in \cite{kr} and \cite{kl}. Namely the reduction procedure
  of a quaternion $\ka$ metric leads automatically to a SDE metric when
  the quotient has real dimension 4.
The classical theorem by
 N. Hitchin  in \cite{hi} reduced to $S^4$ and $\mathbb{CP}^2$
 the only compact SDE smooth Riemannian manifolds with positive scalar
curvature $s>0$. However, there are several examples of  $4-$dimensional
SDE orbifolds $\mathcal{O}^4$ with $s>0$.
 In order to construct these examples,  quaternion $\ka$ and
  $3-$Sasakian reduction techniques are currently the most powerful
   tools. \newline
For any quaternion $\ka$ manifold $\mathcal{M}^{4n}$ of positive scalar
 curvature
$s>0$, the following diagram of fibrations is consistent with the
 reduction construction:

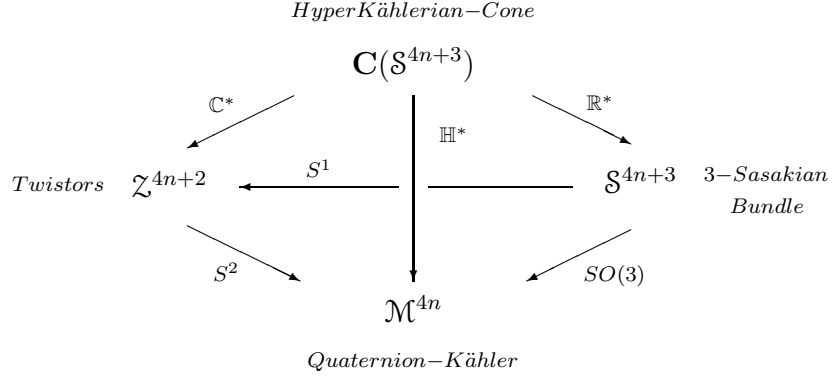
\[
   \dgARROWLENGTH=0.1\dgARROWLENGTH
   \begin{diagram}
   \node[2]{\overset{ \left.\begin{array}{c} {}^{ HyperK\ddot{a}hlerian-Cone } \\ \end{array}\right.}{\mathbf{C}(\mathcal{S}^{4n+3})}}\arrow[2]{s,r}{ \left.\begin{array}{c}\overset{\mathbb{H}^*}{} \\ \\ \\ \end{array}\right.}
      \arrow{sw,t}{ \mathbb{C}^*}  \arrow{se,t}{\mathbb{R}^*} \\
   \node{\left. \begin{array}{cc} \overset{Twistors}{} & \mathcal{Z}^{4n+2}\\ \end{array} \right.  } \arrow{se,b}{S^2}  \node{}
   \arrow{w,t}{S^1}     \node{\left. \begin{array}{cc}\mathcal{S}^{4n+3} &  \underset{{Bundle} }{{}^{3-Sasakian}} \end{array} \right.} \arrow{w,-} \arrow{sw,b}{SO(3)} \\
                           \node[2]{\underset{ \left.\begin{array}{c}  {}_{Quaternion-K\ddot{a}hler} \\ \end{array} \right. }{\mathcal{M}^{4n} } }   \\
   \end{diagram}
\]

 Consider an isometric action of a group $G$   on the $3-$Sasakian manifold
 $\mathcal{S}^{4n+3}$,
 then we can define the $3-$Sasakian moment map
 $\mu_{\mathcal{S}}\, : \mathcal{S}^{4n+3}\rightarrow (Lie\, G)^*\times \mathbb{R}^3\cong
 \mathfrak{g}^*\times\mathbb{R}^3$.
  If we assume
 both that $0$ is a regular value of $\mu_{\mathcal{S}}$ and $G$ acts freely on
  $\mu_{\mathcal{S}}^{-1}(0)\subset
  \mathcal{S}^{4n+3} $ then the quotient $\mu_{\mathcal{S}}^{-1}(0)/G$ is a
  smooth $3-$Sasakian manifold
  of dimension $4(n- dim \, \mathfrak{g})+3$. \newline
   Under some technical hypotheses
  we can reproduce this kind of construction on each  space involved in the
  fibrations 
  \begin{equation}         \label{C0 : 1a}
  \mathcal{S}^{4n+3} \rightarrow \mathcal{Z}^{4n+2}\rightarrow \mathcal{M}^{4n} ,
\end{equation}
   then,
  the mentioned reduction construction can be carried out. \newline
  Now, we give an example of $U(1)-$quotient reduction.
  The diagram below describes  all the homogeneous circle quotients given by
   the $U(1)$
   left multiplication
  on $\mathbb{H}^{n+1}$   $($ $U(1)\times \mathbb{H}^{n+1}\ni(\rho, \uno)\mapsto \rho \uno$
  $\in\mathbb{H}^{n+1})$. Then we have 
  \begin{equation}      \label{C0 : 1a}
  \left. \begin{array}{ccccccccccccc}
           &                           &  \mathbb{H}{}^{n+1}\backslash \{0\}                &                          &                   & & \overset{U(1)}{\Longrightarrow}& &             &                                   & \mathcal{N}                                           &                          &              \\
           & {}^{\mathbb{R}^+}\swarrow &                                                    & \searrow ^{\mathbb{C}^*} &                   & &                                & &             &  {}^{t(\mathbb{C})}\swarrow       &                                                       & \searrow^{t(\mathbb{R})} &              \\
  S^{4n+3} &                           &          \Bigg{ \downarrow}{}^{\mathbb{H}^*}       &                          & \mathbb{CP}^{2n+1}& & \overset{U(1)}{\Longrightarrow}& & \mathcal{Z} &                                   & {}\,\,\,\,\,\,\Bigg{ \downarrow}{}^{t(\mathbb{H})}    &                          &\mathcal{S},  \\
           & \searrow                  &                                                    &  \swarrow                &                   & &                                & &             &  \searrow                         &                                                       & \swarrow                 &              \\
           &                           &    \mathbb{HP}^{n}                                 &                          &                   & & \overset{U(1)}{\Longrightarrow}& &             &                                   &         Gr_{2}(\mathbb{C}^{n+1})                      &                          &              \\
  \end{array} \right.
  \end{equation}
 where  $\mathcal{S}= \frac{U(n+1)}{U(n-1)\times U(1)}$, $\mathcal{N}=
 \mathbf{C}(\mathcal{S})$,
$\mathcal{Z}= \frac{U(n+1)}{U(n-1)\times U(1)\times U(1)}$ and $t(\mathbb{F}):=
\frac{F^*}{\mathbb{Z}_2}$,$\mathbb{F}=\mathbb{R}$,$\mathbb{C}$,$\mathbb{H}$. \efl

\bfl Actions of tori $T^{n-1}_{\Theta}\subset Sp(n+1)$, with weight matrix
$\Theta$, on the Hopf
fibration \efl
\begin{equation}       \label{C0 : 2}
 S^{4n+3}\rightarrow \mathbb{CP}^{2n+1}\rightarrow \mathbb{HP}^n,
\end{equation}
 \bfl produce quotients like:\efl
 \begin{equation}            \label{C0 : 3}
 \mathcal{S}^7 \rightarrow \mathcal{Z}^6\rightarrow \mathcal{O}^4,
 \end{equation}
 \bfl where $\mathcal{O}^4$ is a compact SDE orbifold with positive
  scalar curvature $s>0 $. Now, we give some details about this kind of
   actions and
  the related reduction constructions,
  then we state two remarkable results which involve these facts.
 By varying the integer matrix of weights\efl
 \begin{equation}                    \label{C0 : 4}
 \Theta: =
 \left(\begin{array}{cccc}
 a^1_1 & a^1_2 & \ldots & a^1_{n+1}  \\
  a^2_1 & a^2_2 & \ldots & a^2_{n+1}  \\
  \vdots  & \vdots   & \ddots  & \vdots\\
 a^{n-1}_1 & a^{n-1}_2 & \ldots & a^{n-1}_{n+1}  \\
 \end{array} \right) \in M_{n-1\times n+1} (\mathbb{Z}),
 \end{equation}
\bfl we get different actions of tori $T^{n-1}_{\Theta}\cong
 \overbrace{S^1\times S^1\times \ldots S^1 }^{n-1\,\,
 \mathrm{times}}$ on $\mathbb{HP}^n$. In particular,
  one can act with  $T^{n-1}_{\Theta}$ by left multiplication on the
  quaternionic
  vector column
 $\uno= (u_1, u_2, \ldots,$
 $ u_n, u_{n+1})\in \mathbb{H}^{n+1}$ by a  diagonal matrix \efl

\begin{equation}         \label{C0 : 5a}
B(\Theta):=\left( \begin{array}{ccc}
 \prod^{n-1}_{\alpha=1}  \tau_{\alpha}^{a_1^{\alpha}} & \ldots  & 0 \\
\vdots                                             & \ddots  & \vdots \\
0                                                  &  \ldots &
\prod^{n-1}_{\alpha=1}\tau_{\alpha}^{a_{n+1}^{\alpha}} \\
\end{array} \right),
\end{equation}

 \bfl where $(\tau_1, \ldots, \tau_{n-1})\in S^1\times\ldots\times
 S^1\cong T^{n-1}_{\Theta}$. 
  Choices of admissible weight matrices $\Theta$
 produce quotient
of $S^{4n+3}$ that are smooth $7-$dimensional $3-$Sasakian manifolds
$\mathcal{M}^7(\Theta)$.
The corrisponding quaternion $\ka$ quotients give orbifolds
$\mathcal{O}^4(\Theta)$.
This costruction
allows to prove:  \efl



\begin{teor} \emph{(Boyer-Galicki-Mann-Rees, 1998)}
There exist compact, $T^2-$
symmetric, self dual Einstein orbifolds $(\mathcal{O}^4, h)$
of positive scalar curvature
$s>0$
with arbitrary second Betti number.  $\square$
\end{teor}

\bfl The following  inverse statement holds:  \efl

\begin{teor} \emph{(Bielawski, 1999; Calderbank-Singer, 2006)}
Let $(\mathcal{O}^4, g)$ be a compact SDE orbifold with positive
scalar curvature
$s>0$ and two commuting Killing vector fields. Then $(\mathcal{O}^4, g)$
can be obtained
as quotient of $\mathbb{HP}^n$ by the action of a torus
 $T^{n-1}_{\Theta}$.$\square$
\end{teor}

\vspace{0,5 cm}

\bfl Let us consider another interesting weighted action of the torus
$T^{k}_{\Theta}\subset SO(n+1)\subset$
$ Sp(n+1)$ on $\mathbb{H}^{n+1}$, where $k:=[\frac{n+1}{2}]$, defined
by a  $2\times 2$
  block diagonal
 matrix $A(\Theta)$   \efl
 \begin{equation}         \label{C0 : 5}
 \left(\begin{array}{c|c|c|c}
 A(\theta_1) & 0 & \ldots & 0  \\
 \hline
  0 & A(\theta_2) & \ldots & 0  \\
  \hline
  \vdots  & \vdots   & \ddots  & \vdots\\
  \hline
 0 & 0 & \ldots & A(\theta_{k})  \\
 \end{array} \right)  \quad\quad or\quad\quad
\left(\begin{array}{c|c|c|c}
 1 & 0 & \ldots & 0  \\
 \hline
  0 & A(\theta_1) & \ldots & 0 \\
  \hline
  \vdots  & \vdots   & \ddots  & \vdots\\
  \hline
 0 & 0 & \ldots & A(\theta_{k})  \\
 \end{array} \right).
 \end{equation}
\bfl (according to wheter $n$ is even or odd), and
$A_{\alpha}:=A(\theta_{\alpha})=
\left( \begin{array}{cc}
\cost_{\alpha} & \sit_{\alpha}  \\
-\sit_{\alpha} & \cost_{\alpha}  \\
\end{array} \right),$ where $\theta_{\alpha}=
\sum_{l=1}^{n} a^{\alpha}_l t_{l}, $
$t_l\in [ 0, 2\pi),$ and the integer weights $a^{\alpha}_l$ are
contained in $\Theta$.  \efl

\bfl Moreover, we can obtain another natural quotient by using the
action of $Sp(1)$ via left
 multiplication on $S^{4n+3}\subset\mathbb{H}^{n+1} $
and $\mathbb{HP}^n$:
\begin{equation}                  \label{C0 : 6}
 \left. \begin{array}{ccccc}
\mathbb{H}^{n+1}\supseteq S^{4n+3} &    \overset{Sp(1)}{\Longrightarrow} &
 \frac{SO(n+1)}{SO(n-3)\times Sp(1)}   &
 \overset{SO(3))}{\longrightarrow}  &    \overbrace{\frac{SO(n+1)}{SO(n-3)
 \times SO(4)}.}^{Gr_4(
 \mathbb{R}^{n+1}):=} \\
 \end{array} \right.
\end{equation}
As we will see in a moment, these quotients can be used to obtain $4-$dimensional
 SDE orbifolds.
In fact, if we assume  $n=6$,  the $Sp(1)-$quotient is the $12-$dimensional
 $Gr_4(\mathbb{R}^7)$,
and a further quotient is provided by the following statement.\efl

  \begin{teor} \emph{(Kobak-Swann, 1993)}     The diagonal subgroup
  \begin{equation}
  U(1)\subset U(3)\subset SO(7),
  \end{equation}
 acts on  $Gr_4(\mathbb{R}^7)$ and the quaternionic $\ka$ is an
  $8-$dimensional, compact, SDE
  orbifold
 $\mathbb{Z}_3\backslash G_2/ SO(4).$  $\square$
 \end{teor}

\bfl Together with the corrisponding $3-$Sasakian quotient,
this gives \cite{kk}:

\begin{equation}                  \label{C0 : 7}
\left. \begin{array}{ccccc}
S^{27} &    \overset{Sp(1)}{\Longrightarrow} &
 \overbrace{ \frac{SO(7)}{SO(3)\times Sp(1)}}^{ \mathcal{S}^{15}\cong  }   &
 \overset{U(1)}{\Longrightarrow}  &    \mathbb{Z}_3\backslash G_2/ Sp(1). \\
 \quad \downarrow^{_{S^3}}& &{}^{_{SO(3)}}\downarrow\quad& &{}^{_{SO(3)}}\downarrow\quad \\
{\mathbb{HP}}^6 &    \overset{Sp(1)}{\Longrightarrow} &
 \underbrace{ \frac{SO(7)}{SO(3)\times SO(4)}}_{\cong Gr_{4} ({\mathbb{R}}^7) }   &
 \overset{U(1)}{\Longrightarrow}  &    \mathbb{Z}_3\backslash G_2/ SO(4) \\
\end{array}  \right.
\end{equation}

 The above   quotient construction
admits interesting generalizations. In fact,  the construction in (\ref{C0 : 7})
 can be "deformed"
by introducing
weights for the $U(1)-$action in theorem $0.3$ \cite{sk},
   and these weighted actions
provide  new examples of quotient orbifolds.
Namely, consider the
 weighted action of $T^1_{\Omega}\cong $
 $U(1) \subset  SO(7)$
or $T^2_{\Omega}\cong U(1)\times$
$ U(1)\subset SO(7)$ on $S^{27}\subset \mathbb{H}^7$,
where 
$\Omega$ is either:
\begin{equation}          \label{C0 : 8}
\left(\begin{array}{ccc}
p_1 & p_2 & p_3 \\
\end{array}  \right) \in M_{1\times 3}(\mathbb{Z}), \quad
\mathrm{or}  \quad
\left(\begin{array}{ccc}
p_1 & p_2 & p_3 \\
q_1 & q_2 & q_3 \\
\end{array}  \right) \in M_{2\times 3}(\mathbb{Z}).
\end{equation}
Then, according to the definitions given in (\ref{C0 : 4})
and (\ref{C0 : 5}), we have:
\begin{equation}             \label{C0 : 9}
A_{\alpha}:= A(\omega_{\alpha}) =
\left( \begin{array}{cc}
cos\,\omega_{\alpha} & sin\,\omega_{\alpha}  \\
-sin\,\omega_{\alpha}  & cos\,\omega_{\alpha}  \\
\end{array} \right)\in U(1),
\end{equation}
where $\omega_{\alpha}$ is either $p_{\alpha} t$ or $p_{\alpha}t
+ q_{\alpha}s$ with
$t,s\in [0, 2\pi)$, $\alpha\in \{1,2,3\}$,
and the matrix:
\begin{equation}              \label{C0 : 10}
A(\Omega)=\left(\begin{array}{c|c|c|c}
 1 & 0 & 0 & 0  \\
 \hline
  0 & A(\omega_1) &0& 0 \\
  \hline
  0 & 0   & A(\omega_2) & 0\\
  \hline
 0 & 0 & 0 & A(\omega_3).  \\
 \end{array} \right)\in SO(7),
 \end{equation}
 acts on $S^{27}\subset \mathbb{H}^7$ by rotating pairs
 of quaternionic coordinates  $(u_{2\alpha}, u_{2\alpha+1})$, $\alpha\in\{1,$
 $2,3\}$. The following result holds \efl

\vspace{0,5 cm}

 \begin{teor} \emph{(Boyer-Galicki-Piccinni, 2002)} Consider the sphere
   $S^{27}\subset \mathbb{H}^7$,
 then we have:
\bfl $a) $ for any non zero  triple of weights
 $\Omega=( {p}_{1} ,p_{2} ,p_{3} )\in \mathbb{Z}^3$
  there is an isometric action of
  $Sp(1)\times T^1_{\Omega}$  on  the sphere
  $S^{27}\subset \mathbb{H}^7.$  If
 $0<p_1 <p_2 <p_3$ with $gcd( p_1 \pm p_2, p_1 \pm p_3)=$
 $ 1$, then the Hopf fibration
 $ S^{27}  \rightarrow \mathbb{HP}^6$ gives as  quotient  \efl

 \begin{equation}              \label{C0 : 11}
\left. \begin{array}{ccccc}
  S^{27} & \rightarrow &\mathbb{CP}^{13} & \rightarrow  & \mathbb{HP}^6  \\
  \Downarrow  &        &  \Downarrow     &           &  \Downarrow  \\
   \mathcal{M}^{11} & \rightarrow   &  \mathcal{Z}^{10} & \rightarrow  &   \mathcal{O}^8(\Omega) \\
   \end{array} \right.
 \end{equation}

\bfl \emph{  where  $\mathcal{M}^{11} $ is a smooth  $3-$Sasakian manifold
  and  $ \mathcal{O}^8(\Omega) $ its quaternion
  $\ka$ leaf orbifold. }  \efl \vspace{0,5 cm}

   \bfl  $b)$ Let $\Omega\in M_{2\times 3}(\mathbb{Z})$  be
  such that each of its  $2\times 2$ minor determinants $\Delta_{\alpha\beta}$
 does not vanish. Assume their sum does not vanish and that none of these
 three determinants is equal to the sum of the other two.
 Then the quaternionic K$\ddot{a}hler$ quotient of $\mathbb{HP}^6$
 by $G=Sp(1)\times T^2_{\Omega}$ is a compact SDE $4-$dimensional
 orbifold with positive scalar curvature $s>0$ and with a $1-$dimensional
  group of isometries. $\square$  \efl
  \end{teor}

\bfl We denote by $\mathcal{O}^4(\Omega)$ the orbifolds
obtained as  quotient by the action
of $G=Sp(1)\times$
$ T^2_{\Omega}$
on $S^{27}\subset\mathbb{H}^7$. \efl

\bfl  First examples of compact self-dual Einstein $4-$orbifolds of positive
 scalar curvature were constructed by Galicky and Lawson via quotient reduction
 similarly to those obtained as symplectic quotients in \cite{kl}.
 Currently, known examples
  of such orbifolds
 can be divided in three groups:
 \efl

 \begin{itemize}
 \item[(a)] The $SO(3)$-invariant, cohomogeneity one orbifold
 metric on $S^4$ discovered by
  Hitchin \cite{hd}.  \\

  \item[(b)] $T^2-$invariant toric orbifolds metrics constructed by
  Boyer, Galicki, Mann and Rees \cite{mr}
  ( Theorems $0.1$, $0.2$) as quaternionic K$\ddot{\mathrm{a}}$hler
   quotients \\

  \item[(c)] $S^1-$invariant orbifold metrics of: Galicki-Nitta \cite{gn},
  Boyer-Galicki-Piccinni \cite{bg}
  and the ones we describe in this thesis. \\

 \end{itemize}

  \bfl
  We indicate as compact positive toric
 self-dual $4-$orbifolds $(\mathcal{O}, g)$ the orbifolds
  obtained as a
 quaternion-K$\ddot{\mathrm{a}}$hler reduction of
 some $\mathbb{HP}^{n}$ via an action
 of some $(n-1)$-torus subgroup of $Sp(n+1)$. Note that, the orbifolds
 mentioned in theorem $0.1$ and $0.2$ are toric, instead the ones in
 theorem $0.3$ and $0.4$ are not toric.
 Galicki and Nitta \cite{gn} gave the first examples of compact positive
 self-dual $4-$orbifolds which are not toric. In particular,
they constructed one of these examples  for each
   non-Abelian subgroup of $Sp(1)$.
    The examples provided
    by theorem  $0.4$
    are not toric
   but  they are very different from those one of  Galicki and Nitta \cite{gn},
   in fact in the former case all the orbifold uniformazing groups are Abelian.
  \efl

\bfl The starting point of this thesis is to note that actions of
the type described in theorem $0.4$
 can be considered on any quaternion-K$\ddot{\mathrm{a}}$hler
Grassmannian ( look at diagrams (\ref{C0 : 6}) and (\ref{C0 : 7}))  \efl
\begin{equation}         \label{C0 : 12}
Gr_4 (\mathbb{R}^{n+1})\cong\frac{SO(n+1)}{SO(n-3)\times SO(4)},
\end{equation}
 \bfl
whose real dimension is $4(n-3)$.
 Therefore to end up in dimension $4$ a $(n-4)-$
dimensional
torus is needed. Since the maximal torus in $SO(n+1)$ has
 $dim = [\frac{n+1}{2}]$,
in order to be able to consider such an action with weights
we need       \efl
\begin{equation}        \label{C0 : 13}
 n-4 < \bigg[\frac{n+1}{2}\bigg].
\end{equation}
\bfl
This gives, according to whether $n+1$ is even or odd  \efl
\begin{equation}                       \label{C0 : 14}
\begin{aligned}
& n-4 < \frac{n+1}{2} \quad \quad \quad \mathrm{or} \quad
\quad \quad n -4 < \frac{n}{2}, \\
\mathrm{i.\, \, e.} \quad \quad \quad \quad \quad \quad   &   \\
& \left. \begin{array}{l}
           {}\,\,\, n < 9   \\
            (n+1\,\,\mathrm{even} )
            \end{array} \right.
  \quad \quad \,\, \mathrm{or} \quad \quad \quad \quad \left. \begin{array}{l}
           {}\,\, n < 8.   \\
            (n+1\,\, \mathrm{odd} )
            \end{array} \right.  \hspace{10 cm}\\
\end{aligned}
\end{equation}

\bfl Thus the only possible  cases are   \efl
\begin{equation}  \label{C0: 14}
 \left. \begin{array}{lllll}
   Gr_4 (\mathbb{R}^6)\cong Gr_2(\mathbb{C}^4),  &
     {} \hspace{1 cm}  {} & Gr_4 (\mathbb{R}^7),
       &    {}  \hspace{1 cm} \mathrm{and} \quad \quad  {}& Gr_4 (\mathbb{R}^8). \\
       \end{array} \right.
\end{equation}

\bfl The first two cases in (\ref{C0: 14}) have been examined
 in \cite{mr}  and  \cite{bg}, respectively.   \efl
\bfl The present thesis is essentially
devoted to the third case and to its comparison with the orbifolds construction
 in \cite{bg}.  \efl
  \bfl We describe now the main results obtained. Consider the
   Grassmannian $Gr_4 (\mathbb{R}^8)$,
   acted on isometrically by
 a $3-$
 torus
 $T^3_{\Theta}\subset T^4\subset SO(8)\subset Sp(8)$,
 where $\Theta$ is a weight matrix of the type defined in (\ref{C0 : 5}).
Equivalently, consider the action of $G:=Sp(1)\times T^3_{\Theta}$ $\subset Sp(8)$
on the the sphere $S^{31}\subset \mathbb{H}^8$,
  where $Sp(1)$ acts by left multiplication and the
$T^3_{\Theta}-$action is defined via the matrices  \efl
\begin{equation}      \label{C0 : 15}
\begin{aligned}
&
A(\Theta)=\left( \begin{array}{c|c|c|c}
A_1(t,s,r) & 0             & 0            & 0          \\
\hline
0 & A_2 (t,s,r)          & 0 & 0 \\
\hline
0 & 0 & A_3(t,s,r)  & 0\\
 \hline
0 & 0 & 0 & A_4(t,s,r) \\
\end{array} \right) \in T^4,    \\
  \mathrm{where}  \hspace{1 cm}& \\
 &  A_{\alpha}= A_{\alpha}(t,s,r) =
\left( \begin{array}{cc}
cos(p_{\alpha} t + q_{\alpha} s+ l_{\alpha} r ) & sin(p_{\alpha} t + q_{\alpha} s+ l_{\alpha} r )   \\
-sin(p_{\alpha} t + q_{\alpha} s+ l_{\alpha} r )  & cos(p_{\alpha} t + q_{\alpha} s+ l_{\alpha} r ) \\
\end{array} \right).    \hspace{2 cm}
\end{aligned}
\end{equation}
\bfl and the weights  are collected in the matrix    \efl
\begin{equation}           \label{C0 : 16}
\Theta:= \left( \begin{array}{cccc}
p_1 & p_2 & p_3 & p_4\\
q_1 & q_2 & q_3 & q_4\\
l_1 & l_2 & l_3 & l_4\\
\end{array} \right).
\end{equation}

\bfl We prove in this thesis the following statements,
referred to as theorems $A,\ldots,F$.  \efl

\emph{  \bfl  $\mathbf{Theorem\,\, A}$
 Let $\Theta\in M_{3\times 4}(\mathbb{Z})$ be a non null integral matrix,
 such that each of its  $3\times 3$ minor determinants \efl }
 \begin{equation}  \label{C0 : 17}
 \Delta_{\alpha\beta\gamma} :=
\left\vert \begin{array}{ccc}
p_{\alpha} & q_{\alpha} & l_{\alpha} \\
p_{\beta} & q_{\beta} & l_{\beta} \\
p_{\gamma} & q_{\gamma} & l_{\gamma} \\
\end{array} \right\vert , \,\,\,
  1\leq \alpha<\beta<\gamma\leq4 ,
 \end{equation}
 \emph{ \bfl
does not vanish. Moreover, suppose that these minor determinants
  $\Delta_{\alpha\beta\gamma}$ are such that, the further determinants \efl }
\begin{equation}   \label{C0 : 18}
{}^{1\pm 2}\square_{1\pm 3}^{1\pm 4} :=
\left\vert \begin{array}{ccc}
p_1 \pm p_2 & q_1\pm q_2 & l_1 \pm l_2 \\
p_1 \pm p_3 & q_1\pm q_3 & l_1 \pm l_3 \\
p_1 \pm p_4 & q_1 \pm q_4 & l_1 \pm l_4 \\
\end{array} \right\vert ,
\end{equation}
\emph{ \bfl are all non zero. Then for each such a matrix $\Theta$ which
satisfies these properties,
there exists a  compact SDE  $4-$dimensional
 orbifold $\mathcal{O}^4(\Theta )$ with positive scalar curvature and a
 $1-$dimensional group of isometries.$\square$ \efl }

\bfl Recall that the twistor space  $\mathcal{Z}^6(\Theta)$
 is the $ S^2$ bundle over the riemannian orbifold
 $\mathcal{O}^4(\Theta)$ that parametrizes
 its local almost complex structures. Consider the group
  $\tilde{G}:= G\times U(1)= $
 $T^3_{\Theta}\times Sp(1)\times U(1)$,  which is a subgroup, up to
  the central $\mathbb{Z}/2$,
  of $Sp(8)\cdot Sp(1)\subset SO(32)=Isom(S^{31})$.
  Then, look at the action of $\tilde{G}$ on $N(\Theta)= $
 $\mu^{-1}(0) \cap \nu^{-1}(0)$ ( where
$\mu$ and $\nu$ are the moment maps of $Sp(1)$ and $T^3_{\Theta}$,
 respectively) \efl
 \begin{equation}  \label{C0 : 19}
 \begin{aligned}
 &  \Phi:\, T^3_{\Theta}\times Sp(1)\times U(1)\times
  N(\Theta) \longrightarrow N(\Theta),\\
 & \quad \quad\quad \quad\big( (A(\Theta), \lambda, \rho);(\underline{z},
  \underline{w} ) \big)\longmapsto
 \Phi \big( ( A(\Theta), \lambda, \rho)\big)\big((\underline{z},
  \underline{w} )\big), \\
\mathrm{where} \hspace{1 cm} & \\
 & \quad \quad \quad \quad \Phi \big( ( A(\Theta), \lambda, \rho)
 \big)\big((\underline{z}, \underline{w} ) \big)
 :=   A(\Theta) \lambda  \left( \begin{array}{c}
\underline{z}\\
\underline{w} \\
\end{array}\right)
 \rho.   \hspace{5 cm} {}
 \end{aligned}
 \end{equation}
 \bfl Here we have identified $\mathbb{C}^8\times \mathbb{C}^8\cong \mathbb{H}^8$ by using
 the relation $u_{\alpha} = z_{\alpha} + jw_{\alpha}\in \mathbb{C}^2\times \mathbb{C}^2$
  for each ${\alpha}\in \{1,\ldots, 8\}$ and
  $\zw:=\boldsymbol{u}=
 (u_1,u_2,\ldots,u_7,u_8)\in \mathbb{H}^8$.
  The twistor space $\mathcal{Z}^6(\Theta)$
  will be thought as the leaf
space of the $\tilde{G}-$action on the manifold $N(\Theta)$.
 By considering together the singular locus of the $3-Sasakian$ orbifold
 $\mathcal{M}^7(\Theta)$
 and  the singular locus of the twistor
space $\mathcal{Z}^6(\Theta)$, we give the following
 descriptions of the singularities  of the quaternion $\ka$
 orbifold $\mathcal{O}^4(\Theta)$. \efl

\bfl $\mathbf{Theorem\,\, B}$  \emph{
Let $\mathcal{Z}^6(\Theta)$ be the twistor space over the SDE orbifold
 $\mathcal{O}^4(\Theta)$. Then, depending on the minor determinants
  $\Delta_{\alpha\beta\gamma}$
 and $\sq$,  the singular locus $\Sigma(\Theta)$ of $\mathcal{Z}^6(\Theta)$ consists of }\efl

 \bfl \emph{ 1) one isolated $2-$sphere $S^2$, whose isotropy, up to  the non effective subgroup,
 only depends on one of the determinants $\sq$;  }

 \vspace{0,5 cm}

\emph{ 2) eleven disjoint $2-$spheres $S^2$. These are obtained from
  the following $\tilde{G}-$strata
  $($ with respect to the action
 of $\tilde{G}$ on $\mathbb{H}^8$ $)$ }\efl
 \begin{equation}  \label{C0 : 20}
 \begin{aligned}
& a)\,\,\, \quad \overset{=}{S}{}^{123}_4:=\tilde{G}\cdot \Bigg\{
\left( \begin{array}{cc|cc|cc|cc}
z_1 & z_2  & z_3 & z_4   &    z_5 & z_6        & 0    &  0   \\
0   &    0 &  0  & 0     &   0   &  0          & w_7  &  w_8\\
\end{array}\right) \Bigg\},   \\
&  b) \,\,\, \quad \overset{=}{S}{}^{124}_3:=\tilde{G}\cdot \Bigg\{
\left( \begin{array}{cc|cc|cc|cc}
z_1 & z_1  & z_3 & z_4       &   0   &    0          & z_7   & z_8   \\
0   &    0  & 0 &    0     &   w_5 &   w_6        & 0 & 0\\
\end{array}\right) \Bigg\},   \hspace{10 cm} \\
&  c)\,\,\, \quad  \overset{=}{S}{}^{134}_2:=\tilde{G}\cdot \Bigg\{
\left( \begin{array}{cc|cc|cc|cc}
z_1 & z_2  & 0    &  0           &    z_5 & z_6             & z_7   & z_8  \\
0   &    0 & w_3  & w_4   &   0  &  0   &  0    &  0          \\
\end{array}\right) \Bigg\}, \hspace{10 cm }  \\
&  d) \,\,\, \quad \overset{=}{S}{}^{234}_1:=\tilde{G}\cdot \Bigg\{
\left( \begin{array}{cc|cc|cc|cc}
0   &    0      & z_3 & z_4 &    z_5 & z_6             & z_7   & z_8  \\
w_1  & w_2& 0   &    0     &      0  &    0              &  0    &  0\   \
\end{array}\right) \Bigg\}.  \hspace{20 cm}  \\
& e) \,\,\, \quad \overset{=}{S}{}^{12}_{34}:=\tilde{G}\cdot \Bigg\{
\left( \begin{array}{cc|cc|cc|cc}
z_1 & z_2  & z_3 & z_4   &    0  &  0        & 0    &  0   \\
0   &    0 &  0  & 0     &   w_5   & w_6          & w_7  &  w_8\\
\end{array}\right) \Bigg\},   \\
&  f) \,\,\, \quad \overset{=}{S}{}^{14}_{23}:=\tilde{G}\cdot \Bigg\{
\left( \begin{array}{cc|cc|cc|cc}
z_1 & z_1  & 0 & 0       &   0   &    0          & z_7   & z_8   \\
0   &    0  & w_3 &    w_4     &   w_5 &   w_6        & 0 & 0\\
\end{array}\right) \Bigg\},   \hspace{10 cm} \\
&  g)\,\,\, \quad \overset{=}{S}{}^{12}_{34}:=\tilde{G}\cdot \Bigg\{
\left( \begin{array}{cc|cc|cc|cc}
z_1 & z_2  & 0    &  0           &    z_5 & z_6             & 0  & 0  \\
0   &    0 & w_3  & w_4    &   0  &  0   &  w_7   &  w_8          \\
\end{array}\right) \Bigg\}. \hspace{10 cm }  \\
\end{aligned}
 \end{equation}
\bfl \emph{Each of these strata 
intersects $N(\Theta)$ and
$\overset{=}{S}{}^{\alpha\beta\gamma}_{\delta}\cap N(\Theta)$ is
formed by two connected components which gives two
 $\widetilde{S}$ strata for  the $\tilde{G}-$action on $N(\Theta)$.
  Instead, each of the
 $\overset{=}{S}{}^{\alpha\beta}_{\gamma\delta}\cap N(\Theta)$ is
 connected and provides
 exactly one strata $\widetilde{S}$ of $N(\Theta)$.
 Moreover, we have that
$\widetilde{S}/\tilde{G}\cong S^2$ and for each of the $2-spheres$ linked to one of
    $\overset{=}{S}{}^{\alpha\beta\gamma}_{\delta}$
  the isotropy only depends on
one of   the minor determinants $\pm \Delta_{\alpha \beta \gamma}$, for the remaining
cases the isotropy dipends on a linear combination of $\Delta_{\alpha \beta \gamma}$;   }

\vspace{0,5 cm}

\emph{ 3) three sets of  points such that, for each of them,
 the respective points are joined by
 one of the  $2-$spheres
$S^2\cong  \, \overset{=}{S}{}^{\alpha\beta}_{\alpha\gamma}\cap N(\Theta)\big/\tilde{G}$.
 Each of these sets
consists at most of four points.}\efl

  \bfl \emph{In the case when some of the minor determinants
  $\Delta_{\alpha\beta\gamma}$
 and $\sq$ assumes values $\pm 1$ or when the isotropy,
 associated to some
  of the points mentioned in
  $3)$, vanishes, then  the singular locus
 $\Sigma(\Theta)$ is obtained by removing from the above list the  sets
 whose isotropy depends on  the mentioned determinants. $\square$} \efl

\vspace{0,5 cm}

\bfl $\mathbf{Theorem\,\, C}$ \emph{ The singular points of the
$3-$Sasakian orbifold $\mathcal{M}^7(\Theta)$
 come from the singular strata $\widetilde{S}$
 of $N(\Theta)$,
 which are obtained by intersecting $N(\Theta)$ with
 the following $G-$strata $\overline{S}$
  of $\mathbb{H}^8$,
 and the ones  which
 can be obtained by changing $(w'_1,
 w'_2)$ with $(w''_1,  w''_2)$:  }\efl
 \begin{align*}
 & a) \,\,\, \quad  \overline{S}{}^{1234}:= 
  \bigg\{
\left( \begin{array}{cc|cc|cc|cc}
z_1 & z_2  & z_3 & z_4   &    z_5 & z_6        & z_7    &  z_8   \\
w'_1 & w'_2  &  w'_3  & w'_4     &   w'_5   &  w'_6          & w'_7  &  w'_8\\
\end{array}\right) \bigg\},   \\
& b) \,\,\, \quad \overline{S}{}^{123}_{4}:= 
 \bigg\{
\left( \begin{array}{cc|cc|cc|cc}
z_1 & z_2  & z_3 & z_4   &    z_5 & z_6        & z_7    &  z_8   \\
w'_1 & w'_2  &  w'_3  & w'_4     &   w'_5   &  w'_6          & w''_7  &  w''_8\\
\end{array}\right) \bigg\},  \\
&c) \,\,\,\quad \overline{S}{}^{124}_{3}:= 
\bigg\{
\left( \begin{array}{cc|cc|cc|cc}
z_1 & z_2  & z_3 & z_4   &    z_5   & z_6           & z_7   & z_8  \\
w'_1 & w'_2  &  w'_3  & w'_4     &   w''_5  &  w''_6   & w'_7  &  w'_8\\
\end{array}\right) \bigg\},  \hspace{5 cm}  \\
 & d)\,\,\, \quad \overline{S}{}^{134}_{2}:= 
 \bigg\{
\left( \begin{array}{cc|cc|cc|cc}
z_1 & z_2  & z_3 & z_4   &    z_5 & z_6        & z_7    &  z_8   \\
w'_1 & w'_2 & w''_3  & w''_4    &  w'_5  &  w'_6   & w'_7  &  w'_8\\
\end{array}\right) \bigg\},   \\
\end{align*}
\begin{equation} \label{C0 : 21}{
 \begin{aligned}
& e)\,\,\,\quad \overline{S}{}^{12}_{34}:= 
\bigg\{
\left(\begin{array}{cc|cc|cc|cc}
z_1 & z_2  & z_3 & z_4   &    z_5 & z_6        & z_7    &  z_8   \\
w'_1 & w'_2  &  w'_3  & w'_4    &   w''_5  &  w''_6        & w''_7  &  w''_8\\
\end{array}\right) \bigg\},   \\
& f)\,\,\,\quad \overline{S}{}^{13}_{24}:= 
\bigg\{
\left( \begin{array}{cc|cc|cc|cc}
z_1 & z_2  & z_3 & z_4   &    z_5 & z_6        & z_7    &  z_8   \\
w'_1 & w'_2 & w''_3  & w''_4     &  w'_5  &  w'_6          & w''_7  &  w''_8\\
\end{array}\right) \bigg\},   \\
& g)\,\,\,\quad  \overline{S}{}^{14}_{23}:= 
\bigg\{
\left( \begin{array}{cc|cc|cc|cc}
z_1 & z_2  & z_3 & z_4   &    z_5 & z_6        & z_7    &  z_8   \\
w'_1 & w'_2 & w''_3  & w''_4     &   w''_5 & w''_6   & w'_7  &  w'_8\\
\end{array} \right)    \bigg\}, \\
& h)\,\,\,\quad  \overline{S}{}^{1}_{234} := 
\bigg\{
\left( \begin{array}{cc|cc|cc|cc}
z_1 & z_2  & z_3 & z_4   &    z_5 & z_6        & z_7    &  z_8   \\
w'_1 & w'_2 & w''_3  & w''_4     &   w''_5 & w''_6        & w''_7  &  w''_8\\
\end{array}\right) \bigg\},
 \end{aligned}    }    \hspace{4,5 cm}    {}
 \end{equation}
\bfl \emph{where } \efl
\begin{equation} \label{C0 : 22}
\begin{aligned}
&  (w'_{2\alpha-1},
 w'_{2\alpha }):=\frac{e^{i\delta} }{ sin\, \varphi}
 \big( -z_{2\alpha} +iz_{2\alpha -1}\cos\, \varphi, \, z_{2\alpha-1} +
 iz_{2\alpha}\cos\, \varphi  \big), \hspace{10 cm}\\
and \,\,\,\, & \hspace{10 cm}\\
  &(w''_{2\alpha-1},  w''_{2\alpha}):= \frac{e^{i\delta} }{ sin\, \varphi} \big(-\zap -
i\zad\cos\, \varphi,\, \zad -i\zap\cos\, \varphi  \big),  \hspace{ 5 cm}\\
\end{aligned}
\end{equation}
\bfl \emph{ we can suppose $sin\, \varphi\ne 0$. Here,  the element $\lambda=\epsilon + j\sigma\in Sp(1)$
is written as
  $\epsilon =$
  $ \cost +i(\sit cos\, \varphi)$  and
$\sigma =  \sit sin\, \varphi cos\, \delta +
i(\sit sin\, \varphi sin\, \delta ).$ Then, each stratum $\overline{S}$
   in $(\ref{C0 : 21})$ generates a singular point
  $\big(\overline{S}\cap N(\Theta)\big)/ G  \in \mathcal{M}^7(\Theta)$ whose isotropy only
  depends on one of the minor determinants $\sq$. $\square$ }\efl

\bfl By referring to the orbifolds $\mathcal{O}^4 (\Omega)$ constructed in \cite{bg}
we prove also: \efl

\bfl  $\mathbf{Theorem\,\, D}$ \emph{ Let $\mathcal{Z}^6(\Omega)$ be the twistor space over the SDE orbifold
 $\mathcal{O}^4(\Omega)$. Then, depending on the minor determinants $\Delta_{\alpha\beta}$
 and $\square^{1\pm 2}_{1\pm 3}$  the singular locus $\Sigma(\Omega)\subset \mathcal{Z}^6(\Omega)$
  consists of } \efl

 \bfl \emph{ 1) one isolated $2-$sphere $S^2$, whose isotropy, up to divide by the non effective sugroup,
 only depends on one of the determinants $\square^{1\pm 2}_{1\pm 3}$; }    \efl

\bfl \emph{ 2) six points. These come from
  the following strata $($ respect to the action of $\tilde{G}$ on $\mathbb{H}^7)$}    \efl
 \begin{equation} \label{C0 : 23}
 \begin{aligned}
& a) \,\,\, \quad  {}^+\overset{=}{S}{}^{12}_3:=\tilde{G}\cdot \Bigg\{
\left( \begin{array}{c|cc|cc|cc}
0 & z_2  & z_3 & z_4   &    z_5 &  0    &  0   \\
0   &    0 &  0  & 0     &   0             & w_6  &  iw_6\\
\end{array}\right) \Bigg\},  \\
& b) \,\,\, \quad  {}^-\overset{=}{S}{}^{12}_3:=\tilde{G}\cdot \Bigg\{
\left( \begin{array}{c|cc|cc|cc}
0 & z_2  & z_3 & z_4   &   z_5   & 0    &  0   \\
0   &    0 &  0  & 0     &   0     & w_6  &  -iw_6\\
\end{array}\right) \Bigg\},   \\
&c) \,\,\,\quad  {}^+\overset{=}{S}{}^{13}_2:=\tilde{G}\cdot \Bigg\{
\left( \begin{array}{c|cc|cc|cc}
0 & z_2  & z_3    &    0   & 0           & z_6   & z_7  \\
0 &    0 &  0     &   w_4  &  iw_4   &  0    &  0          \\
\end{array}\right) \Bigg\},  \\
& d)\,\,\,\quad   {}^-\overset{=}{S}{}^{13}_2:=\tilde{G}\cdot \Bigg\{
\left( \begin{array}{c|cc|cc|cc}
0 & z_2  & z_3        &   0   &    0          & z_6   & z_7   \\
0 &    0 &     0      &   w_4 &    -iw_4      & 0     & 0\\
\end{array}\right) \Bigg\},   \\
 & e)\,\,\, \quad  {}^+\overset{=}{S}{}^{23}_1:=\tilde{G}\cdot \Bigg\{
\left( \begin{array}{c|cc|cc|cc}
    0 & 0    &  0      &    z_4 & z_5             & z_6   & z_7  \\
    0 & w_2  & iw_2    &   0    &  0              &  0    &  0          \\
\end{array}\right) \Bigg\},   \\
& f)\,\,\, \quad  {}^-\overset{=}{S}{}^{23}_1:=\tilde{G}\cdot \Bigg\{
\left( \begin{array}{c|cc|cc|cc}
0   & 0    &  0           &    z_4 & z_5      & z_6   & z_7  \\
0   & w_2  & -iw_2        &   0    &  0       &  0    &  0          \\
\end{array}\right) \Bigg\}. \hspace{10 cm }  \\
\end{aligned}
\end{equation}
\bfl \emph{Each of these strata intersects the submanifold $N(\Omega)$. Then,
 the $\tilde{G}-$strata of $N(\Omega)$ are given by
 $\widetilde{S}{}^{\alpha\beta}_{\gamma}:= {}^{\pm}\overset{=}{S}{}^{\alpha\beta}_{\gamma}\cap N(\Omega)$
 and each of the quotients
\begin{equation}
\frac{{}^{\pm}\overset{=}{S}{}^{\alpha\beta}_{\gamma}\cap N(\Omega)}{\tilde{G}},
\end{equation}
 is a point $($ here $\tilde{G}= Sp(1)\times T^2_{\Omega}\times  U(1)$ $)$.
   Moreover, for each of these points the isotropy only depends on
one of the  the minor determinants $ \pm \Delta_{\alpha \beta} $. }\efl

\bfl  \emph{ In the case when some of the minor determinants $\Delta_{\alpha\beta}$
 and $\square^{1\pm 2}_{1\pm 3}$ assume values $\pm 1$,  the singular locus
 $\Sigma(\Omega)$ can be obtained by removing from the above list the singular sets
 whose isotropy depends on one of the mentioned determinants. $\square$ }\efl

 \vspace{0,2 cm}

\bfl $\mathbf{Theorem\,\, E}$   \emph{
The singular points of the $3-$Sasakian orbifold $\mathcal{M}^7(\Omega)$
 come from the singular strata $\widetilde{S}$ of $N(\Omega)$,
 which are obtained by intersecting $N(\Omega)$ with
 the following $G-$strata $\overline{S}$ of $\mathbb{H}^7$  $($ here $G=T^2_{\Omega}\times Sp(1)$$)$}\efl
\begin{equation} \label{C0 : 24}  {
 \begin{aligned}
 & a) \,\,\, \quad  \overline{S}{}^{123}:= 
 \bigg\{
\left( \begin{array}{c|cc|cc|cc}
0 & z_2  & z_3 & z_4   &    z_5 & z_6        & z_7       \\
0 & \widetilde{w}_2  &  \widetilde{w}_3  & \widetilde{w}_4     &   \widetilde{w}_5   &  \widetilde{w}_6      & \widetilde{w}_7  \\
\end{array}\right) \bigg\},  \hspace{10 cm}   \\
& b) \,\,\, \quad \overline{S}{}^{12}_3:= 
\bigg\{
\left( \begin{array}{c|cc|cc|cc}
0 & z_2   &   z_3  & z_4      &    z_5   & z_6            & z_7       \\
0 & \widetilde{w}_2  &  \widetilde{w}_3  & \widetilde{w}_4     &   \widetilde{w}_5   &  \overset{\approx}{w}_6          & \overset{\approx}{w}_7  \\
\end{array}\right) \bigg\},   \hspace{10 cm}\\
&c) \,\,\,\quad \overline{S}{}^{13}_2:= 
\bigg\{
\left( \begin{array}{c|cc|cc|cc}
0 & z_2  & z_3      & z_4            &    z_5                & z_6           & z_7    \\
0 & \widetilde{w}_2  &  \widetilde{w}_3  & \overset{\approx}{w}_4     &  \overset{\approx}{w}_5      &  \widetilde{w}_6         & \widetilde{w}_7  \\
\end{array}\right) \bigg\},   \\
 & d)\,\,\, \quad \overline{S}{}^{23}_1:= 
 \bigg\{
\left( \begin{array}{c|cc|cc|cc}
0 & z_2  & z_3 & z_4   &    z_5 & z_6        & z_7       \\
0 &\overset{\approx}{w}_2 & \overset{\approx}{w}_3  & \widetilde{w}_4    &  \widetilde{w}_5  &  \widetilde{w}_6   & \widetilde{w}_7  \\
\end{array}\right) \bigg\},   \\
& e)\,\,\,\quad \overline{S}{}^{1}_{23}:= 
\bigg\{\left( \begin{array}{c|cc|cc|cc}
0 & z_2  & z_3 & z_4   &    z_5 & z_6        & z_7    \\
0  & \widetilde{w}_2  &  \widetilde{w}_3  & \overset{\approx}{w}_4    &   \overset{\approx}{w}_5  & \overset{\approx}{w}_6        & \overset{\approx}{w}_7  \\
\end{array}\right) \bigg\},   \\
& f)\,\,\,\quad \overline{S}{}^2_{13}:= 
 \bigg\{
\left( \begin{array}{c|cc|cc|cc}
0 & z_2   & z_3    & z_4      &    z_5 & z_6        & z_7      \\
0 & \overset{\approx}{w}_2 & \overset{\approx}{w}_3  & \widetilde{w}_4     &  \widetilde{w}_5  &  \overset{\approx}{w}_6          & \overset{\approx}{w}_7  \\
\end{array}\right) \bigg\},   \\
& g)\,\,\,\quad  \overline{S}{}^3_{12}:= 
\bigg\{
\left( \begin{array}{c|cc|cc|cc}
0 & z_2  & z_3 & z_4   &    z_5 & z_6        & z_7      \\
0 & \overset{\approx}{w}_2 & \overset{\approx}{w}_3  & \overset{\approx}{w}_4     &  \overset{\approx}{w}_5 & \widetilde{w}_6   & \widetilde{w}_7  \\
\end{array} \right)    \bigg\}, \\
& h)\,\,\,\quad  \overline{S}{}_{123} :=  
\bigg\{
\left( \begin{array}{c|cc|cc|cc}
0 & z_2  & z_3 & z_4   &    z_5 & z_6        & z_7    \\
0 & \overset{\approx}{w}_2 & \overset{\approx}{w}_3  & \overset{\approx}{w}_4     &   \overset{\approx}{w}_5 & \overset{\approx}{w}_6        & \overset{\approx}{w}_7  \\
\end{array}\right) \bigg\},
 \end{aligned}    }    \hspace{4,5 cm}    {}
 \end{equation}
 \bfl \emph{ where
  \begin{equation}    \label{C0 : 25}
 \begin{aligned}
& (\widetilde{w}_{2\alpha}, \widetilde{w}_{2\alpha +1}):=   \frac{e^{i\delta} }{  sin \, \varphi}
\bigg(-z_{2\alpha+1} +iz_{2\alpha}\cos\, \varphi,
z_{2\alpha } +iz_{2\alpha+1}cos\, \varphi \bigg),   \\
  and & \\
  & (\overset{\approx}{w}_{2\alpha}, \overset{\approx}{w}_{2\alpha +1}):=    \frac{e^{i\delta} }{  sin \, \varphi}
\bigg( -z_{2\alpha+1} -i\zap
cos\, \varphi, \zap -iz_{2\alpha+1} cos\, \varphi \bigg), \hspace{3 cm} \\
\end{aligned}
\end{equation}
with $sin\, \varphi\ne 0$. Moreover,
the isotropy associated to the $G-$strata in $(\ref{C0 : 24} )$ only depends on the
 determinants $\square^{1\pm 2}_{1\pm 3}$. $\square$ }\efl

 \bfl By comparing theorems B,C,D and E we conclude:\efl

\bfl $\mathbf{Theorem\,\, F}$ \emph{
The SDE orbifolds $\mathcal{O}^4(\Omega)$ and
$\mathcal{O}^4(\Theta)$ give rise to distinct  families
at the twistor and $3-$Sasakian level.} $\square$    \efl

 \newpage

 \chapter*{Aknowledgements}
\addcontentsline{toc}{chapter}{Aknowledgements}

\bfl   Before I proceed to expose the present work, I would like to express my gratitude
to all of the people who
played an important role  during the period of my Ph.D.

 \vspace{0,3 cm}

First of all I am deeply grateful to my advisor Paolo Piccinni.
 He
introduced me to the study of $3-$Sasakian and quaternion $\ka$ Geometry and I wish
to thank him for his constant support and his continuous help during
the development of this thesis.

 \vspace{0,3 cm}

Last year I met in Rome  Krzysztof Galicki  and I started  to talk with him about my project.
  This opportunity has been fundamental on the way
 of writing this thesis.
 I want to thank  him
  for his precious advices,
 his valuable comments and
all of the time he dedicated to me.

 \vspace{0,3 cm}

I wish to express my gratitude also to Maurizio Parton
for his assistance and encoragement.

 \vspace{0,3 cm}

Felt thanks go to my fiancée Gabrielle and to all of my friends, in particular Leonardo,
  who shared with
 me the most important moments of this part and, in general,  all of my life.
 \efl

\[
\left. \begin{array}{lcccl}
\mathrm{ Roma    }             &\hspace{5 cm} & &\hspace{4 cm} &    \\ 
\mathrm{settembre\,\, 2006}    &\hspace{5 cm} & &\hspace{1 cm} &   \mathrm{Luca\,\,Bisconti} \\
\end{array} \right.         \hspace{10 cm}  {}
\]

\newpage

\chapter{The Quotient Orbifolds $\mathcal{O}^4 (\Theta)$}
\pagenumbering{arabic}
\addcontentsline{toc}{chapter}{The Quotient Orbifolds $\mathcal{O}^4 (\Theta)$}

\bfl

As mentioned in the introduction, this thesis presents
 new examples of compact orbifolds
$\mathcal{O}^4 (\Theta)$ which admit a self-dual
Einstein metric of positive curvature $s>0$ with a one-dimensional
group of isometries.   We refer to the articles by S.Salamon, C. Boyer and K. Galicki,
C. Le Brun \cite{lw} for all preliminary notions on quaternionic $\ka$, $3-$Sasakian
and $4-$dimensional Einstein Manifolds and for the related quotient constructions
we are going to use throughout the present work. \newline
In this chapter we construct  these examples as quaternionic
K$\ddot{\mathrm{a}}$hler reduction of the Grassmannian
$Gr_4(\mathbb{R}^8)$ by an isometric action of
 a $3-$torus $T^3_{\Theta}\subset T^4\subset SO(8)\subset Sp(8)$
 where $\Theta$ is a weight matrix.
The starting point is to revisit a quotient construction
by Kobak-Swann
in the context of $3-$Sasakian geometry. We are going to examine
the action of
$T^3_{\Theta}\times Sp(1)\subset Sp(8)$ on the sphere $S^{31}\subset \mathbb{H}^8$
via left multiplication.
 Then we get the following $3$-Sasakian and quaternion K$\ddot{\mathrm{a}}$hler
 quotients \efl

\begin{equation}   \label{C1 : 1}
\left. \begin{array}{ccccc}
\,\, S^{31} &    \overset{Sp(1)}{\Longrightarrow} &
 \overbrace{ \frac{SO(8)}{SO(4)\times Sp(1)}}^{ \mathcal{S}^{15}\cong  }   &
 \overset{T^3_{\Theta}}{\Longrightarrow}  & \,\,   \mathcal{M}^{7}(\Theta). \\
  {}^{_{S^3}}\downarrow\,\,{}& &{}^{_{SO(3)}}\downarrow \quad{}& &{}^{_{SO(3)}}\downarrow\quad {}\\
\,\, {\mathbb{HP}}^7 &    \overset{Sp(1)}{\Longrightarrow} &
 \underbrace{ \frac{SO(8)}{SO(4)\times SO(4)}}_{\cong Gr_{4} ({\mathbb{R}}^8) }   &
 \overset{T^3_{\Theta}}{\Longrightarrow}  & \,\,    \mathcal{O}^4(\Theta) \\
\end{array}  \right.
\end{equation}

\bfl
We will prove the following statement ( Theorem A of the introduction).
\efl

\begin{teor}

\bfl
  Let $\Theta\in M_{3\times 4}$  be
 integral $3\times 4$ matrix $\Theta$ such that each of its four $3\times 3$
 minor determinants
 \begin{equation}      \label{C1 : 2}
 \Delta_{\alpha\beta\gamma} :=
\left\vert \begin{array}{ccc}
p_{\alpha} & q_{\alpha} & l_{\alpha} \\
p_{\beta} & q_{\beta} & l_{\beta} \\
p_{\gamma} & q_{\gamma} & l_{\gamma} \\
\end{array} \right\vert , \,\,\,\,\,\,
  1\leq \alpha<\beta<\gamma\leq4 ,
 \end{equation}

  does not vanish.
   In addition we suppose that the
  $\Delta_{\alpha\beta\gamma}$ are such that also the minor determinants
\begin{equation}     \label{C1 : 3}
{}^{1\pm 2}\square_{1\pm 3}^{1\pm 4} :=
\left\vert \begin{array}{ccc}
p_1 \pm p_2 & q_1\pm q_2 & l_1 \pm l_2 \\
p_1 \pm p_3 & q_1\pm q_3 & l_1 \pm l_3 \\
p_1 \pm p_4 & q_1 \pm q_4 & l_1 \pm l_4 \\
\end{array} \right\vert,
\end{equation}

do not vanish. Then
for any such $\Theta$ there exist a compact $4$-dimensional
 orbifold $\mathcal{O}^4(\Theta )$
 which admits a self-dual  Einstein metric
 of positive scalar curvature with a one-dimensional group of isometries.
 Moreover , this metric can be constructed explicitly as
 quaternionic-K$\ddot{a}$hler
reduction of the symmetric quaternion K$\ddot{a}$hler
metric on the Grassmannian $Gr_4 (\mathbb{R}^8)$ by an isometric action of the
$3-$torus $T_{\Theta}^3$ defined by the weights matrix $\Theta$. \efl
\end{teor}

\vspace{0,5 cm}

\bfl As already mentioned in the introduction,
 actions of this type can be considered on any quaternion-K$\ddot{\mathrm{a}}$hler
Grassmannian   \efl
\begin{equation}   \label{C1 : 3a}
Gr_4 (\mathbb{R}^{n+1})\cong\frac{SO(n+1)}{SO(n-3)\times SO(4)},
\end{equation}
 \bfl
whose real dimension is $4(n-3)$. Therefore to end up in dimension $4$ a $(n-4)$-dimension
torus is needed. Since the maximal torus in $SO(n+1)$ has $dim = [\frac{n+1}{2}]$,
in order to be able to consider such an action with weights we need       \efl
\begin{equation}   \label{C1 : 4}
 n-4 < \bigg[\frac{n+1}{2}\bigg].
\end{equation}
\bfl
This gives, according to whether $n+1$ is even or odd  \efl
\begin{equation}        \label{C1 : 5}
\begin{aligned}
& n-4 < \frac{n+1}{2} \quad \quad \quad \mathrm{or} \quad \quad \quad n -4 < \frac{n}{2}, \\
\mathrm{i.\, \, e.} \quad \quad \quad \quad \quad \quad   &   \\
& \left. \begin{array}{l}
           {}\,\,\, n < 9   \\
            (n+1 \,\,\mathrm{even} )
            \end{array} \right.
  \quad \quad \,\, \mathrm{or} \quad \quad \quad \quad \left. \begin{array}{l}
           {}\,\, n < 8.   \\
            (n+1\,\, \mathrm{odd} )
            \end{array} \right.  \hspace{10 cm}\\
\end{aligned}
\end{equation}

\bfl Thus the only possible  cases are   \efl

\begin{equation}  \label{C1: 6}
 \left. \begin{array}{lllll}
   Gr_4 (\mathbb{R}^6)\cong Gr_2(\mathbb{C}^4),  &
     {} \hspace{1 cm}  {} & Gr_4 (\mathbb{R}^7),
       &    {}  \hspace{1 cm} \mathrm{and} \quad \quad  {}& Gr_4 (\mathbb{R}^8). \\
       \end{array} \right.
\end{equation}

\bfl The first two cases in (\ref{C1: 6}) have been examined in \cite{mr}  and  \cite{bg}, respectively.
 The present work is essentially
devoted to the third case and to its comparison with the orbifolds construction
 in \cite{bg}.  \efl

 \vspace{0,5 cm}

\section{The  action of $Sp(1)\times T^3_{\Theta}$ on 
$N(\Theta)$ }

\vspace{0,5 cm}

\begin{flushleft}

Consider the group $G= Sp(1)\times T^3_{\Theta}$ as a
subgroup of the $3$-Sasakian isometries of the sphere $S^{31}$.
For any triple of weight vectors  $\Theta:=(\mathbf {p}, \mathbf{q}, \mathbf{l})\in
\mathbb{Z}_{+}^{4} \oplus\mathbb{Z}_{+}^{4}\oplus \mathbb{Z}_{+}^{4}\cong
 M_{3\times 4}(\mathbb{Z}_+ )$,
  define the action of the 3-torus $T^3_{\Theta}=U(1)\times U(1)\times U(1)$ on
$S^{31}$ in the following way
\begin{equation}     \label{C1: 7}
\begin{aligned}
f: & [0,2\pi )^3 \rightarrow T^4\subset SO(8), \\
&f(t,s,r):=
\left( \begin{array}{c|c|c|c}
A_1 (t,s,r) & 0 & 0 & 0\\
\hline
 0 & A_2 (t,s,r) & 0& 0\\
 \hline
 0& 0& A_3 (t,s,r)& 0\\
 \hline
 0& 0& 0& A_4 (t,s,r ) \\
\end{array} \right) \in T^4 ,\\
\end{aligned}
\end{equation}
\bfl where \efl
\begin{equation}      \label{C1: 8}
  {}  A_{\alpha} (t,s,r) =
\left( \begin{array}{cc}
cos(p_{\alpha} t + q_{\alpha} s + l_{\alpha} r ) & sin(p_{\alpha} t + q_{\alpha} s + l_{\alpha} r )   \\
-sin(p_{\alpha} t + q_{\alpha} s + l_{\alpha} r )  & cos(p_{\alpha} t + q_{\alpha} s + l_{\alpha} r ) \\
\end{array} \right),
\end{equation}

and $\alpha\in \{1,2,3,4\}$. Consider the moment map $\mu : S^{31}
\rightarrow \mathfrak{sp}(1) \otimes \mathbb{R}^3\cong
\mathbb{R}^9$ associated with the diagonal action on the left of $Sp(1)$
on $S^{31}\subset \mathbb{H}^8$
and the  moment map
$ \nu : S^{31} \rightarrow \mathfrak{u}(1)^3\otimes \mathbb{R}^3$
associated with the action of $T^3_{\Theta}\cong U(1)^3$ just described.\\
The moment map $\mu$
 is given by
\begin{equation}           \label{C1: 9}
\mu (\uno)= ( \sum_{\alpha=1}^{8} \overline{u}_{\alpha} i u_{\alpha}, \,\,
\sum_{\alpha=1}^{8} \overline{u}_{\alpha} j u_{\alpha}, \,\,
\sum_{\alpha=1}^{8} \overline{u}_{\alpha} k u_{\alpha} )\in \mathfrak{sp}(1) \otimes \mathbb{R}^3,
\end{equation}

\bfl where $\uno\in S^{31}\subset \mathbb{H}^8$, and    \efl

\begin{equation}             \label{C1: 10}
 \nu (\uno)= \left( \begin{array}{c}
 \sum_{\alpha = 1}^{4}  p_{\alpha}(\overline{u}_{2\alpha -1 } u_{2\alpha} -  \overline{u}_{2\alpha  } u_{2\alpha -1} )  \\
 \sum_{\alpha = 1}^{4}  q_{\alpha}(\overline{u}_{2\alpha -1 } u_{2\alpha} -  \overline{u}_{2\alpha  } u_{2\alpha -1} )  \\
 \sum_{\alpha = 1}^{4}  l_{\alpha}(\overline{u}_{2\alpha -1 } u_{2\alpha} -  \overline{u}_{2\alpha  } u_{2\alpha -1} )  \\
\end{array} \right)\in \mathfrak{u}(1)^3\otimes \mathbb{R}^3.
\end{equation}

\end{flushleft}

\vspace{0,5 cm}

\begin{flushleft}
$\mathbf{Remark\,\, 1.1}$
The submanifold $\mu^{-1}(0)$  can be
identified with the Stiefel manifold
of oriented orthonormal $4$-frames in $\mathbb{R}^8$. Accordingly,
we describe the elements of $N(\Theta)$
$:=\mu^{-1}(0)$
$\cap
\nu^{-1}(0)$ as $4\times 8$ real matrices
$ \uno :=(u_1 , u_2 ,\dots, u_7 , u_8 ) $
whose columns $u_i$  can be thought as the coefficients of a quaternion
respect to the base $\{
1,i,j,k\}$. Of course the matrix $\boldsymbol{u}$
 has $rank = 4$. $\square$ \end{flushleft}

\begin{defin}
 Any pair $\uua\in \mathbb{H}\times \mathbb{H}$, $\alpha\in\{1,2,3\}$ of quaternionic coordinates of
  $\uno\in S^{31}\subset \mathbb{H}^8$ will be breafly called a \rm{quaternionic pair}.
\end{defin}

\vspace{0,5 cm}

 \begin{flushleft}

\begin{lem}
 Suppose that all the minor determinants
\begin{equation}  \label{C1: 10a}
\Delta_{\alpha\beta\gamma}=
\left\vert \begin{array}{ccc}
p_{\alpha} & q_{\alpha} & l_{\alpha} \\
p_{\beta} & q_{\beta} & l_{\beta} \\
p_{\gamma} & q_{\gamma} & l_{\gamma} \\
\end{array} \right\vert , \quad \,\,\,
  1\leq \alpha<\beta<\gamma\leq4 ,
\end{equation}

do not vanish. Then there is no element $\uno$ of $N(\Theta)$   with a null quaternionic pair.
\end{lem}
\end{flushleft}

\begin{flushleft}

$\mathbf{Proof}. $
 It is sufficient to assume  $(u_7 ,u_8 )$ as a null quaternionic pair on some point of $N(\Theta)$.
Let $\Theta$ be the integral $3\times 4$ matrix  of the weights
$ (\mathbf {p}, \mathbf{q}, \mathbf{l} )\in
\mathbb{Z}_{+}^{4} \oplus\mathbb{Z}_{+}^{4}\oplus \mathbb{Z}_{+}^{4}$.
Let $x_{\alpha}:= \overline{u}_{2\alpha -1 } u_{2\alpha} -  \overline{u}_{2\alpha  }
u_{2\alpha -1}, \,\, \alpha\in\{1,2,3,4\} $,
and rewrite $\nu (\uno )$ in the following way
\begin{equation}       \label{C1: 11}
\left( \begin{array}{ccc|c}
  &            &  & p_4 \\
  & \mathbb{A} &  & q_4 \\
  &            &  & l_4 \\
\end{array} \right)  :=
\underbrace{\left( \begin{array}{ccc|c}
p_1  & p_2 & p_3 & p_4 \\
q_1  & q_2 & q_3 & q_4 \\
l_1  & l_2 & l_3 & l_4 \\
\end{array} \right)}_{=\Theta}\,\, \Rightarrow\,\,
\nu (\uno)= \left( \begin{array}{ccc|c}
  &            &  & p_4\\
  & \mathbb{A} &  & q_4 \\
  &            &  & l_4 \\
\end{array} \right)
\left( \begin{array}{c}
x_1 \\
x_2 \\
x_3 \\
x_4 \\
\end{array}   \right).
\end{equation}

 Thus
\begin{equation}        \label{C1: 12}
\nu (\uno ) = \left( \begin{array}{ccc|c}
  &            &  & 0 \\
  & \mathbb{A} &  & 0 \\
  &            &  & 0 \\
\end{array} \right)
\left( \begin{array}{c}
x_1 \\
x_2 \\
x_3 \\
0   \\
\end{array}\right) +
  \left( \begin{array}{c}
 p_4 x_4 \\
  q_4 x_4\\
 l_4 x_4\\
\end{array} \right).
\end{equation}

\bfl Since $det\, \mathbb{A}=det\, {}^T \mathbb{A}= \Delta_{123} $, the equation $\nu (\uno )=0$
can be written as \efl
\begin{equation}     \label{C1: 13}
\left( \begin{array}{c}
x_1 \\
x_2 \\
x_3 \\
\end{array} \right) =-\mathbb{A}^{-1}
 \left( \begin{array}{c}
 p_4 x_4 \\
  q_4 x_4\\
 l_4 x_4\\
\end{array} \right):=
 \left( \begin{array}{c}
 c_1  x_4 \\
 c_2  x_4\\
 c_3  x_4\\
\end{array} \right),   \,\,\,\,\,\, c_{\alpha}\in \mathbb{Q}.
\end{equation}

\bfl Hence the intersection
$N(\Theta)\cap \{ u_7= u_8 = 0\} $ is described by the following equations  \efl
\begin{equation}        \label{C1: 14}
\left\{ \begin{array}{l}
\overline{u}_1 u_2 -\overline{u}_2 u_1 =c_1 (\overline{u}_7 u_8 -\overline{u}_8 u_7)=0,\\
\overline{u}_3 u_4 -\overline{u}_4 u_3 =c_2 (\overline{u}_7 u_8 -\overline{u}_8 u_7)=0,\\
\overline{u}_5 u_6 -\overline{u}_6 u_5 =c_3 (\overline{u}_7 u_8 -\overline{u}_8 u_7)=0,\\
\sum_{\alpha=1}^6 \overline{u}_{\alpha}\sigma u_{\alpha}=0, \\
\end{array} \right.\,\,\,\, \Rightarrow\,\,\,\,
\left\{ \begin{array}{l}
\overline{u}_1 u_2 =\overline{u}_2 u_1 ,\\
\overline{u}_3 u_4 =\overline{u}_3 u_4 ,\\
\overline{u}_5 u_6 =\overline{u}_6 u_5 ,\\
\sum_{\alpha=1}^6 \overline{u}_{\alpha}\sigma u_{\alpha}=0,
\end{array} \right.
\end{equation}

where $\mu (\uno )=(\sum_{\alpha=1}^6 \overline{u}_{\alpha}\sigma u_{\alpha})|_{\{
\sigma = i,j,k\}}$,
is the moment map of the $Sp(1)$ action. From equations (\ref{C1: 14}) we get
$Im (\barun\um )=0,\, \alpha\in\{1,2,3\}$, and hence $\barun\um\in\mathbb{R}$.
Moreover, from equations (\ref{C1: 14}) we see that we can always assume, up to
a scale, that $u_{\alpha} ={}^T(u_{\alpha}^0, u_{\alpha}^1, u_{\alpha}^2, u_{\alpha}^3)$
 belongs to $Sp(1)$.
The maps
\begin{equation}      \label{C1: 15}
\begin{aligned}
\varphi_{\alpha} :  &Sp(1)\times Sp(1)\longrightarrow Sp(1), \quad\,\,\,\, \alpha\in\{1,2,3\}, \\
&\quad (\un , \um )\longmapsto \barun\um   \in \mathbb{R}\\
\end{aligned}
\end{equation}

are well defined by the previous facts.
 Thus we have that the real numbers
$\barun\um$ belong to the sphere $S^3\cong Sp(1)$ so that the only
possibility  to have
$\barun\um =\pm 1 $ is $\un =\pm \um ,\,\,\, \alpha\in\{1,2,3\}.$
It follows that a real $4\times 8 $  matrix $\uno\in N(\Theta)$ cannot
satisfy at the same
time all the equations (\ref{C1: 14}), because the first three equations
give us that the columns $\un={}^T(\un^0,\un^1,\un^2\un^3)$ and $\um=
{}^T(\um^0,\um^1,\um^2\um^3)$
of each quaternionic pair are proportional to each other:
\begin{equation}
{}^T(\un^0,\un^1,\un^2\un^3)=d_{\alpha}\,
{}^T(\um^0,\um^1,\um^2\um^3),\,\,\,\, d_{\alpha}
\in\mathbb{R}\backslash\{0\},\,\forall\, \alpha\in\{1,2,3\}.
\end{equation}
 Then the matrix
$\uno$ has at most $rank = 3$ and it does not describe any $4$-frame
 in $\mathbb{R}^8$. \\
The consequence is that the  $N(\Theta)\cap \{ u_7= u_8 = 0\} $
is empty. $\square $
\end{flushleft}

\vspace{1 cm}

\begin{prop}
\begin{flushleft}
The action of $G=Sp(1)\times T^3_{\Theta}$ on $N(\Theta)$
 is locally free if and only if all
the followings determinants
\end{flushleft}
\begin{equation}  \label{C1: 15a}
\begin{aligned}
& {}\quad \quad \quad \quad\Delta_{\alpha\beta\gamma}=
\left\vert \begin{array}{ccc}
p_{\alpha} & q_{\alpha} & l_{\alpha} \\
p_{\beta} & q_{\beta} & l_{\beta} \\
p_{\gamma} & q_{\gamma} & l_{\gamma} \\
\end{array} \right\vert , \quad \,\,\,
  1\leq \alpha<\beta<\gamma\leq4 ,  \\
 and\quad \quad \hspace{1 cm} {}& \\
& {}\quad \quad \quad \sq  :=
\left\vert \begin{array}{ccc}
p_1 \pm p_2 & q_1\pm q_2 & l_1 \pm l_2 \\
p_1 \pm p_3 & q_1\pm q_3 & l_1 \pm l_3 \\
p_1 \pm p_4 & q_1 \pm q_4 & l_1 \pm l_4 \\
\end{array} \right|,   \\
do\,\, not\,\, vanish. & {} \hspace{15 cm} {} \\
\end{aligned}
\end{equation}
\end{prop}

\bfl

$\mathbf{Proof}$. By lemma $1.1$  the condition $\Delta_{\alpha\beta\gamma}\ne 0$ insures
that $N(\Theta)$ has no points with a null quaternionic pair.  Then the  fixed point
equations can be written as  \efl
\begin{equation}     \label{C1: 16}
A_{\alpha}(p_{\alpha} t + q_{\alpha} s + l_{\alpha} r ) \left( \begin{array}{c}
u_{2{\alpha}-1}  \\
u_{2{\alpha}}  \\
\end{array} \right) =
\left( \begin{array}{cc}
a_{\alpha} & b_{\alpha} \\
-b_{\alpha} & a_{\alpha} \\
\end{array} \right)  \left( \begin{array}{c}
u_{2{\alpha}-1}  \\
u_{2{\alpha}}  \\
\end{array} \right)  =
\lambda \left( \begin{array}{c}
u_{2{\alpha}-1}  \\
u_{2{\alpha}}  \\
\end{array} \right),
\end{equation}

\bfl ${\alpha}\in\{1,2,3,4\}$, where $a_{\alpha} = \cost_{\alpha}, \,\,
  b_{\alpha} = \sit_{\alpha},\,\, \theta_{\alpha}= p_{\alpha} t + q_{\alpha} s + l_{\alpha} r$,
   $ t, s, r\in [0 ,2\pi )$
and $\lambda\in Sp(1)$.
For each ${\alpha}=1,2,3,4$ we get   \efl
\begin{equation}      \label{C1: 17}
\left. \begin{array}{l}
a_{\alpha} \un + b_{\alpha}\um = \lambda \un ,\\
-b_{\alpha} \un + a_{\alpha} \um = \lambda \um ,\\
\end{array}  \right.
\end{equation}

\bfl
and by multipling on the right the first equation by
 ${\overline{u}}_{2\alpha-1}$  and the
second by $\barum$ we have      \efl
\begin{equation}      \label{C1: 18}
\left. \begin{array}{l}
a_{\alpha} |\un |^2 + b_{\alpha}\um \barun= \lambda |\un |^2 ,\\
-b_{\alpha} \un \barum+ a_{\alpha}  |\um |^2 = \lambda |\um |^2 .\\
\end{array} \right.
\end{equation}

 \bfl Then \efl
\begin{equation}     \label{C1: 19}
a_i ( |\un |^2 + |\um |^2 ) +  b_i ( \um \barun -  \un \barum )
= \lambda  (|\un |^2 + |\um |^2 ),\,\, \alpha\in\{1,2,3,4\}  .
\end{equation}
 \bfl
Under our assumptions there are no null quaternionic  pairs,
and the term multiplying $\lambda$ on the right-hand side of this equation never
vanishes.
 This gives    \efl
\begin{equation}     \label{C1: 20}
Re(\lambda ) =a_{\alpha},\,\,\,\,\,\, Im(\lambda ) =
 b_{\alpha}\frac{ ( \um \barun -  \un \barum ) }{(|\un |^2 + |\um |^2 )},
\end{equation}
\bfl
for each $\alpha\in\{1,2,3,4\}$.
The first of these equation gives    \efl
 \begin{equation}     \label{C1: 21}
  a_1 =a_2 = a_3= a_4 ,
 \end{equation}
 \bfl
and combining this  with $a_{\alpha}^{2} +
b_{\alpha}^{2}=1$ we get                  \efl
\begin{equation}        \label{C1: 21a}
b_1=\pm b_2 =\pm b_3=\pm b_4 .
\end{equation}
\bfl
These  conditions  give   \efl
\begin{equation}      \label{C1: 22}
\left\{ \begin{array}{c}
e^{i(p_1 t + q_1 s + l_1 r)}= e^{\pm i( p_2 t + q_2 +l_2 r)},  \\
e^{i(p_1 t + q_1 s + l_1 r)}= e^{\pm i( p_3 t + q_3 +l_3 r)},  \\
 e^{i(p_1 t + q_1 s + l_1 r )}= e^{\pm i( p_4 t+ q_4 s+l_4 r)}. \\
\end{array} \right.
\end{equation}

\bfl Let us rewrite these equations in the following way  \efl
\begin{equation}         \label{C1: 23}
\left\{ \begin{array}{l}
(p_1 \pm p_2 )t + ( q_1 \pm q_2 )s + (l_1 \pm l_2 )r = 2h_{12}^{\pm} \pi,  \\
(p_1 \pm p_3 )t + ( q_1 \pm q_3 )s + (l_1 \pm l_3 )r = 2h_{13}^{\pm}  \pi,  \\
(p_1 \pm p_4 )t + ( q_1 \pm q_4 )s + (l_1 \pm l_4 )r = 2h_{23}^{\pm} \pi , \\
\end{array} \right.
\end{equation}

\bfl
where $ h_{\alpha\beta}\in \mathbb{Z}$.
The equations (\ref{C1: 23}) describe eight systems in the $(t,s,r)$ variables.
We want all eight of them to have  at most discrete solutions.
This requires that all
the eight determinants      \efl
\begin{equation}     \label{C1: 23a}
{}^{1\pm 2}\square_{1\pm 4}^{1\pm 3} =
\left\vert \begin{array}{ccc}
p_1 \pm p_2 & q_1\pm q_2 & l_1 \pm l_2 \\
p_1 \pm p_3 & q_1\pm q_3 & l_1 \pm l_3 \\
p_1 \pm p_4 & q_1 \pm q_4 & l_1 \pm l_4 \\
\end{array} \right\vert ,
\end{equation}
\bfl
do not vanish. Then our assumptions insure that the fixed point equations
(\ref{C1: 16}) admit
at most a finite number of solutions.  Moreover, the equations (\ref{C1: 21})
 and (\ref{C1: 21a}) give  \efl
\begin{equation}     \label{C1: 24}
\frac{  \um \barun -  \un \barum  }{(|\un |^2 + |\um |^2 )}=
\pm \frac{u_2\overline{u_1}- u_1\overline{u}_2}{ |u_1|^2 + |u_2|^2}, \,\,\,
\alpha\in\{1,2,3,4 \}.
\end{equation}
\bfl Now, under the assumption $\sq\ne 0$, it holds that $\sit_1$ assume only
a finite number of values,
hence the equation $Im(\lambda )= \sit_1
\frac{u_2\overline{u_1}- u_1\overline{2}}{ |u_1|^2 + |u_2|^2}$  can be rewritten
as system of three equations
in three non known parameters. Then the fixed point equations (\ref{C1: 17})
admit only discrete solutions.  $\square$ \efl

\vspace{0,5 cm}

\bfl $\mathbf{Example\,\, 1.1.}$ Assume   $m:=\sqp\ne 0$,  we get \efl

\begin{equation}     \label{C1: 25}
\begin{aligned}
& {}\quad \quad t=\frac{2\pi}{m}
 \left\vert \begin{array}{ccc}
h_{12}^{+} & q_1+ q_2 & l_1 + l_2 \\
h_{13 }^{+}& q_1+ q_3 & l_1 +l_3 \\
h_{14}^{+} & q_1+ q_4 & l_1+ l_4\\
\end{array} \right\vert =
\frac{2\pi}{m }\underbrace{( h_{12}^{+}  {}^{ql}\Delta_{3+4}^{2+ 3} -
 h_{23}^{+} {}^{ql}\Delta_{3+ 4}^{1+ 2}  +  h_{34}^{+}  {}^{ql}\Delta_{2+ 3}^{1+ 2} ) }_{\in \mathbb{Z}}, \\
\mathrm{where} &\\
&\quad \quad {}^{ql}\Delta_{\alpha\pm \beta}^{\gamma\pm \delta} :=
\left\vert \begin{array}{cc}
 q_{\alpha}\pm q_{\beta} & l_{\alpha}\pm l_{\beta} \\
 q_{\gamma}\pm q_{\delta} & l_{\gamma} \pm l_{\delta} \\
\end{array} \right\vert  . \hspace{10 cm} \\
\end{aligned}
\end{equation}

     \bfl
Since  $ 0\leq\frac{t}{\pi}=\frac{ ( h_{12}^{+}  {}^{ql}\Delta_{3+4}^{2+ 3} -
 h_{23}^{+} {}^{ql}\Delta_{3+ 4}^{1+ 2}  +  h_{34}^{+}  {}^{ql}\Delta_{2+ 3}^{1+ 2}}m{}<1$
all the possible solutions $t$ of the system (\ref{C1: 23}) belong to $H:=\{0, \frac{\pi}{|m |},\ldots ,
\frac{\pi (|m |-1)}{|m |}\}$.
Similarly also $s,r$ belong to $H$. Then each  solution $(t,s,r)$ for the equation (\ref{C1: 23})
is a triple whose elements are in $H$. We can use a similar argument
 for the other non null determinants
$\sq$ and also in these cases we get discrete solutions.
 Moreover, by using a short argument linked to the $gcd$
of the minor determinants ${}^{ql}\Delta_{\beta}^{\alpha}$ we can prove that
the cardinality of the subgroup of solutions is exactly equal to $m=\sqp$.  $\square$\efl

\vspace{0,5 cm}

\begin{prop}
\bfl There is no weight matrix $\Theta$ such that the action of
 $G=Sp(1)\times T^3_{\Theta}$ on $N(\Theta)$ is free.  \efl
\end{prop}

\bfl
$\mathbf{Proof}.$ From the previous proof we see that there is a unique solution
 for the fixed point equations (\ref{C1: 23})
if and only if $\sqpm = \pm 1.$  We have the identities    \efl

\begin{equation}     \label{C1: 26}
\left\{ \begin{array}{l}
{}^{1+ 2}\square_{1+ 4}^{1+ 3}= \Delta_{123}- \Delta_{124}+ \Delta_{134}+\Delta_{234}  ,\\
{}^{1+ 2}\square_{1+ 4}^{1- 3}= -\Delta_{123}- \Delta_{124}- \Delta_{134}-\Delta_{234} , \\
{}^{1+ 2}\square_{1- 4}^{1+ 3}= \Delta_{123}+ \Delta_{124}-\Delta_{134}-\Delta_{234}  ,\\
{}^{1+ 2}\square_{1- 4}^{1- 3}= -\Delta_{123}+ \Delta_{124}+\Delta_{134}+\Delta_{234}  ,\\
{}^{1- 2}\square_{1+ 4}^{1+ 3}=- \Delta_{123}+ \Delta_{124}+ \Delta_{134}-\Delta_{234} , \\
{}^{1- 2}\square_{1+ 4}^{1- 3}= \Delta_{123}+ \Delta_{124}- \Delta_{134}+\Delta_{234}  ,\\
{}^{1- 2}\square_{1- 4}^{1+ 3}= -\Delta_{123}- \Delta_{124}- \Delta_{134}+\Delta_{234} , \\
{}^{1- 2}\square_{1- 4}^{1- 3}= \Delta_{123}- \Delta_{124}+ \Delta_{134}-\Delta_{234} .  \\
\end{array} \right.
\end{equation}

\bfl
that can be solved with respect to the $\Delta_{\alpha\beta\gamma}.$  In fact note that the system  \efl
\begin{equation}     \label{C1: 27}
\left\{ \begin{array}{l}
X :={}^{1+ 2}\square_{1+ 4}^{1+ 3}= \Delta_{123}- \Delta_{124}+ \Delta_{134}+\Delta_{234},  \\
Y:={}^{1+ 2}\square_{1+ 4}^{1- 3}= -\Delta_{123}- \Delta_{124}- \Delta_{134}-\Delta_{234} , \\
Z:={}^{1- 2}\square_{1- 4}^{1+ 3}= -\Delta_{123}- \Delta_{124}-\Delta_{134}+\Delta_{234}  ,\\
W:= {}^{1+ 2}\square_{1- 4}^{1- 3}= -\Delta_{123}+ \Delta_{124}+\Delta_{134}+\Delta_{234},  \\
\end{array} \right.
\end{equation}
\bfl has  determinant \efl
\begin{equation}     \label{C1: 28}
\begin{aligned}
&  \left| \begin{array}{cccc}
1 & -1 & 1 & 1 \\
-1 & -1 & -1 & -1 \\
-1 & -1 & -1 & 1 \\
-1 & 1 & 1 & 1 \\
\end{array} \right| =8,     \\
\mathrm{and\,\, its\,\, solutions\,\, are}  \hspace{1 cm} \quad&  \\
 &  {} \quad \left\{ \begin{array}{l}
 \Delta_{123} = -\frac{Y +W}{2},  \\
 \Delta_{124} = -\frac{X + Y}{2}, \\
 \Delta_{134}= \frac{X +Y -Z + W}{2}, \\
 \Delta_{234} =\frac{Z- Y}{2} ,
 \end{array}  \right.     \hspace{15 cm}   {}
\end{aligned}
\end{equation}
\bfl that are also solutions of the remaining identities in (\ref{C1: 26})    \efl
\begin{equation}          \label{C1: 29}
\left\{ \begin{array}{l}
\pm 1 = {}^{1+ 2}\square_{1- 4}^{1+ 3}=  -(X + Y + W) ,\\
\pm 1 = {}^{1- 2}\square_{1+ 4}^{1- 3}=  Z - X -2Y -W , \\
\pm 1 = {}^{1- 2}\square_{1+ 4}^{1+ 3}= Y  + Z + W, \\
\pm 1 = {}^{1- 2}\square_{1- 4}^{1- 3}= X + Y - Z. \\
\end{array}\right.
\end{equation}

\bfl  But since $\Delta_{\alpha\beta\gamma}\in \mathbb{Z}\backslash \{0\}$ and $X=\pm 1,\,
Y=\pm 1$, $Z=\pm 1$ and $W=\pm 1$  we would have
 \begin{equation}
X = Y =W =-Z,
 \end{equation}

and hence $(X, Y, Z, W )$ is either $(1,1,-1,1)$ or $(-1, -1, 1, -1)$
and $(\Delta_{123} , \Delta_{124}, \Delta_{134}, $
$\Delta_{234})$ is either
 $(1,1, -2, 1)$ or $(-1, -1, 2, -1)$.
It is now easy to check that neither of the possible two values of $(X,Y,Z,W)$
is a solution
for the system of equations above. Then the action of $G$ on $N(\Theta)$ is never free. $\square$\\  \efl

 \bfl Using all of these facts  we  get  \efl

\begin{teor} \bfl
The action of $Sp(1)\times T^3_{\Theta}$ on $N(\Theta)=\nu^{-1}(0)\cap \mu^{-1}(0)$ is locally free
if and only if the following conditions hold             \efl
\begin{flushleft}
$1)$ all the determinants $\Delta_{\alpha\beta\gamma}\ne 0,\,\,\,\,\forall\,
(\alpha,\beta,\gamma)$ admissible triple
of indices. \\

 \vspace{0,5 cm}
 
$2)$  $\Delta_{\alpha\gamma\delta}$ are such that
 all the minor determinants ${}^{1\pm 2}\square_{1\pm 3}^{1\pm 4} $
   in equations $(\ref{C1: 26})$ do not vanish. \\

 \vspace{0,5 cm}

In such a case the quotient
\begin{equation}       \label{C1: 29a}
\mathcal{M}^7(\Theta)= \frac{N(\Theta)}{Sp(1)\times T^3_{\Theta}},
\end{equation}
 is a compact $7-$dimensional $3-$Sasakian orbifold.  Moreover for any
 integral $3\times 4$ matrix $\Theta$  which satisfies the assumptions
 $1)$ and $2)$
 there exists a $4$-dimensional orbifold $\mathcal{O}^4(\Theta )$ which
 admits a self-dual  Einstein metric
 of positive scalar curvature with a $1-$dimensional group of isometries.
\end{flushleft}
\end{teor}

\bfl The orbifold family $\mathcal{M}^7(\Theta)$ cannot be toric,
this fact can be seen
by having a careful analysis of the associated foliations. Note that,
 $N(\Theta)$ is a compact submanifold of the sphere $S^{31}\subset
  \mathbb{H}^{8}$ given by the zero locus of two
quadratic functions( the moment maps $\mu$ and $\nu$).
Then, by using a similar argument
like that one in \cite{bg} pag 97,  we can generalize the analysis
of \cite{mk},  and conclude that all the isometries
of $N(\Theta)$ come from the restriction of the isometries of $S^{31}$
(or the Euclidean space
$\mathbb{H}^{8}\cong \mathbb{R}^{32}$ in which this sphere is embedded).
Then, we get that the group of isometries associated
to  $\mathcal{O}^4(\Theta )$ is $1-$dimensional. \efl


\begin{flushleft}
 $\mathbf{Example\,\, 1.2.}$
The  $3\times 4$ matrix
 \begin{equation}    \label{C1: 30}
\left( \begin{array}{cccc}
1 & 0 & 1 & 1 \\
0 & 1 & 1 & 1  \\
1 & 1 & 0 & 1  \\
\end{array} \right),
\end{equation}

\end{flushleft}

\bfl satisfies all the requirements of theorem $1.2$. $\square$ \efl

\chapter{The Singular Locus }

\addcontentsline{toc}{chapter}{The Singular Locus}

\bfl Recall that the twistor space  $\mathcal{Z}^6(\Theta)$
 is the $ S^2$ bundle over the riemannian orbifold $\mathcal{O}^4(\Theta)$
 that parametrizes
 its local almost complex structures. Consider the group $\tilde{G}$
 $:=G\times U(1)= $
 $T^3_{\Theta}\times Sp(1)\times U(1)$,  which is a subgroup, up to the central $\mathbb{Z}/2$,
  of $Sp(8)\cdot Sp(1)\subset$
  $ SO(32)=Isom(S^{31} )$.
  Then, we can define the action of $\tilde{G}$ on $N(\Theta)=   $
 $\mu^{-1}(0) \cap \nu^{-1}(0)$ ( where
$\mu$ and $\nu$ are the moment maps of $Sp(1)$ and $T^3_{\Theta}$, respectively) as follows \efl
 \begin{equation}    \label{C2 : 1}
 \begin{aligned}
 &  \Phi:\, T^3_{\Theta}\times Sp(1)\times U(1)\times
  N(\Theta) \longrightarrow N(\Theta),\\
 & \quad \quad\quad \quad\big( (A(\Theta), \lambda, \rho);(\underline{z},
  \underline{w} ) \big)\longmapsto
 \Phi \big( ( A(\Theta), \lambda, \rho)\big)\big((\underline{z},
  \underline{w} )\big), \\
\mathrm{where} \quad \quad& \\
 & \quad \quad \quad \quad \Phi \big( ( A(\Theta), \lambda, \rho)
 \big)\big((\underline{z}, \underline{w} ) \big)
 :=   A(\Theta) \lambda  \left( \begin{array}{c}
\underline{z}\\
\underline{w} \\
\end{array}\right)
 \rho,
  \hspace{5 cm} {}
 \end{aligned}
 \end{equation}
 \bfl and we have identified $\mathbb{C}^8\times \mathbb{C}^8\cong \mathbb{H}^8$ by
  $u_{\alpha} = z_{\alpha} + jw_{\alpha}\in \mathbb{C}\times \mathbb{C}$
  for each ${\alpha}\in \{1,\ldots,$
  $ 8\}$. Thus
  $\zw:=\boldsymbol{u}=
 (u_1,u_2,\ldots,u_7,u_8)\in \mathbb{H}^8$.  In the following we will
  use 
   both  the notations $\uno$ and $\zw$.
  Here,
  the twistor space $\mathcal{Z}^6(\Theta)$
  will be thought as the leaf
space of the $\tilde{G}-$action on the manifold $N(\Theta)$.
There is a natural
stratification of the twistor space $\mathcal{Z}^6(\Theta)$ defined
by the action of $\tilde{G}$, and the twistor space
 $\mathcal{Z}^6(\Theta)$ is 
  the union
of its stata.
\efl

\begin{defin}
\bfl We say that two points $\zw$, $(\underline{z}_1, \underline{w}_1)\in
\mathcal{Z}^6(\Theta)$ are in the same stratum $S$ if their corresponding
 isotropy subgroups $\tilde{G}_{\zw}$
and $\tilde{G}_{(\underline{z}_1, \underline{z}_2)}$   are conjugate.
 \efl
\end{defin}

 \bfl
  The strata $S$, with respect to the  action of $\tilde{G}$ on
   $\mathcal{Z}^6(\Theta)$,
 can be descibed by studying the action of $\tilde{G}$ on $N(\Theta)$
 and the respective stratification. We denote with $\widetilde{S}$
  the strata
 yielded by the $\tilde{G}-$action $\Phi$
 on $N(\Theta)$. The smooth locus $\mathcal{Z}_0^6(\Theta)$
 is the  dense open subset of points in
 $\mathcal{Z}^6(\Theta)$ whose isotropy subgroup is
 the identity.
 Instead, $\Sigma(\Theta)$ denotes the singolar locus of
 $\mathcal{Z}^6(\Theta)$ and we have that
 $\mathcal{Z}^6(\Theta)= \mathcal{Z}_0^6(\Theta) \cup\Sigma(\Theta)$.
In the case of a $4-$dimensional "$good\,\, orbifold$" $\mathcal{O}^4$ , that is
 when the $3-$Sasakian
bundle $\mathcal{M}^7$ is smooth, the singularities of  $\mathcal{O}^4$ come
 from those of its twistor space $\mathcal{Z}^6$ \cite{tw}. However,
 when the $3-$Sasakian bundle $\mathcal{M}^7$ is not smooth, this fact does not hold,
 and in order to study the singularities of $\mathcal{O}^4$  a different approach is needed.
 As we have seen in chapter $1$, the $3-$Sasakian bundle $\mathcal{M}^7(\Theta)$
 necessarely has orbifold singularities which pass to the $SO(3)-$quotient space
 $\mathcal{O}^4(\Theta)$. Then, with the purpose of describing the singular locus
 of the quaternion
 $\ka$ orbifold   $\mathcal{O}^4(\Theta)$, we will study both
 the singular loci of   $\mathcal{M}^7(\Theta)$ and
$\mathcal{Z}^6(\Theta)$. \efl

\bfl In order
 to understand the  $3$-Sasakian and quaternion $\ka$ reduction we are dealing with,
   look at the diagram    \efl

\begin{equation}     \label{C2 : 2}
\dgARROWLENGTH=0.3\dgARROWLENGTH
   \begin{diagram}
     \node[2]{ \frac{SO(8)}{SO(4)\times  Sp(1)}  } \arrow[2]{e,t}{T^3_{\Theta}} \arrow[2]{s,l}{{   {}^ {}_{SO(3)}}}
       \node[2]{\mathcal{M}^7(\Theta)}  \arrow[2]{s,l}{ {}^{{}_{SO(3)}} } \arrow{se,t}{  U(1) } \\
                                                                          \node[5]{\mathcal{Z}^6(\Theta)} \arrow{sw,b}{S^2}  \\
     \node[2]{\grm} \arrow[2]{e,t}{T^3_{\Theta}} \node[2]{\mathcal{O}^4(\Theta) }
   \end{diagram}
\end{equation}

\bfl where $ \mu^{-1}(0)/Sp(1)\cong SO(8)/{SO(4)\times  Sp(1)}$.
 In this chapter we prove the following statements mentioned in the introduction
as Theorems B and C.\efl

\begin{teor}
\bfl Let $\mathcal{Z}^6(\Theta)$ be the twistor space of the SDE orbifold
 $\mathcal{O}^4(\Theta)$. Then, depending on the minor determinants $\Delta_{\alpha\beta\gamma}$
 and $\sq$  the singular locus $\Sigma(\Theta)$ consists of \efl

 \bfl 1) one isolated $2-$sphere $S^2$, whose isotropy
 only depends on one of the determinants $\sq$;

 \vspace{0,5 cm}

 2) eleven disjoint $2-$spheres $S^2$. These are obtained from
  the following strata
  $($ with respect to the action
 of $\tilde{G}$ on $\mathbb{H}^8$ $)$ \efl
\begin{equation}   \label{C2 : 3}
\begin{aligned}
& a)\,\,\, \quad \overset{=}{S}{}^{123}_4:=\tilde{G}\cdot \Bigg\{
\left( \begin{array}{cc|cc|cc|cc}
z_1 & z_2  & z_3 & z_4   &    z_5 & z_6        & 0    &  0   \\
0   &    0 &  0  & 0     &   0   &  0          & w_7  &  w_8\\
\end{array}\right) \Bigg\},   \\
&  b) \,\,\, \quad \overset{=}{S}{}^{124}_3:=\tilde{G}\cdot \Bigg\{
\left( \begin{array}{cc|cc|cc|cc}
z_1 & z_1  & z_3 & z_4       &   0   &    0          & z_7   & z_8   \\
0   &    0  & 0 &    0     &   w_5 &   w_6        & 0 & 0\\
\end{array}\right) \Bigg\},   \hspace{10 cm} \\
&  c)\,\,\, \quad  \overset{=}{S}{}^{134}_2:=\tilde{G}\cdot \Bigg\{
\left( \begin{array}{cc|cc|cc|cc}
z_1 & z_2  & 0    &  0           &    z_5 & z_6             & z_7   & z_8  \\
0   &    0 & w_3  & w_4   &   0  &  0   &  0    &  0          \\
\end{array}\right) \Bigg\}, \hspace{10 cm }  \\
&  d) \,\,\, \quad \overset{=}{S}{}^{234}_1:=\tilde{G}\cdot \Bigg\{
\left( \begin{array}{cc|cc|cc|cc}
0   &    0      & z_3 & z_4 &    z_5 & z_6             & z_7   & z_8  \\
w_1  & w_2& 0   &    0     &      0  &    0              &  0    &  0\   \
\end{array}\right) \Bigg\}.  \hspace{20 cm}  \\
& e) \,\,\, \quad \overset{=}{S}{}^{12}_{34}:=\tilde{G}\cdot \Bigg\{
\left( \begin{array}{cc|cc|cc|cc}
z_1 & z_2  & z_3 & z_4   &    0  &  0        & 0    &  0   \\
0   &    0 &  0  & 0     &   w_5   & w_6          & w_7  &  w_8\\
\end{array}\right) \Bigg\},   \\
&  f) \,\,\, \quad \overset{=}{S}{}^{14}_{23}:=\tilde{G}\cdot \Bigg\{
\left( \begin{array}{cc|cc|cc|cc}
z_1 & z_1  & 0 & 0       &   0   &    0          & z_7   & z_8   \\
0   &    0  & w_3 &    w_4     &   w_5 &   w_6        & 0 & 0\\
\end{array}\right) \Bigg\},   \hspace{10 cm} \\
&  g)\,\,\, \quad \overset{=}{S}{}^{12}_{34}:=\tilde{G}\cdot \Bigg\{
\left( \begin{array}{cc|cc|cc|cc}
z_1 & z_2  & 0    &  0           &    z_5 & z_6             & 0  & 0  \\
0   &    0 & w_3  & w_4    &   0  &  0   &  w_7   &  w_8          \\
\end{array}\right) \Bigg\}. \hspace{10 cm }  \\
\end{aligned}
 \end{equation}
\bfl Each of these strata 
intersects $N(\Theta)$ and
$\overset{=}{S}{}^{\alpha\beta\gamma}_{\delta}\cap N(\Theta)$ is
formed by two connected components which gives two
 $\widetilde{S}$ strata for  the $\tilde{G}-$action on $N(\Theta)$.
  Instead, each of the
 $\overset{=}{S}{}^{\alpha\beta}_{\gamma\delta}\cap N(\Theta)$ is
 connected and provides
 exactly one strata $\widetilde{S}$ of $N(\Theta)$.
 Moreover, we have that
$\widetilde{S}/\tilde{G}\cong S^2$ and for each of the $2-spheres$ linked to one of
    $\overset{=}{S}{}^{\alpha\beta\gamma}_{\delta}$
  the isotropy only depends on
one of   the minor determinants $\pm \Delta_{\alpha \beta \gamma}$, for the remaining
cases the isotropy dipends on a linear combination of $\Delta_{\alpha \beta \gamma}$;

\vspace{0,5 cm}

 3) three sets of  points such that, for each of them, the respective points are joined by
 one of the  $2-$spheres
$S^2\cong  \big(\, \overset{=}{S}{}^{\alpha\beta}_{\alpha\gamma}\cap N(\Theta)\big)^{\circ}/\tilde{G}$.
 Each of these sets
consists at most of four points.\efl

  \bfl In the case when some of the minor determinants $\Delta_{\alpha\beta\gamma}$
 and $\sq$ assumes values $\pm 1$ or when the isotropy, associated to some
  of the points mentioned in
  $3)$, vanishes, then  the singular locus
 $\Sigma(\Theta)$ is obtained by removing from the above list the  sets
 whose isotropy depends on  the mentioned determinants. $\square$ \efl

\end{teor}

 \begin{teor}
\bfl The singular points on the $3-$Sasakian orbifold $\mathcal{M}^7(\Theta)$
 come from the $G-$strata $\widetilde{S}$ 
  of $N(\Theta)$,
 which are obtained by intersectig $N(\Theta)$ with
 the following strata $\overline{S}$ of $\mathbb{H}^8$,
 and the ones  which
 can be obtained by changing $(w'_1, w'_2)$ with $(w''_1,  w''_2)$:  \efl
  \begin{equation} \label{C2 : 5}{
 \begin{aligned}
 & a) \,\,\, \quad  \overline{S}{}^{1234}:= 
  \bigg\{
\left( \begin{array}{cc|cc|cc|cc}
z_1 & z_2  & z_3 & z_4   &    z_5 & z_6        & z_7    &  z_8   \\
w'_1 & w'_2  &  w'_3  & w'_4     &   w'_5   &  w'_6          & w'_7  &  w'_8\\
\end{array}\right) \bigg\},   \\
& b) \,\,\, \quad \overline{S}{}^{123}_{4}:= 
\bigg\{
\left( \begin{array}{cc|cc|cc|cc}
z_1 & z_2  & z_3 & z_4   &    z_5 & z_6        & z_7    &  z_8   \\
w'_1 & w'_2  &  w'_3  & w'_4     &   w'_5   &  w'_6          & w''_7  &  w''_8\\
\end{array}\right) \bigg\},   \hspace{10 cm}\\
&c) \,\,\,\quad \overline{S}{}^{124}_{3}:= 
\bigg\{
\left( \begin{array}{cc|cc|cc|cc}
z_1 & z_2  & z_3 & z_4   &    z_5   & z_6           & z_7   & z_8  \\
w'_1 & w'_2  &  w'_3  & w'_4     &   w''_5  &  w''_6   & w'_7  &  w'_8\\
\end{array}\right) \bigg\},   \\
 &d)\,\,\, \quad \overline{S}{}^{134}_{2}:= 
 \bigg\{
\left( \begin{array}{cc|cc|cc|cc}
z_1 & z_2  & z_3 & z_4   &    z_5 & z_6        & z_7    &  z_8   \\
w'_1 & w'_2 & w''_3  & w''_4    &  w'_5  &  w'_6   & w'_7  &  w'_8\\
\end{array}\right) \bigg\},   \\
& e)\,\,\,\quad \overline{S}{}^{12}_{34}:= 
\bigg\{\left( \begin{array}{cc|cc|cc|cc}
z_1 & z_2  & z_3 & z_4   &    z_5 & z_6        & z_7    &  z_8   \\
w'_1 & w'_2  &  w'_3  & w'_4    &   w''_5  &  w''_6        & w''_7  &  w''_8\\
\end{array}\right) \bigg\},   \\
& f)\,\,\,\quad \overline{S}{}^{13}_{24}:= 
\bigg\{
\left( \begin{array}{cc|cc|cc|cc}
z_1 & z_2  & z_3 & z_4   &    z_5 & z_6        & z_7    &  z_8   \\
w'_1 & w'_2 & w''_3  & w''_4     &  w'_5  &  w'_6          & w''_7  &  w''_8\\
\end{array}\right) \bigg\},   \\
& g)\,\,\,\quad  \overline{S}{}^{14}_{23}:= 
\bigg\{
\left( \begin{array}{cc|cc|cc|cc}
z_1 & z_2  & z_3 & z_4   &    z_5 & z_6        & z_7    &  z_8   \\
w'_1 & w'_2 & w''_3  & w''_4     &   w''_5 & w''_6   & w'_7  &  w'_8\\
\end{array} \right)    \bigg\}, \\
& h)\,\,\,\quad  \overline{S}{}^{1}_{234} :=  
\bigg\{
\left( \begin{array}{cc|cc|cc|cc}
z_1 & z_2  & z_3 & z_4   &    z_5 & z_6        & z_7    &  z_8   \\
w'_1 & w'_2 & w''_3  & w''_4     &   w''_5 & w''_6        & w''_7  &  w''_8\\
\end{array}\right) \bigg\},
 \end{aligned}    }    \hspace{4,5 cm}    {}
 \end{equation}
\bfl where
\begin{equation} \label{C2 : 6}
\begin{aligned}
&  (w'_{2\alpha-1},
 w'_{2\alpha }):=\frac{e^{i\delta} }{ sin\, \varphi}
 \big( -z_{2\alpha} +iz_{2\alpha -1}\cos\, \varphi, \, z_{2\alpha-1} +iz_{2\alpha}\cos\, \varphi  \big), \\
and \,\,\,\, & \\
  &(w''_{2\alpha-1},  w''_{2\alpha})=: \frac{e^{i\delta} }{ sin\, \varphi} \big(-\zap -
i\zad\cos\, \varphi,\, \zad -i\zap\cos\, \varphi  \big).  \hspace{ 5 cm}\\
\end{aligned}
\end{equation}
 with $sin\, \varphi\ne 0$. Here, the element $\lambda=\epsilon + j\sigma\in Sp(1)$
is written as
  $\epsilon = \cost +$ $i(\sit cos\, \varphi)$  and
$\sigma =  \sit sin\, \varphi cos\, \delta +
i(\sit sin\, \varphi sin\, \delta ).$ Then, each $\overline{S}$
   in $(\ref{C2 : 5})$ generates a singular point
  $\big(\overline{S}\cap N(\Theta)\big)/ G  \in \mathcal{M}^7(\Theta)$ whose isotropy only
  depends on one of the minor determinants $\sq$. $\square$ \efl
 \end{teor}

\vspace{0,5 cm}

\section{The
Twistor Space $\mathcal{Z}^6(\Theta)$  }

\vspace{0,5 cm}

 \bfl By looking at the action of $\tilde{G}= Sp(1)\times T^3_{\Theta}\times U(1)$
on $ \mathbb{H}^8$, let $\uno $ be a non null point in
$  \mathbb{H}^8$. Then for each quaternionic pair $\uu$,
 the fixed point equations can be written as:\efl

  \begin{equation} \label{C2a : 1}
 A(\theta_{\alpha})
\left( \begin{array}{c}
u_{2\alpha-1} \\
u_{2\alpha}  \\
 \end{array}   \right) =
\left( \begin{array}{c}
\lambda u_{2\alpha-1} \rho \\
\lambda u_{2\alpha} \rho \\
 \end{array}   \right),
\end{equation}

\bfl  where
$A(\theta_{\alpha}):=\left(\begin{array}{cc}
\cost_{\alpha} & \sit_{\alpha} \\
-\sit_{\alpha}  & \cost_{\alpha} \\
\end{array} \right)$,
$\lambda= \epsilon + j\sigma\in Sp(1)$, $\rho\in U(1)$,
  $(u_{2\alpha-1}, u_{2\alpha})\in $
  $\mathbb{C}^4\times \mathbb{C}^4$
   and $\alpha\in\{1,2,3,4\}$.  \efl

\begin{prop}
Let $\zw$ be a point in $N(\Theta)$.
 Then, up to $\tilde{G}-$conjugation, we have that
 $\tilde{G}_{\zw}\subset T^3_{\Theta}\times \{\lambda\in Sp(1)\,\, |\,\,
 \sigma\equiv 0 \}\times U(1)$.
\end{prop}

\bfl  $\mathbf{Proof}.$
For each quaternionic pair $\uu$ the fixed point equations in (\ref{C2a : 1})
can be rewritten as follows\efl

 \begin{equation} \label{C2a : 2}
\left\{ \begin{array}{l}
 u_{2\alpha-1}\cost_{\alpha} + u_{2\alpha}\sit_{\alpha} =   \lambda u_{2\alpha-1} \rho, \\
 -u_{2\alpha-1}\sit_{\alpha} + u_{2\alpha}\cost_{\alpha} =   \lambda u_{2\alpha} \rho, \\
\end{array} \right.
 \end{equation}
\bfl where $\alpha\in \{1,2,3,4\}$.
 Now, by multiplying the first equation in
(\ref{C2a : 2}) by $\cost_{\alpha}$,
 the second one by $\sit_{\alpha}$, and by considering the sum and the
difference of these two equations, we get \efl

\begin{equation} \label{C2a : 3}
\left\{ \begin{array}{l}
 u_{2\alpha-1}\costd_{\alpha} + u_{2\alpha}\sitd_{\alpha} =   \lambda (u_{2\alpha-1}\cost_{\alpha} +
u_{2\alpha}\sit_{\alpha} )\rho, \\
 u_{2\alpha-1}  =   \lambda (u_{2\alpha-1}\cost_{\alpha} -
u_{2\alpha}\sit_{\alpha} )\rho. \\
\end{array} \right.
 \end{equation}

\bfl Then, by multiplying on the right the first equation in (\ref{C2a : 3}) by
the conjugated of the second one, we obtain  \efl

\begin{equation} \label{C2a : 4}
\left\{ \begin{array}{l}
u_{2\alpha-1}  =   \lambda (u_{2\alpha-1}\cost_{\alpha} -
u_{2\alpha}\sit_{\alpha} )\rho, \\
 |u_{2\alpha-1}|^2\costd_{\alpha} + u_{2\alpha}\overline{u}_{2\alpha-1}\sitd_{\alpha} = \\
 \hspace{3 cm}  =\lambda (u_{2\alpha-1}\cost_{\alpha} +
u_{2\alpha}\sit_{\alpha} )(\overline{u}_{2\alpha-1}\cost_{\alpha} -
\overline{u}_{2\alpha}\sit_{\alpha} )\overline{\lambda}, \\
\end{array} \right.
 \end{equation}
\bfl or equivalently:\efl
\begin{equation} \label{C2a : 5}
\left\{ \begin{array}{l}
u_{2\alpha-1}  =   \lambda (u_{2\alpha-1}\cost_{\alpha} -
u_{2\alpha}\sit_{\alpha} )\rho, \\
 |u_{2\alpha-1}|^2\costd_{\alpha} + u_{2\alpha}\overline{u}_{2\alpha-1}\sitd_{\alpha} = \\
 \hspace{ 3 cm}  =\lambda (|u_{2\alpha-1}|^2\costq_{\alpha}
-|u_{2\alpha}|^2\sitq_{\alpha} + Im(u_{2\alpha}\overline{u}_{2\alpha-1})\sitd_{\alpha}  )\overline{\lambda}, \\
\end{array} \right.   \hspace{10 cm}
 \end{equation}

\bfl or:  \efl

\begin{equation} \label{C2a : 6}
\left\{ \begin{array}{l}
  |u_{2\alpha-1}|^2\costd_{\alpha} + Re(u_{2\alpha}\overline{u}_{2\alpha-1})\sitd_{\alpha}=
   |u_{2\alpha-1}|^2\costq_{\alpha}
-|u_{2\alpha}|^2\sitq_{\alpha}, \\
Im(u_{2\alpha}\overline{u}_{2\alpha-1})\sitd_{\alpha} =
 \lambda Im(u_{2\alpha}\overline{u}_{2\alpha-1})\overline{\lambda} \sitd_{\alpha},  \\
u_{2\alpha-1}  =   \lambda (u_{2\alpha-1}\cost_{\alpha} -
u_{2\alpha}\sit_{\alpha} )\rho. \\
\end{array} \right.   \hspace{10 cm}
 \end{equation}

\bfl Suppose there exists a quaternionic pair $\uu$
such that  $Im(u_{2\alpha}\overline{u}_{2\alpha-1})\ne 0$. If $\sit_{\alpha}=0$, then
we can easily see that the equations  (\ref{C2a : 6}) describe the non effective
subgroup.
Instead,  when
 $\sitd_{\alpha}\ne 0$ the first two equations give\efl

\begin{equation} \label{C2a : 7}
\left\{ \begin{array}{l}
 (|u_{2\alpha}|^2 -|u_{2\alpha-1}|^2)\sitq_{\alpha}  +
  Re(u_{2\alpha}\overline{u}_{2\alpha-1})\sitd_{\alpha}=0, \\
Im(u_{2\alpha}\overline{u}_{2\alpha-1})=
 \lambda Im(u_{2\alpha}\overline{u}_{2\alpha-1})\overline{\lambda}.  \\
\end{array} \right.   \hspace{10 cm}
 \end{equation}
\bfl Since we are interested in points $\uno\in N(\Theta)$,
 we can suppose
$\uno\in $
$\nu^{-1}(0).$  Then, the first system of equations in (\ref{C1: 14}) gives \efl

\begin{equation}   \label{C2a : 7a}
 Im(u_{2\alpha}\overline{u}_{2\alpha-1})= c_{\alpha}Im( \overline{u}_8 u_7), \quad
 \alpha \in \{1,2,3,4\},
 \end{equation}

 \bfl and $c_{\alpha}\in \mathbb{Q}$. As a direct consequence, the second equation
 in (\ref{C2a : 7}) does not depends on the index $\alpha$, that is, it is
  the same for each quaternionic pair $\uu$, $\alpha \in\{1,2,$
   $3,4\}$. In the
case when $Im(u_{2\alpha}\overline{u}_{2\alpha-1})=i$,  the second equation
 in (\ref{C2a : 7}) reads: \efl
\begin{equation}  \label{C2a : 10}
\lambda i\overline{\lambda}=(\epsilon + j\sigma)i (\overline{\epsilon} -j\sigma)=
i(|\epsilon|^2 - |\sigma|^2) - 2k\sigma\overline{\epsilon}= i,
\end{equation}
 \bfl that, together with the condition $|\epsilon |^2 + |\sigma |^2=1$,
 yields $\sigma\equiv0$.  Thus, in this particular case
 we have that  $\tilde{G}_{\zw}\subset T^3\times U(1)^{\epsilon}\times U(1):=
  T^3_{\Theta}\times \{\lambda\in Sp(1)\,\, |\,\,
 \sigma\equiv 0 \}\times U(1)$.
  Note that, for each quaternionic pair $\uu$ there exists $\tau\in Sp(1)$ such that
   $\tau  Im(u_{2\alpha}\overline{u}_{2\alpha-1}) \overline{\tau}=i$.
 Then, in general the isotropy subgroup of $\zw$ is such that
 $\tilde{G}_{\zw}\subset
   T^3_{\Theta}\times \overline{\tau}U(1)^{\epsilon}\tau\times$
   $ U(1)$.
 Instead, when $Im(u_{2\alpha}\overline{u}_{2\alpha-1})= 0$ for each
$\alpha\in\{1,2,3,$
$4\}$, it follows that $u_{2\alpha}\overline{u}_{2\alpha-1}=$
$c_{\alpha}\in \mathbb{R}$.
 First suppose $c_{\alpha}\ne 0$ for each $\alpha\in\{1,2,3,4\}$.
  Then, the  fixed point
   equations  (\ref{C2a : 2}) become \efl
   \begin{equation}  \label{C2a : 11b}
\left\{ \begin{array}{l}
 u_{2\alpha-1}(\cost_{\alpha} + c_{\alpha} \sit_{\alpha}) =
  \lambda u_{2\alpha-1} \rho, \\
 u_{2\alpha-1}(-\sit_{\alpha} + c_{\alpha} \cost_{\alpha}) =
  c_{\alpha}  \lambda u_{2\alpha-1} \rho. \\
\end{array} \right.
\end{equation}
   \bfl As a consequence, we get  \efl
   \begin{equation}  \label{C2a : 12b}
 (-\sit_{\alpha} + c_{\alpha} \cost_{\alpha}) =  c_{\alpha}  ( \cost_{\alpha} +
 c_{\alpha} \sit_{\alpha}),
  \end{equation}

  \bfl if and only if
$ (1 + c^2_{\alpha}) \sit_{\alpha}=0$, and we get $\sit_{\alpha}=0$.
 Then, up to a sign, the fixed point equations become\efl

 \begin{equation}  \label{C2a : 11b}
\left\{ \begin{array}{l}
 u_1 =\lambda u_1 \rho, \\
  u_3= \lambda u_3 \rho, \\
  u_5 = \lambda u_5 \rho, \\
  u_7 =  \lambda u_7 \rho, \\
\end{array} \right.
\quad \Rightarrow  \quad
\left\{ \begin{array}{l}
 u_1 \overline{u}_{\beta} =\lambda u_1 \overline{u}_{\beta}  \overline{\lambda}, \\
 u_{\beta} \overline{u}_1 = \overline{\rho}\,  \overline{u}_{\beta} u_1  \rho, \\
\end{array} \right.
\end{equation}

\bfl where $\beta\in\{3,5,7\}$. Then, we obtain $\sigma\equiv0$, $\epsilon=\rho=\pm 1$,
by substituting these conditions in the fixed point equations  (\ref{C2a : 1})
it follows that the isotropy subgroup coincide with the non effective subgroup.
 Now, if there exists a quaternionic pair such that
$u_{2\alpha}\overline{u}_{2\alpha-1}$ $= c_{\alpha}= 0$,
 the  fixed point
   equations  (\ref{C2a : 2}) give $\sit_{\alpha}=0$ and using
    the same argumet of the
   previous case we have that the isotropy coincide with the ineffectivity.
   Then  we get the conclusion. $\square$  \efl

\bfl Suppose now that $\zw\in N(\Theta)$ has a quaternionic pair $\uu$ such that
$Im(u_{2\alpha}\overline{u}_{2\alpha-1})=$
$i$.
 Then, by using the $T^3_{\Theta}-$moment map
 equations ( look at (\ref{C1: 14}) and (\ref{C2a : 7a})),
 up to a sign, this condition holds
for each quaternionic pair of $\zw$. 
 Moreover, it follows that
$u_{2\alpha}= $
$(Re(u_{2\alpha}\overline{u}_{2\alpha-1}) +i)u_{2\alpha-1}$.
In general, we have that $Im(u_{2\alpha}\overline{u}_{2\alpha-1})=\tau i\overline{\tau}$,
 $\tau\in$ $ Sp(1)$,
and   $u_{2\alpha}\overline{u}_{2\alpha-1}= (Re(u_{2\alpha}\overline{u}_{2\alpha-1}) +
\tau i\overline{\tau})$ if and only if
 $\overline{\tau}u_{2\alpha}=$
 $ (Re(u_{2\alpha}\overline{u}_{2\alpha-1}) +  $
 $i)\overline{\tau}u_{2\alpha-1}$.    That is, for each quaternionic pair $\uu$, up to multiplication on
  the right by an element $\overline{\tau}\in Sp(1)$, we can assume
   $Im(u_{2\alpha}\overline{u}_{2\alpha-1})=
i$ and $G_{\zw}\subset T^3_{\Theta}\times U(1)^{\epsilon} \times$
$ U(1)$ .
 Then, in order to give a complete account
 of the singular $\tilde{G}-$strata of $N(\Theta)$
 it will be sufficient to describe the points
  $\zw\in N(\Theta)$
 with  isotropy subgroup
   $\tilde{G}_{\zw}$
   $\subset T^3_{\Theta}\times U(1)^{\epsilon}\times U(1)$. We then
    consider the orbits through these points.
In the following we will work under the assumption $\sigma\equiv0$. \efl

 \vspace{0,5 cm}

\bfl Note that the points
$(\underline{z}, \underline{w})\in N(\Theta)\subset S^{31}$
can be written as $4\times 8$ complex matrices, namely: \efl
\begin{equation}  \label{C2 : 7}
{}^T(\underline{z}, \underline{w} ):=
\left(\begin{array}{cc|cc|cc|cc}
z_1 & z_2  & z_3 & z_4 & z_5 & z_6 & z_7 & z_8\\
w_1 & w_2 &  w_3 & w_4 & w_5 & w_6 & w_7 & w_8 \\
\end{array}\right),
\end{equation}
\bfl and for each quaternionic pair the fixed point
  equations  become
    \efl
\begin{equation}  \label{C2 : 8}
A(\theta_{\alpha})
 \left(\begin{array}{cc}
z_{2{\alpha}-1} & w_{2{\alpha}-1} \\
z_{2{\alpha}} & w_{2{\alpha}} \\
\end{array} \right) =
{\left[ \begin{array}{c}
\left(\begin{array}{cc}
\epsilon & -\overline{\sigma}  \\
\sigma & \overline{\epsilon} \\
\end{array} \right)
 \left(\begin{array}{cc}
z_{2{\alpha}-1} & z_{2{\alpha}} \\
w_{2{\alpha}-1} & w_{2{\alpha}} \\
\end{array} \right)
 \left(\begin{array}{cc}
\rho& 0 \\
0 & \rho \\
\end{array} \right)
\end{array}\right]}^T ,
\end{equation}

\bfl where $(z_{\alpha}, w_{\alpha})\in \mathbb{C}$ and $\alpha\in\{1,2,3,4\}$.
  The equations in (\ref{C2 : 8}) can be easily rewritten as:\efl
\begin{equation}    \label{C2 : 9}
\overbrace{\left( \begin{array}{cc|cc}
0                          & -\sigma\rho &  -\sit_{\alpha}  & \cost_{\alpha} - \overline{\epsilon}\rho  \\
-\sigma\rho     & 0                      & \cost_{\alpha} - \overline{\epsilon}\rho & \sit_{\alpha} \\
\hline
-\sit_{\alpha}   & \cost_{\alpha} - \epsilon\rho           &    0       &      \overline{\sigma}\rho  \\
\cost_{\alpha} -\epsilon \rho & \sit_{\alpha} &  \overline{\sigma}\rho   &   0    \\
\end{array} \right)}^{M_{\alpha}:=}
\left(\begin{array}{c}
z_{2\alpha-1} \\
z_{2\alpha}\\
w_{2\alpha-1} \\
w_{2\alpha} \\
\end{array} \right) =
\left(\begin{array}{c}
0 \\
0\\
0\\
0 \\
\end{array} \right),
\end{equation}
 \bfl where $\alpha\in \{1,2,3,4\}$.
 Note that none of the matrices $M_{\alpha}$,
  $\alpha\in \{1,2,3,4\}$, has $rank\,\,4 $.
 Otherwise the corrispondent quaternionic
  pair $(z_{2\alpha-1}, w_{2\alpha-1}, z_{2\alpha}, w_{2\alpha})$ would vanish, and
    this contradicts lemma $1.1$.
   Now, we prove the following fact about the action of
$\tilde{G}=T^3_{\Theta}\times Sp(1)\times U(1)$ on $\mathbb{C}^8\times
 \mathbb{C}^8\cong\mathbb{H}^8$.\efl

\begin{lem}
Let $M_{\alpha}$ be the matrix
 appearing in formula $(\ref{C2 : 9})$. The condition
 $det\, \mal=0$, $\alpha\in \{1,2,3,4\},$ is equivalent to
\begin{equation}   \label{C2 : 18b}
\overline{\rho}e^{i\theta_{\alpha}}= Re(\epsilon) \pm  i\sqrt{Im(\epsilon)^2 + |\sigma|^2}
\,\,\,\,\,\, \mathrm{or}\,\,\,\,\,\,
\overline{\rho}e^{-i\theta_{\alpha}}= Re(\epsilon) \pm
 i\sqrt{Im(\epsilon)^2 + |\sigma|^2}.
\end{equation}
\end{lem}

\bfl$\mathbf{Proof.}$
Each matrix $M_{\alpha}$, $\alpha\in\{1,2,3,4\}$, can be written as follows \efl
\begin{equation}   \label{C2 : 10}
M_{\alpha}=\left(\begin{array}{c|c}
A & B \\
\hline
 C & D \\
 \end{array}\right):=
\left( \begin{array}{cc|cc}
0                          & -\sigma\rho &  -\sit_{\alpha}  & \cost_{\alpha} - \overline{\epsilon}\rho  \\
-\sigma\rho     & 0                      & \cost_{\alpha} - \overline{\epsilon}\rho & \sit_{\alpha} \\
\hline
-\sit_{\alpha}   & \cost_{\alpha} - \epsilon\rho           &    0       &      \overline{\sigma}\rho  \\
\cost_{\alpha} -\epsilon \rho & \sit_{\alpha} &  \overline{\sigma}\rho   &   0    \\
\end{array} \right).
\end{equation}
\bfl Since the submatrices $A,D$ are invertible,
 this gives\efl
\begin{equation}   \label{C2 : 11}
M_{\alpha}=  \left(\begin{array}{cc}
 B & Id\\
  D & 0 \\
 \end{array}\right) \left(\begin{array}{cc}
 D^{-1}C &  Id \\
 A- BD^{-1}C  &  0\\
 \end{array}\right).
\end{equation}
\bfl
 Note that $det\, M_{\alpha}=0 $
 if and only if $det\, (A- BD^{-1}C)=0$, and the matrix $A- BD^{-1}C$ is \efl
\begin{equation}     \label{C2 : 12}
\frac{\sigma}{|\sigma|^2}\left(\begin{array}{cc}
\overline{\rho} sin\, 2\theta_{\alpha} -2 \sit_{\alpha} Re(\epsilon) &   ( 2 \cost_{\alpha} Re(\epsilon)- \overline{\rho}cos\, 2\theta_{\alpha}) - \rho \\
( 2 \cost_{\alpha} Re(\epsilon)- \overline{\rho}cos\, 2\theta_{\alpha}) - \rho  &   -\overline{\rho} sin\, 2\theta_{\alpha} +2 \sit_{\alpha} Re(\epsilon) \\
 \end{array} \right).
\end{equation}
\bfl Then, up to the scale $(\frac{\sigma}{|\sigma|^2})^{-1}$,
we can rewrite the matrix in (\ref{C2 : 12}) as follows 
\begin{equation}   \label{C2 : 13}
\left(\begin{array}{cc}
Y & X-\rho\\
X-\rho & -Y  \\
\end{array}\right),
\end{equation}
where $X =  2 \cost_i Re(\epsilon)- \overline{\rho}cos\, 2\theta_{\alpha}$ and
 $Y=\overline{\rho} sin\, 2\theta_{\alpha} -2 \sit_{\alpha} Re(\epsilon)$. Then,
up to column permutations,  $det\, M_{\alpha}$ can be rewritten as
\begin{equation}   \label{C2 : 14}
\begin{aligned}
 det\, \Bigg[ \left(\begin{array}{cc}
 X & Y\\
-Y & X   \\
\end{array}\right) &- \rho Id_{2\times 2}\Bigg] =
det\, \left(\begin{array}{cc}
 X-\rho & Y\\
-Y & X-\rho   \\
\end{array}\right)  = \\
& \rho^2 -2XY +(X^2 + Y^2)= [\rho - (X +iY)][\rho- (X- iY)],
\end{aligned}
\end{equation}
 or equivalently:
\begin{equation}   \label{C2 : 15}
\begin{aligned}
& det\, (A- BD^{-1}C )=  \frac{\sigma^2}{|\sigma |^4}\big[\rho -
 \big( ( 2 \cost_{\alpha} Re(\epsilon)- \overline{\rho}cos\, 2\theta_{\alpha}) +
 i(\overline{\rho} sin\, 2\theta_{\alpha} +\\
&  -2 \sit_{\alpha} Re(\epsilon)) \big) \big]
\big[\rho
 -\big( ( 2 \cost_{\alpha} Re(\epsilon)- \overline{\rho}cos\, 2\theta_{\alpha}) -i(
\overline{\rho} sin\, 2\theta_{\alpha} -2 \sit_{\alpha} Re(\epsilon)) \big) \big].
\end{aligned}
\end{equation}
 Assume $\sigma \ne 0$, since $rank\,M_{\alpha}\leq 4$, for each
 $\alpha\in\{1,2,3,4\}$, it follows that \efl
\begin{equation} \label{C2 : 16}
\begin{aligned}
 0=det\, \mal = \overline{\sigma}^2\rho^2 &det\, (A- BD^{-1}C )= \\
 &=\rho^2\bigg( \rho+
 \overline{\rho} (e^{-i\theta_{\alpha}})^2
- 2 Re(\epsilon) e^{-i\theta_{\alpha}}\bigg)\bigg( \rho+
\overline{\rho} (e^{i\theta_{\alpha}})^2
- 2 Re(\epsilon) e^{i\theta_{\alpha}}\bigg),
\end{aligned}
\end{equation}
\bfl  if and only if  \efl
\begin{equation}  \label{C2 : 17}
 \rho+ \overline{\rho} (e^{-i\theta_{\alpha}})^2 - 2 Re(\epsilon) e^{-i\theta_{\alpha}}=0
  \,\,\,\,\,\, \mathrm{or}\,\,\,\,\,\,
 \rho+ \overline{\rho} (e^{i\theta_{\alpha}})^2- 2 Re(\epsilon) e^{i\theta_{\alpha}}=0,
\end{equation}
\bfl  or equivalently \efl
\begin{equation}   \label{C2 : 18}
\overline{\rho}e^{i\theta_{\alpha}}= Re(\epsilon) \pm  i\sqrt{Im(\epsilon)^2 + |\sigma|^2}
\,\,\,\,\,\, \mathrm{or}\,\,\,\,\,\,
\overline{\rho}e^{-i\theta_{\alpha}}= Re(\epsilon) \pm  i\sqrt{Im(\epsilon)^2 + |\sigma|^2}.
\end{equation}
\bfl Note that, when $\sigma\equiv 0$ the condition $det\, M_{\alpha}=0$ is equivalent to
$\epsilon\rho= e^{\pm i\theta_{\alpha}}$ or  $\overline{\epsilon}\rho= e^{\pm i\theta_{\alpha}}$,
which are particular cases of the relations in (\ref{C2 : 18}).
Then we get the conclusion.  $\square$ \efl

\vspace{0,5 cm}

\bfl $\mathbf{Remark\,\, 2.1.}$  Looking at the relations in (\ref{C2 : 18}), we
can define the map
\begin{equation}   \label{C2 : 23}
 \begin{aligned}
&  \,\, \Psi\, : \,    Sp(1)\longrightarrow U(1).    \\
 &\quad \quad \epsilon +j\sigma \longmapsto  Re(\epsilon) +
 i\sqrt{Im(\epsilon)^2 + |\sigma|^2} \in U(1)
 \end{aligned}
 \end{equation}
\bfl The elements $\lambda= \epsilon +j\sigma \in Sp(1)$ can be parametrized as follows \efl
\begin{equation}     \label{C2 : 24}
\left\{ \begin{array}{l}
\epsilon = \cost +i(\sit cos\, \varphi), \\
\sigma =  \sit sin\, \varphi sin\, \delta +i(\sit
   sin\, \varphi cos\, \delta ). \\
\end{array} \right.
\end{equation}
\bfl Then   \efl
 \begin{equation}     \label{C2 : 24a}
 \Psi(\epsilon + j\sigma)=  \cost +i\sqrt{\sit^2 sin^2\, \varphi +
 \sit^2 cos^2\, \varphi   }  =  \cost + i|\sit |,
\end{equation}
 and the relations in (\ref{C2 : 18}) can be rewritten as follows\efl
\begin{equation}    \label{C2 : 24b}
\overline{\rho}e^{i\theta_{\alpha}}= e^{\pm i\theta}
\,\,\,\,\,\, \mathrm{or}\,\,\,\,\,\,
\overline{\rho}e^{-i\theta_{\alpha}}= e^{\pm i\theta}.
\end{equation}
 \bfl Note that the map $\Psi$ in (\ref{C2 : 23}) decsribes
  a Lie group homomorphism if and only if $\sigma\equiv 0$. $\square$
  \efl

 \vspace{0,5 cm}

\bfl  $\mathbf{Remark\,\,2.2.}$
 By substituting in each block $M_{\alpha}$ the expressions of
 $\rho\epsilon$ and $\rho\overline{\epsilon}$ in (\ref{C2 : 18}), the equations
  (\ref{C2 : 9})
  can be turned into an eigenvalue problem. Then, each of this choices 
  provides non null solutions $\zw\in \mathbb{H}^8$ for the  equations
  (\ref{C2 : 9}) and
   these eigenvectors 
    are characterized by having the same
   isotropy subgroup. Moreover, this isotropy is generated exactly by $(\rho\epsilon, \rho\overline{\epsilon})$.  \efl


\bfl
$1)$\,\,\,  just one of the relations in  (\ref{C2 : 18}) holds $\iff$ $rank\, \mal=3,$ \\
$2)$\,\,\,    two relations in (\ref{C2 : 18}) hold $\iff$  $rank\, \mal=2,$  \\
$3)$\,\,\,  three or four relations in (\ref{C2 : 18})
 are satisfied  $\iff$  $ \mal=0_{4\times 4}.$
\efl

\bfl Actually in the case $3)$ we have that if three relations of those in
in (\ref{C2 : 18}) hold, then also the last one holds. Moreover, when all of the four
 relations  in  (\ref{C2 : 18})
 hold, for each $\alpha\in\{1,2,3,4\}$, they describe
 the effectivity. In the following we will work up to the non effective subgroup,
so we do not have to care about the case $3)$. $\square$  \efl

\vspace{0,5 cm}

\begin{lem}
\bfl Assume $\sigma\equiv 0$ in equations $(\ref{C2 : 9})$.
Let $\mal$ be such that $rank\,M_{\alpha}= 3$.  Then
the solutions of the eigenvalue problem associated to the equations  $(\ref{C2 : 9})$
describe a complex $1-$dimensional
space which is either \efl
 \begin{equation}  \label{C2 : 30}
\begin{aligned}
1)\,\, \hspace{1 cm} & {}^{\pm}V_1^{\alpha}:=\{z_{2\alpha-1}\in \mathbb{C}\, |
\, (z_{2\alpha-1}, \pm iz_{2\alpha-1}, 0,0)
\in \mathbb{C}^4\}, \quad\quad \mathrm{or}  \\
&   \\
2)\,\, \hspace{1 cm} &{}^{\pm}V_2^{\alpha}:=\{w_{2\alpha-1}\in \mathbb{C}\, |
 \, ( 0,0, w_{2\alpha-1}, \pm iw_{2\alpha-1})\in \mathbb{C}^4\}. \hspace{ 10 cm}    {}
\end{aligned}    \hspace{2 cm}
\end{equation}  \\
When $rank\,M_\alpha =2 $, the possible spaces of solutions are the following:
\begin{equation}   \label{C2 : 30b}
\begin{aligned}
3)\,\,\hspace{1 cm} &{}^{(\pm, \pm )}V_3^{\alpha}:=\{(z_{2\alpha-1},w_{2\alpha-1})\in
 \mathbb{C}^2\, |
 \, (z_{2\alpha-1}, \pm iz_{2\alpha-1}, w_{2\alpha-1}, \pm iw_{2\alpha-1} )
\in \mathbb{C}^4\},  \hspace{10 cm}    \\
& \\
4)\,\,\hspace{1 cm} & V_4^{\alpha}:=\{z_{2\alpha-1},z_{2\alpha} \in \mathbb{C}\, |
\, (z_{2\alpha-1}, z_{2\alpha }, 0,0)
\in \mathbb{C}^4\},  \\
& \\
5)\,\, \hspace{1 cm}&V_5^{\alpha}:=\{w_{2\alpha-1 },w_{2\alpha}\in \mathbb{C}\, |
 \, ( 0,0,
w_{2\alpha-1},
w_{2\alpha})\in \mathbb{C}^4\}. \hspace{10 cm}  {}
 \end{aligned}
 \end{equation}
\end{lem}

\bfl $\mathbf{Proof.}$
 Under  the hypothesis $\sigma\equiv 0$, we get \efl
 \begin{equation}  \label{C2 : 30a}
\begin{aligned}
det M= (-\sit_{\alpha}^2- (\cost_{\alpha} - \overline{\epsilon}\rho)^2)&(-\sit_{\alpha}^2-
(\cost_{\alpha} - \epsilon\rho)^2)
= \\
& =((\overline{\epsilon}\rho)^2 -2\overline{\epsilon}\rho\cost_{\alpha} +1  )
((\epsilon\rho)^2 -2\epsilon\rho\cost_{\alpha} +1  )=0, \\
\end{aligned}
\end{equation}
\bfl and hence  \efl
 \begin{equation}  \label{C2 : 31}
   \rho\epsilon= e^{\pm i\theta_{\alpha}} \,\,\quad\quad \mathrm{or}\,\,\quad \quad
\rho\overline{\epsilon} = e^{\pm i\theta_{\alpha}}.
\end{equation}
\bfl By substituting in the equations (\ref{C2 : 9}) the possible values of
$\rho\epsilon$ and $\overline{\rho}\epsilon$, we get the
respective spaces of solutions for the eigenvalue problem associated to
the equations in (\ref{C2 : 9}). Note that $rank\,\, M_{\alpha}= 3$ if and only if
 just one of the four conditions in (\ref{C2 : 31}) is satisfied.
 In this case, a short computation shows that the space of solutions for the equations (\ref{C2 : 9})
 is either ${}^{\pm}V_1^{\alpha}$ or ${}^{\pm}V_2^{\alpha}$.
 Instead,  $rank\,\, M_{\alpha}=2$, if and only if  exactly two conditions  of those
  in (\ref{C2 : 31}) are verified.
 Then 
  we can describe the spaces of solutions,
 which strictly depends on
  the chosen values for
$\epsilon\rho$ and $\overline{\epsilon}\rho$. So we have
\begin{equation} \label{C2 : 32}
\begin{aligned}
3)\,\,\hspace{0,5  cm}&\left\{\begin{array}{l}
\epsilon\rho= e^{\pm \theta_{\alpha}},\\
\overline{\epsilon}\rho= e^{\pm i\theta_{\alpha}}, \\
\end{array} \right.     \,\,\quad \Rightarrow \,\,\quad {}^{(\pm, \pm )}V_3^\alpha,  \quad \quad
 4)\,\,\hspace{0,5  cm} \left\{\begin{array}{l}
\epsilon\rho= e^{i\theta_{\alpha}},\\
\epsilon\rho= e^{-i\theta_{\alpha}},\\
\end{array} \right.     \,\,\quad \Rightarrow \,\,\quad V_4^{\alpha},\\
&  \\
 5)\,\,\hspace{0,5  cm}& \left\{\begin{array}{l}
\overline{\epsilon}\rho= e^{i\theta_{\alpha}},\\
\overline{\epsilon}\rho= e^{-i\theta_{\alpha}},\\
\end{array}\right.     \,\,\quad \Rightarrow \,\, \quad V_5^{\alpha}.   \hspace{10 cm} {}
\end{aligned}
\end{equation}
This concludes the proof.  $\square$  \efl

\vspace{0.5 cm}

\bfl Now, we are ready to prove that the isotropy of
 a generic point $\zw\in N(\Theta)$, is described by a suitable subset of the
 relations in (\ref{C2 : 18}). As we will see in a moment, this fact allows to prove that,
 up to the $\tilde{G}-$action,
 the quaternionic pairs $\uu$ of any element  $\zw\in N(\Theta)$  are described by
 the eigenspaces in (\ref{C2 : 30}) and (\ref{C2 : 30b}).\efl

\begin{prop}
\bfl Let $(\underline{z}, \underline{w})$ be  a  point
in $N(\Theta)\subset S^{31}\subset \mathbb{H}^8$.  Then the fixed
point equations $( \ref{C2a : 2})$ are
equivalent to the conditions $(\ref{C2 : 18})$ given by
$det\, \mal=0$, $\alpha\in \{1,$
$2,3,4\}$. \efl
\end{prop}

\bfl $\mathbf{Proof.} $  Let $\zw$ a point
 in $N(\Theta)\subset\mathbb{H}^8$. Suppose $\zw$ has discrete but not trivial isotropy, then
 at least one of the relations in (\ref{C2 : 18}) has to be satisfied.  Morever,
 by using the proposition $2.1$ we can assume $\sigma\equiv0$ and $\lambda= \epsilon$.
 Suppose exactly one of the relations in $(\ref{C2 : 18})$ is verified, say $\overline{\epsilon}\rho=
  e^{i\theta_{\alpha}}$, then  for each quaternionic pair $(u_{2\alpha-1},$
  $ u_{2\alpha}),$ $\alpha\in\{1,2,3,4\}$,
  the fixed point equations  (\ref{C2a : 6}) become \efl

\begin{equation}  \label{C2 : 40c}
\left\{ \begin{array}{l}
   e^{i\theta_{\alpha}}u_{2\alpha- 1} = \rho (u_{2\alpha-1 }\cost_{\alpha}-
   u_{2\alpha}\sit_{\alpha} ) \rho, \\
 (|u_{2\alpha}|^2 -|u_{2\alpha-1}|^2)\sitq_{\alpha}  +
  Re(u_{2\alpha}\overline{u}_{2\alpha-1})\sitd_{\alpha}=0. \\
  \end{array} \right.
\end{equation}

\bfl Let us rewrite the first equation in (\ref{C2 : 40c}) in the complex coordinates
$u_{\beta}=z_{\beta} + jw_{\beta}$: \efl

\begin{equation}  \label{C2 : 41c}
\left\{ \begin{array}{l}
  z_{2\alpha-1} e^{i\theta_{\alpha}} =  (z_{2\alpha-1} \cost_{\alpha} -z_{2\alpha}\sit_{\alpha})\rho^2,  \\
  (iw_{2\alpha-1}- w_{2\alpha} )\sit_{\alpha}=0.   \\
\end{array}\right.
\end{equation}

\bfl Now, by multiplying on the right the first
 equation in (\ref{C2 : 40c}) by its conjugated, we get:\efl
 \begin{equation}  \label{C2 : 42c}
 (|u_{2\alpha-1}|^2 -|u_{2\alpha}|^2)\sitq_{\alpha}  +
  Re(u_{2\alpha}\overline{u}_{2\alpha-1})\sitd_{\alpha}=0,
 \end{equation}

\bfl that together with the second equation in (\ref{C2 : 40c}) gives \efl

\begin{equation} \label{C2 : 43c}
\begin{aligned}
& \left\{ \begin{array}{l}
(|u_{2\alpha-1}|^2 -|u_{2\alpha}|^2)\sitq_{\alpha}=0, \\
  Re(u_{2\alpha}\overline{u}_{2\alpha-1})\sitd_{\alpha}=0, \\
\end{array} \right. \\
\mathrm{or \,\,equivalently:} \quad & \\
& \left\{ \begin{array}{l}
  (|z_{2\alpha-1}|^2 + |w_{2\alpha-1}|^2 - |z_{2\alpha}|^2 -
  |w_{2\alpha}|^2)  \sitq_{\alpha}=0, \\
  (\zap\bzad + \wap\bwad)\sitd_{\alpha} =0.
\end{array}\right. \hspace{10 cm}  \\
\end{aligned}
\end{equation}

\bfl Then, if $\sit_{\alpha}\ne 0$, the equations in (\ref{C2 : 41c}) and
(\ref{C2 : 43c}) give  \efl

\begin{equation} \label{C2 : 44c}
\left\{\begin{array}{l}
 |z_{2\alpha-1}|^2 = |z_{2\alpha}|^2, \\
  w_{2\alpha} =  iw_{2\alpha-1}, \\
  \zap = -i |\wad |^2 \zad, \\
  \end{array} \right. \quad \Rightarrow \quad
     \left\{\begin{array}{l}
 |z_{2\alpha-1}|^2 = |z_{2\alpha}|^2, \\
 | \zap |^2 =  |\wad |^2 |\zad |^2. \\
  \end{array} \right.
\end{equation}

\bfl As a consequence, these equations either give
\begin{equation}  \label{C2 : 45c}
 |\wad |^2=|\wap |^2=1 \quad \mathrm{or} \quad
 |z_{2\alpha-1}|^2 = |z_{2\alpha}|^2=0.
 \end{equation}
Since $\zw\in N(\Theta)\subset S^{31}\subset \mathbb{H}^8$, it follows that \efl
\begin{equation}     \label{C2 : 45d}
\sum_{\alpha=1}^8 |u_{\alpha}|^2 =
 \sum_{\alpha=1}^8 (|z_{\alpha}|^2 + |w_{\alpha}|^2)=1.
 \end{equation}
 \bfl Let us recall that none of the
 quaternionic pair $\uu$ can vanish, then  we can exclude the first relation
 in (\ref{C2 : 45c}), so we get $z_{2\alpha-1}=z_{2\alpha}=0$ and  $w_{2\alpha} =
 iw_{2\alpha-1}$.
 Note that we can repeat this argument for each quaternionic pair $\uu $. Then,
 when the isotropy subgroup $\tilde{G}_{\zw}$ is different from the non
  effective subgroup, it
 follows that
 $\uu\in \tilde{G}\cdot {}^{\pm}V_2^{\alpha}$.  As we have seen in
 remark $2.2$ the isotropy
 relative to the eigenspaces  ${}^{\pm}V_2^{\alpha}$ is described by
 the conditions
 $\overline{\epsilon}\rho= e^{\pm i\theta_{\alpha}}$. The other cases
 can be treated similarly.
 Suppose now, there exists $\alpha\in\{1,2,3,$
 $4\}$ such that  two relations
  of those in in (\ref{C2 : 18}) are satisfied. Then, up to a sign, we  have one of the
 following possible cases
 \begin{equation}   \label{C2 : 46c}
 \begin{aligned}
 & 1) \quad \quad \left\{\begin{array}{l}
 \overline{\epsilon}\rho =   e^{i\theta_{\alpha}},  \\
 \overline{\epsilon}\rho =  e^{-i\theta_{\alpha}},   \\
 \end{array} \right.
\quad \mathrm{or} \quad \quad
 2) \quad \quad \left\{\begin{array}{l}
 \epsilon\rho =   e^{i\theta_{\alpha}},    \\
 \overline{\epsilon}\rho =  e^{i\theta_{\alpha}}, \\
 \end{array} \right.  \hspace{15 cm} \\
  &  \\
& 3)\quad \quad \left\{\begin{array}{l}
\epsilon\rho =   e^{i\theta_{\alpha}},  \\
 \overline{\epsilon}\rho =  e^{-i\theta_{\alpha}}. \\
 \end{array} \right. \hspace{15 cm}
 \end{aligned}
 \end{equation}
 In the case $1)$ in (\ref{C2 : 46c}) we have that
\begin{equation}   \label{C2 : 46d}
e^{i\theta_{\alpha}}= e^{ik\pi} \quad\quad
\mathrm{and} \quad\quad \epsilon = \rho':= \rho e^{i k\pi},
\end{equation}
 the fixed point equations  (\ref{C2a : 2}) give \efl
\begin{equation} \label{C2 : 47c}
\left\{\begin{array}{l}
  u_{2\alpha-1}= \rho u_{2\alpha -1} \rho, \\
  u_{2\alpha}= \rho u_{2\alpha } \rho, \\
 \end{array} \right.
\end{equation}

 \bfl if and only if $\zad = \zad \rho^2$ and  $\zap = \zap \rho^2$.
 Then, if $\rho\ne \pm 1$
 we have that $\zad =$ $ \zap=0 $ and  that
  $\uu\in \tilde{G}\cdot V_5^{\alpha}$.  Instead, in the case $2)$ we have that
  $ \epsilon\rho =\overline{\epsilon}\rho$ if and only if $\epsilon= e^{ik\pi}$ and
  $\rho':= \rho e^{ik\pi}= e^{i\theta_{\alpha}}$.  Then, the fixed point equations
 (\ref{C2a : 2}) become \efl

\begin{equation} \label{C2 : 48c}
\left\{\begin{array}{l}
  u_{2\alpha-1}e^{-i\theta_{\alpha}}= u_{2\alpha -1}\cost_{\alpha}- u_{2\alpha}\sit_{\alpha}, \\
   (|u_{2\alpha}|^2 -|u_{2\alpha-1}|^2)\sitq_{\alpha}  +
  Re(u_{2\alpha}\overline{u}_{2\alpha-1})\sitd_{\alpha}=0, \\
 \end{array} \right.
\end{equation}

\bfl or equivalently:\efl

\begin{equation} \label{C2 : 49c}
\left\{\begin{array}{l}
  u_{2\alpha-1}i\sit_{\alpha} = -u_{2\alpha}\sit_{\alpha}, \\
   (|u_{2\alpha}|^2 -|u_{2\alpha-1}|^2)\sitq_{\alpha}  +
  Re(u_{2\alpha}\overline{u}_{2\alpha-1})\sitd_{\alpha}=0. \\
 \end{array} \right.
\end{equation}

 \bfl When $\sit_{\alpha}\ne 0$ we have that  $ u_{2\alpha-1}i= -u_{2\alpha}$
 and the second equation   in  (\ref{C2 : 49c}) is automatically satisfied. Then,
$\uu\in \tilde{G}\cdot {}^{(\pm , \pm)}V_3^{\alpha}$. The case $3)$ can
be treated similarly.
Then we get the conclusion. $\square$  \efl

 \bfl As a consequence, we have\efl

\begin{cor}
\bfl Let $\zw$ be a point in $N(\Theta)$ with non trivial isotropy
subgroup $\tilde{G}_{\zw}$.
Then each quaternionic pair $\uu$ of $\zw$
belongs to one of the sets $\tilde{G}\cdot {}^{\pm}V_{1}^{\alpha}$,
$\tilde{G}\cdot {}^{\pm}V_{2}^{\alpha}$, $\tilde{G}\cdot {}^{(\pm, \pm)}V_{3}^{\alpha}$,
$\tilde{G}\cdot V_{3}^{\alpha}$ or $\tilde{G}\cdot V_{4}^{\alpha}$. \efl
\end{cor}

\bfl Now, we have the following \efl

\begin{prop}
\bfl     Let $\zw$ be  a point on a singular $\tilde{G}-$orbit of $N(\Theta)$.
Then at most
one of the blocks $ M_{\alpha}$,  $\alpha\in\{1,2,3,4\}$, has $rank=2$.  \efl
\end{prop}

\bfl $\mathbf{Proof.}$
 Suppose not. Then we have that
$rank\, M_{\alpha}=3$, $\alpha\in\{1,2,3,4\}$. It follows that, for each
$M_{\alpha}$,
just one of the equations in (\ref{C2 : 31}) is satisfied. We are interested in
 the case when
 there exist
two different blocks $M_{\gamma}$ and $M_{\delta}$, such that for $M_{\gamma}$
the first equation in (\ref{C2 : 31}) is satisfied and for $M_{\delta}$
 the second one. Note that,
in this case we have
   $e^{i\theta_{\gamma}}\ne e^{\pm i \theta_{\delta}}$,
otherwise we could easily show
that both $M_{\gamma}$ and $M_{\delta}$ have $rank = 2$.
If such two blocks do not exist, then it follows that
the spaces of solutions relative the equations in (\ref{C2 : 9}), up to the $\tilde{G}-$action,
  are either
\efl
\begin{equation}  \label{C2 : 33}
\begin{aligned}
&\Bigg\{ X\in M_{4\times8}(\mathbb{C})\,\, |\,\, X=\left(\begin{array}{cc|cc|cc|cc}
* & *  & * & * & * & * & * & *\\
0 & 0  & 0 & 0 & 0 & 0 & 0  & 0 \\
\end{array}\right)\Bigg\},  \\
\mathrm{or} \hspace{1 cm}& \\
&\Bigg\{ X\in M_{4\times8}(\mathbb{C})\,\, |\,\, X=\left(\begin{array}{cc|cc|cc|cc}
0 & 0  & 0 & 0 & 0 & 0 & 0  & 0 \\
* & *  & * & * & * & * & * & *\\
\end{array}\right)\Bigg\},    \hspace{6 cm} {}
\end{aligned}
\end{equation}
\bfl where the above $2\times 2$
blocks
 $\left(\begin{array}{cc}
* & *  \\
0 & 0   \\
\end{array}\right) $ and
$\left(\begin{array}{cc}
0 & 0   \\
* & *  \\
\end{array}\right) $ represents an element of ${}^{\pm}V^{\alpha}_1$
or ${}^{\pm}V^{\alpha}_2$. Since the point $(\underline{z}, \underline{w})$ gives
a $8\times 4$ real matrix with $rank\,\, 4$, then it follows that the
 matrices in (\ref{C2 : 33}) are not admissible as point
$(\underline{z}, \underline{w})\in N$.
Instead, when there exist the mentioned two blocks $M_{\gamma}$ and $M_{\delta}$,
by using the proposition $2.2$, we have that
 for each
$\alpha\in\{1,2,3,4\}$ the spaces of solutions of equations (\ref{C2 : 9})
are either ${}^{\pm}V^{\alpha}_1$ or ${}^{\pm}V^{\alpha}_2$. Thus, up to column
permutations, the spaces of solutions
 have necessarely one of the following shapes\efl
\begin{equation}  \label{C2 : 34}
 \begin{aligned}
 1)\quad \quad \quad & \Bigg\{ X\in M_{4\times8}(\mathbb{C})\,\, |\,\, X=\left(\begin{array}{cc|cc|cc|cc}
* & *  & * & * & * & * & 0 & 0\\
0 & 0  & 0 & 0 & 0 & 0 & * & * \\
\end{array}\right)\Bigg\},   \hspace{10 cm} \\
2)\quad \quad \quad &
 \Bigg\{ X\in M_{4\times8}(\mathbb{C})\,\, |\,\, X=\left(\begin{array}{cc|cc|cc|cc}
 0 & 0  & 0 & 0 & 0 & 0 & * & * \\
* & *  & * & * & * & * & 0 & 0\\
\end{array}\right)\Bigg\},  \\
 3)\quad \quad \quad& \Bigg\{ X\in M_{4\times8}(\mathbb{C})\,\, |\,\, X=\left(\begin{array}{cc|cc|cc|cc}
* & *  & 0 & 0 & * & * & 0 & 0\\
0 & 0  & * & * & 0 & 0 & * & * \\
\end{array}\right)\Bigg\}.  \hspace{10 cm}{}
 \end{aligned}
 \end{equation}
\bfl  where the  $2\times 2$
blocks  represent either an elements
of ${}^{\pm}V^{\alpha}_1$ or ${}^{\pm}V^{\alpha}_2$. Consider the case $1)$
 in (\ref{C2 : 34})
when
$\uu\in{}^{+}V_1^{\alpha}$,  $\alpha\in\{ 1,3\}$, $(u_3,u_4)\in{}^{-}V_1^{2}$
and $(u_7,u_8)\in{}^{+}V_2^{4}$.
 By reading the $T^3_{\Theta}-$moment map equations on this set of solutions we get\efl
 \begin{equation}       \label{C2a : 72a}
 T^3_{\Theta}\quad \Rightarrow\quad
 \left\{\begin{array}{l}
p_1|z_1|^2- p_2|z_3|^2 +p_3|z_5|^2 +p_4|z_7|^2=0,\\
 q_1|z_1|^2- q_2|z_3|^2 +q_3|z_5|^2 +q_4|z_7|^2=0,\\
 l_1|z_1|^2- l_2|z_3|^2 +l_3|z_5|^2 +l_4|z_7|^2=0,\\
 \end{array}\right. \hspace{ 3 cm} {}
 \end{equation}
 \bfl and  \efl
 \begin{equation}     \label{C2a : 73a}
  Sp(1) \quad \Rightarrow\quad
 |z_1|^2  +|z_3|^2 +|z_5|^2 -|w_7|^2=0,
     \hspace{3 cm} {}
  \end{equation}
\bfl  Moreover, we have to consider the sphere equation  \efl
 \begin{equation}
 |z_1|^2
 + |z_3|^2 +|z_5|^2 +|w_7|^2=1.
 \end{equation}
\bfl Then we can rewrite these equations as follows  \efl
 \begin{equation}      \label{C2a : 74a}
\left( \begin{array}{cccc}
 p_1 & -p_2 & p_3  & p_4 \\
 q_1 & -q_2 & q_3  & q_4 \\
  l_1 & -l_2 & l_3  & l_4 \\
  1 & 1  & 1 & -1 \\
 \end{array} \right)
 \left( \begin{array}{c}
 |z_1|^2 \\
  |w_3|^2 \\
  |w_5|^2 \\
  |w_7|^2 \\
 \end{array} \right)    =
 \left( \begin{array}{c}
0\\
 0 \\
  0 \\
  0\\
 \end{array} \right),
 \end{equation}
\bfl and the determinant is exactly one of  $\sq\ne 0$. Then we
 get only the trivial solution.  The remaining cases can be
  treated similarly. $\square$ \efl

\vspace{0,5 cm}

 \bfl Another useful fact is the following \efl

\begin{prop}
\bfl Let $(\underline{z}, \underline{w})$ be
 is on a singular $\tilde{G}-$orbit of $N(\Theta)$ and suppose
  $\sigma\equiv0$ in equations
  $(\ref{C2 : 9})$.
 Suppose that one block  $M_{\beta}$ has
  $rank  =2$ and the system of equations $3)$ in $(\ref{C2 : 32})$ holds.
   Then,  for each $\alpha\in\{1,2,3,4\}$ the block
  $M_{\alpha}$ has $rank\leq2$.   \efl
\end{prop}

\bfl $\mathbf{Proof.}$
 Without loss of generality, we can assume $rank\, M_1 =2$
 and
 \begin{equation}    \label{C2 : 38a}
 \left\{\begin{array}{l}
\epsilon \rho= e^{i\theta_1},\\
\overline{\epsilon} \rho= e^{i\theta_1},\\
\end{array} \right.
 \end{equation}
 Now, consider another block $M_{\gamma}$ with $\gamma\ne 1$. Since
$rank\,\ M_{\gamma}\leq 3$, then at least one of the relations in (\ref{C2 : 31}) holds.
So, we get either $\epsilon\rho = e^{\pm i\theta_{\gamma}}$ or
 $\overline{\epsilon}\rho = e^{\pm i\theta_{\gamma}}$. In the first case
 (the other one can be treated in the same way),
 we get either  \efl

 \begin{equation}   \label{C2 : 38}
\left\{ \begin{array}{l}
\epsilon\rho = e^{i\theta_1}, \\
\overline{\epsilon}\rho = e^{ i \theta_1}, \\
 \epsilon\rho = e^{ i\theta_{\gamma}}, \\
\end{array} \right.   \quad \quad \mathrm{or} \quad \quad
\left\{ \begin{array}{l}
\epsilon\rho = e^{i\theta_1}, \\
\overline{\epsilon}\rho = e^{ i \theta_1}, \\
 \epsilon\rho = e^{ -i\theta_{\gamma}}. \\
\end{array} \right.
\end{equation}

 \bfl In particular it we have either $e^{ i \theta_1}=\overline{\epsilon}\rho
 =e^{i\theta_{\gamma}} $ or $e^{ i \theta_1}=\overline{\epsilon}\rho
 =e^{ -i\theta_{\gamma}} $.
 Thus, we also get one of the further conditions $\overline{\epsilon}
 \rho =e^{\pm i\theta_{\gamma}} $, and it follows that $rank\, M_{\alpha} \leq 2$ for each
 $\alpha\in\{1,2,3,4\}$. $\square$\efl

 \vspace{0,5 cm}

\begin{cor}
Let $\zw\in N(\Theta)$ be such that one of its quaternionic pairs
$\uu$ belongs to one of the spaces ${}^{(\pm, \pm )}V_3^\alpha$,
and another one, say $(u_{2\beta-1}, u_{2\beta})$ with $\alpha\ne \beta$, which
either belongs to $V^4_\beta$ or $V^5_\beta$.
Then the isotropy subgroup $\tilde{G}_{\zw}$ is either the identity $\{Id\}$
or the non effective subgroup.
\end{cor}

\bfl $\mathbf{Proof.}$ Suppose $\uu\in{}^{(+, + )}V_3^\alpha$ and
 $(u_{2\beta-1}, u_{2\beta})\in V_4^\beta$, the other cases can be treated similarly.
 Then the fixed point equations relative to these quaternionic pairs gives    \efl
 \begin{equation}    \label{C2 : 39}
 \left\{\begin{array}{l}
 \epsilon\rho= e^{ i\theta_{\alpha}},\\
\overline{\epsilon}\rho= e^{ i\theta_{\alpha}},\\
 \epsilon\rho= e^{ i\theta_{\beta}},\\
 \epsilon\rho= e^{-i\theta_{\beta}},\\
 \end{array}\right. \quad \quad \iff\quad \quad
 \left\{\begin{array}{l}
 e^{ i(\theta_{\alpha}- \theta_{\beta)} }=1,\\
 \epsilon\rho= e^{ i\theta_{\alpha}},\\
 \epsilon=\overline{\epsilon},\\
\rho=\overline{\rho},\\
 \end{array}\right.
 \end{equation}
\bfl In particular, the system of equations (\ref{C2 : 39}) holds for $\alpha=1$
and it follows that   \efl
 \begin{equation}        \label{C2 : 40}
\left\{\begin{array}{l}
 e^{ i(\theta_1- \theta_{\gamma)} }=1,\\
 \epsilon\rho= e^{ i\theta_1},\\
 \epsilon=\overline{\epsilon},\\
\rho=\overline{\rho}.\\
 \end{array}\right.
 \end{equation}
\bfl 
 Note that, the fixed point equations (\ref{C2 : 40})  describe exactly
 the non effective subgroup. In fact, under the assumption $\sigma\equiv 0$,
 the equations (\ref{C2 : 40}) are equivalent to those described at the end
 of lemma$\, 2.1$,
 ( when all of the relations in (\ref{C2 : 18}) are verified).
  Then we get the conclusion. $\square$  \efl

\vspace{0,5 cm}

\begin{cor}
\bfl Let $\zw$ be a point on a singular $\tilde{G}-orbit$ of $N(\Theta)$. Suppose
$\zw$ has a quaternionic pair $\uu$ which is contained in one
of ${}^{(\pm,\pm)}V_3^\alpha$,
then all of the four quaternionic pairs of $\zw$ are contained in one of the spaces
${}^{(\pm ,\pm )}V_3^\alpha$,
 $\alpha\in\{1,2,3,4\}$.   \efl
\end{cor}

\bfl The proof of this corollary is a straightforward consequence of
proposition $2.4$ and
corollary $2.2$.$\square$ \efl

\vspace{0,5 cm}

\bfl  $\mathbf{Remark\,\,2.3.}$
 Lemma $2.2$ allows to describe all of the possible
strata for the action of $\tilde{G}$ on $\mathbb{C}^8\times\mathbb{C}^8$.
 In fact, each point $\zw\in\mathbb{C}^8\times\mathbb{C}^8\cong \mathbb{H}^8$
 can be thought as  a
 $4\times 8$  complex matrix, then if $\zw$ has no trivial isotropy, its own
 quaternionic pairs are described by the eigenspaces listed in (\ref{C2 : 30})
  and (\ref{C2 : 30b}).
  Moreover,
  the singular $\tilde{G}-$stata of $N(\Theta)$ can be obtained
 by intersecting $N(\Theta)$ with the $\tilde{G}-$strata of
  $\mathbb{H}^8$ having discrete isotropy.
In particular,
Corollary $2.3$ allows to distinguish two different families of singular points,
the first
one  is given by points $\zw$ such that all of their quaternionc pairs are contained
in ${}^{(\pm ,\pm )}V_3^\alpha$, $\alpha\in\{1,2,3,4\}$.  Instead, in the second family,
none of the points
has a quaternionic pair contained in ${}^{(\pm ,\pm )}V_3^\alpha$.
In the following paragraph we consider the case of points $(\underline{z}, \underline{w})\in N(\Theta)$
whose quaternionic pairs
 are contained in ${}^{(\pm,\pm)}V_3^\alpha$. $\square$  \efl

 \vspace{0,5 cm}

 \subsection{  $\tilde{G}-$strata $\tilde{S}_{\alpha}$ characterized by
  quaternionic pairs $\uu$
 of type ${}^{(\pm ,\pm )}V_3^\alpha$}

 \vspace{0,5 cm}

 \bfl By using the
 expression of the
elements in ${}^{(\pm,\pm)}V^{\alpha}_3$, it follows that each point in the
first mentioned
family has the following shape\efl
\begin{equation}   \label{C2 : 41}
{}^T(\underline{z}, \underline{w} ):=
\left(\begin{array}{cc|cc|cc|cc}
z_1 & \pm iz_1  & z_3 & \pm iz_3 & z_5 & \pm iz_5 & z_7 & \pm iz_7\\
w_1 & \pm iw_1 &  w_3 & \pm iw_3 & w_5 & \pm iw_5 & w_7 & \pm iw_7 \\
\end{array}\right).
\end{equation}
\bfl Then the fixed point equations become\efl
\begin{equation}       \label{C2 : 42}
\left\{ \begin{array}{l}
e^{i(\theta_1\pm\theta_\alpha)}=1,\\
   \epsilon\rho = e^{\pm i\theta_1},\\
\overline{\epsilon}\rho = e^{\pm i\theta_1}, \\
\end{array} \right.
\end{equation}
\bfl where $\alpha\in\{1,2,3,4\}.$
Now, by fixing a pair of
 signs $(\pm, \pm)$  for the equations
  $\left\{\begin{array}{l}
\epsilon\rho = e^{\pm i\theta_1}, \\
\overline{\epsilon}\rho = e^{\pm i\theta_1}, \\
\end{array}\right.$
 and a  triple of signs
$(\mp, \mp, \mp)$,
relative to    $e^{i\theta_1}= e^{\pm  i\theta_{\alpha}}$, $\alpha\in\{2,3,4\}$,
 (which  determines exactly one of the
the determinats   $\sq$ ) we have that each eigenspaces associated to the
 eingenvalues problem
 equations in (\ref{C2 : 9}) is univocally determinated  by a  $5-$ple of
  signs $  \big( (\pm,\pm ),(\pm,\pm,\pm ) \big)$.  \efl

 \vspace{0,5 cm}

\begin{defin} \bfl In the following, each eigenspace will be
  indicated with its own $5-ple $ of signs
  $\big( (\pm,\pm ),(\pm,\pm,\pm ) \big)$, that is
  $\overset{=}{S}{}^{(\pm,\pm )}_{(\pm,\pm,\pm )}$.\efl
\end{defin}

  \vspace{0.5 cm}

\bfl $\mathbf{Remark\,\, 2.4}.$
  Let $M_{\alpha}$ be such that $rank =2$ and suppose
  $\uu\in {}^{(\pm, \pm)}V^{\alpha}_3$. Then, the systems of equations
   $3)$ in $(\ref{C2 : 32})$, up to a sign, are either represented by \efl
\begin{equation}       \label{C2 : 43}
\begin{aligned}
 A) &\hspace{ 1,5 cm}  \left\{ \begin{array}{l}
\rho \epsilon = e^{i\theta_i},   \\
\rho \overline{\epsilon}= e^{-i\theta_i},   \\
\end{array} \right. \quad \iff \quad
\left\{ \begin{array}{l}
\rho = \pm 1,            \\
\epsilon=\pm e^{-i\theta},   \\
\end{array} \right. \\
\mathrm{or\,\,by} \hspace{0,5 cm} & \\
 B)& \hspace{1,5 cm} \left\{ \begin{array}{l}
\rho \epsilon = e^{i\theta_i},       \\
\rho \overline{\epsilon}=e^{i\theta},  \\
\end{array} \right. \quad \iff \quad
\left\{ \begin{array}{l}
\rho = \pm e^{-i\theta},             \\
\epsilon= \pm 1.    \\
\end{array} \right.              \hspace{7 cm}   {}
\end{aligned}
\end{equation}

\bfl Now, evaluating the matrices $M_{\alpha}$ on the parameters which come
from the case $A)$, as space of solutions we get either \efl

\begin{equation}    \label{C2 : 43a}
\begin{aligned}
&{}^{(+,- )}V_3^{\alpha}:=\{(z_{2\alpha-1},w_{2\alpha-1})\in \mathbb{C}^2\, |
 \, (z_{2\alpha-1},  iz_{2\alpha-1}, z_{2\alpha-1},-iw_{2\alpha-1} ) \},  \\
 \mathrm{ or } \hspace{1 cm} & \\
&{}^{(-,+ )}V_3^{\alpha}:=\{(z_{2\alpha-1},w_{2\alpha-1})\in \mathbb{C}^2\, |
 \, (z_{2\alpha-1},  -iz_{2\alpha-1}, z_{2\alpha-1},iw_{2\alpha-1} ) \}.  \hspace{7 cm}\\
\end{aligned}
\end{equation}

\bfl Instead in the case $B)$ we either have \efl

\begin{equation}   \label{C2 : 43b}
 \begin{aligned}
&{}^{(+,+ )}V_3^{\alpha}:=\{(z_{2\alpha-1},w_{2\alpha-1})\in \mathbb{C}^2\, |
 \, (z_{2\alpha-1},  iz_{2\alpha-1}, z_{2\alpha-1},iw_{2\alpha-1} ) \},  \\
 \mathrm{ or} \hspace{1 cm} & \\
&{}^{(-,- )}V_3^{\alpha}:=\{(z_{2\alpha-1},w_{2\alpha-1})\in \mathbb{C}^2\, |
 \, (z_{2\alpha-1},  -iz_{2\alpha-1}, z_{2\alpha-1},-iw_{2\alpha-1} ) \}.\, \square  \hspace{7 cm}\\
\end{aligned}
\end{equation}

\vspace{0.5 cm }

 \begin{lem}
 \bfl The singular points $(\underline{z}, \underline{w})\in
 \mathbb{C}^8\times
  \mathbb{C}^8$  which solve the eigenvalues problem $(\ref{C2 : 9})$ and such that
  all of their quaternionic pair are contained in ${}^{(\pm , \pm )}V_3^{\alpha}$,
  can be divided into two
  subfamilies.
 The first family only depends on the case described in  $A)$ at $(\ref{C2 : 43})$
  and the second one
 on the case $B)$.\efl
 \end{lem}
\bfl $\mathbf{Proof.}$
Let $(\underline{z}, \underline{w})$ be a point in $N(\Theta)$
 of the type above described. Note that it is impossible that there
 exist two different indices
$\alpha,\beta\in \{1,2,3,4,\}$ such that $M_{\alpha}$ satisfies the condition $A)$
and $M_{\beta}$ the condition $B)$. Otherwise, up to switch the index signs,
we would have
 \efl
 \begin{equation}   \label{C2 : 44}
 \left\{ \begin{array}{l}
\rho \epsilon = e^{i\theta_{\alpha}}, \\
\rho \overline{\epsilon} = e^{-i\theta_{\alpha}}\\
\rho \epsilon =e^{i\theta_{\beta}},   \\
\rho \overline{\epsilon }=e^{i\theta_{\beta}},   \\
\end{array} \right.\quad  \,\, \iff \quad \,\,
\left\{\begin{array}{l}
e^{i\theta_{\alpha}-\theta_{\beta}} = 1, \\
\rho\epsilon =  e^{i\theta_{\alpha}}, \\
\epsilon = \overline{\epsilon}, \\
\rho= \overline{\rho}. \\
\end{array}\right.
\end{equation}
\bfl Now, the argument in corollary $2.9$
gives that the
  the fixed point equations describe the non effective
subgroup,
against our hypotheses.   Then all the quaternionic pair have to be of the same type.
 That is, either
 all of the quaternionic pair are linked to the case $A)$ or all of them
 are linked to the case $B)$.
So, referring to the cases $A)$ and $B)$ we have that the
only possible points which solve the
fixed point equations in (\ref{C2 : 9})
are listed in the following two subfamilies. First we consider \efl

\begin{equation}   \label{C2 : 45}
\begin{aligned}
& 1)\,\,\, \overset{=}{S}{}^{(+,+)}_{(+,+,-)}:= \Bigg\{ \left(\begin{array}{cc|cc|cc|cc}
z_1 & -iz_1  & z_3 & iz_3 & z_5 &  iz_5 & z_7 &  -iz_7\\
w_1 &  -iw_1 &  w_3 &  iw_3 & w_5 &  iw_5 & w_7 &  -iw_7 \\
\end{array}\right)\Bigg\},  \hspace{10 cm} {}  \\
&2)\,\,\,  \overset{=}{S}{}^{(+,+)}_{(+,-,+)}:= \Bigg\{ \left(\begin{array}{cc|cc|cc|cc}
z_1 & -iz_1  & z_3 & iz_3 & z_5 &  -iz_5 & z_7 &  iz_7\\
w_1 &  -iw_1 &  w_3 &  iw_3 & w_5 &  -iw_5 & w_7 &  iw_7 \\
\end{array}\right)\Bigg\}, \\
&3)\,\,\, \overset{=}{S}{}^{(+,+)}_{(-,+,+)}:= \Bigg\{ \left(\begin{array}{cc|cc|cc|cc}
z_1 & -iz_1  & z_3 & -iz_3 & z_5 &  iz_5 & z_7 &  iz_7\\
w_1 &  -iw_1 &  w_3 &  -iw_3 & w_5 &  iw_5 & w_7 &  iw_7 \\
\end{array}\right)\Bigg\},  \\
& 4)\,\,\,  \overset{=}{S}{}^{(+,+)}_{(+,+,+)}\big):= \Bigg\{ \left(\begin{array}{cc|cc|cc|cc}
z_1 & -iz_1  & z_3 & iz_3 & z_5 &  iz_5 & z_7 &  iz_7\\
w_1 &  -iw_1 &  w_3 &  iw_3 & w_5 &  iw_5 & w_7 &  iw_7 \\
\end{array}\right)\Bigg\}, \\
&5)\,\,\,  \overset{=}{S}{}^{(+,+)}_{(+,-,-)}:= \Bigg\{ \left(\begin{array}{cc|cc|cc|cc}
z_1 & -iz_1  & z_3 & iz_3 & z_5 &  -iz_5 & z_7 &  -iz_7\\
w_1 &  -iw_1 &  w_3 &  iw_3 & w_5 &  -iw_5 & w_7 &  -iw_7 \\
\end{array}\right)\Bigg\},   \\
&6)\,\,\,  \overset{=}{S}{}^{(+,+)}_{(-,+,-)}\big):=\Bigg\{ \left(\begin{array}{cc|cc|cc|cc}
z_1 & -iz_1  & z_3 & -iz_3 & z_5 &  iz_5 & z_7 &  -iz_7\\
w_1 &  -iw_1 &  w_3 &  -iw_3 & w_5 &  iw_5 & w_7 &  -iw_7 \\
\end{array}\right)\Bigg\},  \\
&7)\,\,\,  \overset{=}{S}{}^{(+,+)}_{(-,-,+)}:= \Bigg\{ \left(\begin{array}{cc|cc|cc|cc}
z_1 & -iz_1  & z_3 & -iz_3 & z_5 &  -iz_5 & z_7 &  iz_7\\
w_1 &  -iw_1 &  w_3 &  -iw_3 & w_5 &  -iw_5 & w_7 &  iw_7 \\
\end{array}\right)\Bigg\},    \hspace{10 cm }  \\
\end{aligned}
\end{equation}

\bfl and note that the cases $ \overset{=}{S}{}^{(-,-)}_{ (\pm,\pm,\pm)} $
 are completely similar.
The second subfamily to consider is \efl

\begin{equation}    \label{C2 : 46}
\begin{aligned}
& 8)\,\,\, \overset{=}{S}{}^{(+,-)}_{(+,+,+)} :=\Bigg\{ \left(\begin{array}{cc|cc|cc|cc}
z_1 & -iz_1  & z_3 & iz_3 & z_5 &  iz_5 & z_7 &  iz_7\\
w_1 &  iw_1 &  w_3 & - iw_3 & w_5 &  -iw_5 & w_7 &  -iw_7 \\
\end{array}\right)\Bigg\},  \hspace{10 cm }\\
& 9)\,\,\, \overset{=}{S}{}^{(+,-)}_{(+,+,-)}:=  \Bigg\{ \left(\begin{array}{cc|cc|cc|cc}
z_1 & -iz_1  & z_3 & iz_3 & z_5 &  iz_5 & z_7 &  -iz_7\\
w_1 &  iw_1 &  w_3 & - iw_3 & w_5 &  -iw_5 & w_7 &  iw_7 \\
\end{array}\right)\Bigg\}, \hspace{10 cm }\\
&10)\,\,\, \overset{=}{S}{}^{(+,-)}_{(+,-,+)}  :=  \Bigg\{ \left(\begin{array}{cc|cc|cc|cc}
z_1 & -iz_1  & z_3 & iz_3 & z_5 &  -iz_5 & z_7 &  iz_7\\
w_1 & iw_1 &  w_3 &  -iw_3 & w_5 &  iw_5 & w_7 &  -iw_7 \\
\end{array}\right)\Bigg\}, \\
&11)\,\,\, \overset{=}{S}{}^{(+,-)}_{(-,+,+)}  := \Bigg\{ \left(\begin{array}{cc|cc|cc|cc}
z_1 & -iz_1  & z_3 & -iz_3 & z_5 &  iz_5 & z_7 &  iz_7\\
w_1 &  iw_1 &  w_3 &  iw_3 & w_5 &  -iw_5 & w_7 &  -iw_7 \\
\end{array}\right)\Bigg\},  \\
&12)\,\,\, \overset{=}{S}{}^{(+,-)}_{(+,-,-)}  :=  \Bigg\{ \left(\begin{array}{cc|cc|cc|cc}
z_1 & -iz_1  & z_3 & iz_3 & z_5 &  -iz_5 & z_7 &  -iz_7\\
w_1 & iw_1 &  w_3 & -iw_3 & w_5 &  iw_5 & w_7 &  iw_7 \\
\end{array}\right)\Bigg\}, \hspace{10 cm}   \\
&13)\,\,\, \overset{=}{S}{}^{(+,-)}_{(-,+,-)} \big)  := \Bigg\{ \left(\begin{array}{cc|cc|cc|cc}
z_1 & -iz_1  & z_3 & -iz_3 & z_5 &  iz_5 & z_7 &  -iz_7\\
w_1 &  iw_1 &  w_3 &  iw_3 & w_5 &  -iw_5 & w_7 &  iw_7 \\
\end{array}\right)\Bigg\},  \\
&14)\,\,\, \overset{=}{S}{}^{(+,-)}_{(-,-,+)}  := \Bigg\{ \left(\begin{array}{cc|cc|cc|cc}
z_1 & -iz_1  & z_3 & -iz_3 & z_5 &  -iz_5 & z_7 &  iz_7\\
w_1 &  iw_1 &  w_3 &  iw_3 & w_5 &  iw_5 & w_7 &  -iw_7 \\
\end{array}\right)\Bigg\}, \\
&15)\,\,\, \overset{=}{S}{}^{(+,-)}_{(-,-,-)}   := \Bigg\{ \left(\begin{array}{cc|cc|cc|cc}
z_1 & -iz_1  & z_3 & -iz_3 & z_5 &  -iz_5 & z_7 &  iz_7\\
w_1 &  iw_1 &  w_3 &  iw_3 & w_5 &  iw_5 & w_7 &  -w_7 \\
\end{array}\right)\Bigg\}.   \hspace{10 cm} {}
\end{aligned}
\end{equation}
 \bfl The last  cases $\overset{=}{S}{}^{(-,+)}_{(\pm,\pm,\pm)} $
  have a completely similar description. Note the singular strata reltive to the submanifold $N(\Theta)$
  come from the intersection
  of the above eigenspaces with the manifold $N(\Theta)$. $\square$ \efl

\begin{cor}
\bfl The cases
 of $\overset{=}{S}{}^{(-,-)}_{(+,+,+)}$ and $\overset{=}{S}{}^{(+,+)}_{(-,-,-)}$
 give an empty intersection with the zero set $N(\Theta)$.  \efl
\end{cor}

\bfl $\mathbf{Proof.}$
Consider  the case of $\overset{=}{S}{}^{(-,-)}_{(+,+,+)}$
 (the other case can be treated in the same way ).
An element of this set can be
 thought as\efl
\begin{equation}  \label{C2 : 46a}
{}^T(\underline{z}, \underline{w})= \left(\begin{array}{cc|cc|cc|cc}
z_1 & -iz_1  & z_3 & -iz_3 & z_5 &  -iz_5 & z_7 &  -iz_7\\
w_1 &  -iw_1 &  w_3 &  -iw_3 & w_5 &  -iw_5 & w_7 &  -iw_7 \\
\end{array}\right).
\end{equation}
\bfl The moment map, relative to the $3$-torus $T_{\Theta}^3$, in this case gives the
  equations \efl
\begin{equation}      \label{C2 : 47}
\left\{\begin{array}{l}
\sum_{\alpha=1}^4 d_{\alpha} (|z_{2\alpha-1}|^2+ |w_{2\alpha-1}|^2)=0,  \\
\sum_{\alpha=1}^4 d_{\alpha} (z_{2\alpha-1} w_{2\alpha-1}-
z_{2\alpha-1}w_{2\alpha-1})=0, \\
\end{array}\right.
\end{equation}
\bfl where $d_{\alpha}=p_{\alpha},q_{\alpha}$ and $l_{\alpha}$.
  Then it follows that equations in  (\ref{C2 : 47})
 admit only the trivial solution $(\underline{z},
 \underline{w})=(\underline{0},\underline{0})$.   Thus we get the conclusion.$\square$
\efl

\vspace{0,5 cm}

\bfl  $\mathbf{Remark\,\,2.5}$  Note that, each of the eigenspaces listed in
(\ref{C2 : 45}) is $\tilde{G}-$invariant. Then, these eigenspaces actually describe a
 family of $\tilde{G}$-strata in $\mathbb{H}^8$. Instead, the eigenspaces
  in (\ref{C2 : 46})
 are not $\tilde{G}-$invariant, and in the next proposition $2.5$ we will prove that
 they do not intersect $N(\Theta)$.$\square$ \efl

   \vspace{0.5 cm}

\begin{prop}   The families of eigenspaces
$\overset{=}{S}{}^{(+,-)}_{(\pm, \pm, \pm)} $
  and $\overset{=}{S}{}^{(-,+)}_{(\pm, \pm, \pm)} $
do not intersect the manifold $N(\Theta)=\mu^{-1}(0)\cap \nu^{-1}(0)$.
\end{prop}

\bfl $\mathbf{Proof}$  Consider the  family of solutions
 $\overset{=}{S}{}^{(+,-)}_{(\pm, \pm, \pm)}$ (the proof
   in the case of the remaing sets $\overset{=}{S}{}^{(-,+)}_{(\pm, \pm, \pm)}$
is similar).
 Let us recall that the equations provided by the $Sp(1)$ moment map $\mu$
  are the following\efl

  \begin{equation}    \label{C2 : 48}
\left\{\begin{array}{l}
\sum_{\alpha=1}^4 \overline{(z_{\alpha} + jw_{\alpha})}i (z_{\alpha} + jw_{\alpha})=0, \\
\sum_{\alpha=1}^4 \overline{(z_{\alpha} + jw_{\alpha})}j (z_{\alpha} + jw_{\alpha})=0, \\
 \sum_{\alpha=1}^4 \overline{(z_{\alpha} + jw_{\alpha})}k (z_{\alpha} + jw_{\alpha})=0. \\
\end{array}\right.
\end{equation}

\bfl By reading these  equations on the points $\zw$ in
 $\overset{=}{S}{}^{(-,+)}_{(\pm, \pm, \pm)}$, we obtain\efl

\begin{equation}      \label{C2 : 49}
\left\{\begin{array}{l}
\sum_{\alpha=1}^4  (|z_{2{\alpha}-1}|^2 -|w_{2\alpha-1}|^2)=0 ,\\
\sum_{\alpha=1}^4  z_{2\alpha-1}w_{2\alpha-1}=0. \\
\end{array}\right.
\end{equation}

\bfl In fact, for each point $\zw$ which belongs to one of the spaces
$\overset{=}{S}{}^{(+,-)}_{(\pm, \pm, \pm)}$, we have that
 the last two equations in (\ref{C2 : 48}) are automatically  satisfied. Instead,
  the first of these equations
yields the system of equations in (\ref{C2 : 49}). In particular,
the equations in (\ref{C2 : 49}) are the same
 for all the eigenspaces listed in the second subfamily.
Now consider the moment map $\nu$, relative to the $3$-torus $T^3_{\Theta}$.
 In the following
 $d_{\alpha}$ indicate the weights  $p_{\alpha}, q_{\alpha}$ and $l_{\alpha}$,
for $\alpha\in\{1,2,3,4\}$.
Then, the moment map
$\nu$ yields the following equations \efl

\begin{equation}            \label{C2 : 50}
\left\{\begin{array}{l}
\sum_{\alpha=1}^4 d_{\alpha} Im(z_{2\alpha-1}\overline{z}_{2\alpha} +
w_{2\alpha-1}\overline{w}_{2\alpha})=0, \\
\sum_{\alpha=1}^4 d_{\alpha} (z_{2\alpha-1}w_{2\alpha}-  z_{2\alpha} w_{2\alpha-1} ) =0,
\end{array}\right.
\end{equation}

\bfl By evaluating these equations on a point
$(\underline{z}, \underline{w} )$ which belongs to one of the spaces
 $\overset{=}{S}{}^{(+,-)}_{(\pm, \pm, \pm)}$,
we obtain \efl

\begin{equation}         \label{C2 : 51}
\left\{\begin{array}{l}
\sum_{\alpha=1}^4  (-1)^{m_{\alpha}} d_{\alpha}(|z_{2\alpha-1}|^2 -|w_{2\alpha-1}|^2)=0 ,\\
\sum_{\alpha=1}^4 (-1)^{m_{\alpha}} d_{\alpha} z_{2\alpha-1}w_{2\alpha-1}=0, \\
\end{array}\right.
\end{equation}

\bfl the indices $m_{\alpha}$ only depend on the $5-ple$ of signs of the chosen eigenspace.
Let us rewrite these systems of equations as follows, \efl

\begin{equation}         \label{C2 : 52}
 \left\{\begin{array}{l}
\sum_{\alpha=1}^4 (-1)^{m_{\alpha}} p_{\alpha} (|z_{2{\alpha}-1}|^2 -|w_{2{\alpha}-1}|^2)=0 ,\\
\sum_{\alpha=1}^4 (-1)^{m_{\alpha}} q_{\alpha} (|z_{2{\alpha}-1}|^2 -|w_{2{\alpha}-1}|^2)=0 ,\\
\sum_{\alpha=1}^4 (-1)^{m_{\alpha}} l_i (|z_{2\alpha-1}|^2 -|w_{2{\alpha}-1}|^2)=0 ,\\
 \sum_{\alpha=1}^4  (|z_{2\alpha-1}|^2 -|w_{2\alpha-1}|^2)=0 ,\\
\end{array}\right. \,\,  \cup \,\,
\left\{\begin{array}{l}
 \sum_{\alpha=1}^4 (-1)^{m_{\alpha}} p_{\alpha} z_{2\alpha-1}w_{2\alpha-1}=0,\\
\sum_{\alpha=1}^4 (-1)^{m_{\alpha}} q_{\alpha} z_{2\alpha-1}w_{2\alpha-1}=0,\\
\sum_{\alpha=1}^4 (-1)^{m_{\alpha}}l_{\alpha} z_{2\alpha-1}w_{2\alpha-1}=0, \\
\sum_{\alpha=1}^4 z_{2\alpha-1}w_{2\alpha-1}=0.
\end{array}\right.   \hspace{10 cm} {}
\end{equation}

\bfl These two systems  have the same matrix
of coefficients, and hence
 the same determinant.
Note that, this determinant
 is one of  the
 $\sq$, which never vanish under our assumptions.
  Then, from the first system of equations
 in (\ref{C2 : 52}) we get\efl

 \begin{equation}     \label{C2 : 53}
 |z_{2\alpha-1}|^2 =|w_{2\alpha-1}|^2>0,
 \end{equation}

 \bfl for each $\alpha\in\{1,2,3,4\}$. The condition they have to be bigger than zero
 depends on the fact that each quaternionic pair does not vanish. Similarly
  from the second system of equation in (\ref{C2 : 52}) we get the further conditions\efl

 \begin{equation}       \label{C2 : 54}
 z_{2\alpha-1}w_{2\alpha-1}=0,
 \end{equation}

 \bfl for each $\alpha\in\{1,2,3,4\}$. Then it holds either
$z_{2\alpha-1}=0$ or $w_{2\alpha-1}=0 $, $\alpha\in\{1,2,3,4\}$. So the only solution
for these two system of equations is
$ (\underline{z}, \underline{w})= (\underline{0}, \underline{0})$
 which does not belong to $N(\Theta)$. Then, we do not care about the eigenspaces in
 (\ref{C2 : 46}).$\square$ \efl

 \vspace{0.5 cm}

\begin{prop}
\bfl Consider the families of complex vector spaces
 $\overset{=}{S}{}^{(+,+)}_{(\pm, \pm, \pm)}$ and
  $\overset{=}{S}{}^{(-,-)}_{(\pm, \pm, \pm)}$.
 Then, it follows that just one of these complex spaces intersects the zero set
   $N(\Theta) =\mu^{-1}(0)\cap \nu^{-1}(0)$. \efl
\end{prop}

\bfl $\mathbf{Proof.}$ Consider the the eigenspaces
in the first family. By reading
the $Sp(1)-$
moment map equations in (\ref{C2 : 48}) on any point
$ (\underline{z}, \underline{w}),
 $ which belongs to one of the eigenspaces
 $\overset{=}{S}{}^{(+,+)}_{(\pm, \pm, \pm)}$
 we get the following system of equations\efl

 \begin{equation}   \label{C2 : 55}
\left\{\begin{array}{l}
\sum_{\alpha=1}^4  (|z_{2\alpha-1}|^2 -|w_{2\alpha-1}|^2)=0 ,\\
\sum_{\alpha=1}^4  \overline{w}_{2\alpha-1}z_{2\alpha-1}=0, \\
\end{array}\right.
 \end{equation}

\bfl Note that the equations in (\ref{C2 : 55})
describe the Stiefel manifold $\mathcal{S}=U(4)/U(2)$ of the $2$-complex
 frames in $\mathbb{C}^4$. Besides these equations are common for all
of the eigenspaces $\overset{=}{S}{}^{(+,+)}_{(\pm, \pm, \pm)}$.
 Instead  by evaluating
the $T^3_{\Theta} $ moment map equations  on a generic point
 $ (\underline{z}, \underline{w})$ which belongs
 to one of the spaces we are dealing with,
 we obtain\efl

 \begin{equation}     \label{C2 : 56}
\sum_{\alpha=1}^4  (-1)^{m_{\alpha}}d_{\alpha}(|z_{2\alpha-1}|^2 +|w_{2\alpha-1}|^2)=0,
 \end{equation}

 \bfl  for each $d_{\alpha} =p_{\alpha}, q_{\alpha}, l_{\alpha}$ and $\alpha\in\{1,2,3,4\}$.
  Let us rewrite
 the equations in (\ref{C2 : 56}) as\efl

 \begin{equation}    \label{C2 : 57}
 \sum_{\alpha=1}^4  (-1)^{m_{\alpha}}d_{\alpha} |u_{2\alpha-1}|^2=0,
 \end{equation}

 \bfl where $u_{2\alpha-1} =z_{2\alpha-1}+jw_{2\alpha-1}$ for each $\alpha\in\{1,2,3,4\}$.
  Now, by using the notation we have introduced in (\ref{C2 : 57}),
   the intersection of the  sphere $S^{31}$ with each  eigenspace in
$\overset{=}{S}{}^{(+,+)}_{(\pm, \pm, \pm)}$ is described by the equation  \efl

\begin{equation}       \label{C2 : 58}
\sum_{\alpha=1}^8  |u_{\alpha}|^2=\sum_{\alpha=1}^4   2|u_{2\alpha-1}|^2=1.
\end{equation}
\bfl By considering together the equations in (\ref{C2 : 57}) and (\ref{C2 : 58})
 we get the following
system of equations\efl

\begin{equation}         \label{C2 : 59}
\left(\begin{array}{cccc}
(-1)^{m_1}p_1  & (-1)^{m_2}p_2 & (-1)^{m_3} p_3 & (-1)^{m_4} p_4 \\
(-1)^{m_1}q_1  & (-1)^{m_2}q_2 & (-1)^{m_3} q_3 & (-1)^{m_4} q_4 \\
(-1)^{m_1}l_1  & (-1)^{m_2} l_2 & (-1)^{m_3}l_3 & (-1)^{m_4} l_4 \\
 1   &  1      &   1     &   1     \\
\end{array} \right)
\left(\begin{array}{c}
  |u_{1}|^2 \\
   |u_{3}|^2 \\
    |u_{5}|^2\\
     |u_{7}|^2  \\
\end{array}\right) =
\left(\begin{array}{c}
  0 \\
  0  \\
  0\\
  \frac{1}{2} \\
\end{array}\right).
\end{equation}

\bfl  then we obtain\efl

\begin{equation}         \label{C2 : 60}
\left\{ \begin{array}{l}
2|u_{1}|^2 =\frac{\pm \Delta_{234} }{ \sq}>0,\\
2|u_{3}|^2 =\frac{\pm \Delta_{134}}{ \sq}>0, \\
2|u_{5}|^2 =\frac{\pm \Delta_{124}}{ \sq}>0, \\
2|u_{7}|^2  =\frac{\pm \Delta_{123}}{ \sq}>0.\\
\end{array}\right.
\end{equation}

\bfl Note that, after we have fixed the weigth matrix $\Theta$, the minor determinants
$\Delta_{\alpha\beta\gamma}$ are univocally determinated.
Then the system of
equations in (\ref{C2 : 60}) admits a unique solution in terms of $|u_{\alpha}|^2$, and the relations
 (\ref{C1: 26})) give that
$N(\Theta)$ intersects exactly one of
 the eigenspaces $\overset{=}{S}{}^{(+,+)}_{(\pm, \pm, \pm)}$.
 As a consequence, when $\overset{=}{S}{}^{(+,+)}_{(\pm, \pm, \pm)}\cap N(\Theta)\ne \emptyset$,
 it follows that $dim\, \big(\overset{=}{S}{}^{(+,+)}_{(\pm, \pm, \pm)}\cap N(\Theta)\big)=$
 $ 9$
 and $\big(\overset{=}{S}{}^{(+,+)}_{(\pm, \pm, \pm)}\cap N(\Theta)\big)/ \tilde{G}\cong S^2$. $\square$

 \vspace{1 cm}

$\mathbf{Example\,\, 2.1.}$
In order to
have a better description of Proposition $2.5$, as an example
we consider the case when  the moment map
 equations in (\ref{C2 : 60}) become\efl

\begin{equation}            \label{C2 : 61}
\left(\begin{array}{cccc}
p_1  & p_2 & - p_3 & - p_4 \\
q_1  & q_2 & - q_3 & - q_4 \\
l_1  & l_2 & -l_3 & - l_4 \\
 1   &  1      &   1     &   1     \\
\end{array} \right)
\left(\begin{array}{c}
  |u_{1}|^2 \\
   |u_{3}|^2 \\
    |u_{5}|^2\\
     |u_{7}|^2  \\
\end{array}\right) =
\left(\begin{array}{c}
  0 \\
  0  \\
  0\\
  \frac{1}{2} \\
\end{array}\right),
\end{equation}
\bfl and it follows \efl

\begin{equation}      \label{C2 : 62}
\left\{ \begin{array}{l}
2|u_{1}|^2 =\frac{- \Delta_{234} }{ \sqmpp}>0,\\
2|u_{3}|^2 =\frac{ \Delta_{134}}{ \sqmpp}>0, \\
2|u_{5}|^2 =\frac{ \Delta_{124}}{ \sqmpp}>0, \\
2|u_{7}|^2  =\frac{-\Delta_{123}}{ \sqmpp}>0.\\
\end{array}\right.
\end{equation}

\bfl Suppose this system of equations admit solutions. Hence, we would have that\efl
\begin{equation}     \label{C2 : 62e}
 \frac{\Delta_{234}}{|\Delta_{124}|}=\frac{\Delta_{134}}{|\Delta_{134}|} \quad
 \mathrm{and} \quad  \frac{\Delta_{123}}{|\Delta_{123}|}=\frac{\Delta_{234}}{|\Delta_{234}|}.
 \end{equation}
 \bfl Then,
 according to the sign of $ \sqmpp$,
 the identities (\ref{C2 : 62}) ( which involved the determinants
   $\Delta_{\alpha \beta \gamma}$ and
 $\sq$ as showed in (\ref{C1: 26})) admit exactly one solution
 in terms of $\Delta_{\alpha \beta \gamma}$ and
 $\sq$. Then, the only eigenspace which have a non null
 intersection with $N(\Theta)$ is $\overset{=}{S}{}^{(+,+)}_{(+,+,-)}$.$\square$\efl


\vspace{0,5 cm}

 \subsection{ $\tilde{G}-$strata $\tilde{S}_{\alpha}$ characterized by quaternionic pairs $\uu$
 of type  ${}^{\pm}V_1^{\alpha}, \,\, {}^{\pm}V_2^{\alpha},\,\, V_4^{\alpha}$ and
 $V_5^{\alpha}$}

\bfl  Consider now the  points $\zw\in \mathbb{H}^8$
 which have no
 quaternionic pairs of type ${}^{(\pm , \pm)}V_3^{\alpha}$. First, note that
 the complex vector spaces  ${}^{\pm}V_1^{\alpha}, \,\, {}^{\pm}V_2^{\alpha},\,\, V_4^{\alpha}$ and
 $V_5^{\alpha},\,\, \alpha\in$
 $\{1,$
 $2,3,4\}$ are not invariant for the  action
 of the group $\tilde{G}$. Then, a singular stratum $\overset{=}{S}\subset$
 $\mathbb{H}^8$ whose points have
  quaternionic coordinates   $\uu$
  in one of these spaces 
 can be described as  \efl

 \begin{equation}    \label{C2 : 63}
 \tilde{G}\cdot  \overbrace{\Bigg\{ \left(\begin{array}{cc|cc|cc|cc}
u_1 & u_2  & u_3 & u_4 & u_5 &  u_6 & u_7 &  u_8\\
\end{array}\right)\Bigg\}}^{V:=}.
 \end{equation}

 \bfl Here, $V$ is a vector space whose quaternionic pairs are in one of the spaces ${}^{\pm}V_1^{\alpha},$
  ${}^{\pm}V_2^{\alpha},$ $V_4^{\alpha}$ and $V_5^{\alpha} $.   \efl

  \begin{lem}
  Let $\overset{=}{S}=\tilde{G}\cdot V$ be a singular stratum of those in $(\ref{C2 : 63})$. Then
  $dim\, \tilde{G}\cdot V= dim\, V + 2$.
  \end{lem}

\bfl $\mathbf{Proof.}$ Note that $\tilde{H}:= Sp(1)\times U(1)^{\epsilon}\times U(1)$
 is the biggest subgroup
 of $\tilde{G}$ that fixes each vector space $V$. Consider    \efl

 \begin{equation}    \label{C2 : 63a}
\tilde{G}\times_{\tilde{H}} V:=\big\{ [(g, \uno)] \,\, | \,\,
 g\in \tilde{G} \,\, \mathrm{and} \,\,  \uno\in V \big\},
 \end{equation}

\bfl where  $ [(g, \uno)] =  [(g_1, \uno_1)] $ if and only if
  there exists $h\in \tilde{H}$  such that $(g_1, \uno_1)= (gh^{-1},$
   $ h\cdot\uno)$.
 We have that, $\tilde{G}\times_{\tilde{H}} V$   is a $V$-vector bundle over
  $\tilde{G}/\tilde{H}$, and
 we can define the $\tilde{G}-$action $\Lambda\, : \tilde{G}\times \big(\tilde{G}\times_{\tilde{H}} V \big)
 \rightarrow \tilde{G}\times_{\tilde{H}} V $ , given  by
\begin{equation}    \label{C2a : 63b}
 g'\cdot [(g, \uno)]:= [(g'g, \uno)].
 \end{equation}
 Then, by using the $\tilde{G}-$equivariant diffeomorphism \efl
\begin{equation}      \label{C2 : 63b}
\begin{aligned}
 \Gamma\, :\tilde{G}& \cdot V\longrightarrow   \tilde{G}\times_{\tilde{H}} V, \\
& g\cdot \uno     \longmapsto   [(g, \uno)] \\
\end{aligned}
\end{equation}
 \bfl it follows that $\tilde{G}\cdot V\cong
   \tilde{G}\times_{\tilde{H}} V$. Since $\tilde{G}/\tilde{H}\cong S^2$, we have that
$\tilde{G}\cdot V$ is a vector bundle over $S^2$ with fiber $V$:  \efl

\begin{equation}      \label{C2 : 63c}
V\rightarrow \tilde{G}\cdot V \rightarrow \tilde{G}/\tilde{H}\cong S^2.
\end{equation}

\bfl As a direct consequence $dim\, \tilde{G}\cdot V = dim \, V +2$. $\square$ \efl

\bfl  Let us recall that the $\tilde{G}-$strata $\widetilde{S}$ of $N(\Theta)$ can be obtained by
 taking the intersection of the $\tilde{G}-$strata $\overset{=}{S}$, of $\mathbb{H}^8$, with $H(\Theta)$.
Then lemma $2.4$
will be very useful in order to compute the dimension
of the strata $\widetilde{S}= \frac{\tilde{G}\cdot V\cap N(\Theta)}{\tilde{G}}$ of $\mathcal{Z}^6(\Theta)$.
 Then we have the following result  \efl

\vspace{0,5 cm}

\begin{lem}
In order to describe the moment map equations of $T^3_{\Theta}$ and
$Sp(1)$  on the strata $\tilde{G}\cdot V$, it is sufficient
 to compute them on the relative vector
spaces $V$.
\end{lem}

\bfl $\mathbf{Proof.}$ First, note that
\begin{equation}     \label{C2 : 64}
\begin{aligned}
&  1) \quad \quad \tilde{G}\cdot
\left( \begin{array}{cc}
z_{2\alpha-1} & z_{2\alpha}  \\
0   &    0 \\
\end{array}\right) =
\left( \begin{array}{cc}
v_{2\alpha-1}\epsilon &  v_{2\alpha}\epsilon   \\
v_{2\alpha-1}\sigma   &    v_{2\alpha}\sigma   \\
\end{array}\right) ,   \hspace{ 10 cm}  \\
& \\
& 2) \quad \quad \tilde{G}\cdot \left( \begin{array}{cc}
0    &  0   \\
 w_{2\beta-1} & \pm i w_{2\beta-1}\\
\end{array}\right) =
\left( \begin{array}{cc}
 -v_{2\beta-1}\overline{\sigma}   & \mp iv_{2\beta}\overline{\sigma} \\
  v_{2\beta-1}\overline{\epsilon} & \pm iv_{2\beta}\overline{\epsilon} \\
\end{array} \right),     \hspace{ 10 cm}  \\
\end{aligned}
\end{equation}
  where $v_{2\alpha-1}= \big(\cost_{\alpha}z_{2\alpha-1}+ \sit_{\alpha}z_{2\alpha}\big)\rho$,
 $v_{2\alpha}= \big(-\sit_{\alpha}z_{2\alpha-1}+\cost_{\alpha}z_{2\alpha}\big)\rho$ and $v_{2\beta-1}=$
 $ e^{\pm i\theta_{\beta}}w_{2\beta-1}\rho$.
  In order to describe the
  $Sp(1)$ moment map equations , on the sets in (\ref{C2 : 63}), we have to compute the following \efl
 \begin{equation}    \label{C2 : 65}
\begin{aligned}
& 1) \,\, \overline{   (\epsilon + j\sigma)  v_{\alpha}    }\, i\,  (\epsilon + j\sigma) v_{\alpha}=
\overline{v}_{\alpha} (    \overline{\epsilon} - j\sigma   ) i (\epsilon + j\sigma) v_{\alpha}=
i|v_{\alpha}|^2(|\epsilon |^2- |\sigma |^2) +2k\overline{\epsilon}\overline{\sigma} (v_{\alpha})^2, \\
&    \\
& 2)\,\, \overline{    (-\overline{\sigma} + j\overline{\epsilon})v_{\alpha}     }\, i\,
(-\overline{\sigma} + j\overline{\epsilon})v_{\alpha}=
 \overline{(\epsilon + j\sigma)j  v_{\alpha}    }\, i\,  (\epsilon + j\sigma)j v_{\alpha}=
-i|v_{\alpha}|^2(|\epsilon |^2- |\sigma |^2) -\\
  & \quad 2k\overline{\epsilon}\overline{\sigma} (v_{\alpha})^2,  \\
&   \\
& 3)  \,\, \overline{    (-\overline{\sigma} + j\overline{\epsilon})iv_{\alpha}     }\, i\,
(-\overline{\sigma} + j\overline{\epsilon})iv_{\alpha}=
 \overline{(\epsilon + j\sigma)j i v_{\alpha}    }\, i\,  (\epsilon + j\sigma)ji v_{\alpha}=
  \overline{(\epsilon + j\sigma)k  v_{\alpha}    }\, i\,  (\epsilon +  \\
   &\quad  j\sigma)k v_{\alpha} =  -i|v_{\alpha}|^2(|\epsilon |^2- |\sigma |^2)
 +2k\overline{\epsilon}\overline{\sigma} (v_{\alpha})^2.   \hspace{ 10  cm} \\
\end{aligned}
 \end{equation}

\bfl Then we have to
consider

\begin{equation}   \label{C2 : 66}
\begin{aligned}
&4) \,\,   \overline{   (\epsilon + j\sigma)  v_{\alpha}    }\, j\,  (\epsilon + j\sigma) v_{\alpha}=
\overline{v}_{\alpha} (    \overline{\epsilon} - j\sigma   ) j (\epsilon + j\sigma) v_{\alpha}=
j ((\epsilon )^2+ (\sigma )^2)(v_{\alpha})^2 + 2iIm (\epsilon\overline{\sigma})|v_{\alpha}|^2,\\
& \\
& 5) \,\, \overline{    (-\overline{\sigma} + j\overline{\epsilon})v_{\alpha}     }\, j\,
(-\overline{\sigma} + j\overline{\epsilon})v_{\alpha}=
 \overline{(\epsilon + j\sigma)j  v_{\alpha}    }\, j\,  (\epsilon + j\sigma)j v_{\alpha}=
 j ((\overline{\epsilon} )^2+ (\overline{\sigma} )^2)(v_{\alpha})^2 -  \\
 &\quad  2iIm (\epsilon\overline{\sigma})|v_{\alpha}|^2, \hspace{10 cm}\\
 & \\
 & 6) \,\, \overline{    (-\overline{\sigma} + j\overline{\epsilon})iv_{\alpha}     }\, j\,
(-\overline{\sigma} + j\overline{\epsilon})iv_{\alpha}=
 \overline{(\epsilon + j\sigma)k v_{\alpha}    }\, j\,  (\epsilon + j\sigma)k v_{\alpha}=
- j ((\overline{\epsilon} )^2+ (\overline{\sigma} )^2)(v_{\alpha})^2 -  \\
 &\quad  2iIm (\epsilon\overline{\sigma})|v_{\alpha}|^2, \hspace{15 cm}\\
 \end{aligned}
\end{equation}

 and \efl

\begin{equation}   \label{C2 : 67}
\begin{aligned}
&7) \,\,   \overline{   (\epsilon + j\sigma)  v_{\alpha}    }\, k\,  (\epsilon + j\sigma) v_{\alpha}=
\overline{v}_{\alpha} (    \overline{\epsilon} - j\sigma   ) k (\epsilon + j\sigma) v_{\alpha}=
k ((\epsilon )^2- (\sigma )^2)(v_{\alpha})^2 - 2iRe (\epsilon\overline{\sigma})|v_{\alpha}|^2,\\
& \\
& 8) \,\, \overline{    (-\overline{\sigma} + j\overline{\epsilon})v_{\alpha}     }\, k\,
(-\overline{\sigma} + j\overline{\epsilon})v_{\alpha}=
 \overline{(\epsilon + j\sigma)j  v_{\alpha}    }\, k\,  (\epsilon + j\sigma)j v_{\alpha}=
 -k ((\overline{\epsilon} )^2- (\overline{\sigma} )^2)(v_{\alpha})^2 +  \\
 &\quad  2iRe (\epsilon\overline{\sigma})|v_{\alpha}|^2,\\
 & \\
 & 9) \,\, \overline{    (-\overline{\sigma} + j\overline{\epsilon})iv_{\alpha}     }\, k\,
(-\overline{\sigma} + j\overline{\epsilon})iv_{\alpha}=
 \overline{(\epsilon + j\sigma)k v_{\alpha}    }\, k\,  (\epsilon + j\sigma)k v_{\alpha}=
 k ((\overline{\epsilon} )^2- (\overline{\sigma} )^2)(v_{\alpha})^2 +  \\
 &\quad  2iRe (\epsilon\overline{\sigma})|v_{\alpha}|^2.\\
\end{aligned}
\end{equation}

\bfl By using these relations we can describe the $Sp(1)-$momemt map equations. For example consider
the case of\efl
\begin{equation}     \label{C2 : 68}
{}^{(+,+)}\overset{=}{S}{}^{12}_{34}= \tilde{G}\cdot  \overbrace{\Bigg\{ \left(\begin{array}{cc|cc|cc|cc}
z_1 & z_2  & z_3 & z_4 & 0 &  0 & 0 &  0\\
0   &  0   & 0   & 0   & w_5   &  iw_5  &  w_7 & iw_7 \\
\end{array}\right)\Bigg\}}^{{}^{(+,+)}\overset{=}{V}{}^{12}_{34}:=}.
\end{equation}

\bfl Then by reading the $Sp(1)-$moment map given by
(\ref{C2 : 65}), (\ref{C2 : 66}) and (\ref{C2 : 67}), on ${}^{(+,+)}\overset{=}{S}{}^{12}_{34}$ we have: \efl
\begin{equation}     \label{C2 : 69}
 \left\{\begin{array}{l}
(|\epsilon |^2- |\sigma |^2)\big(\sum_{\alpha=1}^4 |v_{\alpha}|^2 - 2|w_5 |^2 - 2|w_7|^2\big)=0 , \\
      Im(\epsilon\overline{\sigma})    \big(\sum_{\alpha=1}^4 |v_{\alpha}|^2 - 2|w_5 |^2 - 2|w_7|^2\big)=0 , \\
    Re(\epsilon\overline{\sigma})    \big(\sum_{\alpha=1}^4 |v_{\alpha}|^2 - 2|w_5 |^2 - 2|w_7|^2\big)=0 , \\
     ((\epsilon )^2+ (\sigma )^2)   \sum_{\alpha=1}^4  (v_{\alpha})^2=0, \\
     ((\epsilon )^2- (\sigma )^2)  \sum_{\alpha=1}^4  (v_{\alpha})^2=0, \\
      \overline{\epsilon}\overline{\sigma} \sum_{\alpha=1}^4  (v_{\alpha})^2=0, \\
 \end{array}\right.    \hspace{2 cm}
\end{equation}

\bfl if and only if  \efl

\begin{equation}      \label{C2 : 70}
 \left\{\begin{array}{l}
\sum_{\alpha=1}^4 |v_{\alpha}|^2 - 2|w_5 |^2 - 2|w_7|^2=0, \\
\sum_{\alpha=1}^4  (v_{\alpha})^2=0. \\
 \end{array}\right.    \hspace{2 cm}
\end{equation}
\bfl A direct computation shows that these equations coincide with the $Sp(1)-$moment
 map equations
on the vector space ${}^{(+,+)}\overset{=}{V}{}^{12}_{34}$.
 A similar computation allows to describe the $T^3_{\Theta}$
  moment map equations on the vector space ${}^{(+,+)}\overset{=}{V}{}^{12}_{34}$.
 For the other sets in (\ref{C2 : 63}) we can use the same argument.
  Then we get the conclusion.$\square$\efl

 \vspace{0,5 cm}

\begin{prop}
\bfl Let $\zw\in \mathbb{H}^8$ be a singular point of those in $(\ref{C2 : 63})$.
 Then $\zw$ has at most two
 quaternionic
pairs which are either contained in  ${}^{\pm}V_1^{\alpha}$ or ${}^{\pm}V_2^{\alpha}$.   \efl
\end{prop}

\bfl $\mathbf{Proof}.$ We can easily exclude the case of an element which have
all of its quaternionic pairs in  ${}^{\pm}V_1^{\alpha}$ or ${}^{\pm}V_2^{\alpha}$.
In fact, up to row and column permuations, these cases are given by  \efl
 \begin{equation}   \label{C2 : 71}
 \begin{aligned}
& \tilde{G}\cdot \Bigg\{
\left( \begin{array}{cc|cc|cc|cc}
z_1 & \pm iz_1  & 0     &       0   &     0 & 0        & 0    &  0   \\
0   &    0  &  w_3  &  \pm iw_3 &  w_5  & \pm iw_5 & w_7  & \pm i w_7\\
\end{array}\right) \Bigg\},  \\
\mathrm{or} \quad & \\
&  \tilde{G}\cdot \Bigg\{
\left( \begin{array}{cc|cc|cc|cc}
z_1 & \pm iz_1  & z_3 & \pm iz_3   &     0 & 0        & 0    &  0   \\
0   &    0  &  0  &  0 &  w_5  & \pm iw_5 & w_7  & \pm i w_7\\
\end{array}\right) \Bigg\}. \hspace{10 cm} \\
\end{aligned}
\end{equation}
\bfl Then, it follows that   the  moment map equations of $T^3_{\Theta}$ and $Sp(1)$,  on the first family
of sets
in (\ref{C2 : 71}), read  \efl

 \begin{equation}       \label{C2 : 72}
 \begin{aligned}
 & T^3_{\Theta}\quad \Rightarrow\quad
 \left\{\begin{array}{l}
p_1|z_1|^2- p_2|w_3|^2 -p_3|w_5|^2 -p_4|w_7|^2=0,\\
 q_1|z_1|^2- q_2|w_3|^2 -q_3|w_5|^2 -q_4|w_7|^2=0,\\
 l_1|z_1|^2- l_2|w_3|^2 -l_3|w_5|^2 -l_4|w_7|^2=0,\\
 \end{array}\right. \hspace{2 cm}  \\
\mathrm{  and} \quad \quad &  \\
 & 
  Sp(1) \quad \Rightarrow\quad
 |z_1|^2  -|w_3|^2 -|w_5|^2 -|w_7|^2=0.  \hspace{10 cm} \\
 \end{aligned}
  \end{equation}

\bfl In the following we will often  use these notation to indicate the moment map equations
relative to $T^3_{\Theta}$ and $Sp(1)$.  Moreover, we have to consider the sphere equation   \efl
 \begin{equation}
 |z_1|^2
 + |w_3|^2 +|w_5|^2 +|w_7|^2=1.
 \end{equation}

 \bfl Then we can rewrite these equations as follows  \efl

 \begin{equation}      \label{C2 : 74}
\left( \begin{array}{cccc}
 p_1 & -p_2 & -p_3  & -p_4 \\
 q_1 & -q_2 & -q_3  & -q_4 \\
  l_1 & -l_2 & -l_3  & -l_4 \\
  1 & -1  & -1 & -1 \\
 \end{array} \right)
 \left( \begin{array}{c}
 |z_1|^2 \\
  |w_3|^2 \\
  |w_5|^2 \\
  |w_7|^2 \\
 \end{array} \right)    =
 \left( \begin{array}{c}
0\\
 0 \\
  0 \\
  0\\
 \end{array} \right),
 \end{equation}
\bfl the determinant associated is exactly $\sqp\ne 0$, then we
 get only the trivial solution.
The other cases can be treated in the same way.
 Now, we study the case
of singular points which have three quaternionic
pairs which either belong to  ${}^{\pm}V_1^{\alpha}$ or ${}^{\pm}V_2^{\alpha}$.
Then, up to columns permutations, we have only the following possibilities\efl
\begin{equation}     \label{C2 : 75}
 \begin{aligned}
 & \,\, \hspace{1 cm}\zw\in V^1_4\times {}^{\pm}V^2_2\times {}^{\pm}V^3_2
 \times {}^{\pm}V^4_2, \\
 \mathrm{or} & \\
 & \,\, \hspace{1 cm}  \zw\in V^1_5\times {}^{\pm}V^2_1\times {}^{\pm}V^3_1
 \times {}^{\pm}V^4_1. \hspace{10 cm} {}
 \end{aligned}
\end{equation}
\bfl
 In the first case of (\ref{C2 : 75}) the fixed point equations become\efl
\begin{equation}     \label{C2 : 76}
 \left\{ \begin{array}{l}
 \epsilon\rho =e^{i\theta_1}, \\
 \epsilon\rho =e^{-i\theta_1}, \\
 \overline{\epsilon}\rho =e^{\pm i\theta_2}, \\
 \overline\epsilon\rho =e^{\pm i\theta_3}, \\
  \overline\epsilon\rho =e^{\pm i\theta_4}, \\
 \end{array}\right.  \quad \iff \quad
 \left\{ \begin{array}{l}
  e^{ i(\theta_2\mp \theta_3)}=1, \\
 e^{ i(\theta_2\mp \theta_4)}=1, \\
 \overline{\epsilon}\rho =e^{\pm i\theta_2}, \\
 \epsilon\rho =e^{i\theta_1}, \\
 \epsilon\rho =e^{-i\theta_1}, \\
 \end{array}\right.
 \end{equation}
 \bfl then, the determinant  associated to the matrix
  of coefficients, up to rows permutations, becomes \efl

 \begin{equation}     \label{C2 : 78}
 \left| \begin{array}{ccc|cc}
 p_2 \mp p_3 & q_2 \mp  q_3  & l_2 \mp l_3 & 0 & 0 \\
 p_2 \mp p_4 & q_2 \mp  q_4  & l_2 \mp l_4 & 0 & 0 \\
 2p_1 &  2q_1    & 2l_1  & 0 &0 \\
 \hline
 p_1 &  q_1    & l_1  & -1 & -1 \\
 p_2 &  q_2    & l_2  & \pm1 & \mp1 \\
 \end{array}\right| = \pm 4\big( \Delta_{124} - \Delta_{123} + \Delta_{134} \big).
 \end{equation}

 \bfl The sets described by the relations in (\ref{C2 : 75}), up to row
  permuatations have the following shape: \efl
\begin{equation}        \label{C2 : 79}
\begin{aligned}
{}^{(\pm,\pm, \pm)}\overset{=}{S}{}^1_{234}:=\tilde{G}\cdot \Bigg\{
&\left( \begin{array}{cc|cc|cc|cc}
z_1 & z_2  & 0     &       0   &     0 & 0        & 0    &  0   \\
0   &    0 &  w_3  &  \pm iw_3 &  w_5  & \pm iw_5 & w_7 & \pm i w_7\\
\end{array}\right) \Bigg\}= \\
 &\Bigg\{
\left( \begin{array}{cc|cc|cc|cc}
v_1\epsilon &  v_2\epsilon  & -v_3\overline{\sigma}   & \mp iv_3\overline{\sigma}  &
   -v_5\overline{\sigma}   & \mp iv_5\overline{\sigma}      & -v_7\overline{\sigma}   & \mp iv_7\overline{\sigma} \\
v_1\sigma   &    v_2\sigma  &   v_3\overline{\epsilon} & \pm iv_3\overline{\epsilon}  &
v_5\overline{\epsilon} & \pm iv_5\overline{\epsilon}  &
 v_7\overline{\epsilon} & \pm iv_7\overline{\epsilon} \\
\end{array}\right) \Bigg\},
\end{aligned}
\end{equation}
\bfl  where $v_{1}= \big(\cost_{1}z_{1}+ \sit_{1}z_{2}\big)\rho$,
 $v_{2}= \big(-\sit_{1}z_{1}+\cost_{1}z_{2}\big)\rho$ and $v_{2\alpha-1}=$
 $ e^{\pm i\theta_{\alpha}}w_{2\alpha-1}\rho$ with $\alpha\in \{2,3,4\}$. We have  that,
 up to column permutations, the stratum ${}^{(\pm,\pm, \pm)}\overset{=}{S}{}^1_{234}$
 describes all the strata whose
 points have three quaternionic
pairs which either belong to  ${}^{\pm}V_1^{\alpha}$ or ${}^{\pm}V_2^{\alpha}$.
  Since the moment map $\mu$ is $T^3_{\Theta}-$equaivariant,
  $\nu$ is $Sp(1)-$
  equivariant and the submanifold $N(\Theta)$ is $G-$invariant, then in
  order to compute the moment map equations on the stratum
   ${}^{(\pm,\pm, \pm)}\overset{=}{S}{}^1_{234}$, it is sufficient to do this  on the vector space \efl
  \begin{equation}         \label{C2 : 80}
  \bigg\{
\left( \begin{array}{cc|cc|cc|cc}
z_1 & z_2  & 0     &       0   &     0 & 0        & 0    &  0   \\
0   &    0 &  w_3  &  \pm iw_3 &  w_5  & \pm iw_5 & w_7 & \pm i w_7\\
\end{array}\right) \bigg\}.
  \end{equation}
  \bfl Then, by using the notation introduced in (\ref{C2 : 72}) and (\ref{C2 : 73}), the
$T^3_{\Theta}$ and $Sp(1)$  moment map equations read \efl

  \begin{equation}        \label{C2 : 81}
  \begin{aligned}
 & T^3_{\Theta}\quad \Rightarrow\quad
 \left\{\begin{array}{l}
p_1Im(z_{1}\overline{z}_2)- p_2|w_3|^2 -p_3|w_5|^2 -p_4|w_7|^2=0,\\
 q_1Im(z_{1}\overline{z}_2)- q_2|w_3|^2 -q_3|w_5|^2 -q_4|w_7|^2=0,\\
 l_1Im(z_{1}\overline{z}_2)- l_2|w_3|^2 -l_3|w_5|^2 -l_4|w_7|^2=0,\\
 \end{array}\right. \hspace{ 3 cm} {}  \\
\mathrm{and}  \quad \quad&   \\
 &   Sp(1) \quad \Rightarrow\quad
  \left\{\begin{array}{l}
 |z_1|^2 + |z_2|^2 -2|w_3|^2 -2|w_5|^2 -2|w_7|^2=0, \\
 (z_1)^2 + (z_2)^2=0,       \\
 \end{array} \right.      \hspace{10 cm} {}   \\
 \end{aligned}
  \end{equation}

 \bfl Moreover, we have to consider the sphere equation \efl
 \begin{equation}    \label{C2 : 83}
 |z_1|^2 + |z_2|^2
 + 2|w_3|^2 +2|w_5|^2 +2|w_7|^2=1.
 \end{equation}
 \bfl  Let us rewrite the $Sp(1)-$moment map equations in ( \ref{C2 : 81}) and the sphere equation ( \ref{C2 : 83})
  as follows \efl
  \begin{equation}         \label{C2 : 84}
  \left\{\begin{array}{l}
 |z_1|^2 + |z_2|^2 =\frac{1}{2}, \\
 |w_3|^2 +|w_5|^2 +|w_7|^2=\frac{1}{4}, \\
 (z_1)^2 + (z_2)^2=0,       \\
 \end{array} \right.
  \end{equation}
  \bfl
Note that, $dim\, {}^{(\pm,\pm, \pm)}\overset{=}{S}{}^1_{234}= 12 $ and the intersection
  ${}^{(\pm,\pm, \pm)}\overset{=}{S}{}^1_{234}\cap N(\Theta)$ would be
 described by the equations in
 ( \ref{C2 : 81}) and ( \ref{C2 : 84}).
   Then, we would have
  that $dim\,\big( {}^{(\pm,\pm, \pm)}\overset{=}{S}{}^1_{234} \cap N(\Theta) \big)= 5$. Since
  $dim\, \tilde{G}=7$, and the action of $\tilde{G}$
  is locally free on $N(\Theta)$, it follows that each of
   ${}^{(\pm,\pm, \pm)}\overset{=}{S}{}^1_{234}$ does not intersect $N(\Theta)$.
  The other cases can be treated analogously.\efl

 \bfl  Now, consider the case of a singular point $\zw\in \mathbb{H}^8$
  which has two quaternionic
 pairs cointained either
in  ${}^{\pm}V_1^{\alpha}$ or ${}^{\pm}V_2^{\alpha}$. Without loss of generality,
 we can assume these two quaternionic
pairs being  $(u_5, u_6)$ and
 $(u_7, u_8)$ respectively. Note that, up to the $\tilde{G}-$action, $\uu$, $\alpha\in\{3,4\}$
  have to be
 contained in the same eigenspace, analogously $\uu$, $\alpha\in\{1,2\}$ have
 to be contained in the same
 eigenspace and in order to avoid $rank\, \zw< 4$ the only possible solutions
 are either   \efl
 \begin{equation}      \label{C2 : 85}
 \begin{aligned}
 & \,\, \hspace{1 cm}\zw\in V^1_4\times V^2_4\times {}^{\pm}V^3_2\times {}^{\pm}V^4_2, \\
 \mathrm{or} & \\
 & \,\, \hspace{1 cm}  \zw\in V^1_5\times V^2_5\times {}^{\pm}V^3_1\times {}^{\pm}V^4_1, \hspace{10 cm} {}
 \end{aligned}
 \end{equation}

\bfl Consider the first
case ( the second one can be studied similarly). Then we have\efl
\begin{equation}        \label{C2 : 86}
\begin{aligned}
{}^{(\pm,\pm, \pm)}\overset{=}{S}{}^{12}_{34}:=\tilde{G}\cdot \Bigg\{
&\left( \begin{array}{cc|cc|cc|cc}
z_1 & z_2  & z_3 & z_4   &   0   & 0        & 0    &  0   \\
0   &    0 &  0  & 0     &  w_5  & \pm iw_5 & w_7  & \pm i w_7\\
\end{array}\right) \Bigg\}= \\
 &\Bigg\{
\left( \begin{array}{cc|cc|cc|cc}
v_1\epsilon &  v_2\epsilon  & v_3\epsilon &  v_4\epsilon  &
   -v_5\overline{\sigma}   & \mp iv_5\overline{\sigma}      & -v_7\overline{\sigma}
     & \mp iv_7\overline{\sigma} \\
v_1\sigma   &    v_2\sigma  &  v_3\sigma   &    v_4\sigma    &
v_5\overline{\epsilon} & \pm iv_5\overline{\epsilon}  &
 v_7\overline{\epsilon} & \pm iv_7\overline{\epsilon} \\
\end{array}\right) \Bigg\}
\end{aligned}
\end{equation}
\bfl  where $v_{2\alpha-1}= \big(\cost_{\alpha}z_{2\alpha-1}+
 \sit_{\alpha}z_{2\alpha}\big)\rho$,
 $v_{2\alpha}= \big(-\sit_{\alpha}z_{2\alpha-1}+
 \cost_{\alpha}z_{2\alpha}\big)\rho$ and $v_{2\beta-1}=$
 $ e^{\pm i\theta_{\beta}}w_{2\beta-1}\rho$ with
  $\alpha\in \{1,2\}$, $\beta\in\{ 3,4\}$.
 Note that, the stratum described in
 (\ref{C2 : 86}) has dimension $14$.
The fixed point equations
become\efl
 \begin{equation} \label{C2 : 87}
 \left\{ \begin{array}{l}
 \epsilon\rho =e^{i\theta_1}, \\
 \epsilon\rho =e^{-i\theta_1}, \\
 \epsilon\rho =e^{i\theta_2}, \\
 \epsilon\rho =e^{-i\theta_2}, \\
 \overline\epsilon\rho =e^{\pm i\theta_3}, \\
  \overline\epsilon\rho =e^{\pm i\theta_4}, \\
 \end{array}\right.  \quad \iff \quad
 \left\{ \begin{array}{l}
  e^{i(\theta_1-\theta_2)}=1, \\
 e^{ i(\theta_3\mp \theta_4)}=1, \\
 \epsilon\rho =e^{i\theta_1}, \\
 \epsilon\rho =e^{-i\theta_1}, \\
 \overline\epsilon\rho =e^{\pm i\theta_3}. \\
 \end{array}\right.
 \end{equation}
\bfl The determinant associated
to the matrix of coefficients,
up to rows permutations, becomes:\efl

 \begin{equation}     \label{C2 : 88}
 \left| \begin{array}{ccc|cc}
 p_1-  p_2 & q_1 -  q_2  & l_1 - l_2 & 0 & 0 \\
 p_3 \mp p_4 & q_3 \mp  q_4  & l_3 \mp l_4 & 0 & 0 \\
 2p_1 &  2q_1    & 2l_1  & 0 &0 \\
 \hline
 p_1 &  q_1    & l_1  & -1 & -1 \\
 p_3 &  q_3    & l_3  & \pm1 & \mp1 \\
 \end{array}\right| = \pm 4\big(  - \Delta_{123} \pm \Delta_{124} \big).
 \end{equation}

\bfl The moment map equations read\efl
 \begin{equation}            \label{C2 : 89}
 T^3_{\Theta}\quad \Rightarrow\quad
 \left\{\begin{array}{l}
p_1Im(z_{1}\overline{z}_2) + p_2Im(z_{3}\overline{z}_4) \mp p_3|w_5|^2 \mp p_4|w_7|^2=0,\\
 q_1Im(z_{1}\overline{z}_2)+ q_2Im(z_{3}\overline{z}_4) \mp q_3|w_5|^2 \mp q_4|w_7|^2=0,\\
 l_1Im(z_{1}\overline{z}_2)+ l_2Im(z_{3}\overline{z}_4) \mp l_3|w_5|^2 \mp l_4|w_7|^2=0,\\
 \end{array}\right. \hspace{ 3 cm} {}
 \end{equation}
 \bfl and \efl
 \begin{equation}        \label{C2 : 90}
  Sp(1) \quad \Rightarrow\quad
  \left\{\begin{array}{l}
 |z_1|^2 + |z_2|^2 +  |z_3|^2 + |z_4|^2  -2|w_5|^2 -2|w_7|^2=0, \\
 (z_1)^2 + (z_2)^2 + (z_3)^2 + (z_4)^2 =0,       \\
 \end{array} \right.      \hspace{3 cm} {}
  \end{equation}
 \bfl besides, we have to consider the sphere equation \efl
 \begin{equation}       \label{C2 : 91}
 |z_1|^2 + |z_2|^2
 + |z_3|^2 + |z_4|^2 +2|w_5|^2 +2|w_7|^2=1.
 \end{equation}
\bfl  When the quaternionic pairs $(w_5,\pm i w_5)$ and $(w_7, \pm i w_7)$
have the same sign,  the equations (\ref{C2 : 89}) give \efl
\begin{equation}        \label{C2 : 91a}
  Im(z_{1}\overline{z}_2) =\pm \frac{\Delta_{234}}{\Delta_{123}}|w_7|^2, \quad
  Im(z_{3}\overline{z}_4) = \mp \frac{\Delta_{134}}{\Delta_{123}}|w_7|^2 \,\, \mathrm{and}\,\,
  |w_5|^2 = -\frac{\Delta_{124}}{\Delta_{123}}|w_7|^2.
\end{equation}
 \bfl Moreover, the equations (\ref{C2 : 90}) and (\ref{C2 : 91}) yields\efl
 \begin{equation}    \label{C2 : 91b}
 |w_5|^2 + |w_7|^2= \frac{1}{4}.
 \end{equation}
\bfl By taking the equation (\ref{C2 : 91b})  and the value of $|w_5|^2$ in (\ref{C2 : 91a}),
 we get $|w_5|^2( \Delta_{124}$
$-\Delta_{123})= \frac{1}{4}\Delta_{124}.$  Thus when  $|w_5|= $
$|w_7|\ne0$, we have that
 $-\Delta_{123}+\Delta_{124}\ne 0$ and the determinat (\ref{C2 : 88}) does not vanish.
In this case
$dim\,\big({}^{(\pm,\pm, \pm)}\overset{=}{S}{}^{12}_{34}\cap N(\Theta)\big)= 7$.
Since $\dim\, \tilde{G}=7$, we have that
each quotient $\big({}^{(\pm,\pm, \pm)}\overset{=}{S}{}^{12}_{34}\cap N(\Theta)
\big)/ \tilde{G}$ on the twistor space $\mathcal{Z}^6(\Theta)$
coincides with a point.  The remaining spaces of this type can be treated similarly and the associated
 $\tilde{G}-$strata
are just the following\efl

 \begin{align*}
& 1) \,\,\, \quad  {}^{(\pm , \pm)}\overset{=}{S}{}^{12}_{34}:=\tilde{G}\cdot \Bigg\{
\left( \begin{array}{cc|cc|cc|cc}
z_1 & z_2  & z_3 & z_4   &   0    &  0                & 0    &  0   \\
0   &    0 &  0  & 0     &    w_5  &  \pm iw_5        & w_7  &  \pm iw_7\\
\end{array}\right) \Bigg\}, \hspace{10 cm}  \\
& 2) \,\,\, \quad {}^{(\pm , \pm)}\overset{=}{S}{}^{13}_{24}:=\tilde{G}\cdot \Bigg\{
\left( \begin{array}{cc|cc|cc|cc}
z_1 & z_2  & 0   & 0             &  z_5  & z_6        & 0    &  0   \\
0   &    0 & w_3  &  \pm iw_3    &   0   &  0          & w_7  &  \pm iw_7\\
\end{array}\right) \Bigg\},  \hspace{10 cm} \\
&3) \,\,\,\quad  {}^{(\pm , \pm)}\overset{=}{S}{}^{14}_{23}:=\tilde{G}\cdot \Bigg\{
\left( \begin{array}{cc|cc|cc|cc}
z_1 & z_2  & 0     &  0            &    0   & 0       & z_7   & z_8  \\
0   &    0 &  w_3  &  \pm iw_3     &   w_5  &  \pm iw_5   &  0    &  0          \\
\end{array}\right) \Bigg\},  \hspace{10 cm} \\
& 4)\,\,\,\quad   {}^{(\pm , \pm)}\overset{=}{S}{}^{24}_{13}:=\tilde{G}\cdot \Bigg\{
\left( \begin{array}{cc|cc|cc|cc}
0    & 0          & z_3 & z_4       &   0   &    0          & z_7   & z_8   \\
w_1  &  \pm iw_1  & 0 &    0     &   w_5 &    \pm iw_5        & 0 & 0\\
\end{array}\right) \Bigg\},   \\
 & 5)\,\,\, \quad  {}^{(\pm , \pm)}\overset{=}{S}{}^{34}_{12}:=\tilde{G}\cdot \Bigg\{
\left( \begin{array}{cc|cc|cc|cc}
0  &  0 & 0    &  0           &    z_5 & z_6             & z_7   & z_8  \\
w_1  &  \pm iw_1   & w_3  & \pm iw_3    &   0  &  0   &  0    &  0          \\
\end{array}\right) \Bigg\},   \\
\end{align*}
\begin{equation}       \label{C2 : 92}
 \begin{aligned}
& 6)\,\,\, \quad  {}^{(\pm , \pm)}\overset{=}{S}{}^{24}_{13}:=\tilde{G}\cdot \Bigg\{
\left( \begin{array}{cc|cc|cc|cc}
0   &    0       &  z_3 & z_4  &    0 & 0           & z_7   & z_8  \\
w_1  &  \pm iw_1 & 0    & 0    &   w_5 &    \pm iw_5   &  0    &  0          \\
\end{array}\right) \Bigg\}. \hspace{10 cm }  \\
\end{aligned}
\end{equation}

\vspace{0,5 cm}

\bfl The last case we have to consider, is that one when there is
just a quaternionic pair which is either contained in ${}^{\pm}V^{\alpha}_{1}$ or
${}^{\pm}V^{\alpha}_{2}.$
 As a consequence the only possibilities are the following \efl
\begin{equation}     \label{C2 : 93}
 \begin{aligned}
 & \,\, \hspace{1 cm}\zw\in V^1_4\times V^2_4\times {}^{\pm}V^3_4\times {}^{\pm}V^4_2, \\
 \mathrm{or} & \\
 &  \,\, \hspace{1 cm}  \zw\in V^1_5\times V^2_5\times {}^{\pm}V^3_5\times {}^{\pm}V^4_1. \hspace{10 cm} {}
 \end{aligned}
 \end{equation}
  \bfl Note that the two cases in (\ref{C2 : 93}), up to multiply on the
  left by $j$, concide, then all of the possible spaces of solutions are the following\efl
 \begin{align*}
& 1) \,\,\, \quad  {}^+\overset{=}{S}{}^{123}_4:=\tilde{G}\cdot \Bigg\{
\left( \begin{array}{cc|cc|cc|cc}
z_1 & z_2  & z_3 & z_4   &    z_5 & z_6        & 0    &  0   \\
0   &    0 &  0  & 0     &   0   &  0          & w_7  &  iw_7\\
\end{array}\right) \Bigg\}, \hspace{10 cm}  \\
& 2) \,\,\, \quad  {}^-\overset{=}{S}{}^{123}_4:=\tilde{G}\cdot \Bigg\{
\left( \begin{array}{cc|cc|cc|cc}
z_1 & z_2  & z_3 & z_4   &    z_5 & z_6        & 0    &  0   \\
0   &    0 &  0  & 0     &   0   &  0          & w_7  &  -iw_7\\
\end{array}\right) \Bigg\},   \\
&3) \,\,\,\quad  {}^+\overset{=}{S}{}^{124}_3:=\tilde{G}\cdot \Bigg\{
\left( \begin{array}{cc|cc|cc|cc}
z_1 & z_2  & z_3 & z_4   &    0   & 0           & z_7   & z_8  \\
0   &    0 &  0  & 0     &   w_5  &  iw_5   &  0    &  0          \\
\end{array}\right) \Bigg\},  \hspace{10 cm} \\
& 4)\,\,\,\quad   {}^-\overset{=}{S}{}^{124}_3:=\tilde{G}\cdot \Bigg\{
\left( \begin{array}{cc|cc|cc|cc}
z_1 & iz_1  & z_3 & z_4       &   0   &    0          & z_7   & z_8   \\
0   &    0  & 0 &    0     &   w_5 &    -iw_5        & 0 & 0\\
\end{array}\right) \Bigg\},   \\
 & 5)\,\,\, \quad  {}^+\overset{=}{S}{}^{134}_2:=\tilde{G}\cdot \Bigg\{
\left( \begin{array}{cc|cc|cc|cc}
z_1 & z_2  & 0    &  0           &    z_5 & z_6             & z_7   & z_8  \\
0   &    0 & w_3  & iw_3    &   0  &  0   &  0    &  0          \\
\end{array}\right) \Bigg\},   \\
& 6)\,\,\, \quad  {}^-\overset{=}{S}{}^{134}_2:=\tilde{G}\cdot \Bigg\{
\left( \begin{array}{cc|cc|cc|cc}
z_1 & z_2  & 0    &  0           &    z_5 & z_6             & z_7   & z_8  \\
0   &    0 & w_3  & -iw_3    &   0  &  0   &  0    &  0          \\
\end{array}\right) \Bigg\}, \hspace{10 cm }  \\
\end{align*}
\begin{equation}  \label{C2 : 94}
 \begin{aligned}
& 7)\,\,\,\quad   {}^+\overset{=}{S}{}^{234}_1:=\tilde{G}\cdot \Bigg\{
\left( \begin{array}{cc|cc|cc|cc}
0   &    0      &  z_3 & z_4          &    z_5 & z_6             & z_7   & z_8   \\
w_1  & -iw_1  & 0   &    0          &      0  &            0   &  0    &  0\\
\end{array}\right) \Bigg\},   \\
& 8)\,\,\,\quad   {}^-\overset{=}{S}{}^{234}_1:=\tilde{G}\cdot \Bigg\{
\left( \begin{array}{cc|cc|cc|cc}
0   &    0      & z_3 & z_4 &    z_5 & z_6             & z_7   & z_8  \\
w_1  & -iw_1& 0   &    0     &      0  &    0              &  0    &  0\   \
\end{array}\right) \Bigg\}.  \hspace{10 cm}  \\
\end{aligned}
\end{equation}
 \bfl  Note that  ${}^{\pm}\overset{=}{S}{}^{234}_1$ is equivalent to   \efl

 \begin{equation}   \label{C2 : 94a}
 \begin{aligned}
& \tilde{G}\cdot \Bigg\{
\left( \begin{array}{cc|cc|cc|cc}
0   &    0      &  z_3 & z_4          &    z_5 & z_6             & z_7   & z_8   \\
w_1  & \pm iw_1  & 0   &    0          &      0  &            0   &  0    &  0\\
\end{array}\right) \Bigg\} =  \\
&\hspace{5 cm}\tilde{G}\cdot \Bigg\{
\left( \begin{array}{cc|cc|cc|cc}
z_1 & \pm iz_1  & 0       & 0           &   0   &     0      & 0    &  0   \\
0   &    0  & w_3     & w_4         &   w_5 & w_6        & w_7  &  w_8\\
\end{array}\right) \Bigg\}.
\end{aligned}
 \end{equation}

\bfl  Now, consider the pair of strata
 ${}^{\pm}\overset{=}{S}{}^{123}_4$. Then,
 the fixed point equations relative to this stratum become
\begin{equation}    \label{C2 : 95}
 \left\{ \begin{array}{l}
  e^{ i(\theta_1- \theta_2)}=1, \\
 e^{ i(\theta_1- \theta_3)}=1, \\
 \epsilon\rho =e^{i\theta_1}, \\
 \epsilon\rho =e^{-i\theta_1}, \\
  \overline\epsilon\rho =e^{\pm i\theta_4}. \\
 \end{array}\right.
 \end{equation}
 Then the determinant of the coefficient matrix, up to a sign, becomes\efl
 \begin{equation}    \label{C2 : 96}
 \left| \begin{array}{ccc|cc}
 p_1 - p_2 & q_1 -  q_2  & l_1- l_2 & 0 & 0 \\
 p_1- p_3 & q_1 -  q_3  & l_1 -l_3 & 0 & 0 \\
 2p_1 &  2q_1    & 2l_1  & 0 &0 \\
 \hline
 p_1 &  q_1    & l_1  & -1 & -1 \\
 p_4 &  q_4    & l_4  & \pm1 & \mp1 \\
 \end{array}\right| = \pm 4\Delta_{123} .
 \end{equation}
 \bfl In the general case we have that the matrix
 of coefficients associated to the
 fixed point equations, up to a column
 permutation for the first quaternionic pair,
 become \efl
 \begin{equation}       \label{C2 : 97}
 \left| \begin{array}{ccc|cc}
 p_1 - p_{\alpha} & q_1 -  q_{\alpha}  & l_1- l_{\alpha} & 0 & 0 \\
 p_1- p_{\beta} & q_1 -  q_{\beta}  & l_1 -l_{\beta} & 0 & 0 \\
 2p_1 &  2q_1    & 2l_1  & 0 &0 \\
 \hline
 p_1 &  q_1    & l_1  & -1 & -1 \\
 p_{\gamma} &  q_{\gamma}    & l_{\gamma}  & \pm1 & \mp1 \\
 \end{array}\right| = \pm 4\Delta_{1\alpha\beta},
 \end{equation}
 \bfl where $ \alpha, \beta, \gamma\in\{2,3,4\} $ and
 $\alpha, \beta, \gamma$ different to each other. Note that the isotropy
 for the strata listed in (\ref{C2 : 94}), up to move the first quaternionic coordinate,
  only depends on the minor determinants
 $\pm \Delta_{\alpha\beta\gamma}$.
 Now, for the moment map equations relative to the case of
  ${}^{\pm}\overset{=}{S}{}^{123}_4$ we obtain  \efl
\begin{equation}         \label{C2 : 98}
 T^3_{\Theta}\quad \Rightarrow\quad
 \left\{\begin{array}{l}
p_1Im(z_{1}\overline{z}_2) + p_2Im(z_{3}\overline{z}_4) +p_3Im(z_{5}\overline{z}_6) -p_4|w_7|^2=0,\\
 q_1Im(z_{1}\overline{z}_2)+ q_2Im(z_{3}\overline{z}_4) +q_3Im(z_{5}\overline{z}_6) -q_4|w_7|^2=0,\\
 l_1Im(z_{1}\overline{z}_2)+ l_2Im(z_{3}\overline{z}_4) +l_3Im(z_{5}\overline{z}_6) -l_4|w_7|^2=0,\\
 \end{array}\right. \hspace{ 3 cm} {}
 \end{equation}
 \bfl and \efl
 \begin{equation}     \label{C2 : 99}
  Sp(1) \quad \Rightarrow\quad
  \left\{\begin{array}{l}
 |z_1|^2 + |z_2|^2 +  |z_3|^2 + |z_4|^2  + |z_5|^2 + |z_6|^2  -2|w_7|^2=0, \\
 (z_1)^2 + (z_2)^2 + (z_3)^2 + (z_4)^2 +(z_5)^2 + (z_6)^2 =0,       \\
 \end{array} \right.      \hspace{3 cm} {}
  \end{equation}
 \bfl and we have to consider the sphere equation \efl
 \begin{equation}        \label{C2 : 100}
 |z_1|^2 + |z_2|^2
 + |z_3|^2 + |z_4|^2 + |z_5|^2 + |z_6|^2 +2|w_7|^2=1.
 \end{equation}
\bfl The equations in (\ref{C2 : 98}), (\ref{C2 : 99}) and (\ref{C2 : 100}) give
$dim\,\big({}^{\pm}\overset{=}{S}{}^{123}_4\cap N(\Theta)\big)= 9$,
 then in follows that
 $dim\,\big( {}^{\pm}\overset{=}{S}{}^{123}_4\cap N(\Theta)/\tilde{G}\big) =2$,
  that is, on the
 twistor level we have exactly a $2-$sphere $S^2$.
  The other cases in ( \ref{C2 : 94}) can be
 treated in the same way. $\square$  \efl

\bfl Notice that, the singular strata (relative to the action of
 $\tilde{G}$ on $\mathbb{H}^8$) whose points have no
  quaternionic pair is contained in
 ${}^{\pm}V_1^{\alpha}$ or ${}^{\pm}V_2^{\alpha}$, up to row permutations,
 are the following  \efl

  \begin{equation}  \label{C2 : 101} {
 \begin{aligned}
& 1) \,\,\, \quad \overset{=}{S}{}^{123}_4:= \tilde{G}\cdot \bigg\{
\left( \begin{array}{cc|cc|cc|cc}
z_1 & z_2  & z_3 & z_4   &    z_5 & z_6        & 0    &  0   \\
0   &    0 &  0  & 0     &   0   &  0          & w_7  &  w_8\\
\end{array}\right) \bigg\},   \\
&2) \,\,\,\quad  \overset{=}{S}{}^{124}_{3}:=\tilde{G}\cdot \bigg\{
\left( \begin{array}{cc|cc|cc|cc}
z_1 & z_2  & z_3 & z_4   &    0   & 0           & z_7   & z_8  \\
0   &    0 &  0  & 0     &   w_5  &  w_6   &  0    &  0          \\
\end{array}\right) \bigg\},   \\
 & 3)\,\,\, \quad \overset{=}{S}{}^{134}_2:=\tilde{G}\cdot \bigg\{
\left( \begin{array}{cc|cc|cc|cc}
z_1 & z_2  & 0    &  0           &    z_3 & z_4             & z_7   & z_8  \\
0   &    0 & w_3  & w_4    &   0  &  0   &  0    &  0          \\
\end{array}\right) \bigg\},   \\
& 4)\,\,\,\quad  \overset{=}{S}{}^{12}_{34}:=\tilde{G}\cdot \bigg\{\left( \begin{array}{cc|cc|cc|cc}
z_1 & z_2  & z_3 & z_4   &    0 & 0        & 0    &  0   \\
0   &    0 &  0  & 0     &   w_5  &  w_6        & w_7  &  w_8\\
\end{array}\right) \bigg\}.   \\
& 5)\,\,\,\quad  \overset{=}{S}{}^{13}_{24}:=\tilde{G}\cdot \bigg\{
\left( \begin{array}{cc|cc|cc|cc}
z_1 & z_2  & 0   & 0 &    z_5 & z_6        & 0    &  0   \\
0   &    0 & w_3  & w_4     &   0   &  0          & w_7  &  w_8\\
\end{array}\right) \bigg\},   \hspace{10 cm}     \\
& 6)\,\,\,\quad  \overset{=}{S}{}^{14}_{23}:=\tilde{G}\cdot \bigg\{
\left( \begin{array}{cc|cc|cc|cc}
z_1 & z_2  & 0   & 0 &   0 & 0      &  z_7 & z_8   \\
0   &    0 & w_3  & w_4     &   w_5 & w_6   & 0 & 0 \\
\end{array} \right)    \bigg\}, \\
& 7)\,\,\,\quad \overset{=}{S}{}^{1}_{234}:=\tilde{G}\cdot \bigg\{
\left( \begin{array}{cc|cc|cc|cc}
z_1 & z_2  & 0   & 0 &   0 & 0      & 0    &  0   \\
0   &    0 & w_3  & w_4     &   w_5 & w_6        & w_7  &  w_8\\
\end{array}\right) \bigg\}.      \\
 \end{aligned}    }    \hspace{4,5 cm}    {}
 \end{equation}

 \bfl Then we have the following result\efl

 \vspace{0,5 cm}

 \begin{teor}
\bfl The strata $\overset{=}{S}{}^{\alpha\beta\gamma}_{\delta}$ and
$\overset{=}{S}{}^{\alpha\beta}_{\gamma\delta}$
 intersect the submanifold $N(\Theta)$ and their intersections are such that
 \begin{equation}    \label{C2 : 102}
 \begin{aligned}
& 1)\,\,\,\, \overset{=}{S}{}^{\alpha\beta\gamma}_{\delta}\cap N(\Theta)\,\,
is\,\, formed \,\, by\,\, the\,\, two\,\, connected \,\, components\,\,
     {}^+\overset{=}{S}{}^{\alpha\beta\gamma}_{\delta}\cap N(\Theta) \,\,
     \mathrm{and}
     \,\, {}^-\overset{=}{S}{}^{\alpha\beta\gamma}_{\delta}\cap N(\Theta), \\
& 2)\,\,\,\,  \overset{=}{S}{}^{\alpha\beta}_{\gamma\delta}\cap N(\Theta)\,\,
   is\,\, formed \,\, by\,\, one\,\, connected \,\, components. \\
 \end{aligned}    \hspace{7 cm}
 \end{equation}
Moreover, each of the strata $\overset{=}{S}{}^{\alpha\beta}_{\gamma\delta}$
 contains at most four substrata of those listed in $(\ref{C2 : 92})$.\efl
 \end{teor}

\bfl $\mathbf{Proof.}$ The cases $1)$, $2)$, $3)$ and $7)$ in (\ref{C2 : 101})
 are similar, then it is
sufficient to study the first one. Analogously, in order to study the cases $4)$,
$5)$ and $6)$ it is enought to consider $\overset{=}{S}{}^{12}_{34}$.
The moment map equations on the stratum $\overset{=}{S}{}^{123}_4$ become \efl

\begin{equation}       \label{C2 : 103}
 T^3_{\Theta}\quad \Rightarrow\quad
 \left\{\begin{array}{l}
p_1Im(z_{1}\overline{z}_2) + p_2Im(z_{3}\overline{z}_4) +p_3Im(z_{5}\overline{z}_6) +p_4Im(w_{7}\overline{w}_8)=0,\\
 q_1Im(z_{1}\overline{z}_2)+ q_2Im(z_{3}\overline{z}_4) +q_3Im(z_{5}\overline{z}_6) +q_4Im(w_{7}\overline{w}_8)=0, \\
 l_1Im(z_{1}\overline{z}_2)+ l_2Im(z_{3}\overline{z}_4) +l_3Im(z_{5}\overline{z}_6) +l_4Im(w_{7}\overline{w}_8)=0, \\
 \end{array}\right. \hspace{ 3 cm} {}
 \end{equation}

 \bfl and \efl

 \begin{equation}      \label{C2 : 104}
  Sp(1) \quad \Rightarrow\quad
  \left\{\begin{array}{l}
 |z_1|^2 + |z_2|^2 +  |z_3|^2 + |z_4|^2  + |z_5|^2 + |z_6|^2  -|w_7|^2 -|w_8|^2=0, \\
 (z_1)^2 + (z_2)^2 + (z_3)^2 + (z_4)^2 +(z_5)^2 + (z_6)^2 =0,       \\
 (w_7)^2 + (w_8)^2 =0.    \\
 \end{array} \right.      \hspace{3 cm} {}
  \end{equation}

\bfl In order to describe the intersection $\overset{=}{S}{}^{123}_4\cap N(\Theta)$, as usual
we have to consider also the
sphere equation $\sum_{\alpha=1}^6|z_{\alpha}|^2   +|w_7|^2 +|w_8|^2=1$.
  By using these equations, we get     \efl

\begin{equation}       \label{C2 : 105}
\left\{ \begin{array}{c}
|w_7|^2 +|w_8|^2=\frac{1}{2},    \\
 (w_7)^2 + (w_8)^2 =0,     \\
\end{array} \right.   \quad \Rightarrow \quad
\left\{ \begin{array}{c}
w_8 = \pm iw_7 , \\
|w_7|^2=\frac{1}{4}.\\
\end{array} \right.
\end{equation}

\bfl Then, we have that ${}^{\pm}\overset{=}{S}{}^{123}_4\supset\overset{=}{S}{}^{123}_4 = 18$,
 from the moment map equations it follows that $dim\, \overset{=}{S}{}^{123}_4 = 18$, and
$dim\, \big(\big(\,\overset{=}{S}{}^{123}_4\cap N(\Theta)\big)/\tilde{G}\big) =2=
 dim\, \big( \big(\,{}^{\pm}\overset{=}{S}{}^{123}_4\cap N(\Theta)\big) /\tilde{G}\big). $
  Then, $\big(\, {}^{\pm}\overset{=}{S}{}^{123}_4\cap N(\Theta)\big)/\tilde{G}\cong
   \big(\, \overset{=}{S}{}^{123}_4\cap N(\Theta)\big)/\tilde{G}$
 is isomorphic to a $2-$sphere $S^2$. Now, consider the moment map equations on
 $\overset{=}{S}{}^{12}_{34}$, it follows that:   \efl

 \begin{equation}             \label{C2 : 106}
 T^3_{\Theta}\quad \Rightarrow\quad
 \left\{\begin{array}{l}
p_1Im(z_{1}\overline{z}_2) + p_2Im(z_{3}\overline{z}_4) +p_3Im(w_{5}\overline{w}_6) +p_4Im(w_{7}\overline{w}_8)=0,\\
 q_1Im(z_{1}\overline{z}_2)+ q_2Im(z_{3}\overline{z}_4) +q_3Im(z_{5}\overline{z}_6) +q_4Im(w_{7}\overline{w}_8)=0, \\
 l_1Im(z_{1}\overline{z}_2)+ l_2Im(z_{3}\overline{z}_4) +l_3Im(z_{5}\overline{z}_6) +l_4Im(w_{7}\overline{w}_8)=0, \\
 \end{array}\right. \hspace{ 3 cm} {}
 \end{equation}

 \bfl and \efl

 \begin{equation}         \label{C2 : 107}
  Sp(1) \quad \Rightarrow\quad
  \left\{\begin{array}{l}
 |z_1|^2 + |z_2|^2 +  |z_3|^2 + |z_4|^2  - |w_5|^2 - |w_6|^2  -|w_7|^2 -|w_8|^2=0, \\
 (z_1)^2 + (z_2)^2 + (z_3)^2 + (z_4)^2  =0,       \\
 (w_5)^2 + (w_6)^2 + (w_7)^2 + (w_8)^2 =0.    \\
 \end{array} \right.      \hspace{3 cm} {}
  \end{equation}

 \bfl Then it follows  that $dim\, \big(\big(\, \overset{=}{S}{}^{12}_{34}\cap N(\Theta)\big)/\tilde{G}\big) =2$ and
  $\big(\, \overset{=}{S}{}^{12}_{34}\cap N(\Theta)\big)/\tilde{G}\cong S^2$.
  Moreover, each of the strata $4)$, $5)$ and $6)$ in (\ref{C2 : 101}) contains four
    substrata of those in (\ref{C2 : 92}). In fact, by looking at the shape of these two families
    of strata we get the following relations   \efl

\begin{equation}   \label{C2 : 108}
{
\begin{aligned}
&1)\quad \quad \overset{=}{S}{}^{12}_{34}\,\, \mathrm{contains}\,\,
{}^{(\pm, \pm)}\overset{=}{S}{}^{12}_{34}\,\, \mathrm{and} \,\,{}^{(\pm, \pm)}\overset{=}{S}{}^{34}_{12},  \\
&2)\quad \quad \overset{=}{S}{}^{13}_{24}\,\, \mathrm{contains}\,\,
{}^{(\pm, \pm)}\overset{=}{S}{}^{13}_{24}\,\, \mathrm{and} \,\,{}^{(\pm, \pm )}\overset{=}{S}{}^{13}_{24},   \\
&3) \quad \quad \overset{=}{S}{}^{14}_{23}\,\, \mathrm{contains}\,\,
{}^{(\pm, \pm)}\overset{=}{S}{}^{14}_{23}\,\, \mathrm{and} \,\,{}^{(\pm, \pm )}\overset{=}{S}{}^{14}_{23}.  \\
\end{aligned}
}    \hspace{ 10  cm}
\end{equation}

\bfl Each of the  strata in (\ref{C2 : 92})
 yields a point on $\mathcal{Z}^6(\Theta)$,
and the relations in (\ref{C2 : 108}) give that these points can be divided
 in three family  and the points of each family are
 joined by a $2-$sphere. $\square$ \efl


 \vspace{0.5 cm}

 \begin{lem}
\bfl The two families of singular points  $(\ref{C2 : 45})$ and $(\ref{C2 : 101})$
 do not intesect each other.  \efl
 \end{lem}

 \bfl $\mathbf{Proof}.$
 Note that, up to columns  permutations, the strata $1)$, $2)$, $3)$ and $7)$
  in (\ref{C2 : 101})
   (by using the action on the group $\tilde{G}$) can be rewritten as follows
 \efl
 \begin{equation}   \label{C2 : 109}
  \bigg\{
\left( \begin{array}{cc|cc|cc|cc}
v_1\epsilon &  v_2\epsilon  & v_3\epsilon  &  v_4\epsilon    & v_5\epsilon            & v_6\epsilon      & -v_7\overline{\sigma}   & -v_8\overline{\sigma} \\
v_1\sigma   &    v_2\sigma  &  v_3\sigma   &    v_4\sigma    & v_5\sigma              & v_6\sigma        & v_7\overline{\epsilon} & v_8\overline{\epsilon} \\
\end{array}\right) \bigg\},
 \end{equation}
 \bfl where $v_{2\alpha-1}= \big(\cost_{\alpha}z_{2\alpha-1}+ \sit_{\alpha}z_{2\alpha}\big)\rho$,
 $v_{2\alpha}= \big(-\sit_{\alpha}z_{2\alpha-1}+\cost_{\alpha}z_{2\alpha}\big)\rho$, $\alpha\in\{1,2,3\}$ and
 $v_7=\big( \cost_4 w_7 + \sit_4 w_8 \big)\rho$, $v_8=\big( -\sit_4 w_7 + \cost_4 w_8 \big)$.   Instead,
 the cases $4)$, $5)$ and $6)$ in (\ref{C2 : 101}), up to columns permutations, become \efl
  \begin{equation}   \label{C2 : 110}
  \bigg\{
\left( \begin{array}{cc|cc|cc|cc}
v_1\epsilon &  v_2\epsilon  & v_3\epsilon  &  v_4\epsilon    &  -v_5\overline{\sigma}   & -v_6\overline{\sigma}      & -v_7\overline{\sigma}   & -v_8\overline{\sigma} \\
v_1\sigma   &    v_2\sigma  &  v_3\sigma   &    v_4\sigma    & v_5\overline{\epsilon} & v_6\overline{\epsilon}      & v_7\overline{\epsilon} & v_8\overline{\epsilon} \\
\end{array}\right) \bigg\},
 \end{equation}
\bfl  where $v_{2\alpha-1}= \big(\cost_{\alpha}z_{2\alpha-1}+ \sit_{\alpha}z_{2\alpha}\big)\rho$,
 $v_{2\alpha}= \big(-\sit_{\alpha}z_{2\alpha-1}+\cost_{\alpha}z_{2\alpha}\big)\rho$, $\alpha\in\{1,2\}$ and
 $v_{2\beta-1}=\big( \cost_{\beta} w_{2\beta-1} + \sit_{\beta} w_{2\beta} \big)\rho$, $v_{2\beta}
 =\big( -\sit_{\beta} w_{2\beta} + \cost_{\beta} w_{2\beta} \big)$, $\beta\in\{3,4\}$.
   In order to study the intersection of the strata in (\ref{C2 : 94}) and (\ref{C2 : 101}) with those one
   in (\ref{C2 : 45}), it is sufficient to compare the respective quaternionic pairs. Without loss of
   generality, we can compare just the first quaternionic pair. We list all of the possible cases:
   all of the possible cases\efl
 \begin{equation}   \label{C2 : 111}
 \begin{aligned}
1)  \quad &
\left( \begin{array}{cc}
z_1 & iz_1 \\
w_1 & iw_1 \\
\end{array} \right)=
\left( \begin{array}{cc}
v_1\epsilon & v_2\epsilon \\
v_1\sigma & v_2\sigma \\
\end{array} \right),  \quad \quad
2) \quad
\left( \begin{array}{cc}
z_1 & -iz_1 \\
w_1 & -iw_1 \\
\end{array} \right)=
\left( \begin{array}{cc}
v_1\epsilon & v_2\epsilon \\
v_1\sigma & v_2\sigma \\
\end{array} \right),     \\
&  \\
3) \quad &
\left( \begin{array}{cc}
z_1 & iz_1 \\
w_1 & iw_1 \\
\end{array} \right)=
\left( \begin{array}{cc}
-v_1\overline{\sigma} & -v_2\overline{\sigma} \\
v_1\overline{\epsilon} & v_2\overline{\epsilon}\\
\end{array} \right),    \quad
4)\quad
\left( \begin{array}{cc}
z_1 & -iz_1 \\
w_1 & -iw_1 \\
\end{array} \right)=
\left( \begin{array}{cc}
-v_1\overline{\sigma} & -v_2\overline{\sigma} \\
v_1\overline{\epsilon} & v_2\overline{\epsilon}\\
\end{array} \right),   \hspace{6 cm }   {}
\end{aligned}
 \end{equation}
\bfl where $v_1 =(\cost_1 z'_1 + \sit_1 z'_2)\rho$ and $v_2=(-\sit_1 z'_1 + \cost_1 z'_2)\rho$.
Let us rewrite the  equation $1)$ in (\ref{C2 : 111}) as follows\efl
\begin{equation}      \label{C2 : 112}
\left\{ \begin{array}{l}
z_1  = (\cost_1 z'_1 + \sit_1 z'_2)\epsilon\rho, \\
i z_1 = (-\sit_1 z'_1 + \cost_1 z'_2)\epsilon\rho, \\
w_1 =  (\cost_1 z'_1 + \sit_1 z'_2)\sigma\rho, \\
iw_1 =  (-\sit_1 z'_1 + \cost_1 z'_2)\sigma\rho, \\
\end{array} \right.   \hspace{3 cm }   {}
\end{equation}
\bfl then we have \efl
\begin{equation}     \label{C2 : 113}
\left\{\begin{array}{l}
(\cost_1 z'_1 + \sit_1 z'_2)\epsilon\rho = i(\sit_1 z'_1 - \cost_1 z'_2)\epsilon\rho, \\
(\cost_1 z'_1 + \sit_1 z'_2)\sigma\rho = i(\sit_1 z'_1 - \cost_1 z'_2)\sigma\rho, \\
 \end{array}\right.
 \quad \iff \quad
z'_1 =z'_2=0.
\end{equation}
\bfl but that is clearly impossible.
The remaing cases in (\ref{C2 : 111}) can be treated similarly.$\square$\efl

 \vspace{ 0,5 cm}

\section{ The $3-$Sasakian Orbifolds $\mathcal{M}^7(\Theta)$}

 \vspace{0,5 cm}

\bfl Now, we are going to study the singular locus on
 the $3-$Sasakian orbifold $\mathcal{M}^7(\Theta)$ introduced in the first chapter.
The fixed point
equations can be obtained from those one in (\ref{C2 : 8})
 by assuming $\rho=1$, then we obtain \efl
\begin{equation}   \label{C2 : 114}
\left( \begin{array}{cc}
\cost_{\alpha} & \sit_{\alpha} \\
-\sit_{\alpha}  & \cost_{\alpha} \\
\end{array} \right)
 \left(\begin{array}{cc}
z_{2{\alpha}-1} & w_{2{\alpha}-1} \\
z_{2{\alpha}} & w_{2{\alpha}} \\
\end{array} \right) =
{\left[ \begin{array}{c}
\left(\begin{array}{cc}
\epsilon & -\overline{\sigma}  \\
\sigma & \overline{\epsilon} \\
\end{array} \right)
 \left(\begin{array}{cc}
z_{2{\alpha}-1} & z_{2{\alpha}} \\
w_{2{\alpha}-1} & w_{2{\alpha}} \\
\end{array} \right)
\end{array}\right]}^T
\end{equation}

\bfl  where
$\lambda= \epsilon + j\sigma\in Sp(1)$,
  $(z_{\alpha}, w_{\alpha})\in \mathbb{C}$ and ${\alpha}\in\{1,2,3,4\}$.
    These equations can be easily rewritten as:\efl
\begin{equation}  \label{C2 : 115}
\overbrace{\left( \begin{array}{cc|cc}
0                          & -\sigma &  -\sit_{\alpha}  & \cost_{\alpha} - \overline{\epsilon}  \\
-\sigma     & 0                      & \cost_{\alpha} - \overline{\epsilon} & \sit_{\alpha} \\
\hline
-\sit_{\alpha}   & \cost_{\alpha} - \epsilon           &    0       &      \overline{\sigma}  \\
\cost_{\alpha} -\epsilon  & \sit_{\alpha} &  \overline{\sigma}   &   0    \\
\end{array} \right)}^{\tilde{M}_{\alpha}:=}
\left(\begin{array}{c}
z_{2\alpha-1} \\
z_{2\alpha}\\
w_{2\alpha-1} \\
w_{2\alpha} \\
\end{array} \right) =
\left(\begin{array}{c}
0 \\
0\\
0\\
0 \\
\end{array} \right),
\end{equation}
 \bfl where $\alpha\in \{1,2,3,4\}$.
Then we have the followig results  \efl

\begin{lem} \bfl For each $\alpha\in \{1,2,3,4\}$, the matrix  $\tilde{M}_{\alpha}$ in $(\ref{C2 : 115})$
is such that $det\, \tilde{M}_{\alpha}=0$ if and only if it holds either
\begin{equation}      \label{C2 : 116}
e^{i\theta_{\alpha}}= Re( \epsilon ) + i\sqrt{ Im( \epsilon)^2 + |\sigma |^2}\,\,\, \mathrm{or}
\,\,\,  e^{-i\theta_{\alpha}}= Re(\epsilon) + i\sqrt{ Im( \epsilon)^2 + |\sigma |^2}.
\end{equation}
\efl
\end{lem}

\vspace{0,3 cm}

\bfl Lemma $2.7$ can be proved like lemma $2.1$. Moreover, the fixed point
 equations in
$(1.19)$ give that the action of $G$ on $N(\Theta)$ is locally free.
 Notice that, if exactly one of the identities in (\ref{C2 : 116})
  is verified, then $rank\, \tilde{M}_{\alpha}=3$.
Instead, when there exists $\alpha\in \{1,2,3,4\}$ such that both the identities
 in (\ref{C2 : 116}) hold, then we have that $\tilde{M}_{\alpha}=0 $
 for each $\alpha$ and this case
  corresponds to the non-effective subgroup. Then,
   the only interesting cases are realized when
  exactly one of the equations in (\ref{C2 : 116}) is verified and
   $rank\, \tilde{M}_{\alpha}=3$. \efl

\bfl
In this case, in order to describe the singular strata $\overset{-}{S}$ for the
$G-$action on
 $N(\Theta)$,  we use a different approach with respect to the case of $\mathcal{Z}^6(\Theta)$.
 In fact, we compute directly the solutions associatetd to the eigenvalue problem
 in (\ref{C2 : 115}), and then we describe the singular strata.
 We state now a proposition, it will be proved after corollary $2.5$.\efl

\begin{prop}
\bfl  The solutions $\zw$ of the eigenvalue problem associated to
the equations  $(\ref{C2 : 115})$
as quaternionic pairs have either \efl

  \begin{equation} \label{C2 : 117} {
\begin{aligned}
 1)\,\, &\quad \mathrm{if} \,\,\,\, e^{i\theta_{\alpha}}= Re( \epsilon ) +
 i\sqrt{ Im^2(\epsilon) + |\sigma |^2},\quad  \mathrm{then\,\, we\,\, have:}  \quad
 \tilde{V}_1^{\alpha}:=\{z_{2\alpha-1}, z_{2\alpha}\in \mathbb{C}\, | \\
& \, (z_{2\alpha-1},  z_{2\alpha},
\frac{e^{i\delta}  \big(-z_{2\alpha} +iz_{2\alpha-1}\cos\, \varphi  \big)}{sin \, \varphi},
\frac{e^{i\delta}  \big(z_{2\alpha-1} +iz_{2\alpha}\cos\, \varphi  \big)}{ sin \, \varphi})
\in \mathbb{C}^4\},  \\
 \mathrm{or} &  \\
 2)\,\, &\quad \mathrm{if} \,\,\,\, e^{-i\theta_{\alpha}}= Re (\epsilon ) +
 i\sqrt{ Im^2(\epsilon) + |\sigma |^2},\quad  \mathrm{then\,\, we\,\, have:}  \quad
 \tilde{V}_2^{\alpha}:=\{z_{2\alpha-1}, z_{2\alpha}\in \mathbb{C}\, | \\
& \, (z_{2\alpha-1},  z_{2\alpha},
-\frac{e^{i\delta}  \big(\zap +i\zad\cos\, \varphi  \big)}{sin \, \varphi},
\frac{e^{i\delta} \big(\zad -i\zap\cos\, \varphi  \big)}{  sin \, \varphi})
\in \mathbb{C}^4\},
\end{aligned}     } \hspace{ 3 cm}
\end{equation}
\bfl where  $\epsilon = \cost +i(\sit cos\, \varphi)$,
$\sigma =  \sit sin \, \varphi cos\, \delta +i(\sit
   sin \, \varphi sin\, \delta )= $
   $ \sit sin \, \varphi e^{i\delta}$ and $ sin \, \varphi\ne 0$.  Moreover,
  we have that $\tilde{V}_1^{\alpha}$ and
   $\tilde{V}_2^{\alpha}$ are $G-invariant$. \efl
\end{prop}

\vspace{0,3 cm}

\bfl In the following, by assuming $ sin \, \varphi\ne 0$, we denote as    \efl
\begin{equation} \label{C2 : 117a}
\begin{aligned}
&  (w'_{2\alpha-1},
 w'_{2\alpha }):=\frac{e^{i\delta} }{ sin\, \varphi}
 \big( -z_{2\alpha} +iz_{2\alpha -1}cos\, \varphi, \, z_{2\alpha-1} +iz_{2\alpha}\cos\, \varphi  \big), \\
and \,\,\,\, & \\
  &(w''_{2\alpha-1},  w''_{2\alpha})=: \frac{e^{i\delta} }{ sin\, \varphi } \big(-\zap -
i\zad\cos\, \varphi,\, \zad -i\zap\cos\, \varphi  \big).  \hspace{ 5 cm}\\
\end{aligned}
\end{equation}

\vspace{0,5 cm}


\bfl Then we have \efl

\begin{cor}\bfl Consider the action of $G$ on $\mathbb{H}^8$.
Under the hypotheses of lemma $2.8$ and proposition $2.7$
 and, up to change
 $(z_1, z_2, w'_1, w'_2)$ with $(z_1, z_2, w''_1, w''_2)$,  the singularities relative to
 $\mathcal{M}^7(\Theta)$
 come from  the intersection of $N(\Theta)$ with the following strata
  relative to the action of $G$ on $\mathbb{H}^8$ \efl
 \begin{align*}
 & 1) \,\,\, \quad  \overline{S}{}^{1234}:= 
 \bigg\{
\left( \begin{array}{cc|cc|cc|cc}
z_1 & z_2  & z_3 & z_4   &    z_5 & z_6        & z_7    &  z_8   \\
w'_1 & w'_2  &  w'_3  & w'_4     &   w'_5   &  w'_6          & w'_7  &  w'_8\\
\end{array}\right) \bigg\},   \\
& 2) \,\,\, \quad \overline{S}{}^{123}_{4}:= 
 \bigg\{
\left( \begin{array}{cc|cc|cc|cc}
z_1 & z_2  & z_3 & z_4   &    z_5 & z_6        & z_7    &  z_8   \\
w'_1 & w'_2  &  w'_3  & w'_4     &   w'_5   &  w'_6          & w''_7  &  w''_8\\
\end{array}\right) \bigg\},   \hspace{10 cm}\\
&3) \,\,\,\quad \overline{S}{}^{124}_{3}:= 
 \bigg\{
\left( \begin{array}{cc|cc|cc|cc}
z_1 & z_2  & z_3 & z_4   &    z_5   & z_6           & z_7   & z_8  \\
w'_1 & w'_2  &  w'_3  & w'_4     &   w''_5  &  w''_6   & w'_7  &  w'_8\\
\end{array}\right) \bigg\},   \\
 & 4)\,\,\, \quad \overline{S}{}^{134}_{2}:= 
 \bigg\{
\left( \begin{array}{cc|cc|cc|cc}
z_1 & z_2  & z_3 & z_4   &    z_5 & z_6        & z_7    &  z_8   \\
w'_1 & w'_2 & w''_3  & w''_4    &  w'_5  &  w'_6   & w'_7  &  w'_8\\
\end{array}\right) \bigg\},   \\
& 5)\,\,\,\quad \overline{S}{}^{12}_{34}:= 
\bigg\{\left( \begin{array}{cc|cc|cc|cc}
z_1 & z_2  & z_3 & z_4   &    z_5 & z_6        & z_7    &  z_8   \\
w'_1 & w'_2  &  w'_3  & w'_4    &   w''_5  &  w''_6        & w''_7  &  w''_8\\
\end{array}\right) \bigg\},   \\
& 6)\,\,\,\quad \overline{S}{}^{13}_{24}:= 
\bigg\{
\left( \begin{array}{cc|cc|cc|cc}
z_1 & z_2  & z_3 & z_4   &    z_5 & z_6        & z_7    &  z_8   \\
w'_1 & w'_2 & w''_3  & w''_4     &  w'_5  &  w'_6          & w''_7  &  w''_8\\
\end{array}\right) \bigg\},   \\
& 7)\,\,\,\quad  \overline{S}{}^{14}_{23}:= 
\bigg\{
\left( \begin{array}{cc|cc|cc|cc}
z_1 & z_2  & z_3 & z_4   &    z_5 & z_6        & z_7    &  z_8   \\
w'_1 & w'_2 & w''_3  & w''_4     &   w''_5 & w''_6   & w'_7  &  w'_8\\
\end{array} \right)    \bigg\}, \\
\end{align*}
 \begin{equation} \label{C2 : 118} {
 \begin{aligned}
& 8)\,\,\,\quad  \overline{S}{}^{1}_{234} :=
\bigg\{
\left( \begin{array}{cc|cc|cc|cc}
z_1 & z_2  & z_3 & z_4   &    z_5 & z_6        & z_7    &  z_8   \\
w'_1 & w'_2 & w''_3  & w''_4     &   w''_5 & w''_6        & w''_7  &  w''_8\\
\end{array}\right) \bigg\}.
 \end{aligned}    }    \hspace{4,5 cm}    {}
 \end{equation}
\end{cor}

\bfl The proof of this corollary is a straightforward consequence of lemma $2.7$ and proposition $2.8$. \efl

\bfl $\mathbf{Proof}(proposition\,\, 2.8).$  Assume
 $sin \, \varphi \ne 0$. 
and consider the fixed point equations in (\ref{C2 : 114}). Then, by using lemma $2.7$,
 it follows that
 $det\, \tilde{M}_{\alpha}=0$ if and only if it either holds
 $e^{i\theta_{\alpha}}= Re( \epsilon ) +
 i\sqrt{ Im( \epsilon)^2 + |\sigma |^2}$ or
 $ e^{-i\theta_{\alpha}}= Re(\epsilon) +
 i\sqrt{ Im( \epsilon)^2 + |\sigma |^2}.$
Note that, when $\sit_{\alpha}=0$, it follows that the fixed point
 equations (\ref{C2 : 114}) decsribe
the non effective subgroup.
 Now, suppose the first relation in (\ref{C2 : 116}) holds, then we have that \efl
  \begin{equation}  \label{C2 : 119}
 \left\{ \begin{array}{l}
 \epsilon = \cost + i\sit cos\, \varphi , \\
 e^{i\theta_{\alpha}}= Re( \epsilon ) + i\sqrt{ Im( \epsilon)^2 + |\sigma |^2}
  = \cost + i\sit, \\
 \end{array} \right.
  \end{equation}
\bfl hence $\epsilon = \cost_{\alpha} + i\sit_{\alpha} cos\, \varphi$, and  by substituting
this relation in $\tilde{M}_{\alpha}$ we obtain\efl
\begin{equation}    \label{C2 : 120}
\tilde{M}_{\alpha}=
\sit_{\alpha}\left( \begin{array}{cc|cc}
0          & -sin \, \varphi e^{i\delta} &  -1  & icos\, \varphi   \\
-sin \, \varphi e^{i\delta}    & 0                      & icos\, \varphi & 1 \\
\hline
-1   &    -icos\, \varphi       &    0       &      sin \, \varphi e^{-i\delta}   \\
-icos\, \varphi  & 1&  sin \, \varphi e^{-i\delta}    &   0    \\
\end{array} \right), \hspace{ 3 cm }
\end{equation}
 \bfl where $\alpha\in \{1,2,3,4\}$. As a consequence the associated eigenvalue
 problem becomes   \efl
\begin{equation}        \label{C2 : 121}
\left( \begin{array}{cc|cc}
0          & -sin \, \varphi e^{i\delta} &  -1  & i cos\, \varphi   \\
-sin \, \varphi e^{i\delta}    & 0                      & icos\, \varphi  & 1 \\
\hline
-1   &    -i cos\, \varphi     &    0       &      sin \, \varphi e^{-i\delta}   \\
-i cos\, \varphi   & 1 &  sin \, \varphi e^{-i\delta}    &   0    \\
\end{array} \right)\left(\begin{array}{c}
z_{2\alpha-1} \\
z_{2\alpha}\\
w_{2\alpha-1} \\
w_{2\alpha} \\
\end{array} \right) =
\left(\begin{array}{c}
0 \\
0\\
0\\
0 \\
\end{array} \right).
\end{equation}
\bfl From the first two equations in (\ref{C2 : 121}) we get\efl
\begin{equation}          \label{C2 : 122}
\left\{ \begin{array}{l}
 -sin \, \varphi e^{i\delta}\zap - \wad +i\cos\, \varphi\wap=0, \\
 -sin \, \varphi e^{i\delta}\zad +i\cos\, \varphi\wad + \wap=0, \\
\end{array} \right.  \Rightarrow
\left\{ \begin{array}{l}
\wad=\frac{e^{i\delta} sin \, \varphi \big(-z_{2\alpha} +iz_{2\alpha -1}\cos\, \varphi  \big)}{ 1- \cos^2 \, \varphi},\\
\wap=\frac{e^{i\delta} sin \, \varphi \big(z_{2\alpha -1} +iz_{2\alpha}\cos\, \varphi  \big)}{ 1- \cos^2 \, \varphi}. \\
\end{array} \right.
\end{equation}
\bfl Instead, when it holds $ e^{-i\theta_{\alpha}}= Re( \epsilon ) +
i\sqrt{ Im( \epsilon)^2 + |\sigma |^2}$, it follows that \efl
 \begin{equation}       \label{C2 : 123}
 \left\{ \begin{array}{l}
 \epsilon = \cost + i\sit cos\, \varphi , \\
 e^{-i\theta_{\alpha}}= \cost + i\sit, \\
 \end{array} \right.
  \end{equation}
 \bfl hence $\cost + i\sit=cos(-\theta_{\alpha})+ i sin (-\theta_{\alpha})$ and
 $\epsilon = \cost_{\alpha} - i\sit_{\alpha}cos\,\varphi$.
  By substituting the relations (\ref{C2 : 123}) in the equations
 (\ref{C2 : 115}) we get\efl
\begin{equation}
\left( \begin{array}{cc|cc}
0          & -sin \, \varphi e^{i\delta} &  -1  & -i\cos\, \varphi   \\
-sin \, \varphi e^{i\delta}    & 0                      & -i\cos\, \varphi & 1 \\
\hline
-1   &    i\cos\, \varphi     &    0       &      sin \, \varphi e^{-i\delta}   \\
i\cos\, \varphi  & 1 &  sin \, \varphi e^{-i\delta}    &   0    \\
\end{array} \right)\left(\begin{array}{c}
z_{2\alpha-1} \\
z_{2\alpha}\\
w_{2\alpha-1} \\
w_{2\alpha} \\
\end{array} \right) =
\left(\begin{array}{c}
0 \\
0\\
0\\
0 \\
\end{array} \right).
\end{equation}
\bfl  and  \efl
\begin{equation}      \label{C2 : 124}
\left\{ \begin{array}{l}
 -sin \, \varphi e^{i\delta}\zap - \wad -i\cos\, \varphi\wap=0, \\
 -sin \, \varphi e^{i\delta}\zad -i\cos\, \varphi\wad + \wap=0, \\
\end{array} \right.  \Rightarrow
\left\{ \begin{array}{l}
\wad=-\frac{e^{i\delta} sin \, \varphi \big(\zap +i\zad\cos\, \varphi  \big)}{ 1- \cos^2 \, \varphi},  \\
\wap= \frac{e^{i\delta} sin \, \varphi \big(\zad -i\zap\cos\, \varphi  \big)}{ 1- \cos^2 \, \varphi}). \\
\end{array} \right.
\end{equation}
\bfl In the case when $\sigma\equiv 0$, that is $sin\, \varphi =0$, the solutions $\zw$ of the eigenvalue problem
(\ref{C2 : 115}) have quaternionic pair of type $(\zad, \pm i\zad, \wad, \mp i\wad)$
like those described in  (\ref{C2 : 46}). Then by using proposition $2.6$ we can easil exclude these points.
As a consequence, we obtain exactly the strata described in (\ref{C2 : 118} ).
Moreover, if we consider the fixed point equations in (\ref{C1: 17}), (\ref{C1: 18}) and (\ref{C1: 20}) it
follows that \efl

 \begin{equation}       \label{C2 : 125}
Im(\lambda)= \sit_1\frac{  u_2\overline{u}_1- u_1\overline{u}_2 }{(|u_1 |^2 + |u_2 |^2 )}=
 \sit_{\alpha} \frac{\um\barun - \un\barum}{ |\un |^2+ |\um |^2},\,\, \alpha\in\{1,2,3,4\}
 \end{equation}
\bfl and the fixed point equations become\efl
 \begin{equation}           \label{C2 : 126}
 \left\{ \begin{array}{l}
 e^{i(\theta_1\pm \theta_2)}=1, \\
 e^{i(\theta_1\pm \theta_3)}=1, \\
 e^{i(\theta_1\pm \theta_4)}=1, \\
  e^{\pm i\theta_1}= Re( \epsilon ) + i\sqrt{ Im( \epsilon)^2 + |\sigma |^2}, \\
  iIm(\epsilon) + j\sigma = \sit_1\frac{  u_2\overline{u}_1- u_1\overline{u}_2 }{(|u_1 |^2 + |u_2 |^2 )}.  \\
  \end{array} \right.
 \end{equation}
\bfl According to the description of the locally free action of $G$ on $N(\Theta)$,
 given in the first chapter,
 the isotropy subgroup described by the system of equations in (\ref{C2 : 126}) only depends on one
of the determinants $\sq$, which is associated to the subsystem of the first three equations in
(\ref{C2 : 126}). $\square$\efl

\vspace{0,5 cm}


 \subsection{Some computations on the singular $G-$strata of $\mathcal{M}^7(\Theta)$ }

\vspace{1 cm}

\bfl Now, we want to study the intersection of $N(\Theta)$ with
 the $G-$strata described in  (\ref{C2 : 118}),
then we will compute the moment map
 equations on these strata. In order to do this,  look at
  the moment map equations of the $Sp(1)$ and $T^3_{\Theta}$ actions.   \efl

\begin{prop}
\bfl Let us consider the $G-$strata $\overline{S}$ of $\mathbb{H}^8$
 listed in $(\ref{C2 : 118})$. Then, each of the singular $G-$strata $\widetilde{S}=
\overline{S}\cap N(\Theta)\subset N(\Theta)$ gives rise
a singular point on the $3-$
Sasakian orbifold $M^7(\Theta )$.   \efl
\end{prop}

\bfl $\mathbf{Proof.}$ First, we describe the case of $\overline{S}{}^{123}_{4}$, then the cases
 $\overline{S}{}^{124}_{3}$, $\overline{S}{}^{134}_{2}$ and $\overline{S}{}^{1}_{234}$
 can be treated analogously.
 Recall the moment map equations
  relative the $3-$torus are: \efl
 \begin{equation}    \label{C2 : 128}
 \left\{ \begin{array}{l}
 \sum_{\alpha=1}^4 d_{\alpha} Im(z_{2\alpha-1}\overline{z_{2\alpha}} +
w_{2\alpha-1}\overline{w_{2i}})=0, \\
\sum_{\alpha=1}^4 d_{\alpha} (z_{2\alpha-1}w_{2\alpha}-  z_{2\alpha} w_{2\alpha-1} ) =0,
\end{array}\right.
 \end{equation}
\bfl where $d_{\alpha},$ $\alpha\in\{1,2,3,4\}$ are the weights.
We are going to compute these equations on the stratum $\overline{S}{}^{123}_{4}$.
 Then, by considering the quaternionic pair in the case $1)$ of (\ref{C2 : 117}) it follows that\efl
\begin{equation}      \label{C2 : 129}
\begin{aligned}
z_{2\alpha - 1}\overline{z}_{2\alpha}  + w_{2\alpha - 1}\overline{w}_{2\alpha}=   &
\frac{1}{1 -\cos^2\varphi}\big( (1 -cos^2\,\varphi)\zad\bzap + \\
& (-\zap +
 i\zad\cos\, \varphi  )(\bzad -i\bzap\cos\, \varphi )\big)   \\
=\frac{ 1}{1 -\cos^2\varphi}  \big( & (\zad\bzap- \bzad\zap)  +
i(|\zad |^2 + |\zap |^2)cos\, \varphi \big).  
\end{aligned}
\end{equation}
\bfl Instead, for  the quaternionic pair in  the case $2)$ of (\ref{C2 : 117}) we get   \efl
\begin{equation}    \label{C2 : 130}
\begin{aligned}
z_{2\alpha - 1}\overline{z}_{2\alpha}  + w_{2\alpha - 1}\overline{w}_{2\alpha}=  &
\frac{1}{1 -\cos^2\varphi}\big(  (1 -\cos^2\varphi)\zad\bzap -  \\
& (\zap +i\zad\cos\, \varphi  )(\bzad +i\bzap\cos\, \varphi )\big)      \\
= \frac{ 1}{1 -\cos^2\varphi}\big( & (\zad\bzap  -\bzad\zap ) -
i(|\zad |^2 + |\zap |^2)cos\, \varphi \big).
\end{aligned}
\end{equation}

\bfl Then, we compute $\zad w_{2\alpha}- \zap\wad$.
Consider the quaternionic pair in the case $1)$ of (\ref{C2 : 117})\efl
 \begin{equation}     \label{C2 : 131}
\begin{aligned}
\zad w_{2\alpha}- \zap\wad =   &
 \frac{sin \, \varphi e^{i\delta}}{1 -\cos^2\varphi}
  \big( \zad( \zad + i\zap cos\, \varphi)-  \\
  \zap( -\zap + & i\zad cos\, \varphi)  \big) =
     \frac{sin \, \varphi e^{i\delta}}{1 -\cos^2\varphi}
  \big( (\zad)^2 + (\zap)^2 \big).  \hspace{1 cm}
\end{aligned}
\end{equation}

\bfl In the case $2)$ of (\ref{C2 : 117}) we get\efl
\begin{equation}           \label{C2 : 132}
\begin{aligned}
\zad w_{2\alpha}- \zap\wad =    &
 \frac{sin \, \varphi e^{i\delta}}{1 -\cos^2\varphi}
  \big( \zad( \zad - i\zap cos\, \varphi)+ \\
  \zap( \zap +  & i\zad cos\, \varphi)  \big) =
     \frac{sin \, \varphi e^{i\delta}}{1 -\cos^2\varphi}
  \big( (\zad)^2 + (\zap)^2 \big).    \hspace{1 cm}
\end{aligned}
\end{equation}

\bfl Now, we are ready to write down the moment map relative to the $3$-torus $T^3_{\Theta}$
in the case  of   $\overline{S}{}^{123}_{4}$. Then we have  \efl
 \begin{equation}        \label{C2 : 133}
 \begin{aligned}
 \left\{\begin{array}{l}
2 \sum_{\alpha=1}^4
 d_{\alpha}  Im(\bzap \zad)   + cos \, \varphi \big( \sum_{\alpha=1}^3 d_{\alpha} ( |\zad |^2 + |\zap |^2 )
   -d_4(|z_7|^2+ |z_8|^2 ) \big) =0, \\
  \sum_{\alpha=1}^4 d_{\alpha}\big( (\zad)^2 + (\zap)^2 \big) \big) =0.
 \end{array} \right.
 \end{aligned}
 \end{equation}

\bfl Similarly we can compute the moment map equations
linked to the subgroup $Sp(1)\subset G$.  Note that \efl
 \begin{equation}               \label{C2 : 134}
\left\{\begin{array}{l}
\sum_{{\alpha}=1}^4 \overline{(z_{\alpha} + jw_{\alpha})}i (z_{\alpha} + jw_{\alpha})=
\sum_{{\alpha}=1}^4 i\big( |z_{\alpha}|^2 -  |w_{\alpha}|^2\big) + 2k z_{\alpha}w_{\alpha}=0, \\
\sum_{{\alpha}=1}^4 \overline{(z_{\alpha} + jw_{\alpha})}j (z_{\alpha} + jw_{\alpha})=
\sum_{{\alpha}=1}^4  j\big( (z_{\alpha})^2 + (w_{\alpha})^2 \big) + \big( \bwa\za -\bza\wa \big)=0, \\
 \sum_{{\alpha}=1}^4 \overline{(z_{\alpha} + jw_{\alpha})}k (z_{\alpha} + jw_{\alpha})=
\sum_{{\alpha}=1}^4  k \big( (z_{\alpha})^2 - (w_{\alpha})^2 \big) \big)-
i\big( \bwa\za +\bza\wa \big)=0, \\
\end{array}\right.
\end{equation}

\bfl We rewrite now, the first equation in (\ref{C2 : 134}). From the quaternionic pair
described at the point $1)$ of (\ref{C2 : 117}) we get\efl
\begin{equation}        \label{C2 : 135}
\begin{aligned}
\zad\wad + \zap\wap =   & \frac{ sin \, \varphi e^{i\delta}  }{1 -\cos^2\varphi}
 \big( \zad( -\zap + i\zad cos\, \varphi)+  \\
 \zap(  \zad &  + i\zap cos\, \varphi)  \big) =
     i\frac{ \sin\, \varphi  cos \, \varphi e^{i\delta}   }{1 -\cos^2\varphi}
  \big( (\zad)^2 + (\zap)^2 \big).
\end{aligned}
\end{equation}

\bfl Analogously, for the case $2)$ of (\ref{C2 : 117}) it holds \efl

\begin{equation}         \label{C2 : 136}
\begin{aligned}
\zad\wad + \zap\wap =   & \frac{ sin \, \varphi e^{i\delta}  }{1 -\cos^2\varphi}
 \big( \zad( -\zap - i\zad cos\, \varphi)+  \\
 \zap(  \zad &  - i\zap cos\, \varphi)  \big) =
     -i\frac{ \sin\, \varphi  cos \, \varphi e^{i\delta}   }{1 -cos^2\varphi}
  \big( (\zad)^2 + (\zap)^2 \big).
\end{aligned}
\end{equation}

\bfl By reading $|\zad |^2 -  |\wad |^2$ on the quaternionic pair in the case
 $1)$ of (\ref{C2 : 117}) we get\efl
\begin{equation}      \label{C2 : 137}
\begin{aligned}
|\zad |^2 -  |\wad |^2 =  \frac{1 }{1 -\cos^2\varphi}\big(
 |\zad |^2 -|\zap |^2 & -2|\zad |^2 cos^2\,\varphi + \\
& i(\bzad\zap - \zad\bzap ) cos\, \varphi\big),
\end{aligned}
\end{equation}  \label{C2 : 138}
\bfl and \efl
\begin{equation}
\begin{aligned}
|\zap |^2 -  |\wap |^2 =  \frac{1}{1 -\cos^2\varphi}\big(
 |\zap |^2 - |\zad |^2 & -2|\zap |^2  cos^2\,\varphi + \\
 & i(\bzad\zap - \zad\bzap ) cos\,\varphi \big).
 \end{aligned}
\end{equation}
 \bfl
 By adding (\ref{C2 : 137}) and $(2.170)$ we obtain\efl

\begin{equation}   \label{C2 : 139}
\begin{aligned}
 \big( |\zad |^2 - & |\wad |^2 \big)  + \big( |\zap |^2 -  |\wap |^2 \big)=  \\
&  \frac{2icos\, \varphi }{1 -\cos^2\varphi}\big(  i\big( |\zad |^2  +  |\zap |^2  \big)cos\, \varphi
+ \big( \bzad\zap -\zad\bzap  \big)  \big), \hspace{2 cm }
\end{aligned}
\end{equation}

 \bfl The case  $2)$ of (\ref{C2 : 117}) can be obtained in the same way, and we get \efl

\begin{equation}       \label{C2 : 140}
\begin{aligned}
 \big( |\zad |^2 - & |\wad |^2 \big)  + \big( |\zap |^2 -  |\wap |^2 \big)=  \\
&  \frac{2icos\, \varphi }{1 -\cos^2\varphi}\big(  i\big( |\zad |^2  +  |\zap |^2  \big)cos\, \varphi
- \big( \bzad\zap -\zad\bzap   \big)  \big), \hspace{2 cm }
\end{aligned}
\end{equation}

\bfl Then, the first equation in (\ref{C2 : 134}) can be rewritten as follows\efl
\begin{equation}           \label{C2 : 141}
\left\{ \begin{array}{l}
cos \, \varphi \sum_{\alpha=1}^4 ( |\zad |^2 + |\zap |^2 ) +
 2 \big( \sum_{\alpha=1}^3
   Im(\zad\bzap )
   - Im(z_7\overline{z}_8) \big)  =0, \\
   \sum_{\alpha=1}^3 \big( (\zad )^2 + (\zap )^2 \big) \big)
    - \big( (z_7)^2 + (z_8)^2\big) =0.
\end{array} \right.
\end{equation}

\bfl For the second equation in (\ref{C2 : 134}),  with respect to the quaternionic
pairs in the case $1)$ of (\ref{C2 : 117}), we have that \efl
\begin{equation}       \label{C2 : 142}
\begin{aligned}
(\zad)^2 + (\wad)^2  = & \frac{1 }{1- cos^2 \, \varphi}\big( (\zad)^2 (1 - \cos^2 \, \varphi ) +  \\
 & e^{2i\delta}\big( (\zap)^2 -2i\zad\zap cos\, \varphi  -(\zad)^2cos^2\, \varphi \big) \big),  \hspace{ 2 cm}
\end{aligned}
\end{equation}
 \bfl and \efl
 \begin{equation}    \label{C2 : 143}
\begin{aligned}
(\zap)^2 + (\wap)^2  = & \frac{1 }{1- cos^2 \, \varphi}\big( (\zap)^2 (1 - \cos^2 \, \varphi ) +  \\
 & e^{2i\delta}\big( (\zad)^2 +2i\zad\zap cos\, \varphi  -(\zap)^2cos^2 \varphi \big) \big),  \hspace{ 2 cm}
\end{aligned}
\end{equation}
\bfl then it follows \efl
\begin{equation}      \label{C2 : 144}
(\zad)^2 + (\wad)^2 +  (\zap)^2 + (\wap)^2 = \big( (\zad)^2 +(\zap)^2 \big)(1 + e^{2i\delta}).    \hspace{ 2 cm}
\end{equation}
\bfl The case linked to the quaternionic pair at the point
 $2)$ in (\ref{C2 : 117}), is exactly equal.
  Then, also in this case it holds\efl
\begin{equation}        \label{C2 : 145}
(\zad)^2 + (\wad)^2 +  (\zap)^2 + (\wap)^2 = \big( (\zad)^2 +(\zap)^2 \big)(1 + e^{2i\delta}).    \hspace{ 2 cm}
\end{equation}
\bfl Then, consider $\bzad \wad - \zad \bwad$, in the case of the
quaternionic pair at $1)$ of  (\ref{C2 : 117}) we have\efl
\begin{equation}     \label{C2 : 146}
\bzad \wad - \zad \bwad=
\frac{sin \, \varphi e^{i\delta} }{1- cos^2 \, \varphi} \big( 2i |\zad |^2 cos\, \varphi +
 \big( \zad\bzap -\bzad\zap  \big) \big),
\end{equation}
\bfl and  \efl
\begin{equation}   \label{C2 : 147}
\bzap \wap - \zap \bwap=
\frac{sin \, \varphi e^{i\delta} }{1- cos^2 \, \varphi} \big( 2i |\zap |^2 cos\, \varphi +
 \big( \zad\bzap -\bzad\zap  \big) \big),
\end{equation}
\bfl then \efl
\begin{equation}   \label{C2 : 148}
\begin{aligned}
( \bzad \wad & - \zad  \bwad ) + ( \bzap \wap - \zap \bwap )=  \\
 &\frac{2sin\, \varphi e^{i\delta} }{1- cos^2 \, \varphi}
 \big( i \big( |\zap |^2+ |\zad |^2 \big) cos\, \varphi + \big( \zad\bzap -\bzad\zap  \big) \big).
 \end{aligned}
\end{equation}
\bfl Instead, in the case of the quaternionic pair at $2)$ of  (\ref{C2 : 117}) we obtain \efl
\begin{equation}      \label{C2 : 149}
\begin{aligned}
( \bzad \wad & - \zad  \bwad ) + ( \bzap \wap - \zap \bwap )=  \\
 &\frac{2sin\, \varphi e^{i\delta} }{1- cos^2 \, \varphi}
 \big( -i \big( |\zap |^2+ |\zad |^2 \big) cos\, \varphi + \big( \zad\bzap -\bzad\zap  \big) \big).
 \end{aligned}
\end{equation}
\bfl Thus the second equation in (\ref{C2 : 134}) can be rewritten as follows \efl
\begin{equation}        \label{C2 : 150}
\left\{ \begin{array}{l}
2 \sum_{\alpha=1}^4
   Im(\bzap \zad)   +cos \, \varphi \big( \sum_{\alpha=1}^3  ( |\zad |^2 + |\zap |^2 )
   -(|z_7|^2+ |z_8|^2 ) \big) =0, \\
   (1 + e^{2i\delta})\sum_{\alpha=1}^4 \big( (\zad )^2 + (\zap )^2 \big) \big) =0.
\end{array} \right.
\end{equation}

 \bfl Similarly we can describe the last equation in (\ref{C2 : 134}), so we obtain\efl
\begin{equation}       \label{C2 : 151}
\left\{ \begin{array}{l}
 \sum_{\alpha=1}^4
 2Im(\bzap \zad)   + cos \, \varphi \big( \sum_{\alpha=1}^3( |\zad |^2 + |\zap |^2 )
   -(|z_7|^2+ |z_8|^2 ) \big) =0, \\
   (1 - e^{2i\delta} )\sum_{\alpha=1}^4 \big( (\zad )^2 + (\zap )^2 \big) =0.
\end{array} \right.
\end{equation}

 \bfl Then the equations (\ref{C2 : 141}), (\ref{C2 : 150}) and (\ref{C2 : 151})
  can be rewritten as follows\efl

 \begin{equation}       \label{C2 : 152}
Sp(1)\,\, \Rightarrow\,\,
\left\{ \begin{array}{l}
2 \sum_{\alpha=1}^3
 Im(\bzap \zad)   + cos \, \varphi \big( \sum_{\alpha=1}^3( |\zad |^2 + |\zap |^2 ) \big) =0,  \\
\sum_{\alpha=1}^3 \big( (\zad )^2 + (\zap )^2 \big) \big) =0, \\
(z_7)^2 + (z_8)^2 =0,     \\
cos\, \varphi = \frac{Im(\overline{z}_8z_7)}{|z_7 |^2 + |z_8|^2}.  \\
 \end{array} \right.
 \end{equation}

 \bfl As a consequence, the moment map equations of the
  $3-$torus $T^3_{\Theta}$ in (\ref{C2 : 128}) become
 \begin{equation}          \label{C2 : 153}
 T^3_{\Theta}\,\,  \Rightarrow\,\,
 \left\{ \begin{array}{l}
 \sum_{\alpha=1}^3 p_{\alpha} \big( 2Im(\bzap\zad) + ( |\zad |^2 + |\zap |^2 )cos\, \varphi  \big)=0,  \\
 \sum_{\alpha=1}^3 q_{\alpha} \big( 2Im(\bzap\zad) + ( |\zad |^2 + |\zap |^2 )cos\, \varphi  \big)=0,  \\
 \sum_{\alpha=1}^3 l_{\alpha} \big( 2Im(\bzap\zad) + ( |\zad |^2 + |\zap |^2 )cos\, \varphi  \big)=0,  \\
  Im(z_7\overline{z}_8)  - ( |z_7|^2 + |z_8|^2 )cos\, \varphi= 0,  \\
 \end{array} \right.
\end{equation}
 and
\begin{equation}         \label{C2 : 154}
  T^3_{\Theta}\,\,\,\, \Rightarrow\,\,
 \left\{ \begin{array}{l}
 \sum_{\alpha=1}^3 p_{\alpha}\big( (\zad )^2 + (\zap )^2 \big)  =0, \\
 \sum_{\alpha=1}^3 q_{\alpha}\big( (\zad )^2 + (\zap )^2 \big)  =0, \\
 \sum_{\alpha=1}^3 l_{\alpha}\big( (\zad )^2 + (\zap )^2 \big)  =0, \\
 (z_7)^2 + (z_8)^2 =0.    \\
 \end{array} \right.  \hspace{3,5 cm} {}
 \end{equation}
  For each of the two system of equations in  (\ref{C2 : 153}) and
 (\ref{C2 : 154}) consider the subsystem formed respectively by the first three
 equations. These two subsystem of equations have the same matrix of coefficients,
 hence they also have the same
 determinant which is $\Delta_{123}\ne 0$. Then it holds
 \begin{equation}         \label{C2 : 155}
  2Im(\bzap\zad) + ( |\zad |^2 + |\zap |^2 )cos\, \varphi =0, \quad   \alpha\in \{1,2,3\},
  \end{equation}
   that is \efl
  \begin{equation}          \label{C2 : 155a}
  \begin{aligned}
 &  cos\, \varphi  = -\frac{2Im(\bzap\zad)}{|\zad |^2 + |\zap |^2},  \quad  \alpha\in \{1,2,3\}, \\
 \mathrm{ and} \quad\quad  & \\
 & (\zad )^2 + (\zap )^2  =0,\quad \alpha\in \{1,2,3,4\}. \hspace{10 cm} \\
 \end{aligned}
 \end{equation}
 \bfl Note that the moment map equations of $Sp(1)$
 in  (\ref{C2 : 152}) can be obtained from those one in (\ref{C2 : 153}) and
 (\ref{C2 : 154}) relative to $T^3_{\Theta}$.
 In order to describe the intersection of $\overline{S}{}^{123}_{4}$ with $N(\Theta)$
  we have
 also to consider the sphere equation of $S^{31}$. In particular it follows that
 $ dim\, \big(\big(\, \overline{S}{}^{123}_{4}\cap N(\Theta)\big)/ G\big)=0$,
 that is,
  the stratum $\overline{S}{}^{123}_{4}$
  yields one singular point on the orbifolds $\mathcal{M}^7(\Theta)$.
  For the moment map equations
  relative to cases $\overline{S}{}^{12}_{34}$, $\overline{S}{}^{13}_{24}$
   and $\overline{S}{}^{14}_{23}$
  it sufficient to study one of them. Consider $\overline{S}{}^{12}_{34}$,
  then by using the
  computations we have done in the previuos case, then for $Sp(1)$ we get
  \begin{equation}       \label{C2 : 157}
  \begin{aligned}
&\left\{ \begin{array}{l}
  2 \sum_{\alpha=1}^4
 Im( \zad\bzap)   + cos \, \varphi \big( \sum_{\alpha=1}^2( |\zad |^2 + |\zap |^2 )
  - \sum_{\alpha= 3}^4( |\zad |^2 + |\zap |^2 )\big) =0,  \\
 cos \, \varphi  \sum_{\alpha=1}^4( |\zad |^2 + |\zap |^2 ) + 2 \big( \sum_{\alpha=1}^2
 Im( \zad\bzap) -  \sum_{\alpha=3}^4 Im( \zad\bzap) \big) =0,   \\
\sum_{\alpha=1}^2 \big( (\zad )^2 + (\zap )^2 \big) \big)- \sum_{\alpha=3}^4 \big(
(\zad )^2 + (\zap )^2 \big) \big) =0. \\
 \sum_{\alpha=1}^4 \big( (\zad )^2 + (\zap )^2 \big) \big) =0, \\
  \end{array} \right.            \\
& \mathrm{ if\,\, and\,\, only\,\, if\,\, } \\
& Sp(1)\,\, \Rightarrow\,\,
\left\{ \begin{array}{l}
2 \sum_{\alpha=1}^2 Im( \zad\bzap)
 + cos \, \varphi \big( \sum_{\alpha=1}^2( |\zad |^2 + |\zap |^2 ) \big) =0,  \\
 2 \sum_{\alpha=3}^4 Im( \zad\bzap)
   - cos \, \varphi \big( \sum_{\alpha=3}^4( |\zad |^2 + |\zap |^2 ) \big) =0,  \\
\sum_{\alpha=1}^2 \big( (\zad )^2 + (\zap )^2 \big) \big) =0, \\
 \sum_{\alpha=3}^4 \big( (\zad )^2 + (\zap )^2 \big) \big) =0. \\
 \end{array} \right.
 \end{aligned}
 \end{equation}
  The moment map equations of the $3-$torus $T^3_{\Theta}$ in (\ref{C2 : 128})
   become
 \begin{equation}          \label{C2 : 158}
 \left\{ \begin{array}{l}
 \sum_{\alpha=1}^4 p_{\alpha} \big( 2Im(\bzap\zad) + \big( \sum_{\alpha=1}^2 ( |\zad |^2 + |\zap |^2 )- \sum_{\alpha=3}^4 ( |\zad |^2 + |\zap |^2 ) \big)cos\, \varphi  \big)=0,  \\
 \sum_{\alpha=1}^4 q_{\alpha} \big( 2Im(\bzap\zad) + \big( \sum_{\alpha=1}^2 ( |\zad |^2 + |\zap |^2 )- \sum_{\alpha=3}^4 ( |\zad |^2 + |\zap |^2 ) \big)cos\, \varphi  \big)=0,  \\
 \sum_{\alpha=1}^4 l_{\alpha} \big( 2Im(\bzap\zad) + \big( \sum_{\alpha=1}^2 ( |\zad |^2 + |\zap |^2 )- \sum_{\alpha=3}^4 ( |\zad |^2 + |\zap |^2 ) \big)cos\, \varphi \big)=0,  \\
 \sum_{\alpha=1}^3 p_{\alpha}\big( (\zad )^2 + (\zap )^2 \big)  =0, \\
 \sum_{\alpha=1}^3 q_{\alpha}\big( (\zad )^2 + (\zap )^2 \big)  =0, \\
 \sum_{\alpha=1}^3 l_{\alpha}\big( (\zad )^2 + (\zap )^2 \big)  =0. \\
 \end{array} \right.  \hspace{3,5 cm} {}
 \end{equation}
 Then, it follows that\efl
 \begin{equation}      \label{C2 : 159}
  \begin{aligned}
 & 2Im(\bzap\zad) + ( |\zad |^2 + |\zap |^2 )cos\, \varphi =0,    \\
   \mathrm{and}  \quad \quad \hspace{1 cm} &  \\
&  2Im(z_{2\beta-1}\overline{z}_{2\beta}) -( |z_{2\beta-1} |^2 +
 |z_{2\beta} |^2 )cos\, \varphi =0, \hspace{7 cm}  \\
   \end{aligned}
  \end{equation}
  \bfl that is \efl
  \begin{equation}            \label{C2 : 160}
   cos\, \varphi  = -\frac{2Im(\bzap\zad)}{|\zad |^2 + |\zap |^2}
    \quad \mathrm{and}\quad
   cos\, \varphi  = \frac{2Im(z_{2\beta-1}\overline{z}_{2\beta})}{|z_{2\beta-1} |^2 + |z_{2\beta} |^2},
 \end{equation}
 \bfl where $\alpha\in\{1,2\}$ and $\beta\in \{3,4\}$. Moreover,
 it holds that \efl
 \begin{equation}            \label{C2 : 161}
(\zad )^2 + (\zap )^2  =0, \quad \alpha\in \{1,2,3,4\}.
 \end{equation}

 \bfl Let us recall that we have
also to consider the sphere equation of $S^{31}$. In particular it follows that
$ dim\, \big(\big(\, \overline{S}{}^{12}_{34}\cap N(\Theta)\big)/ G\big)=0$, that is,
 the stratum $\overline{S}{}^{12}_{34}$
  yields one singular point on the orbifolds $\mathcal{M}^7(\Theta)$.
   Analogously,
  the case of $\overline{S}{}^{1}_{234}$ can be described by using the
   the prevoius computations.
  In this case the moment map equations become
 \begin{equation}           \label{C2 : 162}
Sp(1)\,\, \Rightarrow\,\,
\left\{ \begin{array}{l}
 \sum_{\alpha=1}^4 \big(
 2Im(\bzap \zad)   + cos \, \varphi \big(( |\zad |^2 + |\zap |^2 ) \big) =0,  \\
\sum_{\alpha=1}^3 \big( (\zad )^2 + (\zap )^2 \big) \big)\big) =0, \\
 \end{array} \right.
 \end{equation}
 \bfl and \efl
 \begin{equation}            \label{C2 : 163}
 T^3_{\Theta}\,\,  \Rightarrow\,\,
 \left\{ \begin{array}{l}
 \sum_{\alpha=1}^4 p_{\alpha} \big( 2Im(\bzap\zad) + ( |\zad |^2 + |\zap |^2 )cos\, \varphi  \big)=0,  \\
 \sum_{\alpha=1}^4 q_{\alpha} \big( 2Im(\bzap\zad) + ( |\zad |^2 + |\zap |^2 )cos\, \varphi  \big)=0,  \\
 \sum_{\alpha=1}^4 l_{\alpha} \big( 2Im(\bzap\zad) + ( |\zad |^2 + |\zap |^2 )cos\, \varphi  \big)=0,  \\
 \sum_{\alpha=1}^4 p_{\alpha}\big( (\zad )^2 + (\zap )^2 \big)  =0, \\
 \sum_{\alpha=1}^4 q_{\alpha}\big( (\zad )^2 + (\zap )^2 \big)  =0, \\
 \sum_{\alpha=1}^4 l_{\alpha}\big( (\zad )^2 + (\zap )^2 \big)  =0. \\
 \end{array} \right.
 \end{equation}
  Note that, also in this case $ dim\, \big(\big(\, \overline{S}{}^{1}_{234}\cap N(\Theta)\big)/ G\big)=0$, that
  is, we have a singular point on $\mathcal{M}^7(\Theta)$.
   Then, each stratum in (\ref{C2 : 118}) yields a singular
   point on $\mathcal{M}^7(\Theta)$.  $\square$\efl

\chapter{The Quotient Orbifolds  $\mathcal{O}^4(\Omega)$ }

\addcontentsline{toc}{chapter}{The Quotient Orbifolds  $\mathcal{O}^4(\Omega)$ }

\bfl Let us recall that the action of the group
$G:=Sp(1)\times T^2_{\Omega}\subset SO(7)\subset Sp(7)$ on the sphere $S^{27}\subset
\mathbb{H}^7$ is defined  by using the left multiplication of the subgroup $Sp(1)$
 and for the $2-$torus $T^2_{\Omega}$, the action is defined via the
homomorphism   \efl
\begin{equation}  \label{C3: 1}
\begin{aligned}
f_{\Omega}(t,s): & [0,2\pi )^2 \rightarrow T^3\subset SO(7), \\
&f_{\Omega}(t,s):=
\left( \begin{array}{c|c|c|c}
1    &  0 &  0 & 0  \\
\hline
0  &  A_1 (t,s) & 0 & 0 \\
\hline
0   &    0 & A_2 (t,s) & 0\\
 \hline
0 &    0& 0& A_3 (t,s)\\
\end{array} \right) \in T^3 ,\\
\end{aligned}
\end{equation}
\bfl where  \efl
\begin{equation}  \label{C3: 2}
\Omega:=
\left( \begin{array}{ccc}
p_1 & p_2 & p_3 \\
q_1 & q_2 & q_3 \\
\end{array} \right)\in M_{2\times 3}(\mathbb{Z}),
\end{equation}
\bfl is an integral $2\times 3$ matrix. Then  the respective  moment maps
 $\mu : S^{27} \rightarrow \mathfrak{sp}(1) \otimes \mathbb{R}^3$ $\cong
\mathbb{R}^9$ 
and 
$ \nu : S^{27} \rightarrow \mathfrak{u}(1)^2\otimes \mathbb{R}^3$ read  \efl


\begin{equation}  \label{C3: 3}
\mu (\uno)=
 ( \sum_{\alpha=1}^{7} \overline{u}_{\alpha} i u_{\alpha}, \,\,
\sum_{\alpha=1}^{7} \overline{u}_{\alpha} j u_{\alpha}, \,\,
\sum_{\alpha=1}^{7} \overline{u}_{\alpha} k u_{\alpha} )\in \mathfrak{sp}(1)
\otimes \mathbb{R}^3,
\end{equation}
\bfl where $\uno\in S^{27}\subset\mathbb{H}^7$, and   \efl
\begin{equation}    \label{C3 : 4}
 \nu (\uno)= \left( \begin{array}{c}
 \sum_{\alpha = 1}^{3}  p_{\alpha}(\overline{u}_{2\alpha } u_{2\alpha+ 1} -
  \overline{u}_{2\alpha +1 } u_{2\alpha } )  \\
 \sum_{\alpha = 1}^{3}  q_{\alpha}(\overline{u}_{2\alpha  } u_{2\alpha +1} -
  \overline{u}_{2\alpha + 1 } u_{2\alpha } )  \\
\end{array} \right)\in \mathfrak{u}(1)^2\otimes \mathbb{R}^3.
\end{equation}

\bfl  According to the study in \cite{bg}, we can distinguish the elements
of $N(\Omega)$ in the following $G-$strata
\begin{equation}  \label{C3 : 5}
\begin{aligned}
& S_0 =\{ \uno\in N(\Omega)\, |\, u_1 =0\}, \\
& S_1 =\{ \uno\in N(\Omega)\, |\, u_1 \ne0 \}, \\
& S_2 = \{ \uno\in N(\Omega)\, |\, \mathrm{some\,\,quaternionic \,\,pair}\,\,
(u_{2\alpha}, u_{2\alpha-1})\,\, \mathrm{vanishes} \}, \\
& S_3 = \{ \uno\in N(\Omega)\, |\, \mathrm{none \,\, of\,\, the\,\,quaternionic \,\,pair}\,\,
(u_{2\alpha}, u_{2\alpha-1})\,\, \mathrm{vanishes} \}.  \hspace{3 cm}
\end{aligned}
\end{equation}
We state now some facts provided in \cite{bg}, also to point out a
minor correction to statements: corollary $2.3$ and lemma $3.2$ of \cite{bg}.\efl

\begin{lem}
\bfl  Suppose that all the minor determinants
\begin{equation}  \label{C3 : 6}
\Delta_{\alpha\beta}=
\left\vert \begin{array}{cc}
p_{\alpha} & q_{\alpha}  \\
p_{\beta} & q_{\beta}  \\
\end{array} \right\vert , \quad \,\,\,
  1\leq \alpha<\beta<\leq3 ,
\end{equation}
do not vanish. Then there is no element $\uno$ of $N(\Omega)$ with a null quaternionic pair.   \efl
\end{lem}

\bfl The proof of this lemma parallels that of the lemma $1.1$. As a direct consequence
of the lemma $3.1$ it follows that $S_2= \emptyset$ and $N(\Omega)=S_3= S_0\cup S_1$.
By using lemma $3.1$ we get \efl

\begin{lem}
\bfl The action of $G=Sp(1)\times T^2_{\Omega}$ on
$S_1= N(\Omega)\cap \{u_1\ne 0\}$ is
\begin{equation}  \label{C3 : 7}
\begin{aligned}
&  1)\quad  locally \,\, free\,\, if \,\, and\,\, only\,\, if \,\,
\Delta_{\alpha\beta}\ne 0, \,\, \forall\, 1\leq \alpha < \beta \leq 3,  \\
&  2) \quad  free\,\, if \,\, and \,\, only\,\, if \,\,
gcd(\Delta_{12}, \Delta_{13}, \Delta_{23})= \pm 1.  \hspace{10 cm} {}
\end{aligned}
\end{equation}
\efl
\end{lem}
\bfl Note the correct conditions $2)$ appearing in (\ref{C3 : 7})
is in fact weaker than that one given in  \cite{bg},
lemma $3.2$ point $2)$.\efl

\bfl $ \mathbf{Proof}.$   Since  $ u_1\ne 0$ it follows that $\lambda =1$, then it is sufficient
to consider the $T^2_{\Omega}-$
action. Since none of the quaternionic pairs vanishes
it follows that the fixed point equations read:\efl
\begin{equation}   \label{C3 : 8}
\overbrace{
\left(\begin{array}{cc}
cos\, \omega_{\alpha} & sin \, \omega_{\alpha} \\
- sin\, \omega_{\alpha} & cos\, \omega_{\alpha} \\
\end{array}  \right)  }^{A(\omega_{\alpha}):=}
\left(\begin{array}{c}
u_{2\alpha} \\
u_{2\alpha+1} \\
\end{array} \right)  =
\left(\begin{array}{c}
u_{2\alpha} \\
u_{2\alpha + 1} \\
\end{array} \right), \quad \alpha\in \{ 1,2,3\},
\end{equation}
\bfl where $\omega_{\alpha}= p_{\alpha}t + q_{\alpha}s$, and $t,s \in [ 0, 2\pi)$.
Then, the equations in (\ref{C3 : 8}) are equivalent to\efl
\begin{equation}   \label{C3: 9}
A(\omega_{\alpha})=id_{2\times 2}, \quad \alpha\in \{ 1,2,3\},
\end{equation}
\bfl that is $e^{i(p_{\alpha}t + q_{\alpha}s)}=1,\,\, \alpha\in \{ 1,2,3\}$.
These equations have trivial solutions if and only if
$gcd(\Delta_{12}, \Delta_{13}, \Delta_{23})= \pm 1$. The locally free conditions was
correctly proved in \cite{bg}. $\square$
\efl

\hspace{0,5 cm}

 \bfl Next, we introduce the main results of this chapter.
  Consider an isometric action of  the
 group  $\tilde{G}:= T^2_{\Omega}\times Sp(1)\times U(1)$ on $N(\Omega)= \mu^{-1}(0)\cap \nu^{-1}(0)$.
  This action is
 defined in the same way of that one in $(2.1)$, namely\efl
 \begin{equation}  \label{C3 : 10}
 \begin{aligned}
 &\quad \quad \tilde{\Phi}:\,\big(T^2_{\Omega}\times Sp(1)\times U(1)\big)\times
 N(\Omega) \longrightarrow N(\Omega),\\
 &\quad \quad \quad \big( (A(\Omega), \lambda, \rho);(\underline{z},
  \underline{w} ) \big)\mapsto
 \tilde{\Phi} \big( ( A(\Omega), \lambda, \rho)\big)\big((\underline{z},
  \underline{w} )\big), \\
\mathrm{where} & \\
 & \quad \quad \quad \quad \tilde{\Phi} \big( ( A(\Omega), \lambda, \rho)
 \big)\big((\underline{z}, \underline{w} )
 =   A(\Omega) \lambda  \left( \begin{array}{c}
\underline{z}\\
\underline{w} \\
\end{array}\right)
 \rho,   \hspace{10 cm} {}
 \end{aligned}
 \end{equation}
 \bfl and we have identified $\mathbb{C}^7\times \mathbb{C}^7\cong \mathbb{H}^7$ by using
 the relation $u_{\alpha} = z_{\alpha} + jw_{\alpha}$ for each $\alpha\in \{1,2,3\}$.
 The  $3$-Sasakian and quaternion $\ka$ reduction we are dealing with,
   is described by the following diagram

\begin{equation}    \label{C3 : 11}
\dgARROWLENGTH=0.3\dgARROWLENGTH
   \begin{diagram}
     \node[2]{ \frac{SO(7)}{SO(4)\times  Sp(1)}  } \arrow[2]{e,t}{T^2_{\Omega}} \arrow[2]{s,l}{{   {}^ {}_{SO(3)}}}
       \node[2]{\mathcal{M}^7(\Omega)}  \arrow[2]{s,l}{ {}^{{}_{SO(3)}} } \arrow{se,t}{  U(1) } \\
                                                                          \node[5]{\mathcal{Z}^6(\Omega)} \arrow{sw,b}{S^2}  \\
     \node[2]{\grms} \arrow[2]{e,t}{T^2_{\Omega}} \node[2]{\mathcal{O}^4(\Omega)  }
   \end{diagram}
\end{equation}

  where $\mu^{-1}(0)/Sp(1)\cong SO(7)/{SO(4)\times Sp(1)}$.
 We prove here the following statements, mentioned in the introduction as Theorems D and E.\efl

\begin{teor}
\bfl Let $\mathcal{Z}^6(\Omega)$ be the twistor space over the positive SDE orbifold
 $\mathcal{O}^4(\Omega)$. Then, depending on the minor determinants $\Delta_{\alpha\beta}$
 and $\square^{1\pm 2}_{1\pm 3}$  the singular locus $\Sigma(\Omega)$ consists of \efl

 \bfl 1) one isolated $2-$sphere $S^2$, whose isotropy only depends on one of the
  determinants $\square^{1\pm 2}_{1\pm 3}$;

 \vspace{0,5 cm}

 2) six points, that is, each of
  the following $\tilde{G-}$strata
  $($ relative to the action
 of $\tilde{G}$ on $\mathbb{H}^7$$)$\efl
\begin{equation}      \label{C3 : 12}
 \begin{aligned}
& a) \,\,\, \quad  {}^+\overset{=}{S}{}^{12}_3:=\tilde{G}\cdot \Bigg\{
\left( \begin{array}{c|cc|cc|cc}
0 & z_2  & z_3 & z_4   &    z_5 &  0    &  0   \\
0   &    0 &  0  & 0     &   0             & w_6  &  iw_6\\
\end{array}\right) \Bigg\},   \\
& b) \,\,\, \quad  {}^-\overset{=}{S}{}^{12}_3:=\tilde{G}\cdot \Bigg\{
\left( \begin{array}{c|cc|cc|cc}
0 & z_2  & z_3 & z_4   &   z_5   & 0    &  0   \\
0   &    0 &  0  & 0     &   0     & w_6  &  -iw_6\\
\end{array}\right) \Bigg\},   \\
&c) \,\,\,\quad  {}^+\overset{=}{S}{}^{13}_2:=\tilde{G}\cdot \Bigg\{
\left( \begin{array}{c|cc|cc|cc}
0 & z_2  & z_3    &    0   & 0           & z_6   & z_7  \\
0 &    0 &  0     &   w_4  &  iw_4   &  0    &  0          \\
\end{array}\right) \Bigg\},  \\
& d)\,\,\,\quad   {}^-\overset{=}{S}{}^{13}_2:=\tilde{G}\cdot \Bigg\{
\left( \begin{array}{c|cc|cc|cc}
0 & z_2  & z_3        &   0   &    0          & z_6   & z_7   \\
0 &    0 &     0      &   w_4 &    -iw_4      & 0     & 0\\
\end{array}\right) \Bigg\},   \\
 & e)\,\,\, \quad  {}^+\overset{=}{S}{}^{23}_1:=\tilde{G}\cdot \Bigg\{
\left( \begin{array}{c|cc|cc|cc}
    0 & 0    &  0      &    z_4 & z_5             & z_6   & z_7  \\
    0 & w_2  & iw_2    &   0    &  0              &  0    &  0          \\
\end{array}\right) \Bigg\},   \\
& f)\,\,\, \quad  {}^-\overset{=}{S}{}^{23}_1:=\tilde{G}\cdot \Bigg\{
\left( \begin{array}{c|cc|cc|cc}
0   & 0    &  0           &    z_4 & z_5      & z_6   & z_7  \\
0   & w_2  & -iw_2        &   0    &  0       &  0    &  0          \\
\end{array}\right) \Bigg\}. \hspace{10 cm }  \\
\end{aligned}
\end{equation}    \label{C3 : 13}
\bfl  Each of these strata intersects the submanifold $N(\Omega)$, and the $\tilde{G}-$strata
 of $N(\Omega)$ are given by
 $\widetilde{S}{}^{\alpha\beta}_{\gamma}:=\overset{=}{S}{}^{\alpha\beta}_{\gamma}\cap N(\Omega)$.
 Then, each quotient
\begin{equation}
\frac{\overset{=}{S}{}^{\alpha\beta}_{\gamma}\cap N(\Omega)}{\tilde{G}},
\end{equation}
 is a point on $\mathcal{Z}^6(\Omega)$.  Moreover,
the isotropy associated to the strata in $(\ref{C3 : 12})$ only depends on the
 minor determinants $ \pm \Delta_{\alpha \beta } $. \efl

\bfl In the case when some of the minor determinants $\Delta_{\alpha\beta}$
 and $\square^{1\pm 2}_{1\pm 3}$ are $\pm 1$,  the singular locus
 $\Sigma(\Omega)$ can be obtained by removing from the above list the
  singular sets
 whose isotropy depends on one of the mentioned determinants.  \efl
\end{teor}

 \begin{teor}
\bfl The singular points of  the $3-$Sasakian orbifold $\mathcal{M}^7(\Omega)$
 come from the singular $G-$strata  $\widetilde{S}$ of $N(\Omega)$,
 which are obtained by intersectig $N(\Omega)$ with
 the following $G-$strata $\overline{S}$ of $\mathbb{H}^7$\efl
 \begin{align*}
 & a) \,\,\, \quad  \overline{S}{}^{123}:= 
  \bigg\{
\left( \begin{array}{c|cc|cc|cc}
0 & z_2  & z_3 & z_4   &    z_5 & z_6        & z_7       \\
0 & \widetilde{w}_2  &  \widetilde{w}_3  & \widetilde{w}_4     &   \widetilde{w}_5   &  \widetilde{w}_6      & \widetilde{w}_7  \\
\end{array}\right) \bigg\},   \\
& b) \,\,\, \quad \overline{S}{}^{12}_3:= 
 \bigg\{
\left( \begin{array}{c|cc|cc|cc}
0 & z_2   &   z_3  & z_4      &    z_5   & z_6            & z_7       \\
0 & \widetilde{w}_2  &  \widetilde{w}_3  & \widetilde{w}_4     &   \widetilde{w}_5   &  \overset{\approx}{w}_6          & \overset{\approx}{w}_7  \\
\end{array}\right) \bigg\},   \hspace{10 cm}\\
&c) \,\,\,\quad \overline{S}{}^{13}_2:= 
\bigg\{
\left( \begin{array}{c|cc|cc|cc}
0 & z_2  & z_3      & z_4            &    z_5                & z_6           & z_7    \\
0 & \widetilde{w}_2  &  \widetilde{w}_3  & \overset{\approx}{w}_4     &  \overset{\approx}{w}_5      &  \widetilde{w}_6         & \widetilde{w}_7  \\
\end{array}\right) \bigg\},   \\
 & d)\,\,\, \quad \overline{S}{}^{23}_1:= 
  \bigg\{
\left( \begin{array}{c|cc|cc|cc}
0 & z_2  & z_3 & z_4   &    z_5 & z_6        & z_7       \\
0 &\overset{\approx}{w}_2 & \overset{\approx}{w}_3  & \widetilde{w}_4    &  \widetilde{w}_5  &  \widetilde{w}_6   & \widetilde{w}_7  \\
\end{array}\right) \bigg\},   \\
& e)\,\,\,\quad \overline{S}{}^{1}_{23}:= 
 \bigg\{\left( \begin{array}{c|cc|cc|cc}
0 & z_2  & z_3 & z_4   &    z_5 & z_6        & z_7    \\
0  & \widetilde{w}_2  &  \widetilde{w}_3  & \overset{\approx}{w}_4    &   \overset{\approx}{w}_5  & \overset{\approx}{w}_6        & \overset{\approx}{w}_7  \\
\end{array}\right) \bigg\},   \\
& f)\,\,\,\quad \overline{S}{}^2_{13}:= 
\bigg\{
\left( \begin{array}{c|cc|cc|cc}
0 & z_2   & z_3    & z_4      &    z_5 & z_6        & z_7      \\
0 & \overset{\approx}{w}_2 & \overset{\approx}{w}_3  & \widetilde{w}_4     &  \widetilde{w}_5  &  \overset{\approx}{w}_6          & \overset{\approx}{w}_7  \\
\end{array}\right) \bigg\},   \\
& g)\,\,\,\quad  \overline{S}{}^3_{12}:= 
\bigg\{
\left( \begin{array}{c|cc|cc|cc}
0 & z_2  & z_3 & z_4   &    z_5 & z_6        & z_7      \\
0 & \overset{\approx}{w}_2 & \overset{\approx}{w}_3  & \overset{\approx}{w}_4     &  \overset{\approx}{w}_5 & \widetilde{w}_6   & \widetilde{w}_7  \\
\end{array} \right)    \bigg\}, \\
\end{align*}
\begin{equation}  \label{C3 : 14}{
 \begin{aligned}
& h)\,\,\,\quad  \overline{S}{}_{123} := 
 \bigg\{
\left( \begin{array}{c|cc|cc|cc}
0 & z_2  & z_3 & z_4   &    z_5 & z_6        & z_7    \\
0 & \overset{\approx}{w}_2 & \overset{\approx}{w}_3  & \overset{\approx}{w}_4     &   \overset{\approx}{w}_5 & \overset{\approx}{w}_6        & \overset{\approx}{w}_7  \\
\end{array}\right) \bigg\},
 \end{aligned}    }    \hspace{4,5 cm}    {}
 \end{equation}
 \bfl where    \efl
 \begin{equation}    \label{C3 : 14a}
 \begin{aligned}
& (\widetilde{w}_{2\alpha}, \widetilde{w}_{2\alpha +1}):= \frac{e^{i\delta} }{ sin\, \varphi}
\bigg(-z_{2\alpha+1} +iz_{2\alpha}\cos\, \varphi,
 z_{2\alpha } +iz_{2\alpha+1}\cos\, \varphi \bigg),  \\
  and & \\
  & (\overset{\approx}{w}_{2\alpha}, \overset{\approx}{w}_{2\alpha +1}):=   \frac{e^{i\delta}  }{ sin\, \varphi}
\bigg( -z_{2\alpha+1} -i\zap
cos\, \varphi,   \zap  -iz_{2\alpha+1}\cos\, \varphi  \bigg), \hspace{3 cm} \\
\end{aligned}
\end{equation}
\bfl with $sin\, \varphi\ne 0$. Let us recall that an element $\lambda=\epsilon + j\sigma\in Sp(1)$
is written as
  $\epsilon = \cost +i(\sit cos\, \varphi)$  and
$\sigma =  \sit sin\, \varphi $ $ cos\, \delta + i(\sit sin\, \varphi sin\, \delta ).$ Moreover,
the isotropy associated to the strata in $(\ref{C3 : 14})$ only depends on the
 minor determinants $\square^{1\pm 2}_{1\pm 3}$. $\square$\efl
 \end{teor}

 \hspace{0,5 cm}

\bfl Finally we get ( Theorem F in the introduction): \efl

\begin{teor}
 By using theorems $2.1$, $2.2$, $3.1$ and $3.2$, we can conclude that
the families of SDE orbifolds $\mathcal{O}^4(\Theta )$ and
$\mathcal{O}{}^4(\Omega)$ are distinct to the twistor and $3-$Sasakian level.
\end{teor}

\hspace{0,5 cm}

\section{The
Twistor Space $\mathcal{Z}^6(\Omega)$  }

\hspace{0,5 cm}

\bfl The study of the action of    $\tilde{G}= Sp(1)\times T^2_{\Omega}\times U(1)$
 on $\mathbb{H}^7$ can be divided in two part. First we consider
  the action of
 $\tilde{G}= Sp(1)\times T^2_{\Omega}\times U(1)$ on
  $S_0=\{ \uno\in N(\Omega)\, |$
  $\, u_1 =0\}$, and  then,
   the case linked to action of $\tilde{G}$ on
 $S_1 =\{ \uno\in N(\Omega)\, |\, u_1 \ne0 \}$.
 Now, we introduce some facts
  shared by both the two mentioned cases.
 Let $\zw$ be a point in $N(\Omega)
 \subset S^{27}\subset \mathbb{H}^7$, then
 the fixed point
  equations for each quaternionic pairs $(u_{2\alpha}, u_{2\alpha + 1})$, $\alpha\in
  \{1,2,3\}$, become   \efl
\begin{equation}   \label{C3 : 22}
A(\omega_{\alpha})
 \left(\begin{array}{cc}
z_{2{\alpha}} & w_{2{\alpha}} \\
z_{2{\alpha}+1 } & w_{2{\alpha}+1} \\
\end{array} \right) =
{\left[ \begin{array}{c}
\left(\begin{array}{cc}
\epsilon & -\overline{\sigma}  \\
\sigma & \overline{\epsilon} \\
\end{array} \right)
 \left(\begin{array}{cc}
z_{2\alpha} & z_{2\alpha +1} \\
w_{2\alpha} & w_{2\alpha +1} \\
\end{array} \right)
 \left(\begin{array}{cc}
\rho& 0 \\
0 & \rho \\
\end{array} \right)
\end{array}\right]}^T,
\end{equation}
\bfl where $A(\omega_{\alpha}):=
\left(\begin{array}{cc}
cos\, \omega_{\alpha} & sin\, \omega_{\alpha} \\
-sin\, \omega_{\alpha}  & cos \omega_{\alpha} \\
\end{array} \right)$,  $\lambda= \epsilon + j\sigma\in Sp(1)$, $\rho\in U(1)$
 and $\alpha\in \{1,2,3\}$.
We indicate the points
$\uno=\zw\in \mathbb{C}^7\times\mathbb{C}^7$ as follows  \efl
\begin{equation}  \label{C3 : 23}
{}^T(\underline{z}, \underline{w} ):=
\left(\begin{array}{c|cc|cc|cc}
z_1 & z_2  & z_3 & z_4 & z_5 & z_6 & z_7  \\
w_1 & w_2 &  w_3 & w_4 & w_5 & w_6 & w_7  \\
\end{array}\right),
\end{equation}
\bfl where $(z_{\alpha}, w_{\alpha})\in \mathbb{C}\times\mathbb{C} $, $\alpha\in \{1,2,3,4,5,6,7\}$.
These equations can be rewritten in the following way:\efl
\begin{equation}  \label{C3 : 24}
\overbrace{\left( \begin{array}{cc|cc}
0                          & -\sigma\rho &  -sin\omega_{\alpha}  & cos\, \omega_{\alpha} - \overline{\epsilon}\rho  \\
-\sigma\rho     & 0                      & cos\, \omega_{\alpha} - \overline{\epsilon}\rho & sin\omega_{\alpha} \\
\hline
-sin\omega_{\alpha}   & cos\, \omega_{\alpha} - \epsilon\rho           &    0       &      \overline{\sigma}\rho  \\
cos\, \omega_{\alpha} -\epsilon \rho & sin\omega_{\alpha}  &  \overline{\sigma}\rho   &   0    \\
\end{array} \right)}^{M_{\alpha}:=}
\left(\begin{array}{c}
z_{2\alpha} \\
z_{2\alpha +1}\\
w_{2\alpha} \\
w_{2\alpha +1} \\
\end{array} \right) =
\left(\begin{array}{c}
0 \\
0\\
0\\
0 \\
\end{array} \right),
\end{equation}
 \bfl for each $\alpha\in \{1,2,3\}$. Note that, none of the matrices $M_\alpha$,
  $\alpha\in \{1,2,3\}$, has $rank\,\,4 $.
    Then we have \efl

\begin{prop}
Let $\zw$ be a point in $N(\Omega)$. Then, up to $\tilde{G}-$conjugation, we have that
 $\tilde{G}_{\zw}\subset T^2_{\Omega}\times \{\lambda\in Sp(1)\,\, |\,\,
 \sigma\equiv 0 \}\times U(1)$.
\end{prop}

\bfl The proof is completely similar to
 that of proposition $2.1$. In the following we assume $\sigma\equiv 0$\efl

\subsection{The  action of $\tilde{G}$ on $S_0$ }

\vspace{0,5 cm}

\bfl When $u_1\ne 0$, we get the following \efl

\begin{lem}    \bfl
Condider the action of $Sp(1)\times T^2_{\Omega}\times U(1)$ on $\mathbb{C}^7\times\mathbb{C}^7$, and
let $\uno=(\underline{z}, \underline{w})\in $ \\
$\mathbb{C}^7\times\mathbb{C}^7 $ be a non null point.
Then, for the first quaternionic coordinate $u_1$, the fixed point equations give either
\begin{equation}     \label{C3 : 15}
\rho= Re(\epsilon) +i\sqrt{Im(\epsilon)^2 + |\sigma|^2}\,\,\quad  \mathrm{or}\,\,
\quad \rho= Re(\epsilon) -i\sqrt{Im(\epsilon)^2 + |\sigma|^2}.
\end{equation}
\efl
\end{lem}

\bfl $\mathbf{Proof.}$
For the first quaternionic coordinate $u_1$, the fixed point equations give
\begin{equation}    \label{C3 : 16}
u_1 = \lambda u_1 \rho,
\end{equation}
 where $\lambda = \epsilon + j\sigma\in Sp(1)$ and $\rho\in U(1)$. Then
 a straightforward computation yields
 \begin{equation}   \label{C3 : 17}
\left\{ \begin{array}{l}
z_1 =(\epsilon z_1- \overline{\sigma}w_1)\rho, \\
w_1 =(\overline{\epsilon}w_1 + \sigma z_1 )\rho, \\
\end{array}\right.    \,\, \quad \iff \,\,\quad
\left( \begin{array}{cc}
1-\epsilon\rho & \overline{\sigma}\rho   \\
 -\sigma\rho &  1- \overline{\epsilon}\rho  \\
   \end{array} \right)
\left(
 \begin{array}{c}
 z_1 \\
 w_1\\
 \end{array}\right)   =
 \left(
 \begin{array}{c}
 0\\
 0\\
 \end{array}\right).
 \end{equation}
\bfl Then, in order to get non trivial solution we assume \efl

 \begin{equation}   \label{C3 : 18}
\left| \begin{array}{cc}
1-\epsilon\rho & \overline{\sigma}\rho   \\
 -\sigma\rho &  1- \overline{\epsilon}\rho  \\
   \end{array} \right| = \rho^2 |\sigma |^2 +   \rho^2 |\epsilon |^2 + 1 -2\rho Re(\epsilon)=
  \rho^2 + 1 -2\rho Re(\epsilon)=0  ,
 \end{equation}

 \bfl if and only  \efl

 \begin{equation}   \label{C3 : 19}
 \rho= Re(\epsilon) +i\sqrt{Im(\epsilon)^2 + |\sigma|^2} \,\,\,\,\,\, \mathrm{or}\,\,\,\,\,\,
 \rho= Re(\epsilon) -i\sqrt{Im(\epsilon)^2 + |\sigma|^2}.  \,\,\,\, \square
 \end{equation}
\efl

\hspace{0,5 cm}

 \bfl $\mathbf{Remark\,\, 3.1.}$ Note that, when $u_1\ne 0$, if
 $z_1=0$ the fixed point equation in $(\ref{C3 : 17})$ becomes
 \begin{equation}    \label{C3 : 20}
 \left\{ \begin{array}{l}
\overline{\sigma}w_1\rho=0, \\
w_1 =\overline{\epsilon}w_1 \rho, \\
\end{array}\right.   \,\, \quad \iff \,\,\quad \sigma\equiv0,\,\, \overline{\epsilon}\rho=1.
 \end{equation}
  Analogously, when $u_1\ne 0$ and $w_1=0$, it follows that\efl
  \begin{equation}    \label{C3 : 21}
 \left\{ \begin{array}{l}
z_1 =\epsilon z_1\rho, \\
\overline{\epsilon}w_1 \rho=0, \\
\end{array}\right.   \,\, \quad \iff \,\,\quad \sigma\equiv0,\,\, \epsilon\rho=1. \,\, \square
 \end{equation}

\vspace{0,5 cm}

\begin{lem}
\bfl Let $\zw$ be a in $N(\Omega)$.
 Then, the fixed point
equations in $(\ref{C3 : 22})$ are equivalent to the conditions given by $det\, M_{\alpha}=0.$
In particular $det\, M_{\alpha}=0$ if and only if it holds either
\begin{equation}    \label{C3 : 25}
 \rho+ \overline{\rho} (e^{-i\theta_i})^2 - 2 Re(\epsilon) e^{-i\theta_i}=0\,\,\,\,\,\, \mathrm{or}\,\,\,\,\,\,
  \rho+
\overline{\rho} (e^{i\theta_i})^2
- 2 Re(\epsilon) e^{i\theta_i} =0,
 \end{equation}
  if and only if
\begin{equation}    \label{C3 : 26}
 \overline{\rho}e^{i\theta_{\alpha}}= Re(\epsilon)
\pm  i\sqrt{Im(\epsilon)^2 + |\sigma|^2}
\,\,\,\,\,\, \mathrm{or}\,\,\,\,\,\,
\overline{\rho}e^{-i\theta_{\alpha}}= Re(\epsilon) \pm  i\sqrt{Im(\epsilon)^2 +
 |\sigma|^2}.
 \end{equation}
 \efl
\end{lem}

\bfl The proof of this lemma is exactly equal to that of lemma $2.1$. \efl

\vspace{0,5 cm}

\begin{prop}
\bfl Let $(\underline{z}, \underline{w})$ be  a point
in $N(\Omega)\subset \mathbb{C}^7\times \mathbb{C}^7$.  Then the fixed
point equations $( \ref{C3 : 22})$ are
equivalent to the conditions $(\ref{C3 : 26})$ given by
$det\, \mal=0$, $\alpha\in \{1,2,3\}$. \efl
\end{prop}

\bfl The proof of this proposition parallels that
of proposition $2.2$. \efl

\vspace{0,5 cm}

\bfl $\mathbf{Remark\,\,3.2.}$ Note that, when  both the equations
in (\ref{C3 : 19}) are satisfied
 it follows that
 \begin{equation}     \label{C3 : 27}
 \left\{ \begin{array}{l}
  \rho= Re(\epsilon) +i\sqrt{Im(\epsilon)^2 + |\sigma|^2}, \\
 \rho= Re(\epsilon) -i\sqrt{Im(\epsilon)^2 + |\sigma|^2}, \\
 \end{array} \right.  \,\, \iff \,\,
\left\{\begin{array}{l}
  Re(\epsilon)=\rho= \pm 1, \\
  \sigma\equiv 0, \\
  Im(\epsilon)=0. \\
   \end{array} \right.
 \end{equation}
 By substituting the conditions (\ref{C3 : 27}) in the equations (\ref{C3 : 24}),
 we get
 $M_{\alpha}=0$ for each $\alpha\in \{1,2,3\}$. Then, when the conditions in
 (\ref{C3 : 27}) hold,
  the fixed point equations describe the non effective subgroup.  $\square$ \efl

 \vspace{0.5 cm}

 \bfl $\mathbf{Remark\,\, 3.3.}$ Since $\mathcal{M}^7(\Omega)$ is an $U(1)-$bundle over
  $\mathcal{Z}^6(\Omega)$ and $G$ acts locally free on
$N(\Omega)$, it follows that also the action of
$\tilde{G}$ on $N(\Omega)$ has to be locally free.
  Now, assume $\sigma\equiv0$
 in the equations (\ref{C3 : 22}).
 Then, the fixed point equations relative to the first quaternionic
  coordinate
 become $\epsilon u_1 \rho = u_1$ , that is
 \begin{equation}  \label{C3 : 34}
 \epsilon z_1 \rho + j\overline{\epsilon}w_1 \rho = z_1 +jw_1.
 \end{equation}   \label{C3 : 34}
As a consequence, when $z_1\ne 0$ and $w_1\ne 0$ then
\begin{equation}     \label{C3 : 35}
 \left\{\begin{array}{l}
  \epsilon \rho=1, \\
 \overline{\epsilon} \rho=1, \\
 \end{array}\right.\,\,\quad  \iff\,\,\quad
  \left\{\begin{array}{l}
  \epsilon\rho=1, \\
  \epsilon=\rho=\pm 1. \\
 \end{array}\right.
\end{equation}
Suppose the conditions in (\ref{C3 : 35}) hold. By keeping in mind the fixed point equations
 (\ref{C3 : 24}), the conditions in (\ref{C3 : 35})
give that  the stabilizer $\tilde{G}_{(\underline{z}, \underline{w})}$ is
 contained in the non effective subgroup.  Then we do not care about this case.
 Then we consider points  $(\underline{z}, \underline{w})$ such that
  $u_1\ne 0$ and $z_1 =0$ or $w_1=0$.  $\square$ \efl

 \bfl We prove now the following \efl

\begin{lem}    \bfl
Suppose $u_1\ne0$ and $\sigma\equiv 0$ in the  equations  $(\ref{C3 : 24})$.
The  spaces of solutions for the eigenvalue problem
asociated  to the equations $(\ref{C3 : 24})$,  $\alpha\in \{1,2,3\}$,
are the following
\begin{equation}  \label{C3 : 36}
\begin{aligned}
1)\,\, \hspace{1 cm} &{}^{\pm}\tilde{V}_1^{\alpha}:=\{z_{2\alpha}\in \mathbb{C}\, |\, (z_{2\alpha}, \pm iz_{2\alpha}, 0,0)\in \mathbb{C}^4\},    \\
 &  \\
2)\,\, \hspace{1 cm}&{}^{\pm}\tilde{V}_2^{\alpha}:=\{w_{2\alpha}\in \mathbb{C}\, |
 \, ( 0,0,w_{2\alpha}, \pm iw_{2\alpha})\in \mathbb{C}^4\}, \\
&    \\
3)\,\,\hspace{1 cm} & \tilde{V}_3^{\alpha}:=\{(z_{2\alpha},z_{2\alpha+1}) \in \mathbb{C}\times \mathbb{C}\, |\, (z_{2\alpha}, z_{2\alpha + 1}, 0,0)
\in \mathbb{C}^4\},  \\
&    \\
4)\,\, \hspace{1 cm}&\tilde{V}_4^{\alpha}:=\{(w_{2\alpha},w_{2\alpha+1})\in \mathbb{C}\times \mathbb{C}\, |
 \, ( 0,0,w_{2\alpha}, w_{2\alpha + 1})\in \mathbb{C}^4\}. \hspace{10 cm}\\
 \end{aligned}
 \end{equation}
 \efl
\end{lem}

 \bfl $\mathbf{Proof.}$
  Under our hypotheses,  the equations in (\ref{C3 : 26}) become
\begin{equation}   \label{C3 : 37}
  \rho \epsilon = e^{\pm\theta_1}
\,\,\,\,\,\mathrm{  and}\,\,\,\,\, \rho \overline{\epsilon}= e^{\pm i\theta_1}.
\end{equation}
  Let $(\underline{z}, \underline{w} )\in \mathbb{C}^7\times \mathbb{C}^7$ on a singular point
  and such that $u_1\ne 0$ . Then according to whether
  $z_1= 0$ or $w_1= 0$ it holds either
  $\rho\epsilon=1$ or $\rho\overline{\epsilon}=1$. Suppose $\rho\epsilon=1$,
  the other case can be treated
  in the same way. Let $M_{\alpha}$ such that two relations of those
   one in
   (\ref{C3 : 37}) are satisfied,
  then, for this block we get either one of
   the following system of fixed point equations  \efl
\begin{equation}    \label{C3 : 38}
\begin{aligned}
& 1)\quad \quad  \left\{\begin{array}{l}
\rho \epsilon= e^{i\theta_{\alpha}}, \\
\rho\overline{\epsilon}= e^{i\theta_{\alpha}}, \\
 \rho\epsilon=1, \\
\end{array} \right.  \quad \quad 2) \quad  \quad
 \left\{\begin{array}{l}
\rho \epsilon= e^{i\theta_{\alpha}},     \\
\rho\overline{\epsilon}= e^{-i\theta_{\alpha}}, \\
 \rho\epsilon=1. \\
\end{array} \right. \\
&  \\
& 3)\quad \quad \left\{\begin{array}{l}
\rho \overline{\epsilon}= e^{i\theta_{\alpha}}, \\
\rho\overline{\epsilon}= e^{-i\theta_{\alpha}}, \\
 \rho\epsilon=1. \\
\end{array} \right.  \hspace{ 10 cm} {}
\end{aligned}
  \end{equation}
 \bfl A short computation shows that the system of equations
 $1)$ and $2)$ in
 (\ref{C3 : 38}) give $\epsilon \rho=1$ and $\overline{\epsilon}\rho=1$, then
they describes the non
  effective subgroup. Instead in the case $3)$, the equation in (\ref{C3 : 24}) becomes    \efl
  \begin{equation}       \label{C3 : 39}
 \left( \begin{array}{cc|cc}
0                          & 0 &   0  & 0  \\
0     & 0                      & 0 & 0 \\
\hline
0   & \pm 1 - 1           &    0       &      0 \\
\pm 1 - 1 & 0  &               0        &   0    \\
\end{array} \right).
\left( \begin{array}{c}
z_{2\alpha} \\
z_{2\alpha +1} \\
w_{2\alpha} \\
w_{2\alpha +1}  \\
\end{array}\right) =
 \left( \begin{array}{c}
 0 \\
 0 \\
 0 \\
 0 \\
 \end{array} \right).
  \end{equation}
\bfl Hence, when the block $M_{\alpha}\ne 0$ the space of solutions is exactly
 $\tilde{V}_4^{\alpha}$.
   Now suppose, $M_{\alpha}$
   is such that just one of the equations in
   (\ref{C3 : 37}) is satisfied, then we have the following
   possibilities     \efl

 \begin{align*}
&  4) \quad \quad\left\{\begin{array}{l}
\rho \overline{\epsilon}= e^{\pm i\theta_{\alpha}}, \\
 \rho\epsilon=1, \\
\end{array} \right. \quad \iff \quad
\left\{\begin{array}{l}
\rho^2= e^{\pm i\theta_{\alpha}}, \\
 \rho\epsilon=1, \\
\end{array} \right.\\
& \mathrm{or}      \hspace{ 15 cm }  {} \\
 \end{align*}
 \begin{equation}     \label{C3 : 40}
 \begin{aligned}
& 5)\quad \quad  \left\{\begin{array}{l}
\rho \epsilon= e^{\pm i\theta_{\alpha}},     \\
 \rho\epsilon=1, \\
 \end{array} \right.  \quad \iff \quad
 \left\{\begin{array}{l}
 e^{\pm i\theta_{\alpha}}=1, \\
 \rho\epsilon=e^{\pm i\theta_{\alpha}}. \\
\end{array} \right.   \hspace{ 8 cm }  {}
 \end{aligned}
\end{equation}
\bfl In the case of the relations $4)$ in (\ref{C3 : 40}) the equations (\ref{C3 : 24}) become \efl
\begin{equation}    \label{C3 : 41}
 \left( \begin{array}{cc|cc}
0                          & 0 &  -\sit_{\alpha}  & \pm i\sit_{\alpha}  \\
0     & 0                      & \pm i\sit_{\alpha} & \sit_{\alpha} \\
\hline
-\sit_{\alpha}   & \cost_{\alpha} - 1           &    0       &      0 \\
\cost_{\alpha} - 1 & \sit_{\alpha}  &               0        &   0    \\
\end{array} \right)\left( \begin{array}{c}
z_{2\alpha} \\
z_{2\alpha +1} \\
w_{2\alpha} \\
w_{2\alpha +1}  \\
\end{array}\right) =
 \left( \begin{array}{c}
 0 \\
 0 \\
 0 \\
 0 \\
 \end{array} \right).
\end{equation}
\bfl Then, the space of solutions, relative to this system of equations
  is given by  ${}^{\pm}\tilde{V}_2^{\alpha}$. Instead, in the case
of the relations $5)$ in (\ref{C3 : 40}) the eigevalue problem in (\ref{C3 : 24}) can be
 rewritten as follows  \efl
 \begin{equation}      \label{C3 : 42}
 \left( \begin{array}{cc|cc}
0                          &       0       &  0            & 1- \rho^2  \\
0                          &       0       & 1- \rho^2     &     0\\
\hline
0   & 0                    &       0       &       0         \\
0   & 0                    &       0       &        0       \\
\end{array} \right)
\left( \begin{array}{c}
z_{2\alpha} \\
z_{2\alpha +1} \\
w_{2\alpha} \\
w_{2\alpha +1}  \\
\end{array}\right) =
 \left( \begin{array}{c}
 0 \\
 0 \\
 0 \\
 0 \\
 \end{array} \right).
 \end{equation}
\bfl In this case $rank\, M_{\alpha}=2$ and the space of solutions is
 $\tilde{V}_3^{\alpha}$.
  Analogously, when  $\rho \overline{\epsilon}=1$,   we can prove that
   the only spaces of solutions
  for the eigenvalue problem in (\ref{C3 : 24})
  are given by ${}^{\pm}\tilde{V}_1^{\alpha}$,  $\tilde{V}_3^{\alpha}$ and $\tilde{V}_4^{\alpha}$. $\square$\efl

\hspace{0,5 cm}

\begin{cor}
\bfl Let $\zw$ be a point in $N(\Omega)$ with non trivial isotropy subgroup $\tilde{G}_{\zw}$.
Then each quaternionic pair $\uu$ of $\zw$
belongs to one of the sets $\tilde{G}\cdot {}^{\pm}\tilde{V}{}_{1}{}^{\alpha}$,
$\tilde{G}\cdot {}^{\pm}\tilde{V}{}_{2}^{\alpha}$, $\tilde{G}\cdot \tilde{V}{}_{3}^{\alpha}$
 or $\tilde{G}\cdot \tilde{V}{}_{4}^{\alpha}$. \efl
\end{cor}

\vspace{0,5 cm}

 \begin{prop}
\bfl There are no singular points $(\underline{z}, \underline{w} )\in N(\Omega)$
 such that their first quaternionic pair
$u_1\ne0$. \efl
 \end{prop}

\bfl  $\mathbf{Proof.}$
Under the hypotheses of lemma $3.5$ it follows that the possible spaces of solutions
relative to the eigenvalues problem in  (\ref{C3 : 24}) are just the following
 \begin{equation}     \label{C3 : 43}
V:= \Bigg\{ \left(\begin{array}{c|cc|cc|cc}
u_1 & u_2  & u_3 & u_4 & u_5 &  u_6 & u_7 \\
\end{array}\right)\Bigg\},
 \end{equation}
where $u_1\ne 0$, $z_1=0$ or $w_1=0$ and each quaternionic coordinates
  $(u_{2\alpha}, u_{2\alpha-1})$ is
 contained in one of the the complex vector spaces
  ${}^{\pm}\tilde{V}_1^{\alpha}$, ${}^{\pm}\tilde{V}_2^{\alpha}$,
  $\tilde{V}_3^{\alpha}$, $\tilde{V}_4^{\alpha}$, $\alpha\in\{1,2,3\}$.
Then, by substituting in the equations (\ref{C3 : 24}) a suitable
subset of the relations in (\ref{C3 : 37}),  the associated  eigenvalue problem admit solutions
  $\zw$ which form
a vector space $V$ of those in (\ref{C3 : 43}).  Moreover
all of the points $\zw\in V$ have exactly the same isotropy.
Let us recall that each stratum $\overset{=}{S}$
 is formed by all of the points in $\mathbb{H}^7$ which have, up to conjugation,
 the same isotropy subgroup. Note that,  the spaces  in  (\ref{C3 : 36}) are
  not $\tilde{G}-$invariant.
 If a point $\zw\in V$ belongs to $\overset{=}{S}$,
  it follows that
$ V\subseteq \overset{=}{S}$ and  $\tilde{G}\cdot V\subseteq
 \overset{=}{S}$.  By using corollary $3.1$ we have that
 $\tilde{G}\cdot V\cong \overset{=}{S}$ and the intersection
 $\tilde{G}\cdot V\cap$ $  N(\Omega)\cong \tilde{G}\cdot (V \cap N(\Omega))$
gives rise a  stratum  $\widetilde{S}\subset N(\Omega)$ associated to the $\tilde{G}-$action on $N(\Omega)$.
 Actually, all of the strata
$\widetilde{S}$ relative to the
submanifold $N(\Omega)$ can be obtained in this way. Suppose a stratum
 $\widetilde{S}$ contains a point $\zw$ such that the conditions
 $\{ u_1\ne 0, z_1\ne 0\}$ or
$\{ u_1\ne 0, w_1 \ne 0\}$ hold.  Then
 the orbit $\tilde{G}\cdot\zw$ is contained in  $\widetilde{S}$,
  and by using the $\tilde{G}-$action,
  there exists a point $(\underline{z}', \underline{w}')\in \tilde{G}\cdot\zw$
with $u'_1\ne 0$,  $z'_1\ne 0$ and  $w'_1\ne 0$.
Then the argument in remark $3.3$  gives that the isotropy of  $\widetilde{S}$
coincide with the non effective subgroup.
 Then we can conclude. $\square$
\efl

\vspace{0,5 cm}

\subsection{The  action of $\tilde{G}$ on $S_1$ }

\vspace{0,5 cm}

\bfl Now, we can study the case of the points $\zw\in \mathbb{H}^7$
such that $u_1=0$. This case can be treated as a subcase of that one
studied in chapter 2. In fact, the action of the group $Sp(1)\times
T^3_{\Theta}\times
U(1)$ on $\mathbb{H}^8$, coincide with that one we are dealing with
as soon as we assume $r=0$ ( the third parameter of $T^3_{\Theta}$)  and the group
acts on the subspace $\{ \zw\in \mathbb{H}^8\,\, |\,\, (u_1, u_2)=(0,0)\}$.
   In order
  to describe the singular strata relative to the action of $\tilde{G}$ on $\mathbb{H}^7$, we can study
  the spaces of solutions relative to the eigenvalue problem in (\ref{C3 : 24}).
  Then \efl

\begin{lem}
\bfl Assume $\sigma\equiv 0$ in equations $(\ref{C3 : 24})$.
Let $\mal$ be such that $rank\,M_{\alpha}= 3$.  Then
the solutions of the eigenvalue problem $(\ref{C3 : 24})$
describe a complex $1-$dimensional
space which is either \efl
 \begin{equation}   \label{C3 : 44}
\begin{aligned}
1)\,\, \hspace{1 cm} & {}^{\pm}\tilde{V}_1^{\alpha}:=\{z_{2\alpha}\in \mathbb{C}\,
|\, (z_{2\alpha}, \pm iz_{2\alpha}, 0,0)
\in \mathbb{C}^4\}, \quad \quad  \mathrm{or} \\
 &  \\
2)\,\, \hspace{1 cm} &{}^{\pm}\tilde{V}_2^{\alpha}:=\{w_{2\alpha}\in \mathbb{C}\, |
 \, ( 0,0, w_{2\alpha}, \pm iw_{2\alpha})\in \mathbb{C}^4\}. \hspace{ 10 cm}    {}
\end{aligned}    \hspace{2 cm}
\end{equation}  \\
When $rank\,M_\alpha =2 $, the possible spaces of solutions are the following:
\begin{equation}   \label{C3 : 45}
\begin{aligned}
3)\,\,\hspace{1 cm} &{}^{(\pm, \pm )}\tilde{V}_3^{\alpha}:=\{(z_{2\alpha},w_{2\alpha})\in
 \mathbb{C}^2\, |
 \, (z_{2\alpha}, \pm iz_{2\alpha}, w_{2\alpha}, \pm iw_{2\alpha} )
\in \mathbb{C}^4\},      \\
& \\
4)\,\,\hspace{1 cm} & \tilde{V}_4^{\alpha}:=\{z_{2\alpha},z_{2\alpha+1}
\in \mathbb{C}\, |\, (z_{2\alpha}, z_{2\alpha+1 }, 0,0)
\in \mathbb{C}^4\},  \\
& \\
5)\,\, \hspace{1 cm}&\tilde{V}_5^{\alpha}:=\{w_{2\alpha },w_{2\alpha+1}\in \mathbb{C}\, |
 \, ( 0,0,
w_{2\alpha},
w_{2\alpha + 1})\in \mathbb{C}^4\}. \hspace{10 cm}  {}
 \end{aligned}
 \end{equation}
\end{lem}
\bfl The proof is a straightforward consequence of that of lemma
$2.4$. By using the proof of proposition $2.2$ and corollary $2.1$ we get
the following results \efl

\begin{prop}
\bfl Let $(\underline{z}, \underline{w})$ be  a point
in $N(\Omega)\subset \mathbb{C}^7\times \mathbb{C}^7$.  Then the fixed
point equations $( \ref{C3 : 22})$ are
equivalent to the conditions $(\ref{C3 : 26})$ given by
$det\, \mal=0$, $\alpha\in \{1,2,3\}$. \efl
\end{prop}

\begin{cor}
\bfl Let $\zw$ be a point in $N(\Omega)$ with non trivial isotropy subgroup $\tilde{G}_{\zw}$.
Then each quaternionic pair $\uu$ of $\zw$
belongs to one of the sets $\tilde{G}\cdot {}^{\pm}\tilde{V}_{1}^{\alpha}$,
$\tilde{G}\cdot {}^{\pm}\tilde{V}_{2}^{\alpha}$, $\tilde{G}\cdot {}^{(\pm, \pm)}\tilde{V}_{3}^{\alpha}$,
$\tilde{G}\cdot\tilde{V}_{4}^{\alpha}$ or $\tilde{G}\cdot \tilde{V}_{5}^{\alpha}$. \efl
\end{cor}

 \bfl Now, by using the corollary $2.1$, remark $2.3$ and lemma $3.6$ we can
 divide the strata (relative to the action of $Sp(1)\times T^2_{\Omega}\times U(1)$ on
 $\mathbb{H}^7$) which intersect the submanifold $N(\Omega)$ in two families depending
 on the vector spaces listed in (\ref{C3 : 44}) and (\ref{C3 : 45}) respectively.
  First, consider the vector spaces \efl
 \begin{equation}   \label{C3 : 46}
\begin{aligned}
& 1)\,\,\, \overset{=}{S}_{(+,+,-)}:= \Bigg\{ \left(\begin{array}{c|cc|cc|cc}
0  & z_3 & iz_3 & z_5 &  iz_5 & z_7 &  -iz_7\\
0  &  w_3 &  iw_3 & w_5 &  iw_5 & w_7 &  -iw_7 \\
\end{array}\right)\Bigg\},  \hspace{10 cm} {}  \\
&2)\,\,\,  \overset{=}{S}_{(+,-,+)}:= \Bigg\{ \left(\begin{array}{c|cc|cc|cc}
0   & z_3 & iz_3 & z_5 &  -iz_5 & z_7 &  iz_7\\
0   &  w_3 &  iw_3 & w_5 &  -iw_5 & w_7 &  iw_7 \\
\end{array}\right)\Bigg\}, \\
&3)\,\,\, \overset{=}{S}_{(-,+,+)}:= \Bigg\{ \left(\begin{array}{c|cc|cc|cc}
0 & z_3 & -iz_3 & z_5 &  iz_5 & z_7 &  iz_7\\
0  &  w_3 &  -iw_3 & w_5 &  iw_5 & w_7 &  iw_7 \\
\end{array}\right)\Bigg\},  \\
& 4)\,\,\,  \overset{=}{S}_{(+,+,+)}\big):= \Bigg\{ \left(\begin{array}{c|cc|cc|cc}
0  & z_3 & iz_3 & z_5 &  iz_5 & z_7 &  iz_7\\
0 &  w_3 &  iw_3 & w_5 &  iw_5 & w_7 &  iw_7 \\
\end{array}\right)\Bigg\}, \\
&5)\,\,\,  \overset{=}{S}_{(+,-,-)}:= \Bigg\{ \left(\begin{array}{c|cc|cc|cc}
0  & z_3 & iz_3 & z_5 &  -iz_5 & z_7 &  -iz_7\\
0 &  w_3 &  iw_3 & w_5 &  -iw_5 & w_7 &  -iw_7 \\
\end{array}\right)\Bigg\},   \\
&6)\,\,\,  \overset{=}{S}_{(-,+,-)}:=\Bigg\{ \left(\begin{array}{c|cc|cc|cc}
0  & z_3 & -iz_3 & z_5 &  iz_5 & z_7 &  -iz_7\\
0 &  w_3 &  -iw_3 & w_5 &  iw_5 & w_7 &  -iw_7 \\
\end{array}\right)\Bigg\},  \\
&7)\,\,\,  \overset{=}{S}_{(-,-,+)}:= \Bigg\{ \left(\begin{array}{c|cc|cc|cc}
0 & z_3 & -iz_3 & z_5 &  -iz_5 & z_7 &  iz_7\\
0 &  w_3 &  -iw_3 & w_5 &  -iw_5 & w_7 &  iw_7 \\
\end{array}\right)\Bigg\},  \\
&8)\,\,\,  \overset{=}{S}_{(-,-,-)}:= \Bigg\{ \left(\begin{array}{c|cc|cc|cc}
0 & z_3 & -iz_3 & z_5 &  -iz_5 & z_7 &  -iz_7\\
0 &  w_3 &  -iw_3 & w_5 &  -iw_5 & w_7 &  -iw_7 \\
\end{array}\right)\Bigg\}.   \hspace{10 cm }  \\
\end{aligned}
\end{equation}
\bfl A short computation shows that these spaces are $\tilde{G}-$invariant. Moreover,
all of the points $\zw$  contained in the same space
have the same isotropy subgroup. Then, each of the vector spaces
$ \overset{=}{S}_{(\pm,\pm,\pm)}$ represents a stratum relative to the action of
$\tilde{G}$ on $\mathbb{H}^7$.\efl

\begin{lem}
\bfl The strata
 of $\overset{=}{S}_{(+,+,+)}$ and $\overset{=}{S}_{(-,-,-)}$
 do not intersect the submanifold $N(\Omega)$.  \efl
\end{lem}

\bfl $\mathbf{Proof.}$
First, note that the moment map equations of the $2$-torus $T^2_{\Omega}$
in the coordinate $u_{\alpha}= z_{\alpha} + jw_{\alpha}$ becomes \efl
\begin{equation} \label{C3 : 47}
\left\{\begin{array}{l}
\sum_{\alpha=1}^3 d_{\alpha} Im(z_{2\alpha}\overline{z}_{2\alpha+1} +
w_{2\alpha}\overline{w}_{2\alpha+1})=0, \\
\sum_{\alpha=1}^3 d_{\alpha} (z_{2\alpha}w_{2\alpha+1 }-  z_{2\alpha+1} w_{2\alpha} ) =0,
\end{array}\right.
\end{equation}
\bfl
 where $d_{\alpha}$ indicates the weights  $p_{\alpha}$ and $q_{\alpha}$
for $\alpha\in\{1,2,3\}$. Now,
 consider  the case of $\overset{=}{S}_{(+,+,+)}$, the other one can be treated in the same way.
 The equations in (\ref{C3 : 47}) can be rewritten as\efl
\begin{equation}   \label{C3 : 48}
\left\{\begin{array}{l}
\sum_{\alpha=1}^3 d_{\alpha} (|z_{2\alpha}|^2+ |w_{2\alpha}|^2)=0,  \\
\sum_{\alpha=1}^3 d_{\alpha} (z_{2\alpha} w_{2\alpha}-
z_{2\alpha+1}w_{2\alpha+1})=0. \\
\end{array}\right.
\end{equation}
\bfl  Then it follows that equations in  (\ref{C3 : 48})
 admit only the trivial solution $(\underline{z},
 \underline{w})=(\underline{0},\underline{0})$.   Thus we get the conclusion.$\square$
\efl

\begin{prop}
\bfl Only one of the strata $\overset{=}{S}_{(\pm,\pm,\pm)}$ intersects
the submanifold
   $N(\Omega)$. Moreover, the intersection $\overset{=}{S}_{(\pm,\pm,\pm)}\cap N(\Omega)$
   yields exatly one isolated point on the twistor space
   $\mathcal{Z}^6(\Omega)$. \efl
\end{prop}

\bfl $\mathbf{Proof.}$  By reading
the $Sp(1)$ moment map equations in (\ref{C3: 3}) on any point
$ (\underline{z}, \underline{w}),$
 which belongs to one of the eigenspaces
 $\overset{=}{S}_{(\pm, \pm, \pm)}$
 we get the following system of equations\efl
 \begin{equation}     \label{C3 : 49}
\left\{\begin{array}{l}
\sum_{\alpha=1}^3  (|z_{2\alpha}|^2 -|w_{2\alpha}|^2)=0 ,\\
\sum_{\alpha=1}^3  \overline{w}_{2\alpha}z_{2\alpha}=0. \\
\end{array}\right.
 \end{equation}
\bfl  These equations are the same for all
 the  $\overset{=}{S}_{(\pm, \pm, \pm)}\big) $.
 Instead,  by evaluating
the $T^2_{\Omega} $ moment map equations  on a generic point
 $ (\underline{z}, \underline{w})$ which belongs
 to one of the spaces we are dealing with,
 we obtain\efl
 \begin{equation}     \label{C3 : 50}
\sum_{\alpha=1}^3  (-1)^{m_{\alpha}}d_{\alpha}(|z_{2\alpha}|^2 +|w_{2\alpha}|^2)=0,
 \end{equation}
 \bfl  for each $d_{\alpha} =p_{\alpha}, q_{\alpha}$ and $\alpha\in\{1,2,3\}$.
  Let us rewrite
 the equations in (\ref{C3 : 50}) as\efl
 \begin{equation}     \label{C3 : 51}
 \sum_{\alpha=1}^3  (-1)^{m_{\alpha}}d_{\alpha} |u_{2\alpha}|^2=0,
 \end{equation}
 \bfl where $u_{2\alpha} =z_{2\alpha}+jw_{2\alpha}$ for each $\alpha\in\{1,2,3\}$.
The intersection of the  sphere $S^{31}$ with one of the eigenspaces in
$\overset{=}{S}_{(\pm, \pm, \pm)}$ is described by the equation  \efl
\begin{equation}        \label{C3 : 52}
\sum_{\alpha=1}^7  |u_{\alpha}|^2=  \sum_{\alpha=1}^3   2|u_{2\alpha}|^2=1.
\end{equation}
\bfl Considering together the equations in (\ref{C3 : 51}) and (\ref{C3 : 52}) we get the following
system of equations\efl
\begin{equation}    \label{C3 : 53}
\left(\begin{array}{ccc}
(-1)^{m_1}p_1  & (-1)^{m_2}p_2 & (-1)^{m_3} p_3  \\
(-1)^{m_1}q_1  & (-1)^{m_2}q_2 & (-1)^{m_3} q_3  \\
 1             &  1             &   1         \\
\end{array} \right)
\left(\begin{array}{c}
  |u_{2}|^2 \\
   |u_{4}|^2 \\
    |u_{6}|^2\\
\end{array}\right) =
\left(\begin{array}{c}
  0  \\
  0\\
  \frac{1}{2} \\
\end{array}\right).
\end{equation}
\bfl  then we obtain\efl
\begin{equation}     \label{C3 : 54}
\left\{ \begin{array}{l}
2|u_{2}|^2 =\frac{\pm \Delta_{23} }{ \square^{1\pm 2}_{1\pm 3}}>0,\\
2|u_{4}|^2 =\frac{\pm \Delta_{13}}{ \square^{1\pm 2}_{1\pm 3}}>0, \\
2|u_{6}|^2 =\frac{\pm \Delta_{12}}{ \square^{1\pm 2}_{1\pm 3}}>0. \\
\end{array}\right.
\end{equation}
\bfl Note that, after we have fixed the weigth matrix $\Omega$, the minor determinants
$\Delta_{\alpha\beta}$ are univocally determinated.
Thus the system of
equations in (\ref{C3 : 54}) admits a unique solution. Then it follows that
$N(\Omega)$ intersects exactly one of
 the eigenspaces $\overset{=}{S}_{(\pm, \pm, \pm)}$.  $\square$  \efl

\bfl Now, consider the following sets:  \efl
 \begin{equation}      \label{C3 : 55}
 \tilde{G}\cdot  \overbrace{\Bigg\{ \left(\begin{array}{c|cc|cc|cc}
0 & u_2  & u_3 & u_4 & u_5 &  u_6 & u_7 \\
\end{array}\right)\Bigg\}}^{V:=},
 \end{equation}
\bfl such that each of the quaternionic coordeinate belongs to
one of the spaces ${}^{\pm}\tilde{V}_1^{\alpha}, \,\, {}^{\pm}\tilde{V}_2^{\alpha},\,\, \tilde{V}_4^{\alpha}$ and
 $\tilde{V}_5^{\alpha},\,\, \alpha\in\{1,2,3\}$.  Then, by using the  argument
in lemma $(2.7)$, we get that the possible strata in (\ref{C3 : 55}), which intersect
$N(\Omega)$, are just those one listed in the following two families

 \begin{equation}  \label{C3 : 56}
 \begin{aligned}
& 1) \,\,\, \quad  {}^+\overset{=}{S}{}^{12}_3:=\tilde{G}\cdot \Bigg\{
\left( \begin{array}{c|cc|cc|cc}
0 & z_2  & z_3 & z_4   &    z_5 &  0    &  0   \\
0   &    0 &  0  & 0     &   0             & w_6  &  iw_6\\
\end{array}\right) \Bigg\},   \\
& 2) \,\,\, \quad  {}^-\overset{=}{S}{}^{12}_3:=\tilde{G}\cdot \Bigg\{
\left( \begin{array}{c|cc|cc|cc}
0 & z_2  & z_3 & z_4   &   z_5   & 0    &  0   \\
0   &    0 &  0  & 0     &   0     & w_6  &  -iw_6\\
\end{array}\right) \Bigg\},   \\
&3) \,\,\,\quad  {}^+\overset{=}{S}{}^{13}_2:=\tilde{G}\cdot \Bigg\{
\left( \begin{array}{c|cc|cc|cc}
0 & z_2  & z_3    &    0   & 0           & z_6   & z_7  \\
0 &    0 &  0     &   w_4  &  iw_4   &  0    &  0          \\
\end{array}\right) \Bigg\},  \\
& 4)\,\,\,\quad   {}^-\overset{=}{S}{}^{13}_2:=\tilde{G}\cdot \Bigg\{
\left( \begin{array}{c|cc|cc|cc}
0 & z_2  & z_3        &   0   &    0          & z_6   & z_7   \\
0 &    0 &     0      &   w_4 &    -iw_4      & 0     & 0\\
\end{array}\right) \Bigg\},   \\
 & 5)\,\,\, \quad  {}^+\overset{=}{S}{}^{23}_1:=\tilde{G}\cdot \Bigg\{
\left( \begin{array}{c|cc|cc|cc}
    0 & 0    &  0      &    z_4 & z_5             & z_6   & z_7  \\
    0 & w_2  & iw_2    &   0    &  0              &  0    &  0          \\
\end{array}\right) \Bigg\},   \\
& 6)\,\,\, \quad  {}^-\overset{=}{S}{}^{23}_1:=\tilde{G}\cdot \Bigg\{
\left( \begin{array}{c|cc|cc|cc}
0   & 0    &  0           &    z_4 & z_5      & z_6   & z_7  \\
0   & w_2  & -iw_2        &   0    &  0       &  0    &  0          \\
\end{array}\right) \Bigg\}, \hspace{10 cm }  \\
\end{aligned}
\end{equation}
 and
\begin{equation}     \label{C3 : 57}
\begin{aligned}
& 7) \,\,\, \quad  \overset{=}{S}{}^{12}_3:=\tilde{G}\cdot \Bigg\{
\left( \begin{array}{c|cc|cc|cc}
0 & z_2  & z_3 & z_4   &    z_5 &  0    &  0   \\
0   &    0 &  0  & 0     &   0             & w_6  &  w_7\\
\end{array}\right) \Bigg\},   \\
& 8) \,\,\,\quad  \overset{=}{S}{}^{13}_2:=\tilde{G}\cdot \Bigg\{
\left( \begin{array}{c|cc|cc|cc}
0 & z_2  & z_3    &    0   & 0           & z_6   & z_7  \\
0 &    0 &  0     &   w_4  &  w_5   &  0    &  0          \\
\end{array}\right) \Bigg\},  \\
& 9)\,\,\, \quad  \overset{=}{S}{}^{23}_1:=\tilde{G}\cdot \Bigg\{
\left( \begin{array}{c|cc|cc|cc}
    0 & 0    &  0      &    z_4 & z_5             & z_6   & z_7  \\
    0 & w_2  & w_3    &   0    &  0              &  0    &  0          \\
\end{array}\right) \Bigg\}.  \hspace{10 cm }  \\
\end{aligned}
\end{equation}
\efl

 \begin{teor}
\bfl The strata listed in $(\ref{C3 : 56})$ and $(\ref{C3 : 57})$
 intersect the submanifold $N(\Omega)$ and their intersections are such that,
 $dim\,  \big(\overset{=}{S}{}^{\alpha \beta}_{\gamma}\cap N(\Omega)\big)= 0 $. Moreover,
  each of the quotients
 $\big(\, \overset{=}{S}{}^{\alpha \beta}_{\gamma}\cap N(\Omega)\big) / \tilde{G}$ describes a point on
 $\mathcal{Z}^6(\Omega)$ and we have
 \begin{equation}     \label{C3 : 58}
 \begin{aligned}
& 1) \quad \quad \overset{=}{S}{}^{12}_3\cap N(\Omega)\,\, \mathrm{  is\,\, formed\,\, by }\,\,
     {}^+\overset{=}{S}{}^{12}_3\cap N(\Omega) \,\, \mathrm{and} \,\, {}^-\overset{=}{S}{}^{12}_3\cap N(\Omega), \\
& 2) \quad \quad \overset{=}{S}{}^{13}_2\cap N(\Omega)\,\, \mathrm{  is\,\, formed\,\, by }\,\,
     {}^+\overset{=}{S}{}^{13}_2\cap N(\Omega) \,\, \mathrm{and} \,\, {}^-\overset{=}{S}{}^{13}_2\cap N(\Omega), \\
& 3) \quad \quad \overset{=}{S}{}^{23}_1\cap N(\Omega)\,\, \mathrm{  is\,\, formed\,\, by }\,\,
     {}^+\overset{=}{S}{}^{23}_1\cap N(\Omega) \,\, \mathrm{and} \,\, {}^-\overset{=}{S}{}^{23}_1\cap N(\Omega). \\
 \end{aligned}    \hspace{5 cm}
 \end{equation}
  \efl
\end{teor}

\bfl $\mathbf{Proof.}$
The moment maps equations, relative to $Sp(1)$ and $T^2_{\Omega}$,
 on the strata listed in (\ref{C3 : 57}), up to a sign,
become \efl
\begin{equation}     \label{C3 : 59}
\begin{aligned}
& T^2_{\Omega}\quad \Rightarrow\quad
 \left\{\begin{array}{l}
p_1Im(z_{2}\overline{z}_3) + p_2Im(z_{4}\overline{z}_5) +p_3Im(w_{6}\overline{w}_7)=0,\\
 q_1Im(z_{2}\overline{z}_3)+ q_2Im(z_{4}\overline{z}_5) +q_3Im(w_{6}\overline{w}_7)=0, \\
 \end{array}\right.      \\
 \mathrm{and}   & \\
&  Sp(1) \quad \Rightarrow\quad
  \left\{\begin{array}{l}
  |z_2|^2 +  |z_3|^2 + |z_4|^2  + |z_5|^2 - |w_6|^2  -|w_7|^2=0, \\
  (z_2)^2 + (z_3)^2 + (z_4)^2 +(z_5)^2  =0,       \\
 (w_6)^2 + (w_7)^2 =0,    \\
 \end{array} \right.      \hspace{8 cm} {}
 \end{aligned}
  \end{equation}
\bfl and the sphere equation has to be added
\begin{equation}   \label{C3 : 60}
 |z_2|^2 +  |z_3|^2 + |z_4|^2  + |z_5|^2 + |w_6|^2  +|w_7|^2 =1.
\end{equation}
In the following we use often the notation introduced in (\ref{C3 : 59}), for the
$Sp(1)$ and $T^2_{\Omega}$ moment map equations.
By combining the equations  (\ref{C3 : 59}) and (\ref{C3 : 60}) we get that
\begin{equation} \label{C3 : 61}
 \left\{ \begin{array}{c}
 |w_6|^2  +|w_7|^2 =1,   \\
  (w_6)^2 + (w_7)^2 =0,   \\
 \end{array} \right.
\end{equation}
 that is, $w_7= \pm iw_6$. As a direct consequence, we can rewrite the equations
in (\ref{C3 : 59}) as \efl
\begin{equation}  \label{C3 : 62}
\begin{aligned}
& \quad  T^2_{\Omega}\quad \Rightarrow\quad
 \left\{\begin{array}{l}
p_1Im(z_{2}\overline{z}_3) + p_2Im(z_{4}\overline{z}_5) \pm 2p_3|w_{6}|^2=0,\\
 q_1Im(z_{2}\overline{z}_3)+ q_2Im(z_{4}\overline{z}_5) \pm 2q_3|w_{6}|^2=0, \\
 \end{array}\right.             \\
  \mathrm{and} &                          \\
&  Sp(1) \quad \Rightarrow\quad
  \left\{\begin{array}{l}
  |z_2|^2 +  |z_3|^2 + |z_4|^2  + |z_5|^2 - 2|w_6|^2  =0, \\
  (z_2)^2 + (z_3)^2 + (z_4)^2 +(z_5)^2  =0.      \\
 \end{array} \right.      \hspace{15 cm} {}
 \end{aligned}
  \end{equation}

 \bfl Then, a direct computation shows that these equations describe the same submanifold
  generated   by the intersection of the strata in (\ref{C3 : 56}) with $N(\Omega)$.
  We get
 \begin{equation}    \label{C3 : 62a}
 dim\, \big(\overset{=}{S}{}^{\alpha\beta}_{\gamma}\cap N(\Omega)\big)
 = dim\,\big( {}^{\pm}\overset{=}{S}{}^{\alpha\beta}_{\gamma}\cap N(\Omega) \big)= 6.
 \end{equation}
 Hence
 $dim\, \big( \big(\,\overset{=}{S}{}^{\alpha\beta}_{\gamma}\cap N(\Omega) \big)/\tilde{G} \big)
 = dim\, \big({}^{\pm}\overset{=}{S}{}^{\alpha\beta}_{\gamma}\cap N(\Omega)/\tilde{G} \big)=0$, that
 is, each of the strata in (\ref{C3 : 56}) and (\ref{C3 : 57}) describes exactly a point on the twistor level.
  Since the strata in (\ref{C3 : 56}) are naturally included
 in those in (\ref{C3 : 57}), it follows that they produce the same points on
 $\mathcal{Z}^6(\Omega)$. $\square$\efl

\section{
The $3-$Sasakian Orbifolds $\mathcal{M}^7(\Omega)$}

\bfl As already mentioned in chapter two, in order to describe
the singularities relative
to the quaternion-$\ka$ orbifold $\mathcal{O}^4(\Omega)$, we
  study the singular locus of the $3-$Sasakian orbifold
  $\mathcal{M}^7(\Omega)$ and those one of
 $\mathcal{Z}^6(\Omega)$. If we assume
the hyphoteses of lemma $3.2$, in particular
 $gcd(\Delta_{12}, \Delta_{13}, \Delta_{23})= \pm 1$, we have that
 the action of $Sp(1)\times $
 $T^2_{\Omega}$ on $N(\Omega)\cap \{u_1\ne 0\}$ is free,
 then the singularities on  $\mathcal{M}^7(\Omega)$ come from the action
 of $Sp(1)\times T^2_{\Omega}$ on $N(\Omega)\cap \{u_1= 0\}$. Let us rewrite the
 fixed point equation (\ref{C3 : 8})  as follows:\efl
\begin{equation}  \label{C2a : 62b}
\overbrace{\left( \begin{array}{cc|cc}
0                          & -\sigma &  -sin\, \omega_{\alpha}  & cos\,\omega_{\alpha} - \overline{\epsilon}  \\
-\sigma     & 0                      & cos\,\omega_{\alpha} - \overline{\epsilon} & sin\, \omega_{\alpha}  \\
\hline
-sin\, \omega_{\alpha}    & cos\,\omega_{\alpha} - \epsilon           &    0       &      \overline{\sigma}  \\
cos\,\omega_{\alpha} -\epsilon  & sin\, \omega_{\alpha}  &  \overline{\sigma}   &   0    \\
\end{array} \right)}^{\tilde{M}_{\alpha}:=}
\left(\begin{array}{c}
z_{2\alpha} \\
z_{2\alpha+1}\\
w_{2\alpha} \\
w_{2\alpha+1} \\
\end{array} \right) =
\left(\begin{array}{c}
0 \\
0\\
0\\
0 \\
\end{array} \right),
\end{equation}
 \bfl where $\alpha\in \{1,2,3\}$.
  The following results
 hold  \efl

\begin{lem} \bfl For each $\alpha\in \{1,2,3\}$, the matrix  $\tilde{M}_{\alpha}$
in  $(\ref{C2a : 62b})$
is such that $det\, \tilde{M}_{\alpha}=0$ if and only if it holds either
\begin{equation} \label{C3 : 63}
e^{i\theta_{\alpha}}= Re( \epsilon ) + i\sqrt{ Im( \epsilon)^2 +
|\sigma |^2}\,\,\,\, \mathrm{or}
\,\,\,\,  e^{-i\theta_{\alpha}}= Re(\epsilon) + i\sqrt{ Im( \epsilon)^2 + |\sigma |^2}.
\end{equation}
\efl
\end{lem}

 \vspace{0,5 cm}

\begin{prop}
\bfl The solutions $\zw$ of the eigenvalues problem in $(\ref{C3 : 63})$
are such that their quaternionic pairs are either contained in one of the following
  $G$-invariant sets $\tilde{V}_1^{\alpha}$ and $\tilde{V}_2^{\alpha}$.\efl
 \begin{equation} \label{C3 : 64} {
\begin{aligned}
 1)\,\, &\quad \mathrm{If} \,\,\,\, e^{i\theta_{\alpha}}= Re( \epsilon ) +
 i\sqrt{ Im^2(\epsilon) + |\sigma |^2},\quad  \mathrm{then\,\, we\,\, have:}  \quad
 \tilde{V}_1^{\alpha}:=\{z_{2\alpha}, z_{2\alpha+1}\in \mathbb{C}\, | \\
& \, (z_{2\alpha},  z_{2\alpha+1},
\frac{e^{i\delta}  \big(-z_{2\alpha+1} +iz_{2\alpha}\cos\, \varphi  \big)}{ sin \, \varphi},
\frac{e^{i\delta}  \big(z_{2\alpha } +iz_{2\alpha+1}\cos\, \varphi  \big)}{ sin \, \varphi})
\in \mathbb{C}^4\},  \\
 &  \\
 2)\,\, &\quad \mathrm{If} \,\,\,\, e^{-i\theta_{\alpha}}= Re (\epsilon ) +
 i\sqrt{ Im^2(\epsilon) + |\sigma |^2},\quad  \mathrm{then\,\, we\,\, have:}  \quad
 \tilde{V}_2^{\alpha}:=\{z_{2\alpha}, z_{2\alpha+1}\in \mathbb{C}\, | \\
& \, (z_{2\alpha},  z_{2\alpha+1},
-\frac{e^{i\delta}  \big(z_{2\alpha+1} +i\zap\cos\, \varphi  \big)}{ sin \, \varphi},
\frac{e^{i\delta}  \big(\zap -iz_{2\alpha+1}\cos\, \varphi  \big)}{ sin \, \varphi})
\in \mathbb{C}^4\},
\end{aligned}     } \hspace{ 3 cm}
\end{equation}
\bfl where  $\epsilon = \cost +i(\sit cos\, \varphi)$,
$\sigma =  \sit sin \, \varphi cos\, \delta +i(\sit
   sin \, \varphi sin\, \delta )= $  \\
   $ \sit sin \, \varphi e^{i\delta} $ and $sin\varphi\ne 0$. \efl
\end{prop}

\bfl The proofs parallel those of lemma $2.6$ and proposition $2.9$.
Moreover,  we have the next proposition  \efl

\begin{prop}\bfl Consider the action of $G$ on $\mathbb{H}^7$.
Then  it holds that  the singularities relative to
 $\mathcal{M}^7(\Omega)$
 comes from  the intersection of $N(\Omega)$ with the following  strata
 on $\mathbb{H}^7$:\efl

 \begin{align*}
 & 1) \,\,\, \quad  \overline{S}{}^{123}:= 
  \bigg\{
\left( \begin{array}{c|cc|cc|cc}
0 & z_2              & z_3               & z_4                 &    z_5              & z_6                   &        z_7       \\
0 & \widetilde{w}_2  &  \widetilde{w}_3  & \widetilde{w}_4     &   \widetilde{w}_5   &  \widetilde{w}_6      & \widetilde{w}_7  \\
\end{array}\right) \bigg\},   \\
& 2) \,\,\, \quad \overline{S}{}^{12}_3:= 
\bigg\{
\left( \begin{array}{c|cc|cc|cc}
0 & z_2   &   z_3  & z_4      &    z_5   & z_6            & z_7       \\
0 & \widetilde{w}_2  &  \widetilde{w}_3  & \widetilde{w}_4     &   \widetilde{w}_5   &  \overset{\approx}{w}_6          & \overset{\approx}{w}_7  \\
\end{array}\right) \bigg\},   \hspace{10 cm}\\
&3) \,\,\,\quad \overline{S}{}^{13}_2:= 
 \bigg\{
\left( \begin{array}{c|cc|cc|cc}
0 & z_2  & z_3      & z_4            &    z_5                & z_6           & z_7    \\
0 & \widetilde{w}_2  &  \widetilde{w}_3  & \overset{\approx}{w}_4     &  \overset{\approx}{w}_5      &  \widetilde{w}_6         & \widetilde{w}_7  \\
\end{array}\right) \bigg\},   \\
 & 4)\,\,\, \quad \overline{S}{}^{23}_1:= 
  \bigg\{
\left( \begin{array}{c|cc|cc|cc}
0 & z_2  & z_3 & z_4   &    z_5 & z_6        & z_7       \\
0 &\overset{\approx}{w}_2 & \overset{\approx}{w}_3  & \widetilde{w}_4    &  \widetilde{w}_5  &  \widetilde{w}_6   & \widetilde{w}_7  \\
\end{array}\right) \bigg\},   \\
& 5)\,\,\,\quad \overline{S}{}^{1}_{23}:= 
\bigg\{\left( \begin{array}{c|cc|cc|cc}
0 & z_2  & z_3 & z_4   &    z_5 & z_6        & z_7    \\
0  & \widetilde{w}_2  &  \widetilde{w}_3  & \overset{\approx}{w}_4    &   \overset{\approx}{w}_5  & \overset{\approx}{w}_6        & \overset{\approx}{w}_7  \\
\end{array}\right) \bigg\},   \\
& 6)\,\,\,\quad \overline{S}{}^2_{13}:= 
 \bigg\{
\left( \begin{array}{c|cc|cc|cc}
0 & z_2   & z_3    & z_4      &    z_5 & z_6        & z_7      \\
0 & \overset{\approx}{w}_2 & \overset{\approx}{w}_3  & \widetilde{w}_4     &  \widetilde{w}_5  &  \overset{\approx}{w}_6          & \overset{\approx}{w}_7  \\
\end{array}\right) \bigg\}, \hspace{10 cm}  \\
\end{align*}
\begin{equation}   \label{C3 : 65} {
 \begin{aligned}
& 7)\,\,\,\quad  \overline{S}{}^3_{12}:= 
\bigg\{
\left( \begin{array}{c|cc|cc|cc}
0 & z_2  & z_3 & z_4   &    z_5 & z_6        & z_7      \\
0 & \overset{\approx}{w}_2 & \overset{\approx}{w}_3  & \overset{\approx}{w}_4     &  \overset{\approx}{w}_5 & \widetilde{w}_6   & \widetilde{w}_7  \\
\end{array} \right)    \bigg\}, \\
& 8)\,\,\,\quad  \overline{S}{}_{123} := 
 \bigg\{
\left( \begin{array}{c|cc|cc|cc}
0 & z_2  & z_3 & z_4   &    z_5 & z_6        & z_7    \\
0 & \overset{\approx}{w}_2 & \overset{\approx}{w}_3  & \overset{\approx}{w}_4     &   \overset{\approx}{w}_5 & \overset{\approx}{w}_6        & \overset{\approx}{w}_7  \\
\end{array}\right) \bigg\}.
 \end{aligned}    }    \hspace{4,5 cm}    {}
 \end{equation}
 \bfl where $(\widetilde{w}_{2\alpha}, \widetilde{w}_{2\alpha +1}):=
(\frac{e^{i\delta}  \big(-z_{2\alpha+1} +iz_{2\alpha}\cos\, \varphi
  \big)}{ sin \, \varphi},
\frac{e^{i\delta}  \big(z_{2\alpha } +iz_{2\alpha+1}\cos\, \varphi
 \big)}{ sin \, \varphi})\}$ and $(\overset{\approx}{w}_{2\alpha},$
 $ \overset{\approx}{w}_{2\alpha +1}):=
( \frac{e^{i\delta}  \big(-z_{2\alpha+1} -i\zap\cos\, \varphi  \big)}{ sin \, \varphi},
\frac{e^{i\delta}  \big(\zap -iz_{2\alpha+1}\cos\, \varphi  \big)}{ sin \, \varphi})\}$ with $sin\varphi\ne 0$.  \efl
\end{prop}

 \vspace{0,5 cm}

\bfl $\mathbf{Remark\, 3.4}$ By substituting the condition
 $e^{i\theta_{\alpha}}= Re( \epsilon ) +  i\sqrt{ Im( \epsilon)^2 +
|\sigma |^2}$ in  equations (\ref{C2a : 62b}), as  solutions we get exactly
 $ \{(z_{2\alpha}, z_{2\alpha +1},  \widetilde{w}_{2\alpha}, \widetilde{w}_{2\alpha +1})\}.$
 Instead, by substituting in the equations (\ref{C2a : 62b}) the condition
  $e^{-i\theta_{\alpha}}= Re( \epsilon ) +  i\sqrt{ Im( \epsilon)^2 +
|\sigma |^2}$, as  solutions we obtain
 $ \{(z_{2\alpha}, z_{2\alpha +1}, \overset{\approx}{w}_{2\alpha}, \overset{\approx}{w}_{2\alpha +1})\}.$
Then, by using together the relations (\ref{C3 : 63})
and the fixed point equations
$(2.11)$ in \cite{bg}, it follows that the fixed point equations for the points  in
 the strata (\ref{C3 : 65}) become
\begin{equation}   \label{C3 : 66}
\left\{ \begin{array}{l}
e^{i(\theta_1\pm \theta_2)}=1, \\
e^{i(\theta_1\pm \theta_3)}=1, \\
e^{\pm i\theta_1}= Re(\epsilon) + i\sqrt{Im(\epsilon)^2  + |\sigma |^2},  \\
iRe(\epsilon) + j\sigma = \frac{\sit_1 ( u_3\overline{u}_2 - u_{2}\overline{u}_3)}{|u_2|^2 + |u_3|^2}, \\
\end{array} \right.
\end{equation}
and the determinat associated to the subsystem of the first two equations coincides
with one of the determinants $\square^{1\pm 2}_{1\pm 3}$.  $\square$\efl

\vspace{0,5 cm}

\bfl Now, by using the computations we have done in chapter 2  with respect to
 the $3-$Sasakian manifold
$\mathcal{M}^7(\Theta)$, we can compute the moment map equations for the subgroups
$T^2_{\Omega}$ and $Sp(1)$ on the strata listed in (\ref{C3 : 65}).
Then we get: \efl

\begin{prop}
\bfl Let us consider the $G-$strata $\overline{S}$ of $\mathbb{H}^7$
 listed in $(\ref{C3 : 65})$. Then, each of the singular $G-$strata $\widetilde{S}=
\overline{S}\cap N(\Omega)\subset N(\Omega)$ gives rise
a singular point on the $3-$
Sasakian orbifold $M^7(\Omega )$.   \efl
\end{prop}

\bfl $\mathbf{Proof.}$ Note that
the strata $\overline{S}{}^{12}_3$, $\overline{S}{}^{13}_2$, $\overline{S}{}^{23}_4$, $\overline{S}{}^{1}_{23}$,
   $\overline{S}{}^{2}_{13}$ and
    $\overline{S}{}^{3}_{12}$
 can be treated in the same way. Then, without loss of generality, we deal with
  the case of
 $\overline{S}{}^{12}_3$ and we get \efl
 \begin{equation}  \label{C3 : 67}
Sp(1)\,\, \Rightarrow\,\,
\left\{ \begin{array}{l}
2 \sum_{\alpha=1}^2
 Im( \zap\bzed)   + cos \, \varphi \big( \sum_{\alpha=1}^2( |\zap |^2 + |\zed |^2 ) \big) =0,  \\
\sum_{\alpha=1}^2 \big( (\zap )^2 + (\zed )^2 \big) \big) =0, \\
(z_6)^2 + (z_7)^2 =0.     \\
cos\, \varphi = \frac{Im(z_6\overline{z}_7)}{|z_6 |^2 + |z_7|^2}.  \\
 \end{array} \right.
 \end{equation}
 \bfl The moment map equations of the $2-$torus $T^2_{\Omega}$  become\efl
 \begin{equation}   \label{C3 : 68}
 T^2_{\Omega}\,\,  \Rightarrow\,\,
 \left\{ \begin{array}{l}
 \sum_{\alpha=1}^2 p_{\alpha} \big( 2Im(\bzed\zap) + ( |\zap |^2 + |\zed |^2 )cos\, \varphi  \big)=0,  \\
 \sum_{\alpha=1}^2 q_{\alpha} \big( 2Im(\bzed\zap) + ( |\zad |^2 + |\zed |^2 )cos\, \varphi  \big)=0,  \\
  Im( z_6\overline{z}_7)  + ( |z_6 |^2 + |z_7|^2 )cos\, \varphi= 0,  \\
 \end{array} \right.
\end{equation}
\bfl and  \efl
\begin{equation}   \label{C3 : 69}
  T^2_{\Omega}\,\,\,\, \Rightarrow\,\,
 \left\{ \begin{array}{l}
 \sum_{\alpha=1}^2 p_{\alpha}\big( (\zap )^2 + (\zed )^2 \big)  =0, \\
 \sum_{\alpha=1}^2 q_{\alpha}\big( (\zap )^2 + (\zed )^2 \big)  =0, \\
 (z_6)^2 + (z_7)^2 =0.    \\
 \end{array} \right.  \hspace{4,5 cm} {}
 \end{equation}
 \bfl For each of the two system of equations in
 (\ref{C3 : 68}) and (\ref{C3 : 69}), consider the subsystem given by the first three
 equations. These two subsystems have the same matrix of coefficients,
 whose
 determinant is $\Delta_{12}\ne 0$. Then \efl
 \begin{equation}    \label{C3 : 70}
  2Im(\zap\bzed) + ( |\zap |^2 + |\zed |^2 )cos\, \varphi =0, \quad   \alpha\in \{1,2,3\}
  \end{equation}     \label{C3 : 71}
 \bfl that is     \efl
 \begin{equation}
   cos\, \varphi  = \frac{2Im(\zap\bzed)}{|\zap |^2 + |\zed |^2},
 \end{equation}
 \bfl and \efl
 \begin{equation}    \label{C3 : 72}
(\zap )^2 + (\zed )^2  =0,
 \end{equation}
 \bfl where $\alpha\in \{1,2,3\}$. In particular it holds that
$ dim\, \big(\big(\,\overline{S}{}^{12}_{3}\cap N(\Omega)\big)/ G\big)=0$, that is,
 the stratum $\overline{S}{}^{12}_{3}$
  yields one singular point on the orbifolds $\mathcal{M}^7(\Omega)$.
  In order to study
   the cases $\overline{S}{}^{123}$,
   $\overline{S}{}_{123}$, we can consider just one of them. Consider $\overline{S}{}^{123}$.
   Then the fixed point equations
    become \efl
     \begin{equation} \label{C3 : 73}
Sp(1)\,\, \Rightarrow\,\,
\left\{ \begin{array}{l}
 \sum_{\alpha=1}^3 \big(
 2Im( \zap\bzed)   + cos \, \varphi ( |\zap |^2 + |\zed |^2 ) \big) =0,  \\
\sum_{\alpha=1}^3 \big( (\zap )^2 + (\zed )^2 \big) \big) =0.    \\
 \end{array} \right.
 \end{equation}
 \bfl For the $T^3_{\Omega}$ in $(2.118)$ becomes\efl
 \begin{equation}    \label{C3 : 74}
 T^3_{\Omega}\,\,  \Rightarrow\,\,
 \left\{ \begin{array}{l}
 \sum_{\alpha=1}^3 p_{\alpha} \big( 2Im(\bzed\zap) + ( |\zap |^2 + |\zed |^2 )cos\, \varphi  \big)=0,  \\
 \sum_{\alpha=1}^3 q_{\alpha} \big( 2Im(\bzed\zap) + ( |\zad |^2 + |\zed |^2 )cos\, \varphi  \big)=0,  \\
 \sum_{\alpha=1}^3 p_{\alpha}\big( (\zap )^2 + (\zed )^2 \big)  =0, \\
 \sum_{\alpha=1}^3 q_{\alpha}\big( (\zap )^2 + (\zed )^2 \big)  =0, \\
 \end{array} \right.
 \end{equation}
\bfl Also in this case
$ dim\, \big(\big(\,\overline{S}{}^{123}\cap N(\Omega)\big)/ G\big)=0$, that is,
each of the strata $\overline{S}{}^{123}$ and $\overline{S}{}_{123}$
  yields a singular point on the orbifolds $\mathcal{M}^7(\Omega)$. $\square$ \efl

\section{  Directions of future research}

\bfl  Here, we give a list of the next natural steps in order to improve this work:  \efl

 \begin{itemize}
  \item[(a)] First, it would be interesting to compute some topological invariants, \\

 \item[(b)] describe locally some of these metrics as solutions of the $SU(\infty)$
 Toda lattice equation  \cite{to}, \\

 \item[(c)] give a classification of this new orbifold's family as the one in the toric case. \\

\end{itemize}

\end{document}